\documentclass[a4paper,12pt]{amsart}
\usepackage[T1]{fontenc}
\usepackage{palatino}
\usepackage{upgreek}

\usepackage{upgreek}
\usepackage{amssymb}
\usepackage{amsmath}

\usepackage{esint}

\usepackage[dvipsnames]{xcolor}
\usepackage{tikz}
\usepackage{enumerate}
\usepackage{comment}
\usepackage{bm}

\newtheorem*{restatablethm}{Restatement of Theorem 1.2}
\newcounter{restatement}

\usepackage{hyperref} 

\hypersetup{
    colorlinks,
    linkcolor=blue,
    citecolor=blue,
    urlcolor=blue
}
\usetikzlibrary{arrows.meta}

\colorlet{FadedBrickRed}{BrickRed!50!white}

\usepackage{graphicx}

\usepackage{caption}
\usepackage{subcaption}
\numberwithin{equation}{section}

\usepackage[leftcaption]{sidecap}
\sidecaptionvpos{figure}{m}

\usepackage{subcaption}
\usepackage[bottom]{footmisc}

\usepackage{accents}

\usepackage[
  a4paper,
  left=23mm,
  right=23mm,
  top=25mm,
  bottom=27mm,
    bindingoffset=0mm
]{geometry}

\makeatletter
\renewcommand{\boxed}[1]{\text{\fboxsep=.2em\fbox{\m@th$\displaystyle#1$}}}
\makeatother

\DeclareRobustCommand{\righttriangle}{%
\begin{tikzpicture}[scale=0.1,baseline=0ex]
\draw (0,0) -- (3,0) -- (3,3) -- cycle;
\end{tikzpicture}}

\newcommand{\hollow}{%
\begin{tikzpicture}[scale=0.1,baseline=0ex]
\draw (0,0) -- (3,0) -- (3,3) -- cycle;
\draw (1,0.5) -- (2.5,0.5) -- (2.5,2) -- cycle;
\end{tikzpicture}%
}

\newcommand{\tre}{%
\begin{tikzpicture}[scale=0.1,baseline=0ex]
\draw (0,0) -- (3,0) -- (3,3) -- cycle;
\draw (1,0) -- (3,2);
\end{tikzpicture}%
}

\newcommand{\urt}{%
\begin{tikzpicture}[scale=0.1,baseline=0ex]
\draw (0,0) -- (0,3) -- (3,3) -- cycle;
\end{tikzpicture}%
}
\DeclareRobustCommand{\rt}{\righttriangle}
\usepackage{fancyhdr}
\newcommand{\interior}[1]{%
  {\kern0pt#1}^{\mathrm{o}}%
}
\definecolor{dgrey}{rgb}{0.7, 0.75, 0.71}
\definecolor{lgrey}{rgb}{0.8, 0.85, 0.81}

\newtheorem{thm}{Theorem}[section]

\newtheorem{cor}[thm]{Corollary}
\newtheorem{lemma}[thm]{Lemma}
\newtheorem{df}[thm]{Definition}
\newtheorem{proposition}[thm]{Proposition}

\newtheorem{conj}[thm]{Conjecture}

\newtheorem{thmalpha}{Theorem}

{
\theoremstyle{definition}
\newtheorem{example}[thm]{Example}
\newtheorem{rem}[thm]{Remark}
\newtheorem{que}[thm]{Question}
}
\usepackage{nomencl}
\makenomenclature
\author{Samuel G.\ G.\ Johnston and Joscha Prochno}

\begin{document}

\title[The macroscopic shape of Gelfand--Tsetlin patterns and free probability]{The macroscopic shape of \\Gelfand--Tsetlin patterns and free probability}

\maketitle

\begin{abstract}

A compression is a function $F:\mathbb{R}\times[0,1]\to[0,1]$ such that each $F(\cdot,\tau)$ is the distribution function of a measure of mass $\tau$, while each $F(x,\cdot)$ is increasing and $1$-Lipschitz. Compressions are continuum analogues of Gelfand--Tsetlin patterns: if $(t_{k,j})_{ 0 \leq j \leq k \leq n}$ is a Gelfand--Tsetlin pattern, setting $F(t_{k,j},k/n) = j/n$ and interpolating creates a compression. For a differentiable compression $F$, we define the compression entropy
\begin{align*}
\mathcal{H}[F]
:=\int_{-\infty}^\infty\int_0^1 F_x
\left\{-\log F_x+\log\sin(\pi F_\tau)+1-\log\pi\right\}
\mathrm{d}\tau\mathrm{d}x.
\end{align*}
If $\mu$ is absolutely continuous and compactly supported, we prove
\begin{equation*}
\sup\left\{\mathcal{H}[F]:F \text{ compression}, F(\cdot,1)\text{ is the distribution function of }\mu\right\}
=\chi[\mu],
\end{equation*}
where $\chi[\mu]$ is Voiculescu's free entropy. By identifying the Euler--Lagrange equations for $\mathcal{H}[F]$ with a Burgers equation for Cauchy transforms, we show that the supremum is attained uniquely by the free compression of free probability theory.

We also view Gelfand--Tsetlin patterns as Ginzburg--Landau $\nabla\phi$-interface models with a hard-core interaction, and compute the surface tension:
\begin{equation*}
\sigma(u_1,u_2)
=-\log(u_1+u_2)-\log\sin\left(\pi\frac{u_1}{u_1+u_2}\right)-1+\log\pi.
\end{equation*}
Finally, we prove that uniform $n$-dimensional Gelfand--Tsetlin patterns with deterministic bottom rows converging to $\mu$ satisfy a large deviation principle with speed $n^2$ and rate function
\begin{equation*}
I_\mu[F]=-\mathcal{H}[F]+\chi[\mu].
\end{equation*}
These results resolve a conjecture of Shlyakhtenko and Tao stating that the Euler--Lagrange equations for free compression arise from the statistical mechanics of interlacing point processes.

\end{abstract}

\noindent\textbf{MSC2020:}
Primary: 82B41, 82B20, 46L54. Secondary: 15A52, 60G55, 60F10, 49Q20

\vspace{0.5em}
\noindent\textbf{Keywords:}
Gelfand--Tsetlin pattern, free probability, free entropy, free compression, Burgers equation, Euler--Lagrange equations, stochastic interface model, Ginzburg--Landau $\nabla \phi$-interface model, surface tension, bead process,  random matrices, log-concave functions.

\tableofcontents

%%%%%%%%%%%%%%%%%%%%%%%%%%%%%%%%%%%%%%%%%%%%%%
\section{Introduction and motivation}
%%%%%%%%%%%%%%%%%%%%%%%%%%%%%%%%%%%%%%%%%%%%%%
\subsection{Gelfand--Tsetlin patterns}
\begin{figure}[ht]
    \centering
    \includegraphics[width=0.69\textwidth]{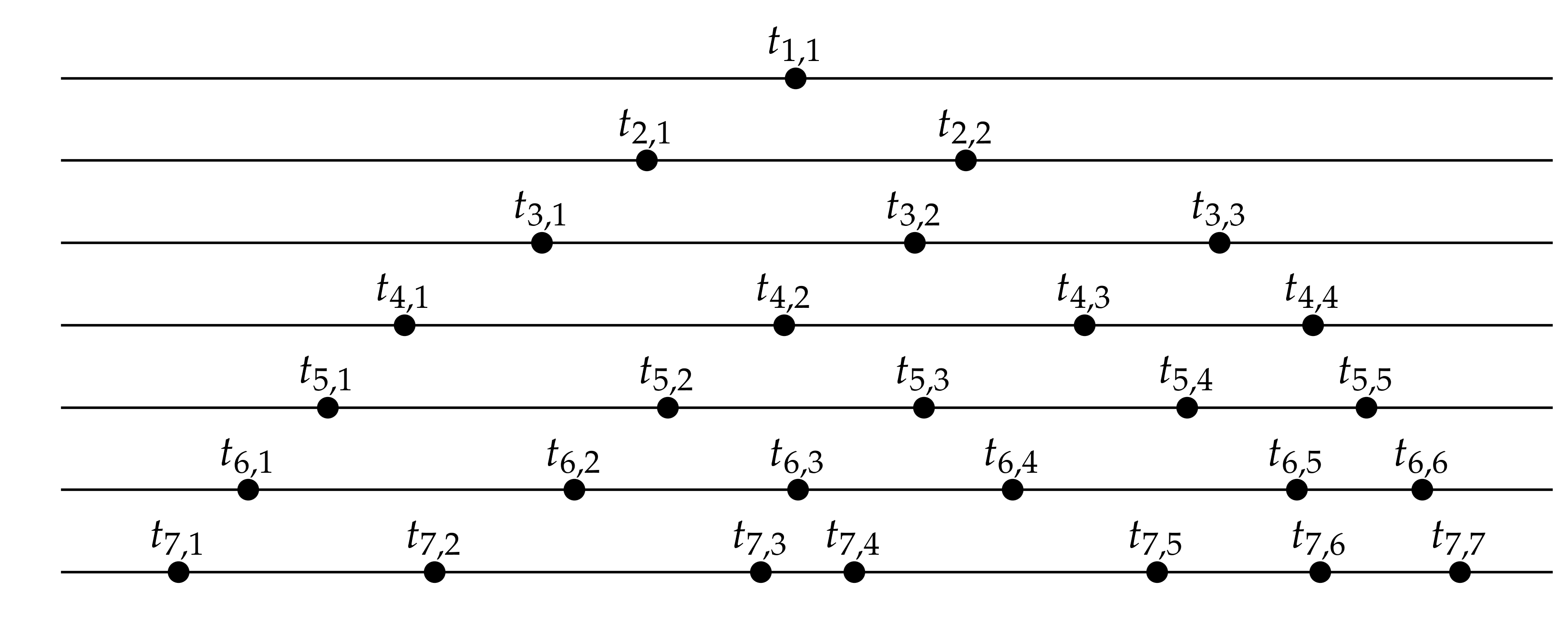}
    \caption{A Gelfand--Tsetlin pattern $(t_{k,j})_{1 \leq j \leq k \leq n}$ with $n = 7$.}
    \label{fig:GT7}
\end{figure}

Gelfand--Tsetlin patterns, introduced by I.~Gelfand and M.~Tsetlin in their study of representations of the unitary group \cite{GT1950}, are fundamental objects in representation theory, algebraic combinatorics, and random matrix theory. 
An ($n$-dimensional) \textbf{Gelfand--Tsetlin pattern} is a triangular array $(t_{k,j})_{1 \leq j \leq k \leq n}$ of real numbers satisfying the interlacing property
\begin{equation} \label{eq:GTp}
t_{k+1,j} \leq t_{k,j} \leq t_{k+1,j+1} \qquad \text{ for all } 1 \leq j \leq k \leq n-1.
\end{equation}
An example of a Gelfand--Tsetlin pattern is given in Figure \ref{fig:GT7}.

A natural way in which Gelfand--Tsetlin patterns arise is through the eigenvalues of minors of matrices. Indeed, for an $n$-by-$n$ Hermitian matrix $A$ with (real) eigenvalues $s_1 \leq \ldots \leq s_n$, define
\begin{align*}
t_{k,j} := \text{$j^{\text{th}}$ smallest eigenvalue of the $k\times k$ principal minor $(A_{i,j})_{1 \leq i,j \leq k}$ of $A$}.
\end{align*}
Then as a consequence of Cauchy's interlacing theorem, $(t_{k,j})_{1 \leq j \leq k \leq n}$ is a Gelfand--Tsetlin pattern with bottom row $t_{n,j} = s_j$, $1\leq j\leq n$. We call this array the eigenvalue process of the matrix $A$. 

The set $\mathrm{GT}(s_1,\ldots,s_n)$ of Gelfand--Tsetlin patterns with bottom row $s_1 \leq \cdots \leq s_n$ may be associated with a compact polytope in $\mathbb{R}^{n(n-1)/2}$ which has positive Lebesgue measure whenever the $s_i$ are distinct. By the Weyl dimension formula, its volume (in terms of $n(n-1)/2$-dimensional Lebesgue measure) is given by
\begin{align} \label{eq:wdf}
\mathrm{Leb}_{n(n-1)/2} \Big(\mathrm{GT}(s_1,\ldots,s_n)\Big) = G(n+1)^{-1} \prod_{ 1 \leq j < k \leq n } (s_k - s_j),
\end{align}
where $G(n) := \prod_{j =1}^{n-2} j!$ is the Barnes G-function (see, e.g., \cite[Lemma 1.12]{baryshnikov}). By associating the set $\mathrm{GT}(s_1,\ldots,s_n)$ with a compact subset of $\mathbb{R}^{n(n-1)/2}$ of finite volume, we can make sense of a \emph{uniformly chosen Gelfand--Tsetlin pattern with bottom row $s_1 < \ldots < s_n$}. It is also possible to make sense of this uniform measure when two or more of the $s_i$ are the same. 

We are particularly interested in \emph{random} Gelfand--Tsetlin patterns arising as the eigenvalue processes of random Hermitian matrices.
% We say that a random matrix $A$ is unitarily invariant if $A$ has the same law as $U^*AU$, where $U$ is a Haar unitary random matrix (i.e., distributed uniformly according to the unique, translation-invariant Haar measure on the unitary group) and $U^*$ is its conjugate transpose.
Here we have the following result due to Baryshnikov \cite[Proposition 4.7]{baryshnikov}. 

\begin{thm}[Baryshnikov \cite{baryshnikov}] \label{thm:Aar}
Let $A$ be an $n$-by-$n$ unitarily invariant random Hermitian matrix with eigenvalues $s_1 \leq \cdots \leq s_n$.
Then the eigenvalue process of $A$ is uniformly distributed on the set $\mathrm{GT}(s_1,\ldots,s_n)$ of Gelfand--Tsetlin patterns with bottom row $s_1 \leq \cdots \leq s_n$.  
\end{thm}

\subsection{Free probability} \label{sec:introfree}

Free probability was originally introduced by D.~Voiculescu to tackle problems in operator algebras, though it later proved to be a powerful tool in random matrix theory, where it describes the asymptotic behaviour of eigenvalues of large random matrices under the basic matrix operations: addition, multiplication, and taking minors \cite{Vcon1, Vcon2, Vcon3, Vent1, Vent2, Vent3}. We will be particularly interested in an operation on probability measures known as \textbf{free compression}, which describes how the empirical spectra of large random matrices behave under taking minors. Here we have the following result:

\begin{thm}[Voiculescu \cite{Vcon3, VIMRN}, Nica and Speicher \cite{NS}] \label{thm:compression}
For each probability measure $\mu$ and $\tau \in [0,1]$, there is a subprobability measure $\mu_\tau$ with the following property. 
Let $(A_n)_{n \geq 1}$ be a sequence of unitarily invariant Hermitian random matrices with deterministic empirical spectra converging weakly to $\mu$, and for $1 \leq j \leq k \leq n$ write $t^n_{k,j}$ for the $j^{\text{th}}$ smallest eigenvalue of the $k \times k$ principal minor of $A_n$. Then, for any sequence $(k_n)_{n \geq 1}$ with $1 \leq k_n \leq n$ and $k_n/n \to \tau \in (0,1)$, we have convergence in probability in the weak topology
\begin{align*}
\frac{1}{n} \sum_{j=1}^{k_n} \delta_{t^n_{k_n,j}} \stackrel{\mathbb P}{\longrightarrow} \mu_\tau, \quad{\text{as } n\to\infty}.
\end{align*}
\end{thm}

The measures $(\mu_\tau)_{\tau \in [0,1]}$ are sometimes called the free compressions of the probability measure $\mu$. In our atypical formulation, $\mu_\tau$ is not a probability measure but a subprobability measure of total mass $\tau$. 

Theorem \ref{thm:compression} is conventionally proved using asymptotic freeness: a unitarily invariant matrix becomes asymptotically free from the deterministic projection onto its first $k_n$ coordinates, and the limiting spectral measure of the corresponding compression is then identified through its moments, or equivalently its free cumulants.
However, in light of Theorem \ref{thm:Aar}, Theorem \ref{thm:compression} may instead be regarded purely as a statement about the macroscopic behaviour of large random Gelfand--Tsetlin patterns.
Our starting point is the following restatement of Theorem \ref{thm:compression} that makes no mention of Hermitian matrices or their eigenvalues:

\begin{restatablethm} \label{thm:restatement}
For each probability measure $\mu$ and $\tau \in [0,1]$, there is a subprobability measure $\mu_\tau$ with the following property. 
Let $(s^n)_{n \geq 1}$ be a sequence of vectors $s^n = (s_1^n,\ldots,s_n^n) \in \mathbb{R}^n$ with increasing coordinates and whose empirical coordinates converge weakly to $\mu$. Let $t^n := (t^n_{k,j})_{ 1 \leq j \leq k \leq n }$ be a sequence of random Gelfand--Tsetlin patterns such that $t^n$ is uniformly distributed on the set of Gelfand--Tsetlin patterns with bottom row $s^n$. Then, for any sequence $(k_n)_{n \geq 1}$ with $1 \leq k_n \leq n$ and $k_n/n \to \tau \in (0,1)$, we have convergence in probability in the weak topology
\begin{align*}
\frac{1}{n} \sum_{j=1}^{k_n} \delta_{t^n_{k_n,j}} \stackrel{\mathbb P}{\longrightarrow} \mu_\tau, \quad{\text{as } n\to\infty}.
\end{align*}
\end{restatablethm}

A central motivating goal of the article at hand is to give a proof of the Restatement of Theorem \ref{thm:restatement} directly from the entropic properties of Gelfand--Tsetlin patterns, without reference to random matrices or their moments.
Aside from its intrinsic interest as a problem in statistical mechanics, our goal is motivated by a desire to explain two notions arising 
 in free probability: Voiculescu's free entropy formula, and Shlyakhtenko and Tao's Euler--Lagrange equations for free compression. Beginning with the first of these notions, in a series of papers in the early 1990s \cite{Vent1, Vent2, Vent3}, Voiculescu introduced the notion of \textbf{free entropy} of a probability measure $\mu$, which is given by 
\begin{align} \label{eq:freeent}
\chi[\mu] = \frac{1}{2} \int_{-\infty}^\infty\int_{-\infty}^\infty \log|y-x| \mu(\mathrm{d}y)\mu(\mathrm{d}x) + \frac{3}{4},
\end{align}
and quantifies how favourably the mass of the measure is distributed with respect to logarithmic interactions. The logarithmic kernel $\log|y-x|$ assigns a larger contribution to pairs of points that are well separated, while points that are close together contribute negatively. Consequently, probability measures whose mass is more evenly spread out tend to have larger free entropy, whereas highly concentrated measures have smaller free entropy.

As observed in \cite{JO}, the free entropy of a probability measure may be regarded as a macroscopic analogue of the Weyl dimension formula \eqref{eq:wdf} for the volume of the Gelfand--Tsetlin polytope with boundary conditions imposed by $\mu$. Namely, let $Q_\mu:[0,1] \to \mathbb{R}$ be the quantile function of a probability measure $\mu$ with bounded support.
Then, recalling \eqref{eq:wdf}, the rescaled asymptotic volume of $(n+1)$-dimensional Gelfand--Tsetlin patterns with bottom row $nQ_\mu(0/n),\ldots,nQ_\mu(n/n)$ satisfies 
\begin{align} \label{eq:chi2}
\lim_{n \to \infty} \frac{1}{n(n+1)} \log \mathrm{Leb}_{n(n+1)/2} \mathrm{GT}\left( n Q_\mu(0/n), n Q_\mu(1/n), \ldots,nQ_\mu(n/n) \right) = \chi[\mu],
\end{align}
where we have used the asymptotics $\frac{1}{n^2} \log G(n) =  (1/2) \log n - 3/4+ o(1)$. There is some inconsistency in the literature with respect to the constants in \eqref{eq:freeent}; our selections were chosen to match the limit in \eqref{eq:chi2}. 

Another interpretation of Voiculescu's free entropy (the so-called microstates free entropy) is that for large $n$, it provides an asymptotic appraisal of the Lebesgue measure of the set $M_n(\mathbb{C})_{\mathrm{sa}}$ of $n$-by-$n$ Hermitian matrices whose empirical spectrum is close to $\mu$:
\begin{align} \label{eq:lebby}
\mathrm{Leb}_{M_n(\mathbb{C})_{\mathrm{sa}}}( \text{$A$ has spectrum close to $\mu$}) \approx \exp \left\{2  n^2 (\chi[\mu] + c_n  ) \right\},
\end{align}
where $c_n := -\frac{1}{4}\log n - \frac{3}{8}+\frac{1}{4}\log(2\pi)$; see \cite[Section 7]{MS} for a precise statement.

\medskip
%\subsubsection{Euler--Lagrange equations for free compression}
More recently, Shlyakhtenko and Tao \cite[Theorem 1.7]{ST} studied monotonicity properties of free entropy $\chi[\mu]$ under free convolution powers. In the process, they observed in formal calculations in free probability that the free compression measures $(\mu_\tau)_{\tau \in [0,1]}$ satisfy certain partial differential equations. Specifically, they construct a surface (see Section \ref{sec:thmA}) from the quantile functions of these measures which satisfies a partial differential equation coinciding with the Euler--Lagrange equations associated with a Lagrangian of the form
\begin{align} \label{eq:lagrangian}
L(u_1,u_2) := \log(u_1+u_2) + \log \sin\left( \frac{\pi u_1}{u_1+u_2} \right).
\end{align}
% We will construct this surface explicitly in Section \ref{sec:thmA}; note that we have used a slightly different coordinate system. 
Shlyakhtenko and Tao \cite{ST} conjectured that the source of this variational behaviour was the statistical mechanics properties of interlacing patterns: 

\begin{conj}[Shlyakhtenko and Tao \cite{ST}] \label{conj:ST}
The Lagrangian density $L(u_1,u_2)$ is proportional to the entropy of the bead process.
\end{conj}

\begin{figure}[ht]
    \centering
    \includegraphics[width=0.49\textwidth]{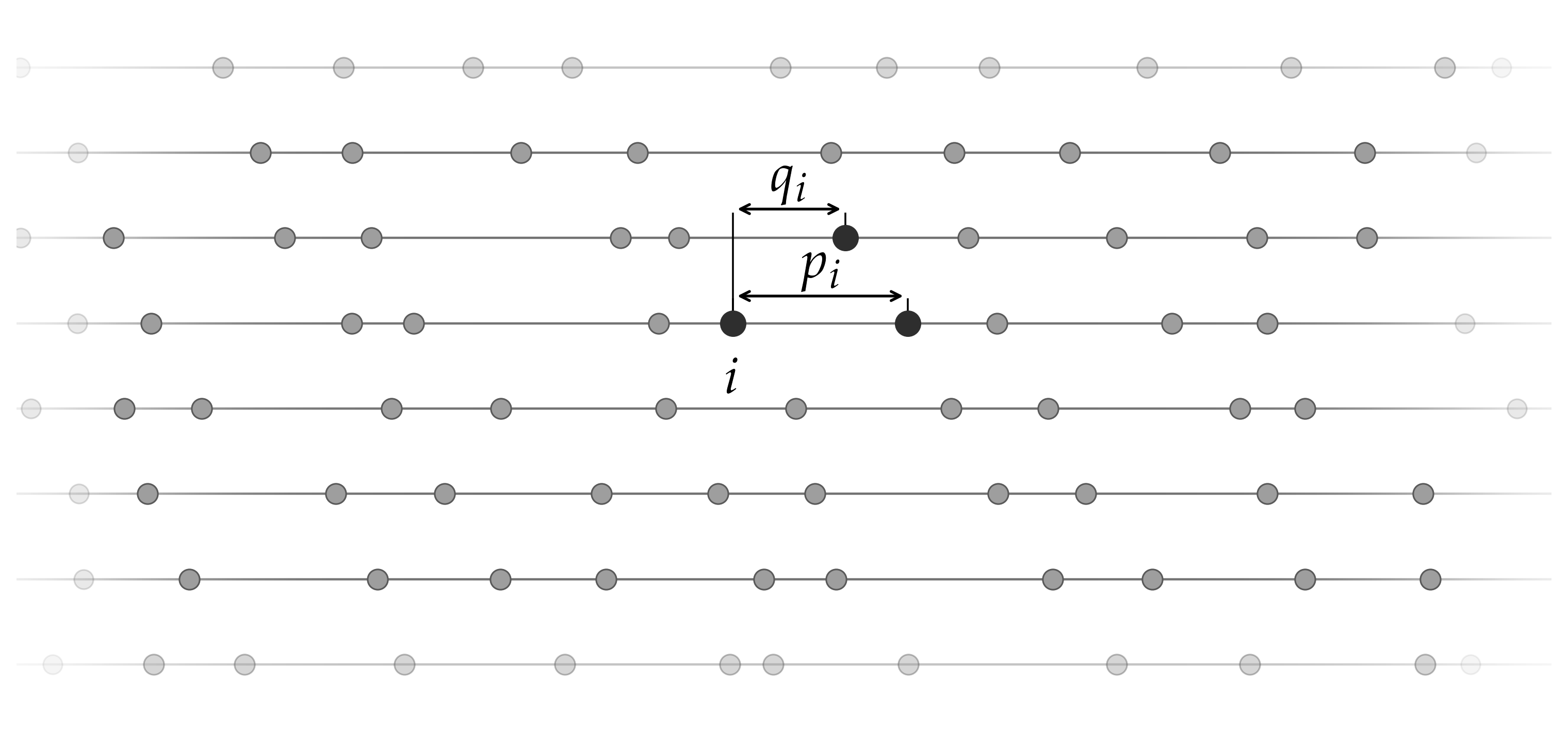}
    \caption{A finite patch within the bead process. The density of the process is the inverse of the average of the $p_i$'s. The tilt of the process is the quotient of the average of the $q_i$'s with the average of the $p_i$'s.}
    \label{fig:bead}
\end{figure}

This conjecture is inspired in part by work of Metcalfe \cite{metcalfe}, who proved that the asymptotic correlation behaviour of points inside the bulk of a large random Gelfand--Tsetlin pattern converges 
to that of a semi-discrete model in statistical mechanics called the bead process, introduced by Boutillier in \cite{boutillier}. Roughly speaking, the bead process is a random point 
process on $\mathbb{R} \times \mathbb{Z}$ chosen uniformly among point processes satisfying the interlacing property: between every two consecutive points $(s,k), (s',k)$ on $\mathbb{R} \times \{k\}$, there exists a unique $t \in [s,s')$ such that $(t,k+1)$ is a point of the process; see Figure \ref{fig:bead}. For each $k \in \mathbb{Z}$, we think of $\mathbb{R} \times \{k\}$ as a \emph{string}, and refer to these points as \emph{beads} (threaded along the string).

Boutillier \cite{boutillier} observed that there are two parameters governing the asymptotic behaviour of bead configurations: one controlling the density of beads per string and the other controlling the inter-string interaction. To define these, given a bead $i$ lying on some string, let $p_i > 0$ denote the distance to the next bead on the same string, and let $q_i \in (0,p_i)$ denote the distance to the next bead on the string above. Define parameters $u_1,u_2 > 0$ by 
\begin{align*}
u_1 &:= \text{Average of the $q_i$},\\
u_1 + u_2 &:= \text{Average of the $p_i$},
\end{align*}
where the averages are taken asymptotically over larger and larger regions of beads. Bead processes on $\mathbb{R} \times \mathbb{Z}$ have their behaviour governed by these two parameters: $u_1 +u_2$ is inversely proportional to the density of beads per string, and the ratio $\gamma := u_1/(u_1+u_2)$ (called the `\emph{tilt}' of the process \cite{johnston}) dictates the average position of neighbouring beads. 
Conjecture \ref{conj:ST} predicts that $L(u_1,u_2)$ somehow measures the average volume of possibilities for the position of a given bead in a large configuration with these parameters. 

As Shlyakhtenko and Tao noted, Conjecture \ref{conj:ST} is well supported by evidence in the statistical physics literature. Sun \cite{sun} used Legendre duality to argue that the formula $L(u_1,u_2)$ appears naturally through a deformation of free energies of lozenge tilings, the first author and O'Connell \cite{JO} reverse-engineered this same formula using the stability of the semicircle law under free compression, and Gordenko \cite{gordenko} recently obtained a similar formula based on scaling limits of Young diagrams. 
In Theorem \ref{thm:B} of the present article, we confirm Conjecture \ref{conj:ST} in the rigorous framework of surface tension for random surfaces. More precisely, we identify the surface tension of the corresponding hard-core $\nabla\phi$-interface model explicitly, and later use this formula to describe the macroscopic behaviour of large random Gelfand--Tsetlin patterns.

% Shlyakhtenko and Tao \cite{ST} note that $\sigma_0(u_1,u_2)$ also appears in work by Sun \cite{sun} on scaling limits of the dimer model, non-rigorous work \cite{JO} by the first author and Neil O'Connell on connections between Gelfand--Tsetlin patterns and the semicircle law, and more recent work by Gordenko on scaling limits of Young diagrams \cite{gordenko}. 

%%%%%%%%%%%%%%%%%%%%%%%%%%%%%%%%%%%%%%%%%%%%%%%%%%%%%%%%%%%%%%
\subsection{Summary of our main results}
%%%%%%%%%%%%%%%%%%%%%%%%%%%%%%%%%%%%%%%%%%%%%%%%%%%%%%%%%%%%%%

With a view to understanding the macroscopic shape of Gelfand--Tsetlin patterns, we introduce an object called a \textbf{compression}, which may be understood as a continuum analogue of a Gelfand--Tsetlin pattern. 

\begin{df}[Compression] \label{df:compression}
A compression $F$ is a function $F:\mathbb{R} \times [0,1] \to \mathbb{R}$ that is (weakly) increasing and right-continuous in $x$ with limits $\lim_{x \to -\infty} F(x,\tau) = 0$ and $\lim_{x \to +\infty}F(x,\tau) = \tau$, and satisfies the continuum interlacing inequalities
\begin{align} \label{eq:continuum}
F(x,\tau) \leq F(x,\tau' ) \leq F(x,\tau) + (\tau' - \tau) \qquad \text{for all $x \in \mathbb{R}$, $0 \leq \tau \leq \tau' \leq 1$}.
\end{align}
 
\end{df}

The continuum interlacing inequalities \eqref{eq:continuum} are equivalent to 
$F(x,\tau)$ being increasing and $1$-Lipschitz in $\tau$, for all $x \in \mathbb{R}$. 
 For each $\tau \in [0,1]$, $F(\cdot,\tau)$ may be thought of as the distribution function 
 \begin{align*}
 \mu_\tau((-\infty,x]) := F(x,\tau), \qquad x \in \mathbb{R},
 \end{align*}
of some measure $\mu_\tau$ of total mass $\tau$. 

 \begin{df}[A compression of a probability measure $\mu$] \label{df:compressionmu}
We say that $F$ is a compression of $\mu$ if $F$ is a compression such that $F(\cdot,1)$ is the distribution function of $\mu$, i.e.
\begin{align*}
F(x,1) = \mu( (-\infty,x]), \qquad x \in \mathbb{R}.
\end{align*}
 \end{df}

We can treat Gelfand--Tsetlin patterns of arbitrary dimensions as compressions. From here on, it will be more convenient for our purposes to index our Gelfand--Tsetlin patterns from $0 \leq j \leq k \leq n$ rather than $1 \leq j \leq k \leq n$. 

\begin{figure}[ht]
    \centering
    \includegraphics[width=1.01\textwidth]{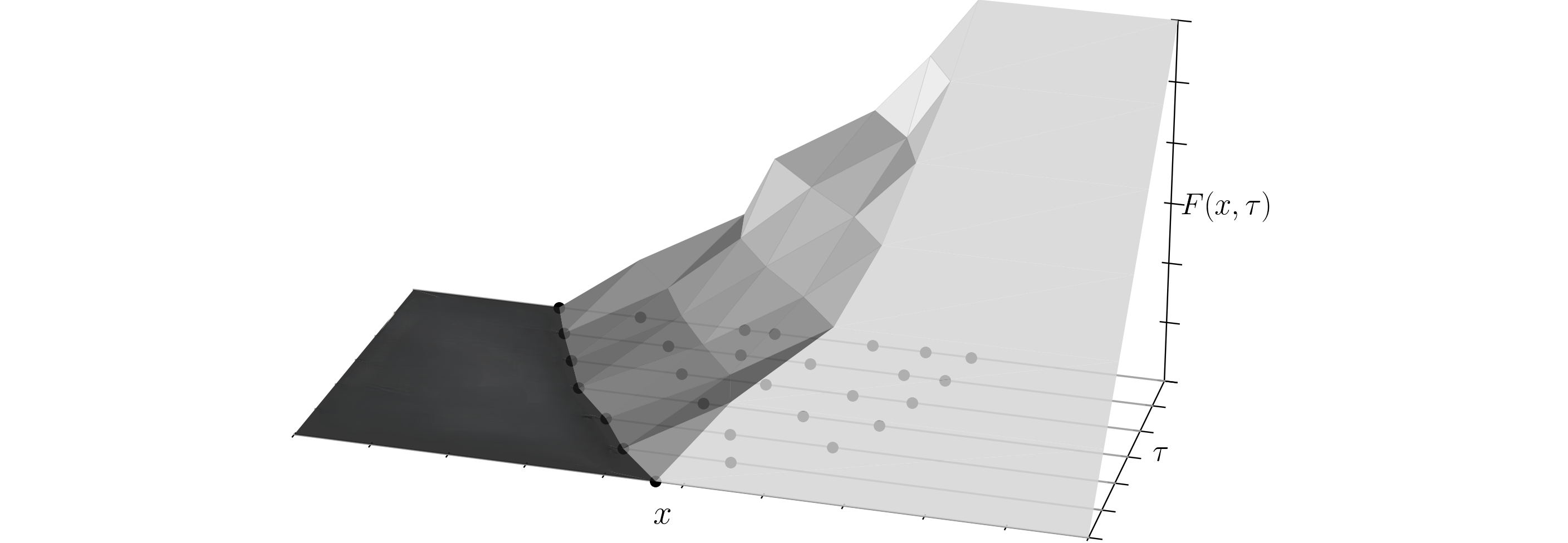}
    \caption{A Gelfand--Tsetlin pattern $(t_{k,j})_{0 \leq j \leq k \leq n}$ gives rise to a compression $\{F(x,\tau) : x \in \mathbb{R}, \tau \in [0,1] \}$ by letting $F(t_{k,j},k/n) := j/n$ and linearly interpolating.}
    \label{fig:gtF}
\end{figure}

\begin{df}[The compression associated with a Gelfand--Tsetlin pattern] \label{df:GTcomp}
The compression associated with a Gelfand--Tsetlin pattern $(t_{k,j})_{0 \leq j \leq k \leq n}$ with strictly increasing rows is defined by setting 
\begin{align*}
F(t_{k,j},k/n) := j/n, \qquad 0 \leq j \leq k \leq n,
\end{align*} linearly interpolating on triangles, and extending the definition naturally for all sufficiently small or large $x$. 
\end{df}

By linearly interpolating on triangles, we mean that we linearly interpolate the value of $F$ on all possible triangles taking the form $ABC$ or $AB'C$, where $A,B,B',C$ take the form $A = (t_{k,j},k/n), B = (t_{k-1,j}, (k-1)/n), B' = (t_{k+1,j+1},(k+1)/n), C = (t_{k,j+1}, k/n)$. We extend the definition naturally outside the union of these triangles by setting
$F(x,\tau)=0$ to the left of the leftmost boundary and $F(x,\tau)=\tau$ to the
right of the rightmost boundary. More specifically, if $\tau=(k+\lambda)/n$ for
some $0\leq k\leq n-1$ and $\lambda\in[0,1]$, then we set $F(x,\tau)=0$ whenever
$x\leq (1-\lambda)t_{k,0}+\lambda t_{k+1,0}$, and $F(x,\tau)=\tau$ whenever
$x\geq (1-\lambda)t_{k,k}+\lambda t_{k+1,k+1}$.
See Figure \ref{fig:gtF}.

We can extend the definition to Gelfand--Tsetlin patterns without strictly increasing rows as follows: if $x = t_{k,j} = \cdots = t_{k,j+\ell}$ we set $F(x,k/n) = (j+\ell)/n$. 

It is a brief exercise (which Figure \ref{fig:gtF} may assist) to verify that \eqref{eq:GTp} implies that \eqref{eq:continuum} is satisfied by the associated compression, so that Definition \ref{df:GTcomp} indeed constructs a compression. In fact, asymptotically speaking, the relationship goes both ways:

\begin{lemma}
\begin{enumerate}
\item Every Gelfand--Tsetlin pattern gives rise to a compression via Definition \ref{df:GTcomp}. 
\item Every compression arises as a limit of Gelfand--Tsetlin patterns in the Wasserstein metric (\eqref{eq:metric}). 
\end{enumerate}
\end{lemma}

See Proposition \ref{prop:closure} for a proof of the second part.

While in general there are many possible compressions of a probability measure $\mu$, there is a canonical compression of $\mu$ from the perspective of free probability. We will give an extrinsic definition here:

\begin{df}[Extrinsic definition of the free compression of $\mu$]  \label{def:freecompression}
 Given a probability measure $\mu$, let $(\mu_\tau)_{\tau \in [0,1]}$ be the free compression measures associated with $\mu$ (as in Theorem \ref{thm:compression}). The \textbf{free compression} $F^*$ is the compression given by 
 \begin{align} \label{eq:freecomp}
 F^*(x,\tau) := \mu_\tau ((-\infty,x]), \qquad x \in \mathbb{R}, \tau \in [0,1].
 \end{align}
\end{df}

We will give an intrinsic definition of free compression in terms of $R$-transforms in Section \ref{sec:Rtransform}, and prove there that this definition \eqref{eq:freecomp} indeed creates a compression.

\medskip
Recall that Lebesgue's monotone differentiation theorem says that a monotone function is differentiable almost everywhere. In particular, a compression $F$, which is increasing in both the $x$ and $\tau$ variable, has first-order derivatives $F_x$ and $F_\tau$ defined for almost all $(x,\tau) \in \mathbb{R} \times [0,1]$. We say that a compression $F$ is absolutely continuous if there exists a nonnegative measurable function $f:\mathbb{R} \times [0,1] \to [0,\infty)$ such that
\begin{align*}
F(x,\tau) := \int_{-\infty}^x f(y,\tau) \mathrm{d}y, \qquad x \in \mathbb{R},\tau \in [0,1].
\end{align*}
We call $f$ the density associated with the compression. We have $F_x= f$ in this case.
It is not necessary to make any regularity assumptions for the $\tau$ variable: the derivative $F_\tau$ exists in $[0,1]$ almost everywhere by monotonicity and Lipschitz continuity.
With this in mind, the following notion of entropy for a compression plays a central role in this article:

\begin{df}[The compression entropy] \label{df:compressionentropy}
If $F$ is an absolutely continuous compression, we define its compression entropy $\mathcal{H}[F] \in [-\infty,+\infty]$ by
\begin{align*}
\mathcal{H}[F] := \int_{-\infty}^\infty \int_0^1 F_x \left\{  - \log F_x + \log \sin \pi F_\tau + 1 - \log \pi \right\}   \mathrm{d}\tau \mathrm{d}x,
\end{align*}
where we are using the shorthand $F_x := \frac{\partial F}{\partial x}(x,\tau)$ and $F_\tau := \frac{\partial F}{\partial \tau}(x,\tau)$ for partial derivatives.

If $F$ is not absolutely continuous, we set $\mathcal{H}[F] = -\infty$. 
\end{df}

The compression entropy quantifies how favourable a compression is from the viewpoint of free probability. It combines two effects: a classical entropy term, which rewards spatially spread-out mass, and a correction which favours a balanced evolution in the $\tau$-direction. Thus, a compression has high compression entropy if it is neither overly concentrated in space nor excessively skewed in its evolution.

% We will see later that the integrand $F_x ( - \log F_x + \log \sin \pi F_\tau + \log \pi - 1 )$ has a statistical physics interpretation in terms of the surface tension of the bead model. 

With these definitions at hand, we are now ready to summarise our three main results. The first of these lies in free probability, the second in the rigorous statistical mechanics of random surfaces, and the last in the theory of large deviations:

\begin{itemize}
\item \textbf{Theorem \ref{thm:A}: A variational identity for compression entropy.}
We begin with a simple variational observation that will provide the foundation for the results that follow. Theorem \ref{thm:A} relates 
 Voiculescu's free entropy $\chi[\mu]$ of a measure $\mu$ \eqref{eq:freeent} to the compression entropy $\mathcal{H}[F]$ of a compression $F$ through the relation
\begin{align} \label{eq:relation}
\chi[\mu] = \sup_{ F: F \text{ compression of $\mu$} } \mathcal{H}[F] = \mathcal{H}[F^*],
\end{align}
where the supremum is taken over all compressions $F$ of $\mu$ (i.e., $F(\cdot,1)$ is the distribution function of $\mu$), and the latter equality says that the free compression $F^*$ is the unique compression attaining the supremum. Although the proof is short, this identity provides the foundation for understanding our subsequent large-deviation result.
\item \textbf{Theorem \ref{thm:B}: The surface tension of the bead model.} 
We show that Gelfand--Tsetlin patterns can be reformulated as two-dimensional Ginzburg--Landau $\nabla\phi$-interface models $\phi:\mathbb{Z}^2 \to \mathbb{R}$ with a hard-core 
interaction potential $V(\phi_y-\phi_x) = +\infty \mathrm{1}\{ \phi_x > \phi_y\}$; we call these random Gelfand--Tsetlin surfaces. 
These surfaces have a notion of surface tension, which measures the asymptotic size of the partition function of a large square patch of surface forced to lie at a certain gradient. 
Our second main result, Theorem \ref{thm:B}, states that the surface tension of Gelfand--Tsetlin surfaces is given by
\begin{align} \label{eq:sigma}
\sigma(u_1,u_2) = -\log(u_1+u_2) - \log \sin (\pi u_1/(u_1+u_2)) - 1 + \log \pi,
\end{align}
thereby providing a rigorous resolution of Conjecture \ref{conj:ST} with explicit constants. We note that to our knowledge, this is the only example of an explicit surface tension formula for a $\nabla\phi$-interface model with non-Gaussian interaction potential in two or more dimensions. 

\item \textbf{Theorem \ref{thm:C}: A large deviation principle for random Gelfand--Tsetlin patterns.} Our third result, Theorem \ref{thm:C}, is a large deviation principle for large random Gelfand--Tsetlin patterns linking together Theorem \ref{thm:A} and Theorem \ref{thm:B}. Roughly speaking, it says that if $(t^n_{k,j})_{0 \leq j \leq k \leq n}$ is a large random Gelfand--Tsetlin pattern with bottom row $(t^n_{n,j})_{ 0 \leq j \leq n}$ approximately empirically distributed like $\mu$, then for a compression $F$ of $\mu$ we have 
\begin{align} \label{eq:X}
\mathbf{P} \left( (t^n_{k,j})_{0 \leq j \leq k \leq n} \approx F \Big| (t^n_{n,j})_{ 0 \leq j \leq n} \approx \mu \right) \approx \exp \left\{ - n^2 I_\mu[F] \right\}, \qquad I_\mu[F] = - \mathcal{H}[F]+\chi[\mu ],
\end{align}
where we are using 
Definition \ref{df:GTcomp} to embed Gelfand--Tsetlin patterns of varying dimensions into the set of compressions. 
Equivalently, after conditioning a large unitarily invariant Hermitian matrix to have empirical spectrum close to $\mu$, the probability that its eigenvalue process is close to $F$ is governed by \eqref{eq:X}.
\end{itemize}

Theorem \ref{thm:A}, Theorem \ref{thm:B} and Theorem \ref{thm:C}, taken together, make sense of free compression as the maximiser of a statistical mechanics problem involving asymptotic Gelfand--Tsetlin patterns. Voiculescu's free entropy encodes the maximal volume of Gelfand--Tsetlin patterns with a fixed bottom row approximating $\mu$, and Shlyakhtenko and Tao's Euler--Lagrange equations codify the fact that free compression characterises the maximal shape. Let us emphasise that 
in light of \eqref{eq:relation} and the fact that the supremum in \eqref{eq:relation} is attained uniquely by the free compression, the rate functional $I_\mu$ in \eqref{eq:X} satisfies, for compressions $F$ of $\mu$, 
\begin{align*}
I_\mu[F] \geq 0 \qquad \text{ with equality if and only if $F$ is the free compression of $\mu$}.
\end{align*}
In particular, \eqref{eq:X} guarantees the convergence of the uniform Gelfand--Tsetlin pattern with bottom row approximating $\mu$ to the free compression of $\mu$. Thus, Theorem \ref{thm:C} encodes a significant generalisation of the free probability result, Theorem \ref{thm:compression}. We also emphasise that \eqref{eq:X} gives rise to the extension
\begin{align} \label{eq:lebby2}
\mathrm{Leb}_{M_n(\mathbb{C})_{\mathrm{sa}}}( \text{$A$ has eigenvalue process close to $F$}) \approx \exp \left\{2  n^2 (\mathcal{H}[F] + c_n  ) \right\}
\end{align}
of \eqref{eq:lebby}. To our knowledge, this formula offers the first explicit free entropy calculation beyond the spectrum of a single matrix itself as in \eqref{eq:lebby}.

%%%%%%%%%%%%%%%%%%%%%%%%%%%%%%%%%%%%%%%%%%%%%%%%%%%%%%%%%%%%%%%%%%%%%%%%%%%%%%%%%%%%%%%%%%%%%%%%%%%%%%%%%
\subsection{Overview}
%%%%%%%%%%%%%%%%%%%%%%%%%%%%%%%%%%%%%%%%%%%%%%%%%%%%%%%%%%%%%%%%%%%%%%%%%%%%%%%%%%%%%%%%%%%%%%%%%%%%%%%%%

The remainder of the paper is organised as follows.

\begin{itemize}
\item In Section \ref{sec:results} we give precise statements of our main results.
\item In Section \ref{sec:discussion} we discuss further related literature, and touch on some of the tools developed in our approach.
\end{itemize}

The rest of the paper is dedicated to proofs, and is divided into three parts.
The first is concerned with the proof of Theorem \ref{thm:A}:
\begin{itemize}
\item In Section \ref{sec:free} we show that free compression $F^*$ satisfies the Euler--Lagrange equations for the compression entropy, and moreover prove the relation $\mathcal{H}[F^*] = \chi[\mu]$. Strictly speaking, we do not provide here a fully rigorous proof that $F^*$ maximises $\mathcal{H}$ over all compressions of $\mu$. This latter fact will follow from our work in Section \ref{sec:Cproof}. 
\end{itemize}

The second part is devoted to the proof of Theorem \ref{thm:B} and related results:
\begin{itemize}
\item In Section \ref{sec:randomsurfaces} we study some basic properties of random surfaces, and prove a pair of Pr\'ekopa--Leindler-type inequalities for random surfaces with convex potentials.
\item In Section \ref{sec:import} we work towards the proof of Theorem \ref{thm:B}, in fact proving a more detailed version with explicit bounds on the second-order terms. We do this by importing several results from \cite{johnston} on the bead model on the torus, and giving a construction that unravels toric bead configurations to create Gelfand--Tsetlin surfaces. We ultimately prove that if $T_n(u_1,u_2)$ denotes the partition function of Gelfand--Tsetlin surfaces on a triangle with side length $n$ and boundary given by the linear function $\varphi^u_x := u_1x_1 + u_2x_2$, then
\begin{align*}
\frac{1}{n^2/2} \log T_n(u_1,u_2) = - \sigma(u_1,u_2) + O_u(\log(n)/n),
\end{align*}
where $\sigma$ is as in \eqref{eq:sigma}.
\item In Section \ref{sec:ST2proof} we extend this result further, showing that if $T_n(\varphi)$ denotes the partition function of a Gelfand--Tsetlin surface on a triangle with side length $n$ and boundary given by $\varphi$, then provided $|\varphi_x - \varphi^u_x| \leq \delta n$ and $\varphi$ is sufficiently well spaced, we have
\begin{align*}
\frac{1}{n^2/2} \log T_n(\varphi) = - \sigma(u_1,u_2) + o_{\delta,n}(1),
\end{align*}
as $\delta \downarrow 0, n \uparrow \infty$. With this machinery at hand, we conclude with our proof of Theorem \ref{thm:B}.
\end{itemize}

The final part of the paper is devoted to the proof of our large deviation principle, Theorem \ref{thm:C}:
\begin{itemize}
\item In Section \ref{sec:metric}, as preparation we study properties of the space of compressions as a metric space, proving the upper semicontinuity of the entropy functional.
\item In the final section, Section \ref{sec:Cproof}, we complete the proof of Theorem \ref{thm:C}.
\end{itemize}

Throughout the article $C > 0$ will denote a sufficiently large universal constant that may change from appearance to appearance.
When a function is referred to as `increasing', this is always taken weakly. A function defined on a subset of $\mathbb{Z}^2$ or $\mathbb{R}^2$ is said to be increasing if it is weakly increasing in both coordinates.

% \subsection*{Acknowledgements}
% The authors are grateful to Octavio Arizmendi, Colin McSwiggen, Hariharan Narayanan and Jenny Pi for their valuable comments.
% \subsection*{AI Disclosure}
% AI was used to assist the generation of images in this article.

%%%%%%%%%%%%%%%%%%%%%%%%%%%%%%%%%%%%%%%%%%%%%%
\section{Main results} \label{sec:results}
%%%%%%%%%%%%%%%%%%%%%%%%%%%%%%%%%%%%%%%%%%%%%%

\subsection{Gelfand--Tsetlin patterns and compressions as surfaces} \label{sec:surfaces}

In order to prepare and further contextualise the statements of our main results, in this section we describe how both Gelfand--Tsetlin patterns and their continuum analogues, compressions, can be regarded as surfaces. In particular, Gelfand--Tsetlin patterns can be recast in the framework of random surfaces \cite{sheffield} (which are also called stochastic interface models \cite{funaki} or Ginzburg--Landau $\nabla \phi$-interface models \cite{FS} in the wider literature). Passing to the continuum limit through this construction, we obtain a natural formulation of compressions as increasing functions defined on a triangular subset of $\mathbb{R}^2$. 

\subsubsection{Uniform random Gelfand--Tsetlin patterns as random surfaces} \label{sec:GTsurface}

Given a finite and graph-theoretically connected subset $C$ of the lattice $\mathbb{Z}^2$, let $C^*$ denote the collection of directed nearest-neighbour edges $\langle x,y \rangle$ in $C$. These edges are directed in that for all such edges either
\begin{align*}
y = x + \mathbf{e}_1 \qquad \text{or}\qquad y = x + \mathbf{e}_2.
\end{align*}
Suppose $C$ is decomposed as a disjoint union $C = A \sqcup B$. Typically $B$ will be a boundary set. A \textbf{random surface} with boundary conditions $\{ \varphi_x \,:\, x \in B \}$ is a random real-valued function on $C$ distributed according to a probability measure on $\mathbb{R}^C$ taking the form
\begin{align} \label{eq:rs}
\frac{1}{Z_A(\varphi)} \exp \left\{ - \sum_{ \langle x,y \rangle \in C^*} V(\phi_y - \phi_x) \right\} \prod_{x \in B} \delta_{\varphi_x}(\mathrm{d}\phi_x) \prod_{x \in A} \mathrm{d}\phi_x,
\end{align}
where $V:\mathbb{R} \to (-\infty,+\infty]$ is a suitably regular function which we call the \textbf{interaction potential}, and where $Z_A(\varphi)$ is the normalising constant, or \textbf{partition function}, given by the integral
\begin{align} \label{eq:pf}
Z_A(\varphi) = Z_A(\varphi_x : x \in B) := \int_{\mathbb{R}^C} \exp \left\{ - \sum_{ \langle x,y \rangle \in C^*} V(\phi_y - \phi_x) \right\} \prod_{x \in B} \delta_{\varphi_x}(\mathrm{d}\phi_x) \prod_{x \in A} \mathrm{d}\phi_x.
\end{align}
The measure \eqref{eq:rs} on the space of height functions $\mathbb{R}^C$ is a Gibbs probability measure. The product of Dirac measures $\prod_{x \in B} \delta_{\varphi_x}(\mathrm{d}\phi_x)$ fixes the boundary values to coincide with $\varphi_x$ for $x \in B$, whereas $\prod_{x \in A} \mathrm{d}\phi_x$ 
denotes Lebesgue measure on the unconstrained interior heights. Consequently, the measure is supported on all surfaces that satisfy the prescribed boundary conditions, with probabilities determined by their interaction energy.

There is a vast literature on random surfaces \cite{funaki, sheffield, DGI}.
Most of this literature operates under the rather strong condition that the interaction potential $V$ satisfies $0 < c_1 \leq V''(x) \leq c_2$. We will instead be interested in an extended-valued interaction potential $V$ that fails even to be continuous.
In this direction, for any subset $C$ of $\mathbb{Z}^2$, we define a \textbf{Gelfand--Tsetlin surface} on $C$ to be a function $\phi:C \to \mathbb{R}$ that is increasing on $C$, i.e.
\begin{align} \label{eq:global}
\phi_x \leq \phi_y, \qquad \text{whenever $x,y \in C$ with $x \leq y$},
\end{align}
where, for $x = (x_1,x_2)$ and $y = (y_1,y_2)$ in $\mathbb{Z}^2$, we write $x \leq y$ if $x_1 \leq y_1$ and $x_2 \leq y_2$. 

If $C$ is \textbf{path-comparable}, meaning that for every pair $x \leq y \in C$ there is a path $x = z_0 \leq z_1 \leq \ldots \leq z_\ell = y$ such that each $\langle z_{i-1},z_i \rangle \in C^*$, then \eqref{eq:global} is equivalent to the statement
\begin{align} \label{eq:local}
\phi_x \leq \phi_y \qquad \text{whenever $\langle x , y \rangle \in C^*$}.
\end{align}
For a general set $C$ that is not path-comparable, \eqref{eq:local} is weaker than \eqref{eq:global}. For the majority of the article, we will have a path-comparable set $C$ in mind. 

Given $\phi:C \to \mathbb{R}$, define
\begin{align} \label{eq:globalind}
\mathrm{1}_{\mathrm{GT}}^C(\phi) := \prod_{x \leq y \in C} \mathrm{1}\{ \phi_x \leq \phi_y \}.
\end{align}
If $C = A \sqcup B$ is finite, we define a \textbf{random Gelfand--Tsetlin surface} with boundary conditions $\varphi:B \to \mathbb{R}$ to be a random function $\phi:C \to \mathbb{R}$ chosen uniformly among the set of increasing functions satisfying the boundary conditions imposed by $\varphi$. In other words,
this is a function distributed according to a probability measure of the form
\begin{align} \label{eq:jok1}
\frac{1}{Z_A(\varphi)} \mathrm{1}_{\mathrm{GT}}^C(\phi) \prod_{x \in B} \delta_{\varphi_x}(\mathrm{d}\phi_x)\prod_{x \in A} \mathrm{d}\phi_x,
\end{align}
provided the Gelfand--Tsetlin partition function
\begin{align} \label{eq:jok2}
Z_A(\varphi) := \int_{\mathbb{R}^C}  \mathrm{1}_{\mathrm{GT}}^C(\phi) \prod_{x \in B} \delta_{\varphi_x}(\mathrm{d}\phi_x)\prod_{x \in A} \mathrm{d}\phi_x
\end{align}
takes a value in $(0,\infty)$. 

If $C$ is path-comparable, the equations \eqref{eq:jok1} and \eqref{eq:jok2} coincide precisely with \eqref{eq:rs} and \eqref{eq:pf} in the special case where $V$ is the hard-core interaction potential
\begin{align} \label{eq:hc}
V_{\mathrm{GT}}(\phi) := 
\begin{cases}
0 \qquad &: \, \phi \geq 0,\\
+\infty \qquad &: \phi < 0.
\end{cases}
\end{align}
Thus random Gelfand--Tsetlin surfaces (up to the issue of path-comparability) may be considered as a special case of random surfaces.
\medskip

To see the connection with Gelfand--Tsetlin patterns, define
\begin{align*}
\rt_n := \{ (x_1,x_2) \in \mathbb{Z}^2 : 0 \leq x_2 \leq x_1\leq n \}.
\end{align*}
Given a Gelfand--Tsetlin pattern $(t_{k,j})_{0 \leq j \leq k \leq n}$ with bottom row $t_{n,j} = s_j$, we define a function $\phi:\rt_n \to \mathbb{R}$ by setting
\begin{align} \label{eq:rec1} 
\phi_{x_1,x_2} := t_{n-x_1+x_2,x_2}, \qquad (x_1,x_2) \in \rt_n.
\end{align}
Note that $\phi_{j,j} = t_{n,j} = s_j$ for $0 \leq j \leq n$. 
The interlacing inequalities \eqref{eq:GTp} are equivalent to $\phi_{x_1,x_2}$ being an increasing surface, i.e., an increasing function of $x_1$ and $x_2$.

If the Gelfand--Tsetlin pattern is chosen uniformly from the set of Gelfand--Tsetlin patterns with bottom row $s_0,\ldots,s_n$, then the associated random function $\phi$ is a Gelfand--Tsetlin random surface on $
C = \rt_n$
with boundary conditions $\{ \varphi_{j,j} = s_{j} : 0 \leq j \leq n\}$ on the diagonal $
B = \{ (j,j) : 0 \leq j \leq n \}$.
In this triangular case the associated partition function measures the volume of the set of Gelfand--Tsetlin surfaces on $A := C \setminus B$ whose extension to all of $C$ with boundary values $\{ \varphi_x : x \in B \}$ is still a Gelfand--Tsetlin surface. The formula for this partition function in this case is given by \eqref{eq:wdf}. 

In summary, our central task of describing the macroscopic behaviour of uniform random Gelfand--Tsetlin patterns can be reformulated as describing the macroscopic behaviour of uniformly chosen increasing functions $\phi:\rt_n \to \mathbb{R}$ with boundary conditions on the diagonal $B$. 

\subsubsection{Compressions as surfaces}

Recall Definition \ref{df:compression}. Inspired by the previous section, which related Gelfand--Tsetlin patterns to increasing functions on triangular subsets of $\mathbb{Z}^2$, here we relate compressions to increasing functions on triangular subsets of $\mathbb{R}^2$. Define 
\begin{align*}
\rt := \{ (s,t) \in \mathbb{R}^2 : 0 \leq t \leq s \leq 1 \}.
\end{align*}

\begin{df}[Increasing surface]
An increasing surface $\psi:\rt \to \mathbb{R}$ is an upper-semicontinuous real-valued function that is (weakly) increasing in both variables.
\end{df}

We make no requirement that an increasing surface is continuous.

It transpires that increasing surfaces are essentially equivalent to compressions, and can be constructed from compressions using their mixed quantile functions. The quantile functions $\{ Q(r,\tau) : 0 \leq r \leq \tau \leq 1\}$ associated with a compression $\{ F(x,\tau) : x \in \mathbb{R}, \tau \in [0,1] \}$ are the unique upper semicontinuous functions that are increasing in $r$ and satisfy 
\begin{align*}
\int_{-\infty}^\infty f(x) F(\mathrm{d}x,\tau) = \int_0^\tau f(Q(r,\tau))\mathrm{d}r
\end{align*}
for all bounded and measurable $f:\mathbb{R} \to \mathbb{R}$. 

\begin{proposition} \label{prop:bijection}
There is a one-to-one correspondence between compressions $F$ and increasing surfaces $\psi$ defined by setting
\begin{align*}
\psi(s,t) := Q(t,1-s+t), \qquad (s,t) \in \rt := \{(s,t) \in \mathbb{R}^2 : 0 \leq t \leq s \leq 1 \},
\end{align*}
where $Q$ is the mixed quantile function of $F$. 
\end{proposition}

In other words, the diagonal $\psi(1-\tau+r,r)$ encodes the quantile function of the measure associated with $F(\cdot,\tau)$. In particular, if $F$ is a compression of $\mu$, and $\psi$ is the associated increasing surface, then the diagonal $\{ \psi(s,s) : s \in [0,1]\}$ is the quantile function $\{Q_\mu(s): s \in [0,1] \}$ of $\mu$.

The continuum interlacing inequality \eqref{eq:continuum} implies that $Q(r,\tau+h)$ is decreasing in $h$, and that $Q(r+h,\tau+h)$ is increasing in $h$. Accordingly:

\begin{rem}[Interlacing inequalities]
The continuum interlacing inequalities \eqref{eq:continuum} for $F$ are equivalent to $\psi$ being increasing in both variables.
\end{rem}

The compression associated with a Gelfand--Tsetlin pattern (Definition \ref{df:GTcomp}) is more transparent when cast in terms of an increasing surface:

\begin{df}[The increasing surface associated with a Gelfand--Tsetlin pattern] \label{df:GTcomp2}
We define the increasing surface $\xi:\rt \to \mathbb{R}$ associated with a Gelfand--Tsetlin pattern $(t_{k,j})_{0 \leq j \leq k \leq n }$ by  
\begin{align*}
\xi(x_1/n,x_2/n) := t_{n-x_1+x_2,x_2}, \qquad 0 \leq x_2 \leq x_1 \leq n, \quad x_1,x_2 \text{ integers},
\end{align*}
and then linearly interpolating $\xi$ on triangles between these points.
\end{df}

By linearly interpolating on triangles, we mean that we linearly interpolate $\xi$ over all triangles of the form $ABC$ or $AB'C$, where
\begin{align*}
A = (k/n,j/n), \qquad B = ((k+1)/n,j/n), \qquad B' = (k/n,(j+1)/n), \qquad C = ((k+1)/n,(j+1)/n).
\end{align*}
One checks that for $0 \leq r \leq n$, $\{ \xi( (n-r+j)/n, j/n) : j = 0,\ldots, r \}$ corresponds to the $r^{\text{th}}$ row $\{ t_{r,j} : j = 0,\ldots,r \}$ of the Gelfand--Tsetlin pattern. 

\begin{figure}[ht]
    \centering
    \includegraphics[width=0.33\textwidth]{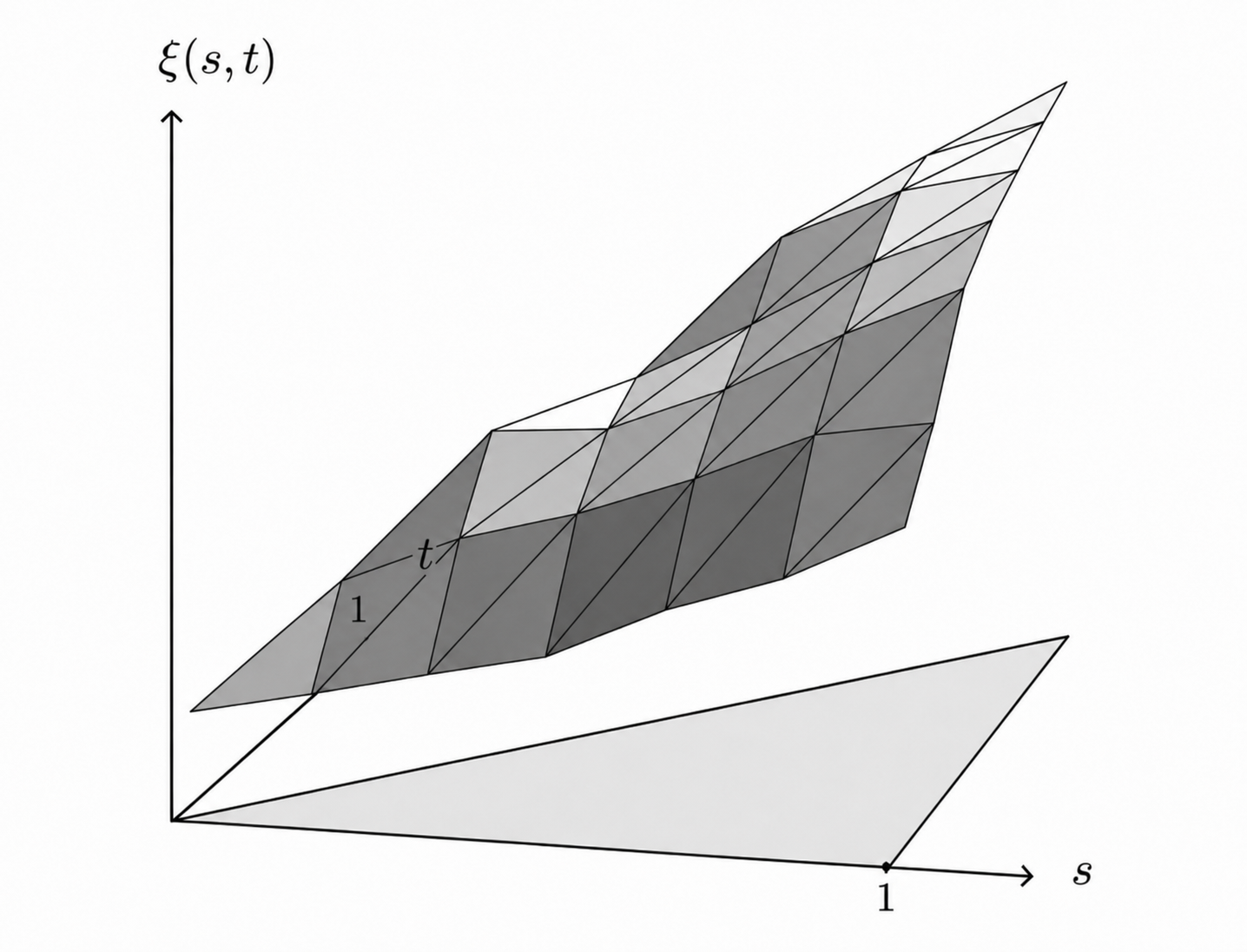}
     \caption{The increasing surface associated with the Gelfand--Tsetlin pattern in Figure \ref{fig:GT7} (after relabelling so that $(t_{k,j})_{ 0\leq j \leq k \leq 6}$).}
    \label{fig:gtsurface2}
\end{figure}

Definition \ref{df:GTcomp} for the compression associated with a Gelfand--Tsetlin pattern and Definition \ref{df:GTcomp2} for the increasing surface associated with a Gelfand--Tsetlin pattern are compatible with the bijection in Proposition \ref{prop:bijection}.

We also observe that the discrete interlacing inequalities \eqref{eq:GTp} for Gelfand--Tsetlin patterns read
\begin{align*}
\xi(x_1/n,x_2/n) \leq \xi((x_1+1)/n,x_2/n) \quad \text{and}  \quad \xi(x_1/n,x_2/n) \leq \xi( x_1/n, (x_2+1)/n),
\end{align*}
for suitable $0 \leq x_2 \leq x_1 \leq n$. In particular, when we linearly interpolate between these points we indeed obtain an increasing surface.

Finally, let us emphasise that an increasing surface $\psi:\rt \to \mathbb{R}$ need not be continuous, but by Lebesgue's monotone differentiation theorem, it is differentiable almost everywhere on $\rt$. Discontinuities in $x$ of $F(x,\tau)$ correspond to flat regions for $\psi$. 

%%%%%%%%%%%%%%%%%%%%%%%%%%%%%%%%%%%%%%%%%%%%%%%%%%%%%%%%%%%%%%%%%%%%%%%%%%%%%%%%%%%%%%%%%%%
\subsection{Theorem \ref{thm:A}: Compressions and entropy}
\label{sec:thmA}
%%%%%%%%%%%%%%%%%%%%%%%%%%%%%%%%%%%%%%%%%%%%%%%%%%%%%%%%%%%%%%%%%%%%%%%%%%%%%%%%%%%%%%%%%%%

In order to state our main results, we topologise the space of compressions. In this direction, we say that a compression $F$ of $\mu$ is integrable if $\mu(\mathrm{d}x) := F(\mathrm{d}x,1)$ has finite first moment.
We can metrize the collection of integrable compressions with an $L^1$-type metric as follows. Let $\psi,\psi'$ be increasing surfaces associated with integrable compressions $F,F'$. Then 
\begin{align} \label{eq:metric}
\mathrm{d}(F, F') := \int_0^1 \int_{-\infty}^\infty |F'(x,\tau) - F(x,\tau)| \mathrm{d}x\mathrm{d}\tau = \int_{\rt} |\psi'(s,t) - \psi(s,t)| \mathrm{d}s\mathrm{d}t =: \mathrm{d}(\psi,\psi'),
\end{align} 
where the fact that the two integrals coincide is a consequence of the well-known equality between $L^1$-Wasserstein distances and $L^1$-distances between quantile functions (see e.g., \cite[Proposition 2.17]{santa}).

% With the topology induced by the metric in \eqref{eq:metric}, one can show that every compression is a limit of Gelfand--Tsetlin patterns (Proposition \ref{prop:closure}). 

Recall that we defined the entropy $\mathcal{H}[F]$ of a compression in Definition \ref{df:compressionentropy} by the formula
\begin{align} \label{eq:chiagain}
\mathcal{H}[F] := \int_{-\infty}^\infty \int_0^1 (-\Sigma)(F_x,F_\tau) \mathrm{d}\tau \mathrm{d}x, \qquad \Sigma(F_x,F_\tau) := F_x \left\{  \log F_x - \log \sin \pi F_\tau - 1 + \log \pi \right\}.
\end{align}
If $\psi$ is the surface associated with a compression $F$, we will abuse notation and write $\mathcal{H}[\psi]$ for $\mathcal{H}[F]$. A change of variable calculation tells us that $\mathcal{H}[\psi]$ is given by the formula
\begin{align} \label{eq:entropysurface}
\mathcal{H}[\psi] = \int_{\rt} (-\sigma)(\nabla \psi) \mathrm{d}s\mathrm{d}t \qquad \sigma(\nabla\psi) := - \log( \psi_s + \psi_t) - \log \sin \left( \frac{ \pi \psi_s}{ \psi_s+\psi_t} \right) - 1 + \log \pi,
\end{align}
where $\psi_s$ and $\psi_t$ denote partial derivatives with respect to $s$ and $t$. 
% The definition of $\mathcal{H}[\psi]$ in \eqref{eq:entropysurface} is valid for all increasing surfaces $\psi:\rt \to \mathbb{R}$, regardless of any continuity requirements, with the convention that $\log 0 = -\infty$.

We will prove in Theorem \ref{thm:lowersemi} that the compression entropy is upper semicontinuous at every compression with bounded support.

% Note that if $\psi,\psi'$ are increasing surfaces, then so is $\lambda \psi + (1-\lambda) \psi'$. We also note that $\psi \mapsto \sigma(\nabla \psi)$ is convex in $\psi$. Accordingly, $\psi \mapsto \mathcal{H}[\psi]$ is a concave function of an increasing surface in the sense that
% \begin{align*}
% \mathcal{H}[ \lambda \psi + (1-\lambda) \psi' ] \geq \lambda \mathcal{H}[\psi] + (1-\lambda) \mathcal{H}[\psi'].
% \end{align*}
% In particular, any sufficiently regular stationary point of $\mathcal{H}[\psi]$ among increasing surfaces with fixed diagonal boundary data is necessarily a maximiser.

We turn to structural properties of maximisers of $\mathcal{H}[\cdot]$. The Euler--Lagrange equations for the functional in question read 
\begin{align} \label{eq:EL}
0
=
\frac{\partial}{\partial x} \frac{\partial \Sigma}{\partial F_x}
+
\frac{\partial}{\partial \tau} \frac{\partial \Sigma}{\partial F_\tau}
=
\frac{\partial}{\partial x} \left( \log F_x - \log \sin(\pi F_\tau) \right)
-
\frac{\partial}{\partial \tau}\left( \pi F_x \frac{\cos(\pi F_\tau)}{\sin(\pi F_\tau)} \right).
\end{align}
While \eqref{eq:EL} looks untenable, we can interpret it as the \emph{real part} of a complex Burgers equation. In this direction, let $\mathbb{C}^+$ and $\mathbb{C}^-$ denote the strict upper and lower half-planes in $\mathbb{C}$. We define the \textbf{joint Cauchy transform} of a compression $\{ F(x,\tau) : x \in \mathbb{R}, \tau \in [0,1]\}$ to be the function $\{ G(z,\tau) : z \in \mathbb{C}^+, \tau \in [0,1]\}$ taking values in $\mathbb{C}^-$ and defined by
\begin{align*}
G(z,\tau) := \int_{-\infty}^\infty \frac{1}{z-x}\, F(\mathrm{d}x,\tau).
\end{align*}

Recall our extrinsic Definition \ref{def:freecompression} of the free compression measures. In Section \ref{sec:Rtransform} we give a more intrinsic definition of the free compression measures using $R$-transforms, and show that the associated Cauchy transforms satisfy the Burgers equation 
\begin{align*}
\partial_\tau G + \frac{\partial_z G}{G} = 0 \qquad z \in \mathbb{C}^+, \tau \in (0,1).
\end{align*}
% We will verify directly there using the Burgers equation that the free compression measures do indeed give rise to a compression in the sense of Definition \ref{df:compression}.

Given a probability measure $\mu$, let
\begin{align*}
\mathcal{A}_\mu := \{ F : \text{$F$ is a compression of $\mu$} \}.
\end{align*}
Observe that $\mathcal{A}_\mu$ is a convex set: if $F$ and $\tilde{F}$ are compressions of $\mu$, so is $\lambda F + (1-\lambda)\tilde{F}$. We also observe that $\Sigma(F_x,F_\tau)$ is a convex function of $(F_x,F_\tau)$, and accordingly, $\mathcal{H}:\mathcal{A}_\mu \to [-\infty,\infty)$ is concave in that
\begin{align*}
\mathcal{H}[\lambda F + (1-\lambda)\tilde{F} ] \geq \lambda  \mathcal{H}[F ] + (1-\lambda) \mathcal{H} [ \tilde{F}].
\end{align*}
Since $\Sigma(F_x,F_\tau)$ is in fact strictly convex, any maximiser of $\mathcal{H}$ in $\mathcal{A}_\mu$ is necessarily unique. 
We are now ready to state our first main result:

\begin{thmalpha}(Free entropy and free compression) \label{thm:A}
Let $\mu$ be a compactly supported absolutely continuous probability measure with $\chi[\mu] > - \infty$.
\begin{enumerate}
\item Free compression $F^*$ satisfies the Euler--Lagrange equations for $\mathcal{H}[F]$, and by the strict concavity of $\mathcal{H}[ \cdot]$ is the unique compression of $\mu$ maximising $\mathcal{H}[F]$.  
\item The compression entropy and the free entropy of a probability measure are related by 
\begin{align} \label{eq:entropy}
\mathcal{H}[F^*] = \sup_{ F \in \mathcal{A}_\mu } \mathcal{H}[F] = \chi[\mu].
\end{align}
\end{enumerate}
\end{thmalpha}

The first part of Theorem \ref{thm:A} is in essence a recapitulation of Shlyakhtenko and Tao's observation that free compression satisfies certain Euler--Lagrange equations. We will see in fact that it manifests algebraically through the relation
\begin{align*}
\mathrm{Re}\left( \partial_\tau G + \frac{\partial_x G}{G} \right)  = \frac{\partial}{\partial x} \frac{\partial \Sigma}{\partial F_x}
+
\frac{\partial}{\partial \tau} \frac{\partial \Sigma}{\partial F_\tau},
\end{align*}
identifying the Euler--Lagrange equations for compression entropy as coinciding with the real part of the Burgers equation for the joint Cauchy transform. 

Our approach to proving both parts of Theorem \ref{thm:A} is reminiscent of work by Kenyon and Okounkov \cite{KO} (see also Kenyon and Prause \cite{KP,KP2}), which relates limit shapes minimising surface tension to solutions of Burgers equations.

%%%%%%%%%%%%%%%%%%%%%%%%%%%%%%%%%%%%%%%%%%%%%%%%%%%%%%%%%%%%%%
\subsection{Theorem \ref{thm:B}: The surface tension of the bead model} \label{sec:thmB}
%%%%%%%%%%%%%%%%%%%%%%%%%%%%%%%%%%%%%%%%%%%%%%%%%%%%%%%%%%%%%%

In this section we state our next main result, Theorem \ref{thm:B}.
Recall from Section \ref{sec:surfaces} the definition of a random surface with interaction potential $V$. Given an interaction potential $V$, the associated \textbf{surface tension} $\sigma_V:\mathbb{R}^2 \to [-\infty,+\infty]$ is an asymptotic measure of the cost of constraining a random surface to have gradient $(u_1,u_2)$. More precisely, let $\Box_n := \{0,1,\ldots,n\}^2$, and let $\partial \Box_n$ and $\interior{\Box}_n$ denote its boundary and interior respectively. 
The square partition function with boundary slope $(u_1,u_2)$ is given by
\begin{align*}
S_n(u_1,u_2) := \int_{\mathbb{R}^{\Box_n} } \exp \Big\{ - \sum_{ \langle x,y \rangle \in \Box_n^*} V(\phi_y - \phi_x) \Big\}  \prod_{x \in \partial \Box_n} \delta_{u_1 x_1 + u_2 x_2 }(\mathrm{d}\phi_x) \prod_{x \in \interior{\Box}_n} \mathrm{d}\phi_x.
\end{align*}
Note that with the notation of \eqref{eq:pf}, we have 
\begin{align} \label{eq:squarebox}
S_n(u_1,u_2) := Z_{\interior{\Box}_n}( \varphi^u_x : x \in \partial \Box_n ), \qquad \varphi^u_x = u_1 x_1 + u_2 x_2.
\end{align}
In Funaki and Spohn \cite{FS}, it is shown under the condition $0 < c_1 \leq V''(x) \leq c_2$ that the limit 
\begin{align} \label{eq:st0}
\sigma_V(u_1,u_2) := - \lim_{n \to \infty} \frac{1}{n^2} \log S_n(u_1,u_2),
\end{align}
exists and takes values in $(-\infty,+\infty)$. The value $\sigma_V(u_1,u_2)$ is known as the surface tension at gradient $(u_1,u_2)$ of the interaction potential $V$. 

There is an analogous definition of surface tension for stochastic interface models in $\mathbb{Z}^d$ involving a limit of partition functions taken over $\{0,1,\ldots,n\}^d$. 
The surface tension of one-dimensional interface models can be calculated explicitly for essentially any interaction potential $V$ \cite[Section 5.5]{funaki}. To our knowledge, however, in the existing literature on Ginzburg--Landau $\nabla\phi$-interface models in $d \geq 2$ dimensions, there are no explicit formulas for the surface tension $\sigma_V$ with the exception of the case where $V$ is the Gaussian interaction potential $V(u) = \frac{a}{2} u^2$ (the surface tension is also a quadratic function; see \cite[Proposition 5.2]{funaki} for a formula). 
However, in our second main result, we are able to offer an explicit formula for the two-dimensional surface tension associated with the hard-core interaction potential $V_{\mathrm{GT}}$ in \eqref{eq:hc}: 
% This result is a rigorous resolution of Conjecture \ref{conj:ST}:

\begin{thmalpha}(Surface tension of the bead model) \label{thm:B}
The surface tension $\sigma(u_1,u_2)$ associated with the hard-core interaction potential given in \eqref{eq:hc} exists in $(-\infty,+\infty]$ and is given by
\begin{align} \label{eq:beadtension}
\sigma(u_1,u_2) =
\begin{cases}
- \log(u_1 + u_2) - \log \sin \left( \pi \frac{u_1}{u_1 + u_2} \right) + \log \pi - 1  &: u_1 > 0 \text{ and } u_2 > 0,\\
+ \infty  &: \text{otherwise}. 
\end{cases}
\end{align}
\end{thmalpha}

Up to the additive normalisation constant $\log\pi-1$, this surface tension is the negative of the Lagrangian density in \eqref{eq:lagrangian}. 
% Let us explain how Theorem \ref{thm:B} is related to the conjectured entropy of the bead process described in Section \ref{sec:21}.

As noted above, Theorem \ref{thm:B} gives a rigorous formulation of related formulas appearing in numerous places, such as \cite{ST, sun, JO, gordenko}. Our proof of Theorem \ref{thm:B} relies on explicit determinantal partition function and correlation function formulas from recent work of the first author \cite{johnston}.

%%%%%%%%%%%%%%%%%%%%%%%%%%%%%%%%%%%%%%%%%%%%%%%%%%%%%%%%%%%%%%
\subsection{Theorem \ref{thm:C}: A large deviation principle for Gelfand--Tsetlin patterns}
\label{sec:thmC}
%%%%%%%%%%%%%%%%%%%%%%%%%%%%%%%%%%%%%%%%%%%%%%%%%%%%%%%%%%%%%%

In this section we state our large deviation principle for Gelfand--Tsetlin patterns.
First we recall the basic tenets of large deviation theory \cite{DZ}. 
Let $\mathcal X$ be a Hausdorff topological space. A function $I : \mathcal{X} \to [0,\infty]$ is called a rate function if it is lower semicontinuous, i.e., if all the sublevel sets are closed. 
A sequence of $\mathcal{X}$-valued random variables $(X_n)_{n \geq 1}$ is said to satisfy a \textbf{large deviation principle} (LDP) with speed $(a_n)_{n \geq 1}$, where $a_n \to \infty$, and rate function $I : \mathcal{X} \to [0,\infty]$, if for every closed set $F \subseteq \mathcal{X}$ and every open set $G \subseteq \mathcal{X}$ one has
\begin{align*}
\limsup_{n \to \infty} \frac{1}{a_n} \log \mathbf{P}(X_n \in F) &\leq - \inf_{x \in F} I(x), \\
\liminf_{n \to \infty} \frac{1}{a_n} \log \mathbf{P}(X_n \in G) &\geq - \inf_{x \in G} I(x).
\end{align*}
Informally, one may write $\mathbf{P}( X_n \approx x) \approx e^{ - a_n I(x)}$ as $n \to \infty$.
\begin{thmalpha} \label{thm:C}
Let $\mu$ be an absolutely continuous compactly supported probability measure on $\mathbb{R}$ satisfying $\chi[\mu] > - \infty$, and let $Q_\mu:[0,1] \to \mathbb{R}$ be its quantile function. Let $(t^n)_{n \geq 1}$ be a sequence of Gelfand--Tsetlin patterns such that for each $n\in\mathbb N$, $t^n := (t^n_{k,j})_{0 \leq j \leq k \leq n}$ is a uniform random Gelfand--Tsetlin pattern with bottom row
\begin{align*}
t^n_{n,j} := Q_\mu(j/n), \qquad 0 \leq j \leq n.
\end{align*}
Let $\xi_n:\rt \to \mathbb{R}$ be the increasing surface associated with the Gelfand--Tsetlin pattern $t^n$. Then the sequence $(\xi_n)_{n \geq 1}$ satisfies a large deviation principle on the metric space of integrable increasing surfaces with speed $n^2$ and rate function given by
\begin{align*}
I_\mu[\psi] :=
\begin{cases}
- \mathcal{H}[\psi] + \chi[\mu] \quad &\text{if $\psi(s,s) = Q_\mu(s)$ for all $s \in [0,1]$,}\\
+\infty \quad &\text{otherwise}.
\end{cases}
\end{align*}
Here $\mathcal{H}[\psi]$ is as in \eqref{eq:entropysurface}, while $\chi[\mu]$ denotes Voiculescu's free entropy as in \eqref{eq:freeent}.
\end{thmalpha}

An equivalent statement in the setting of Theorem \ref{thm:C} is that if $F^n$ is the compression associated with $t^n$, then the sequence $(F^n)_{n \geq 1}$ of compressions satisfies a large deviation principle with speed $n^2$ and rate function $I_\mu[F] = - \mathcal{H}[F]+\chi[\mu]$ for compressions $F$ of $\mu$, thereby making rigorous sense of \eqref{eq:X}.

It follows from Theorem \ref{thm:A} that $\inf_\psi I_\mu[\psi] = 0$, and that the infimum is attained by the surface $\psi^*$ associated with the free compression of $\mu$. We also note from Baryshnikov's result, Theorem \ref{thm:Aar}, that a uniform Gelfand--Tsetlin pattern whose bottom row approximates a probability measure $\mu$ is equivalent to the eigenvalue process of a random matrix whose empirical spectrum approximates $\mu$. Thus, Theorem \ref{thm:C} describes the eigenvalue processes of large random matrices. In other words, as a corollary of Theorem \ref{thm:C}, we have the following.

\begin{cor} \label{cor:total}
In the setting of Theorem \ref{thm:C}, the sequence $(\xi_n)_{n \geq 1}$ of increasing surfaces converges in probability in the metric \eqref{eq:metric} to the increasing surface $\psi^*:\rt \to \mathbb{R}$ associated with the free compression $F^*$ of $\mu$.
\end{cor}

Note in particular that, under the hypotheses of Theorem \ref{thm:C}, Corollary \ref{cor:total} implies Theorem \ref{thm:compression}.

\medskip
Thus, Theorems \ref{thm:A}, \ref{thm:B}, and \ref{thm:C} provide a unifying picture explaining free compression (Theorem \ref{thm:compression}), free entropy \eqref{eq:freeent}, and Shlyakhtenko and Tao's Euler--Lagrange equations for free compression. Theorem \ref{thm:B} says that, as conjectured by Shlyakhtenko and Tao, the Lagrangian in \eqref{eq:lagrangian} arises naturally as a free energy in the statistical mechanics of Gelfand--Tsetlin random surfaces. Theorems \ref{thm:A} and \ref{thm:C} together describe how Baryshnikov's result, Theorem \ref{thm:Aar}, gives rise macroscopically to free compression and Voiculescu's free entropy.

That completes the statement of our main results.
%%%%%%%%%%%%%%%%%%%%%%%%%%%%%%%%%%%%%%%%%%%%%%
%%%%%%%%%%%%%%%%%%%%%%%%%%%%%%%%%%%%%%%%%%%%%%
%%%%%%%%%%%%%%%%%%%%%%%%%%%%%%%%%%%%%%%%%%%%%%
\section{Further discussion and related work} \label{sec:discussion}
%%%%%%%%%%%%%%%%%%%%%%%%%%%%%%%%%%%%%%%%%%%%%%
%%%%%%%%%%%%%%%%%%%%%%%%%%%%%%%%%%%%%%%%%%%%%%
%%%%%%%%%%%%%%%%%%%%%%%%%%%%%%%%%%%%%%%%%%%%%%

In this section we collect together further consequences of our main results, tools developed in our proofs that might be of broader application, further discussion of related literature, and conjectures and outlines for future research directions.

\subsection{Further discussion of related literature}

\subsubsection{Large random Gelfand--Tsetlin patterns}

The macroscopic behaviour of random Gelfand--Tsetlin patterns has previously been studied through both lozenge tilings and random matrix theory. For discrete patterns, Petrov \cite{petrov} and Duse and Metcalfe \cite{DM1,DM2} analysed uniformly random Gelfand--Tsetlin patterns with general fixed boundary data, describing their limit shapes, liquid regions and local asymptotics. In the continuous setting, Metcalfe \cite{metcalfe} studied uniform Gelfand--Tsetlin patterns with a fixed top row and related their macroscopic row distributions to free compression. Collins and Metcalfe \cite{CM} subsequently obtained estimates controlling deviations of extremal particles from this limiting shape.

\subsubsection{Random surfaces and their large deviations}
As mentioned above, much of the classical analytic literature on continuous
$\nabla\phi$-interface models assumes that the interaction potential $V$ is
uniformly strictly convex. In particular, Deuschel et al.\ \cite{DGI}
established a large deviation principle for the rescaled height functions in
the $L^2$ topology, under the condition $0<c_1\leq V''\leq c_2$; their rate
function is given by an integral of the surface tension over the weak gradient
and is finite only on the appropriate Sobolev space. See also \cite{FN} for
the dynamical setting.

Our hard-core interaction is highly singular and lies outside this framework.
Correspondingly, the rate functional
\begin{align*}
I_\mu[\psi] := - \mathcal{H}[\psi] + \chi[\mu], \qquad
\text{with} \quad
\mathcal{H}[\psi] = - \int_{\rt} \sigma(\nabla \psi)
\mathrm{d}s\mathrm{d}t,
\end{align*}
can be finite for monotone surfaces of much lower regularity, including some
surfaces with jump discontinuities. Indeed, by Lebesgue's theorem on the
almost-everywhere differentiability of monotone functions, applied to
one-dimensional coordinate sections, an increasing function
$\psi:\rt\to\mathbb{R}$ has first-order partial derivatives almost everywhere.
These almost-everywhere derivatives need not form a weak gradient, since the
distributional derivatives of $\psi$ may also contain singular parts.
Nevertheless, they suffice to define the integral above, giving our rate
functional a substantially larger effective domain than in the uniformly
convex setting.

Speaking more broadly about the literature on random surfaces, there is a large body of work, much of which is captured in research monographs by Sheffield \cite{sheffield} and Funaki \cite{funaki}. In terms of specific examples, the most commonly studied random surfaces are those with Gaussian interaction potentials $V(u) = c u^2$. The corresponding model is the discrete Gaussian free field \cite{berestycki}. The discrete Gaussian free field has the advantageous property that the marginal law of any given coordinate has a Gaussian distribution. Under suitable scaling, it may be shown that scaling limits of the discrete Gaussian free field give rise to the continuum Gaussian free field \cite{berestycki}. Quartic interfaces with interaction potentials of the form $V(u) = (u^2-a)^2$, known as the symmetric double-well potential, are occasionally studied as a continuous-height analogue of the Ising model \cite{CGW} and are of considerable interest in quantum mechanics and quantum field theory. Another case that appears occasionally in the literature is the Toda lattice, which in our language is the random surface with the exponential interaction potential $V(u) = e^u$ \cite{JO, oconnell, OSZ}. Various other models fall into analogous frameworks: height functions of domino and lozenge tilings provide important integer-valued examples of random surfaces \cite{gorin}, and the Ising model can be regarded as a random surface taking values in $\{-1,1\}$. Let us highlight in particular the famous work of Cohn, Kenyon and Propp \cite{CKP}, who established a variational principle for domino tilings.

\subsubsection{Large deviations for random matrix eigenvalues}

Ben Arous and Guionnet \cite{BAG} established a large deviation principle, with speed $n^2$, for the empirical spectral measure of the GUE, with rate function given by the logarithmic Coulomb energy.
We expect that by adapting the tools of the present article, it should be possible to show that the compression arising from the eigenvalue process associated with a large GUE random matrix satisfies a large deviation principle with rate function of the form
\begin{align*}
I_{\mathrm{GUE}}[F] = -\mathcal{H}[F] + \frac{1}{2}\int_{-\infty}^\infty x^2 F(\mathrm{d}x,1) + \mathrm{const}.
\end{align*}
This theory has since been extended to broad classes of invariant matrix ensembles and related orbital models; see, for example, \cite{AGZ,BGH}. These results concern the empirical spectrum of a single matrix, whereas our large deviation principle describes the entire surface formed by the spectra of all principal minors, conditional on the limiting spectrum of the full matrix. In joint work of Octavio Arizmendi and the first author \cite{AJ}, the authors applied a large deviation principle to the quadrature formula of Marcus, Spielman and Srivastava \cite{MSS} to formulate the operations of free probability in terms of entropic optimal transport.

\subsubsection{Zeroes of polynomials}

The zeroes of the successive derivatives of a real-rooted polynomial form a Gelfand--Tsetlin pattern by Rolle's theorem. Kabluchko \cite{kabluchko} proved that taking a fixed number of derivatives does not change the limiting empirical distribution of the zeroes, while Steinerberger \cite{steinerberger,steinerbergerfree} and Hoskins and Kabluchko \cite{HK} studied the regime in which the number of derivatives is proportional to the degree. In this regime, the limiting distribution evolves according to free compression. Thus, the deterministic Gelfand--Tsetlin pattern generated by repeated differentiation has the same row-wise macroscopic limit as a uniform random Gelfand--Tsetlin pattern with the same boundary data.

%%%%%%%%%%%%%%%%%%%%%%%%%%%%%%%%%%%%%%%%%%%%%%%%%%%%
\subsection{Convex geometry tools} \label{sec:convex}
%%%%%%%%%%%%%%%%%%%%%%%%%%%%%%%%%%%%%%%%%%%%%%%%%%%%

While our hard-core interaction potential is not twice differentiable,
and hence the associated tools of the $\mathcal{C}^2$ theory for random surfaces (such as the Helffer-Sj\"ostrand representation \cite{HS}) are not available, $V_{\mathrm{GT}}$ is convex, 
and much of our proof of Theorem \ref{thm:B} relies on tools from
convex geometry. We were inspired in part by the work of
Narayanan, Sheffield, and Tao \cite{NST} on augmented hives, which uses
high-dimensional concentration results for log-concave measures to prove
concentration of augmented hives. 
As byproducts of our proof, we establish several facts about random surfaces
with convex interaction potentials, which need not be
differentiable, continuous or even finite. These include a Pr\'ekopa--Leindler inequality
for partition functions, the log-concavity of marginal densities, and a
derivation of the convexity of surface tension. We are optimistic that this
robust approach may be of broader interest to researchers working on random
surfaces.

To explain these ideas, recall that a density function
$f:\mathbb{R}^k \to [0,\infty)$ of a random vector is said to be
\textbf{log-concave} if it takes the form
$f(s_1,\ldots,s_k)=e^{-W(s)}$ for some extended-valued convex
function $W:\mathbb{R}^k\to(-\infty,\infty]$. Log-concave measures exhibit a
range of useful concentration and stability properties.

Our primary convex-geometric tools are two consequences of the
Pr\'ekopa--Leindler inequality for random surfaces. The following proposition,
which we prove in Section \ref{sec:randomsurfaces}, follows fairly quickly from
standard forms of the Pr\'ekopa--Leindler inequality and plays a central role
throughout the article.
\begin{proposition}[The Pr\'ekopa--Leindler inequalities for random surfaces] \label{prop:pl}
Let $C = A \sqcup B$ be finite. Let $\phi = (\phi_x)_{x \in C}$ be a random surface with probability density function taking the form \eqref{eq:rs} for a convex interaction potential $V$. Then we have the following:
\begin{enumerate}
\item Let $\varphi, \tilde{\varphi}:B \to \mathbb{R}$ and $\lambda \in [0,1]$. Then with $Z_A(\varphi)$ as in \eqref{eq:pf} we have
\begin{align} \label{eq:pfpl}
Z_A(\lambda \varphi + (1-\lambda)\tilde{\varphi}) \geq Z_A(\varphi)^\lambda Z_A(\tilde{\varphi})^{1-\lambda}.
\end{align}
\item Let $D$ be any nonempty subset of $A$. Then the marginal density in $\mathbb{R}^{|D|}$ of the random variables $\{\phi_x : x \in D\}$ is log-concave.
\end{enumerate}
\end{proposition}

The first part of Proposition \ref{prop:pl} not only plays an integral role in our proof of Theorem \ref{thm:C}, but also has ramifications for surface tension (\eqref{eq:st0}). Funaki and Spohn show in \cite{FS} that the surface tension associated with an interaction potential $V$ satisfying the strong condition $0 < c_1 \leq V''(x) \leq c_2$ is convex. The following generalisation of this result (operating under the much weaker condition of convexity alone) is an immediate corollary of Proposition \ref{prop:pl} and the definition of surface tension \eqref{eq:st0}:

\begin{cor} \label{cor:STconvex}
Suppose the surface tension $\sigma_V:U \to (-\infty,+\infty]$ (defined in \eqref{eq:st0}) exists for all $(u_1,u_2)$ in a convex subset $U$ of $\mathbb{R}^2$. Then $\sigma_V$ is a convex function on $U$. 
\end{cor}
\begin{proof}
Let $u = (u_1,u_2)$ and $\tilde{u} = (\tilde{u}_1,\tilde{u}_2)$ be slopes in $U$.
Set $\varphi = \varphi^u$ and $\tilde{\varphi} = \varphi^{\tilde{u}}$. Then \eqref{eq:squarebox}, \eqref{eq:st0} and \eqref{eq:pfpl} imply that
\begin{align*}
\sigma_V(\lambda u_1 + (1-\lambda)\tilde{u}_1 , \lambda u_2 + (1-\lambda) \tilde{u}_2)
\leq
\lambda \sigma_V(u_1,u_2) + (1-\lambda)\sigma_V(\tilde{u}_1,\tilde{u}_2),
\end{align*}
as required.
\end{proof}

The second part of Proposition \ref{prop:pl} has been observed in several places \cite{sheffield, MP}. 
It furnishes a key step in obtaining Theorem \ref{thm:B} from the volume formula given in \cite{johnston}. We combine this marginal log-concavity with a concentration result of Lov\'asz--Vempala \cite{LV} to show that sub-configurations within larger toric bead configurations are reasonably likely to occupy their 
expected positions. 
% The Lov\'asz--Vempala result in question states that if $f:\mathbb{R}^k \to [0,\infty)$ is a log-concave density function for a random vector with mean vector $(\mu_1,\ldots,\mu_k)$, we have the absolute bound
% \begin{align} \label{eq:LV}
% f(\mu_1,\ldots,\mu_k) \leq \sup_{s_1,\ldots,s_k} f(s_1,\ldots,s_k) \leq 2^{2k} f(\mu_1,\ldots,\mu_k).
% \end{align}
% In other words, the density at the mean of a log-concave random vector is within a universal exponential factor of its maximum.

\medskip
Recall from above that the bulk behaviour of large uniform Gelfand--Tsetlin patterns was shown by Metcalfe \cite{metcalfe} to converge to Boutillier's bead process \cite{boutillier}. It was shown by Boutillier that regardless of the underlying tilt $\gamma := u_1/(u_1+u_2)$, the marginal law of the point process on any string (say the string at level $\ell \in \mathbb{Z}$) is given by the sine process, whose law is characterised by the determinantal correlation formula:
\begin{align} \label{eq:sine}
\mathbf{P}^\gamma \left( \text{Locations $(\mathrm{d}x_1,\ell),\ldots,(\mathrm{d}x_m,\ell)$ contain beads} \right) = \det_{i,j=1}^m \left( \frac{ \sin(\pi (x_j-x_i) ) }{ \pi (x_j - x_i) } \right) \mathrm{d}x_1 \cdots \mathrm{d}x_m.
\end{align}

From the alternative perspective of Proposition \ref{prop:pl}, the marginal density of any sub-collection of points of a uniform Gelfand--Tsetlin surface is log-concave. Combining Boutillier's observation \eqref{eq:sine}, Metcalfe's proof of the bulk convergence \cite{metcalfe}, and Proposition \ref{prop:pl} is therefore strong evidence that the joint distribution of the gaps between successive points in a sine process is log-concave.
% \begin{conj}[The sine process has log-concave marginals] \label{conj:sine}
% Take a determinantal point process on $\mathbb{R}$ with correlation kernel
% \begin{align*}
% K(x,y) = \frac{\sin(\pi(y-x))}{\pi(y-x)},
% \end{align*}
% with the usual diagonal convention $K(x,x)=1$, and label the points of the process $\cdots < Z_{-1} < Z_0 < Z_1 < Z_2 < \cdots$ so that $Z_{-1} < 0 \leq Z_0$.  Then the marginal density on $\mathbb{R}^k$ of any finite sub-collection $(Z_{i_1},\ldots,Z_{i_k})$ of points is log-concave.
% \end{conj}
This idea is pursued in further detail in a recent preprint \cite{johnstonGUE} by the first author, where the log-concavity of gaps between GUE eigenvalues is established, and is used to sharpen bounds obtained in Tao \cite{taoGUE}.

\subsection{A free probability interpretation of compressions}

A compression as in Definition \ref{df:compression} has the following operator-algebraic interpretation.

\begin{rem}[Operator-algebraic interpretation of compression]
Let $(\mathcal{A},\phi)$ be a tracial noncommutative probability space. We say that an element $p$ of $\mathcal{A}$ is a projection if $p^2=p=p^*$. The trace $\phi(p)$ of a projection takes values in $[0,1]$. We say that $\{p_\tau:\tau\in[0,1]\}$ is a filtered collection of projections if $\phi(p_\tau)=\tau$ for each $\tau\in[0,1]$, and if $p_{\tau_2}-p_{\tau_1}$ is a projection whenever $\tau_1\leq \tau_2$.

Let $a=a^*$ be a self-adjoint element of $\mathcal{A}$ and let $\{p_\tau:\tau\in[0,1]\}$ be a filtered collection of projections. For each $\tau\in[0,1]$, the compressed operator $p_\tau a p_\tau$ has spectral distribution of the form $(1-\tau)\delta_0 + \mu_\tau$, where $\mu_\tau$ is a subprobability measure of total mass $\tau$. 
Thus, a self-adjoint element of $\mathcal{A}$ together with a filtered family of projections gives rise to a family of subprobability measures $(\mu_\tau)_{\tau \in [0,1]}$, or equivalently to a family of distribution functions
\begin{align*}
F(x,\tau):=\mu_\tau((-\infty,x]).
\end{align*}
The operator-theoretic interlacing of spectra guarantees that this family of distribution functions is a compression. 

The free compression is the special case in which the filtered family of projections $\{p_\tau:\tau\in[0,1]\}$  is free from $a$.
\end{rem}

The essence of the compression definition is present in earlier works in the free probability literature, see e.g., Bercovici and Voiculescu \cite[Section 3]{BV}.

%%%%%%%%%%%%%%%%%%%%%%%%%%%%
\subsection{Further examples of compressions}
%%%%%%%%%%%%%%%%%%%%%%%%%%%%

Given a probability measure $\mu$ with distribution function
$F_\mu:\mathbb{R} \to [0,1]$, we can define the \textbf{classical compression}
of $\mu$ by setting
\begin{align*}
F^{\mathrm{classical}}(x,\tau) := \tau F_\mu(x).
\end{align*}
Suppose that $\mu$ has a density
$f_\mu:\mathbb{R} \to [0,\infty)$ with respect to Lebesgue measure and that
the relevant integrals are finite. A direct calculation (using $\int_{-\infty}^\infty f_\mu \log \sin( \pi F_\mu) \mathrm{d}x = - \log (2)$) gives
\begin{align*}
\mathcal{H}[F^{\mathrm{classical}}]
=
-\frac{1}{2}\int_{-\infty}^\infty
f_\mu(x)\log f_\mu(x)\mathrm{d}x
+\frac{3}{4}-\frac{1}{2}\log(2\pi).
\end{align*}
Thus, up to a factor of $1/2$ and an additive universal constant, the
compression entropy of the classical compression is the classical
differential entropy of $f_\mu$. By Theorem \ref{thm:A}, we deduce that
\begin{align*}
-\int_{-\infty}^\infty f_\mu(x)\log f_\mu(x)\mathrm{d}x
-\log(2\pi)
\leq
\int_{-\infty}^\infty\int_{-\infty}^\infty
\log|y-x|\,\mu(\mathrm{d}y)\mu(\mathrm{d}x).
\end{align*}
This is precisely the one-dimensional logarithmic
Hardy--Littlewood--Sobolev inequality; see, for example, Carlen and Loss
\cite{CL}.

Other examples of compressions of $\mu$ are given by the stochastic maximum and minimum compressions,
defined by
\begin{align*}
F^{\mathrm{max}}(x,\tau)
:=
\max\{F_\mu(x)+\tau-1,0\},
\qquad
F^{\mathrm{min}}(x,\tau)
:=
\min\{F_\mu(x),\tau\}.
\end{align*}
For each $\tau$, these place the mass $\tau$ as far to the right and as far
to the left, respectively, as is possible while respecting the continuum interlacing inequalities.

%%%%%%%%%%%%%%%%%%%%%%%%%%%%
\subsection{Atoms} \label{sec:atoms}
%%%%%%%%%%%%%%%%%%%%%%%%%%%%
In this section we outline how the discussion and results of the introduction
might extend to probability measures $\mu$ with atoms. Let $(\mu_\tau)_{\tau \in [0,1]}$ be the measures associated with a compression $F$. The size of the atom at $x$ of $\mu_\tau$ (if any) is given by
\begin{align} \label{eq:dcon}
\mu_\tau(\{x\}) := F(x,\tau) - \lim_{y \uparrow x} F(y,\tau).
\end{align}
We have the following consequence of \eqref{eq:dcon} and the continuum interlacing inequalities \eqref{eq:continuum}:

\begin{rem} \label{rem:atoms}
Let $(\mu_\tau)_{\tau \in [0,1]}$ be the measures associated with a
compression. Then for each $x \in \mathbb{R}$, $\mu_\tau(\{x\})$ is a $1$-Lipschitz function of $\tau \in [0,1]$. 
\end{rem}

We now turn to the behaviour of atoms under \emph{free} compression. The compression apparatus gives a neat way of reformulating a result of Belinschi and Bercovici \cite{BBatoms}:

\begin{thm}[Reformulation of Theorem 3.1 of Belinschi and Bercovici \cite{BBatoms}] \label{thm:BB}
Among compressions of $\mu$, free compression has the smallest atoms possible while adhering to Remark \ref{rem:atoms}. Namely, if $(\mu_\tau)_{\tau \in [0,1]}$ is the free compression of $\mu$, then 
\begin{align*}
\mu_{1-r}(\{x\}) = (\mu(\{x\}) - r)_+, \qquad \text{for all $x \in \mathbb{R}$ and $r \in [0,1]$},
\end{align*}
where $y_+ := \max\{y,0\}$. 
% If $\mu$ has an atom at $x$ of size $a>0$, then, for
% $0 \leq r < a$, the measure $\mu_{1-r}$ has an atom of size $a-r$ at $x$,
% whereas, for $a \leq r \leq 1$, the measure $\mu_{1-r}$ has no atom at $x$.
\end{thm}
 We would now like to give a heuristic variational explanation for this result using Theorem \ref{thm:A}.
Thinking of a compression $F$ in terms of its associated increasing surface $\psi:\rt \to \mathbb{R}$, these atoms manifest geometrically as regions of $\rt$ on which $\psi$ is constant. On such regions we have $\psi_s = \psi_t = 0$, and hence the integrand satisfies $\sigma(\psi_s,\psi_t) = +\infty$. We saw in Theorem \ref{thm:A} that the increasing surface associated with free compression endeavours to minimise the integral
\begin{align*}
\psi \mapsto \int_{\rt} \sigma(\nabla \psi) \mathrm{d}s\mathrm{d}t.
\end{align*}
Heuristically, this variational principle assigns infinite local cost to flat two-dimensional regions, so the optimiser should create only the flat regions forced by the terminal atoms. The minimal flat regions predicted by this principle are precisely those described in Theorem \ref{thm:BB}.

\begin{figure}[ht]
\centering
\includegraphics[width=0.33\textwidth]{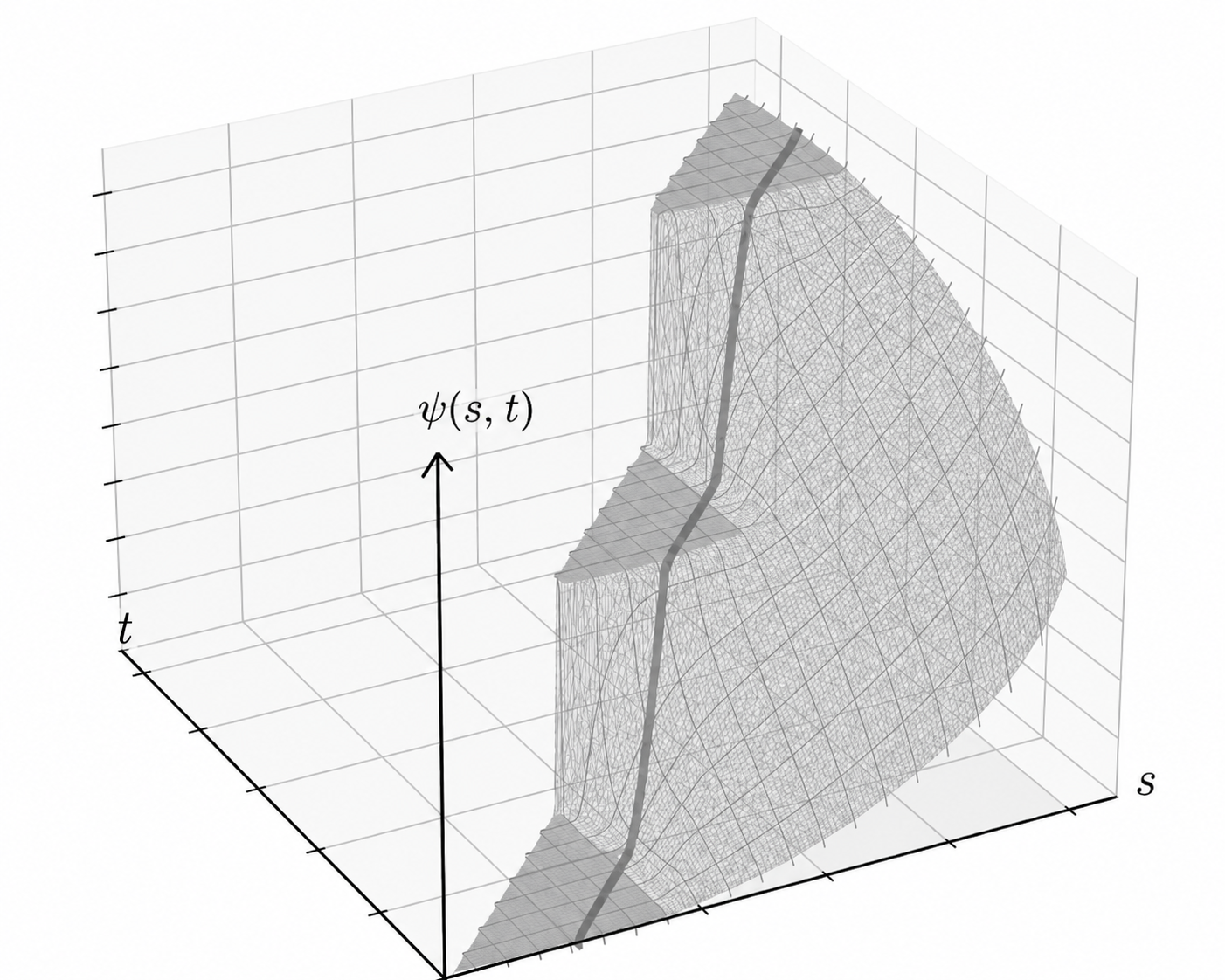}
\caption{The increasing surface associated with the free compression $(\mu_\tau)_{\tau \in [0,1]}$ of the measure $\mu_1 = \frac{1}{3} \delta_0 + \frac{1}{3} \delta_1 + \frac{1}{3}\delta_2$. For all $2/3 < \tau \leq 1$, the diagonal $\{ (1-\tau+r,r) : 0 \leq r \leq \tau \}$ has flat segments, corresponding to $\mu_\tau$ having atoms. The diagonal $\tau = 0.8$ is highlighted.}
\label{fig:atomic-free-compression}
\end{figure}

Using the tools developed in the present article, we expect that it should be possible to establish a large deviation principle for Gelfand--Tsetlin patterns whose bottom-row empirical measures converge to a probability measure $\mu$ that may have atoms. To avoid the additional technical complications arising in this setting, we do not pursue this extension here.

\subsection{The entropy of a Gelfand--Tsetlin pattern}

The compression entropy of the compression $F$ associated with a Gelfand--Tsetlin pattern $(t_{k,j})_{ 0 \leq j \leq k \leq n}$ can be written directly in terms of the $(t_{k,j})$ by the formula 
\begin{equation}
\begin{aligned}
\mathcal{H}[F] 
&= 
\frac{1}{2n^2}
\sum_{k=1}^{n}\sum_{j=0}^{k-1}
\left[
\log\!\bigl(n\Delta_{k,j}\bigr)
+
\log\sin\!\left(
\pi\frac{t_{k-1,j}-t_{k,j}}{\Delta_{k,j}}
\right) + (1 - \log \pi)
\right]
\\
&+
\frac{1}{2n^2}
\sum_{k=1}^{n-1}\sum_{j=0}^{k-1}
\left[
\log\!\bigl(n\Delta_{k,j}\bigr)
+
\log\sin\!\left(
\pi\frac{t_{k,j+1}-t_{k+1,j+1}}{\Delta_{k,j}}
\right) + (1 - \log \pi)
\right] ,
\end{aligned}
\end{equation}
where $\Delta_{k,j} := t_{k,j+1}-t_{k,j}$.
Here the entropy is understood to equal $-\infty$ whenever one of the
arguments of a logarithm or sine vanishes. If $\xi:\rt \to \mathbb{R}$ is the increasing surface associated with the Gelfand--Tsetlin pattern, then each term in the sum corresponds to an integral $- \int_{U} \sigma(\nabla \xi)$ on a subtriangle $U$ of $\rt$ with side-length $1/n$. 

It would be interesting to relate this formula to notions of finite free information \cite{GSS}.

\subsection{The semicircle law and the free central limit theorem}

The reader may wonder about the special role played by the semicircle law in free probability. The standard semicircle law is the probability distribution $\mu^{\mathrm{sc}}$ given by
\begin{align} \label{eq:semicircle}
\mu^{\mathrm{sc}}(\mathrm{d}x) = \frac{1}{2\pi} \sqrt{4-x^2} \mathrm{1}\{ x \in [-2,2]\} \mathrm{d}x.
\end{align}

Up to affine rescaling, the semicircle law is the unique non-degenerate probability distribution with finite variance whose free compressions are affine rescalings of the law itself. Indeed, if $\mu^{\mathrm{sc}}_\tau$ is its $\tau^{\text{th}}$ free compression, then $\mu^{\mathrm{sc}}_\tau$ satisfies the scaling property
\begin{align*}
\mu^{\mathrm{sc}}_\tau((-\infty,\sqrt{\tau}x ]) = \tau \mu^{\mathrm{sc}}((-\infty,x]).
\end{align*}
Equivalently, the increasing surface associated with the semicircle law has a self-similarity property. Let $Q^{\mathrm{sc}}$ be the quantile function associated with $\mu^{\mathrm{sc}}$. If $\psi^{\mathrm{sc},*}:\rt \to \mathbb{R}$ is the minimiser of $\int_{\rt} \sigma(\nabla \psi)$ subject to boundary condition $\psi^{\mathrm{sc},*}(s,s) = Q^{\mathrm{sc}}(s)$, then $\psi^{\mathrm{sc},*}$ takes the form
\begin{align} \label{eq:psisemi}
\psi^{\mathrm{sc},*}(s,t) = \sqrt{1-s+t} Q^{\mathrm{sc}}(t/(1-s+t)) \qquad (s,t) \in \rt.
\end{align}
It is possible, though the calculation is quite involved, to prove directly that $\psi^{\mathrm{sc},*}$ in \eqref{eq:psisemi} satisfies the Euler--Lagrange equations associated with maximising $\mathcal{H}[\psi]$ using \eqref{eq:semicircle} and \eqref{eq:psisemi}. The reader may consult the final section of \cite{JO} for a related calculation. 
In summary, up to affine rescalings, the semicircle law is the unique non-degenerate measure with finite variance whose minimal surface enjoys the self-similarity property in \eqref{eq:psisemi}.

Let us touch on the free central limit theorem and its ramifications for minimal surfaces. If $\mu$ is a probability measure with mean $m$ and variance $\sigma^2 \in (0,\infty)$, and $\mu_\tau$ is its $\tau^{\text{th}}$ free compression, then it is a consequence of the free central limit theorem that
\begin{align*}
\lim_{\tau \downarrow 0} \frac{1}{\tau}\mu_\tau\big((-\infty,m+\sigma\sqrt{\tau}x]\big)
=\mu^{\mathrm{sc}}((-\infty,x]),
\end{align*}
see, e.g., \cite{ST}. 
% Consequently, if $\psi^*$ is the minimal surface associated with the quantile function of $\mu$, then for every $x,y>0$,
% \begin{align} \label{eq:tipsemicircle}
% \lim_{\varepsilon\downarrow 0}
% \frac{\psi^*(1-\varepsilon x,\varepsilon y)-m}
% {\sigma\sqrt{\varepsilon}}
% =
% \sqrt{x+y}\,
% Q^{\mathrm{sc}}\left(\frac{y}{x+y}\right).
% \end{align}
% Indeed, at the point $(1-\varepsilon x,\varepsilon y)$, the compression parameter is $\tau=\varepsilon(x+y)$, while the corresponding quantile parameter is $y/(x+y)$. 
Consequently, if $\psi^*$ is the minimal surface associated with the quantile function of $\mu$ (or equivalently, the increasing surface associated with the free compression of $\mu$), 
as $(s,t)\to(1,0)$ with $t/(1-s+t)$ bounded away from $0$ and $1$, we have
\begin{align*}
\psi^*(s,t)
=
m+\sigma\psi^{\mathrm{sc},*}(s,t)
+o\left(\sqrt{1-s+t}\right).
\end{align*}
Thus, after centring and rescaling, the minimal surface associated with a finite variance measure $\mu$ has the semicircle surface as its universal asymptotic shape near the tip $(1,0)$. 

There is an analogous heavy-tailed picture: if $\mu$ lies in the domain of attraction of a non-Gaussian free stable law, then suitably centred and rescaled free compressions converge near the tip to that law, yielding the corresponding free-stable surface as the universal asymptotic shape; see Bercovici and Pata \cite{BP}. Strictly speaking, however, the results of the present article are proved only for compactly supported measures.

\subsection{The hive model} \label{sec:hive}
%%%%%%%%%%%%%%%%%%%%%%%%%%%%

In this section we discuss the hive model, which is a certain type of random surface which plays a role analogous to the Gelfand--Tsetlin pattern for the eigenvalues of sums --- as opposed to minors --- of Hermitian matrices. Whereas this article solidifies the connection between the macroscopic shape of Gelfand--Tsetlin patterns and free compression, it is an outstanding open problem to establish such a connection between the macroscopic shape of hives and free convolution, a free probability operation describing the asymptotic behaviour of the spectra of sums of large random matrices.

We outline this research program here. In 1962, Horn \cite{horn} raised the problem of characterizing the possible eigenvalues $c_1 \leq \cdots \leq c_n$ of a sum $A+B$ of Hermitian matrices $A$ and $B$ with respective eigenvalues $a_1 \leq \cdots \leq a_n$ and $b_1 \leq \cdots \leq b_n$. Horn conjectured that the possible eigenvalues are described by a polytope of inequalities involving the $(a_i)$ and $(b_i)$; this conjecture was subsequently resolved by Klyachko \cite{klyachko} and Knutson and Tao \cite{KT} in the 1990s.

One consequence of Knutson and Tao's approach is a reframing of Horn's problem in terms of a combinatorial object called a \textbf{hive}. To give a brief description of a hive here, consider tiling a unit equilateral triangle $T$ in $\mathbb{R}^2$ with subtriangles of side length $1/n$, and let $T_n \subseteq \mathbb{R}^2$ denote the set of vertices of these subtriangles. A hive is simply a rhombus-concave function $\phi:T_n \to \mathbb{R}$: that is, it has the property that if $ABC$ and $BCD$ ($A,B,C,D \in T_n$) are neighbouring triangles sharing an edge $BC$, then $\phi(B)+\phi(C) \geq \phi(A)+\phi(D)$. 
One can label the boundary height changes of a hive with reals $a_1 \leq \ldots \leq a_n$, $b_1 \leq \ldots \leq b_n$ and $c_1 \leq \ldots \leq c_n$; see Figure \ref{fig:hive}. Knutson and Tao \cite{KT} proved that there are Hermitian matrices $A, B, C$ satisfying $A+B+C = 0$ with respective eigenvalues $(a_i),(b_i),(c_i)$ if and only if there exists a hive with boundary height changes given by the $(a_i),(b_i),(c_i)$.

\begin{figure}[h!]
\begin{minipage}[c]{0.5\textwidth}
\centering
\begin{tikzpicture}[x={(1.8cm,0cm)}, z={(0cm,1.3cm)}, y={(0.7cm,1.1cm)}]
%%%%%%%%% BOTTOM LAYER
%outer triangle

\draw[->] (0,0,0) -- (0,0,3.5) node[right] {Height of $\phi$};

\draw[gray] (0,0,0) -- (2,3.464,0) -- (4,0,0) -- (0,0,0);
\draw[gray] (0.5,0.866,0) -- (3.5,0.866,0);
\draw[gray] (1, 1.732,0) -- (3, 1.732,0);
\draw[gray] (1.5, 2.598,0) -- (2.5, 2.598,0);
% Horizontal lines
\draw[gray] (0.5,0.866) -- (3.5,0.866);
\draw[gray] (1, 1.732) -- (3, 1.732);
\draw[gray] (1.5, 2.598) -- (2.5, 2.598);
% Diagonal lines from bottom-left to top-right
%\draw[gray] (0,0) -- (2,1.732);
\draw[gray] (1,0) -- (2.5,2.598);
\draw[gray] (2,0) -- (3,1.732);
\draw[gray] (3,0) -- (3.5,0.866);
% Diagonal lines from bottom-right to top-left
%\draw[gray] (4,0) -- (2,1.732);
\draw[gray] (3,0) -- (1.5,2.598);
\draw[gray] (2,0) -- (1,1.732);
\draw[gray] (1,0) -- (0.5,0.866);

%bottom boundary spikes
\draw[thick] (1,0,0) -- (1,0,1.2);
\draw[thick] (2,0,0) -- (2,0,1.9);
\draw[thick] (3,0,0) -- (3,0,1.562);
\draw[thick] (4,0,0) -- (4,0,0.92);
\draw[very thick, Tan, <-] (0,0,0)-- node[midway,below] {$c_1$} (1,0,1.2);
\draw[very thick, Tan, <-] (1,0,1.2) -- node[midway,below] {$c_2$} (2,0,1.9);
\draw[very thick, Tan, <-] (2,0,1.9) -- node[midway,below] {$c_3$} (3,0,1.562);
\draw[very thick, Tan, <-] (3,0,1.562) -- node[midway,below] {$c_4$} (4,0,0.92);

%row 1
\draw[thick] (0.5, 0.866, 0) -- (0.5, 0.866, 1.4);
\draw[thick] (1.5, 0.866, 0) -- (1.5, 0.866, 2.1);
\draw[thick] (2.5, 0.866, 0) -- (2.5, 0.866, 2.2);
\draw[thick] (3.5, 0.866, 0) -- (3.5, 0.866, 1.8);
\draw[OliveGreen] (0.5,0.866,1.4) -- (1.5,0.866,2.1) -- (2.5,0.866,2.2) -- (3.5,0.866,1.8);
\draw[OliveGreen] (0,0,0) -- (0.5,0.866,1.4) -- (1,0,1.2) -- (1.5,0.866, 2.1) -- (2,0,1.9) -- (2.5,0.866,2.2) -- (3,0,1.562) -- (3.5,0.866,1.8) -- (4,0,0.92);

\draw[very thick, Bittersweet, ->] (0,0,0)-- node[midway,left] {$a_4$} (0.5,0.866,1.4);

\draw[very thick, BrickRed, ->]  (3.5,0.866,1.8) -- node[midway,right] {$b_1$}  (4,0,0.92);

%row 2
\draw[thick] (1,1.732,0) -- (1,1.732,1.524);
\draw[thick] (2,1.732,0) -- (2,1.732, 1.811);
\draw[thick] (3,1.732,0) -- (3,1.732, 1.621);
\draw[OliveGreen] (1,1.732,1.524) -- (2,1.732,1.811) -- (3,1.732,1.621);
\draw[OliveGreen] (0.5,0.866,1.4) -- (1,1.732,1.524)  -- (1.5,0.866, 2.1) -- (2,1.732,1.811) -- (2.5,0.866,2.2) --(3,1.732,1.621)-- (3.5,0.866,1.8);

\draw[very thick, Bittersweet, ->] (0.5,0.866,1.4)-- node[midway,left] {$a_3$}  (1,1.732,1.524);
\draw[very thick, Bittersweet, ->] (1,1.732,1.524) -- node[midway,above] {$a_2$}  (1.5, 2.598, 1.071);
\draw[very thick, Bittersweet, ->] (1.5, 2.598, 1.071) -- node[midway,above] {$a_1$}  (2, 3.464, 0.33);

\draw[very thick, BrickRed, ->]  (3,1.732,1.621) -- node[midway,right] {$b_2$}  (3.5,0.866,1.8);

\draw[very thick, BrickRed, ->]  (2, 3.464, 0.33) -- node[midway, above] {$b_4$}  (2.5, 2.598, 1.107);

\draw[very thick, BrickRed, ->]  (2.5, 2.598, 1.107) -- node[midway,above right] {$b_3$}  (3,1.732,1.621);

%row 3
\draw[thick, gray] (1.5, 2.598, 0) -- (1.5, 2.598, 1.071);
\draw[thick, gray] (2.5, 2.598, 0) -- (2.5, 2.598, 1.107);
\draw[OliveGreen] (2.5, 2.598, 1.107) -- (1.5, 2.598, 1.071) -- (2, 3.464, 0.33) -- (2.5, 2.598, 1.107);
\draw[OliveGreen] (1,1.732,1.524)  --  (1.5, 2.598, 1.071)-- (2,1.732,1.811) --(2.5, 2.598, 1.107) --(3,1.732,1.621);

%row 4 = top spike
\draw[thick, gray] (2, 3.464, 0) -- (2, 3.464, 0.33);

\end{tikzpicture}
\end{minipage}\hfill
\begin{minipage}[c]{0.5\textwidth}
\caption{
A hive is a rhombus-concave function defined on a triangular lattice. Knutson and Tao \cite{KT} showed that the height changes on the boundary of a hive may be associated with the possible eigenvalues of Hermitian matrices satisfying $A+B+C=0$.
} \label{fig:hive}
\end{minipage}
\end{figure}

Now on the one hand, Coquereaux and Zuber \cite{CZ} (see also \cite{CMZ}) recently showed that if one chooses a hive randomly from the set of hives with two of their boundaries fixed by $a_1 \leq \cdots \leq a_n$ and $b_1 \leq \cdots \leq b_n$, and the remaining boundary biased by the Vandermonde determinant $\Delta(c) := \prod_{1 \leq i < j \leq n } (c_j - c_i)$, then the random free boundaries $\gamma_1 \leq \cdots \leq \gamma_n$ have the same law as the eigenvalues of $- (A+B)$, where $A$ and $B$ are independent unitarily invariant Hermitian matrices with respective eigenvalues $(a_i)$ and $(b_i)$. On the other hand, one of the central results in free probability states that if $(A_n)_{n \geq 1}$ and $(B_n)_{n\geq 1}$ are independent sequences of unitarily invariant random matrices with empirical spectra converging to $\mu$ and $\nu$, then the empirical spectrum of $A_n+B_n$ converges almost surely to a probability measure $\mu \boxplus \nu$ called the free convolution of $\mu$ and $\nu$; see e.g.\ \cite{Vcon1, Vcon3}.

We saw in our main results that we could recover free compression from the statistical mechanics of Gelfand--Tsetlin patterns. The previous paragraph then raises the question, can we recover additive free convolution from the statistical mechanics of random hives?
On this front, there has been recent progress by Narayanan and Sheffield \cite{NarS} and Narayanan, Sheffield and Tao \cite{NST}.
In \cite{NarS}, it is shown that the free edge of a large random hive, which characterises the empirical spectrum of $A_n + B_n$, satisfies a large deviation type bound.
The article \cite{NST} considers \textbf{augmented hives}, which are hives with a Gelfand--Tsetlin pattern attached in order to swallow the Vandermonde factor $\Delta(c)$. The authors of \cite{NST} show that augmented hives exhibit concentration when the two free boundaries are chosen randomly according to the eigenvalues of independent GUE random matrices. See further work by Narayanan and coauthors \cite{hari1, hari2, hari3, hari4}, as well as blog posts by Tao  \cite{taoblog,taoblog2}, for further discussion of this circle of problems.

\subsection{Gaussian free field fluctuations for Gelfand--Tsetlin patterns}

At the macroscopic level, Ginzburg--Landau $\nabla \phi$-interface models in two dimensions with fixed boundary conditions on the $n$-by-$n$ square have been shown, at least in the case $0 < c_1 \leq V''(x) \leq c_2$, to have fluctuations of order $\sqrt{\log n}$ with Gaussian free field correlations \cite{miller}.
The case $d=2$ of \cite[Theorem 1.4]{MP} echoes this order of fluctuations under weaker conditions on the interaction potential.
% Further still, let us mention that Magazinov and Peled \cite{MP} prove concentration properties for random surfaces with interaction potentials satisfying the weaker condition that the second derivative of $V$ exists and is positive almost everywhere. The case $d=2$ of \cite[Theorem 1.4]{MP} appears to be strong evidence that in a large random Gelfand--Tsetlin pattern with bottom row $s_1 < \ldots < s_n$, where $s_{i+1} - s_i = O(1)$, the variance of a bulk point is $O(\lo)$. 
We conjecture that this behaviour continues to hold for Gelfand--Tsetlin patterns:

\begin{conj} \label{conj:GFF}
Gelfand--Tsetlin random surfaces on sufficiently nice large subsets of $\mathbb{Z}^2$ of diameter $n$ have Gaussian free field fluctuations of order $\sqrt{\log n}$ times the order of the gaps on the boundary. In particular, in the setting of Theorem \ref{thm:C}, in the bulk the variables $(t^n_{k,j})_{0 \leq j \leq k \leq n}$ have Gaussian free field correlations with fluctuations of order $\sqrt{\log n}/n$. 
\end{conj}

Let us touch on how this conjecture, which is based on literature on random surfaces, is also compatible with work in random matrix theory. Consider a Ginzburg--Landau $\nabla\phi$-interface model on $S_n := \{1,\ldots,n\}^2$ with hard-core interactions, no boundary conditions, and a Gaussian external field along the diagonal. That is, consider a random function $\phi:S_n \to \mathbb{R}$ sampled from the probability measure
\begin{align*}
\frac{1}{Z_n} \mathrm{1}_{\mathrm{GT}}^{S_n}(\phi) \exp \left\{ -\frac{1}{2} \sum_{i=1}^n \phi_{i,i}^2 \right\} \prod_{x\in S_n} \mathrm{d}\phi_x,
\end{align*}
where as in \eqref{eq:globalind}, $\mathrm{1}_{\mathrm{GT}}^{S_n}(\phi)$ is the indicator function that $\phi:S_n \to \mathbb{R}$ is increasing. 
It was observed in \cite{JO} that the marginal law of the diagonal $(\phi_{1,1},\ldots,\phi_{n,n})$ of such a surface coincides, up to scaling, with the law of the $n$ eigenvalues of a GUE random matrix. With this observation in mind, Conjecture \ref{conj:GFF} accords with Gustavsson's central limit theorem \cite{gustavsson} for the position of a bulk 
eigenvalue of a GUE random matrix, which states that for $k \in (\varepsilon n, (1-\varepsilon)n)$, $\phi_{k,k}$ is approximately Gaussian distributed with fluctuations of the order $O(\sqrt{\log (n)/n})$.

It would also be interesting to understand whether these conjectural Gaussian free field fluctuations admit an interpretation within second-order free probability \cite{sfree1,sfree2,sfree3}, or more broadly within the theory of fluctuations of Hermitian random matrices \cite{johansson}.
We mention that global fluctuations of several classes of interlacing arrays and random tilings have been shown to converge to the Gaussian free field; see, for example, Borodin and Gorin \cite{BG} and Bufetov and Gorin \cite{BuG}.

\subsection{Higher dimensional Gelfand--Tsetlin surfaces}

There is a natural generalisation of Gelfand--Tsetlin surfaces to dimensions other than two, which we define as follows. In parallel to our discussion in Section \ref{sec:surfaces}, one can define a $d$-dimensional Gelfand--Tsetlin surface to be a function $\phi:C \to \mathbb{R}$ defined on a subset $C$ of the $d$-dimensional lattice $\mathbb{Z}^d$ and satisfying $\phi_x \leq \phi_y$ whenever $x \leq y$ in the natural partially ordered set structure of $\mathbb{Z}^d$. One can define the surface tension in $d$ dimensions in analogy to how it was defined for $d=2$ in Section \ref{sec:thmB}.

When $d=1$, it is an easy calculation using the fact that
\begin{align*}
\int_{0 \leq \phi_1 \leq \ldots \leq \phi_{n-1} \leq u_1 n } \mathrm{d}\phi_1 \cdots \mathrm{d}\phi_{n-1} = \frac{(u_1n)^{n-1}}{(n-1)!}
\end{align*}
to show that the one-dimensional Gelfand--Tsetlin surface tension is given by
\begin{align*}
\sigma(u_1) =
\begin{cases}
-\log u_1 - 1 \quad &\text{if $u_1 > 0$,}\\
+ \infty, \quad &\text{otherwise},
\end{cases}
\end{align*}
and of course Theorem \ref{thm:B} tells us that the Gelfand--Tsetlin surface tension in two dimensions is given by
\begin{align*}
\sigma(u_1,u_2) = 
\begin{cases} - \log(u_1 + u_2) - \log \sin \pi (u_1/(u_1+u_2)) - 1 + \log \pi, \quad &\text{if $u_1,u_2 > 0$}\\
+ \infty, \quad &\text{otherwise}.
\end{cases}
\end{align*}
This leaves the following natural open problem:

\begin{que}
Is there an explicit formula for the surface tension of Gelfand--Tsetlin surfaces in dimension $d$ for $d \geq 3$?
\end{que}

In fact, let us observe that given any partially ordered set $(P,\preceq)$, one can define a Gelfand--Tsetlin surface on $P$ to be a function $\phi:P \to \mathbb{R}$ such that $\phi_x \leq \phi_y$ whenever $x \preceq y$. One can then of course ask about statistical mechanics properties of Gelfand--Tsetlin surfaces on large sets with a partially ordered set structure. The case where $P$ is a tree and the relation is ancestry is closely related to work on accessibility percolation; see e.g., \cite{mattdiana, BBS}.

%%%%%%%%%%%%%%%%%%%%%%%%%%%%%%%%%%%%%%%%%%%%%%
%%%%%%%%%%%%%%%%%%%%%%%%%%%%%%%%%%%%%%%%%%%%%%
%%%%%%%%%%%%%%%%%%%%%%%%%%%%%%%%%%%%%%%%%%%%%%
\section{Proof of Theorem A} \label{sec:free}
%%%%%%%%%%%%%%%%%%%%%%%%%%%%%%%%%%%%%%%%%%%%%%
%%%%%%%%%%%%%%%%%%%%%%%%%%%%%%%%%%%%%%%%%%%%%%
%%%%%%%%%%%%%%%%%%%%%%%%%%%%%%%%%%%%%%%%%%%%%%

In this section we begin the proofs of our main results and essentially prove Theorem \ref{thm:A}. More precisely, we prove Theorem \ref{thm:freeEL}, which states that the free compression $F^*$ satisfies the Euler--Lagrange equations for $\mathcal{H}[\cdot]$, and Theorem \ref{thm:compressionachieve}, which states that $\mathcal{H}[F^*] = \chi[\mu]$.
We do not, however, prove here with full rigour that $F^*$ achieves the supremum of $\mathcal{H}[F]$ over all compressions $F$ of $\mu$. This will follow from Theorem \ref{thm:C} and be completed in Section \ref{sec:Cproof}.

\medskip
Throughout Section \ref{sec:free}, all measures we will consider will have compact support and a density with respect to Lebesgue measure.

%%%%%%%%%%%%%%%%%%%%%%%%%%%%%%%%%%%%
\subsection{Cauchy and log-transforms} \label{sec:cauchy}
%%%%%%%%%%%%%%%%%%%%%%%%%%%%%%%%%%%%
We briefly collect together some first properties of Cauchy transforms and log transforms of finite measures. Let us define the upper and lower half-plane by 
\begin{align*}
\mathbb{C}^\pm := \{ x \pm iy : x \in \mathbb{R}, y > 0 \}.
\end{align*}
Let $\mu$ be any finite measure on the real line, not necessarily a probability measure. Let $\log:\mathbb{C}^+ \to \mathbb{R} \times i(0,\pi)$ be the principal branch of the complex logarithm, and define the Cauchy transform and log-potential of $\mu$ to be the analytic functions $G_\mu:\mathbb{C}^+ \to \mathbb{C}^-$ and $H_\mu:\mathbb{C}^+ \to \mathbb{C}^+$ given by 
\begin{align*}
G_\mu(z) := \int_{-\infty}^\infty \frac{1}{z-x} \mu(\mathrm{d}x) \quad \text{and} \quad H_\mu(z) := \int_{-\infty}^\infty \log (z - x) \mu ( \mathrm{d}x).
\end{align*}
Note that $H_\mu'(z) = G_\mu(z)$.
We will separate the real and imaginary parts of $G_\mu(z)$ and $H_\mu(z)$ by setting
\begin{align} 
G_\mu(z) &= \pi (u_\mu(z) - i v_\mu(z)), \label{eq:IR1} \\
H_\mu(z) &= \pi (U_\mu(z) - i V_\mu(z)).\label{eq:IR2}
\end{align}
Note $v_\mu(z) > 0$ and $V_\mu(z) < 0$ for all $z \in \mathbb{C}^+$. 

If $\mu$ has a density with respect to Lebesgue measure, then it may be verified that $G_\mu(z)$ and $H_\mu(z)$ have limits as they approach the real line at almost all points.
For $G_\mu$ we have the Plemelj formula
\begin{align} \label{eq:plemelj}
\lim_{ \varepsilon \downarrow 0} G_\mu( x + i \varepsilon ) =: \pi (u_\mu(x) - i v_\mu(x)),
\end{align}
where $v_\mu(x)$ coincides precisely with the density of $\mu$, and $u_\mu$ is the principal value integral
\begin{align} \label{eq:plemelj2}
u_\mu(x) = \frac{1}{\pi} \mathrm{p.v.} \int_{-\infty}^\infty \frac{1}{x-y} v_\mu(y) \,\mathrm{d}y := \frac{1}{\pi} \lim_{\varepsilon \downarrow 0} \int_{\mathbb{R}\setminus(x-\varepsilon,x+\varepsilon)} \frac{1}{x-y} v_\mu(y) \,\mathrm{d}y.
\end{align}
Likewise, for $H_\mu$ we have the integrated Plemelj formula
\begin{align*}
\lim_{ \varepsilon \downarrow 0} H_\mu( x + i \varepsilon ) =: \pi (U_\mu(x) - i V_\mu(x)),
\end{align*}
where
\begin{align} \label{eq:VUdef}
V_\mu(x) = - \int_x^\infty v_\mu(y) \,\mathrm{d}y \qquad \text{and} \qquad U_\mu(x) = \frac{1}{\pi} \int_{-\infty}^\infty \log | x-y| \,v_\mu(y) \,\mathrm{d}y,
\end{align}
provided the latter integral exists.
Note that for derivatives with respect to the real variable $x$ we have 
\begin{align} \label{eq:vudef}
V_\mu'(x) = v_\mu(x) \qquad \text{and} \qquad U_\mu'(x) = u_\mu(x).
\end{align}

\subsection{The $R$-transform construction of free compression and the Burgers equation} \label{sec:Rtransform}

In this section we construct free compression using the $R$-transform and show its Cauchy transform satisfies a Burgers equation.
% We close the section on free probability by briefly explaining how our construction of free compression in terms of the Burgers equation coincides with the usual definition of free compression in terms of $R$-transforms.

We define the $R$-transform $R_\mu$ of a probability measure $\mu$ to be the complex-valued function defined in a neighbourhood of the origin in $\mathbb{C}$ and satisfying
\begin{align} \label{eq:Rtran}
R_\mu(G_\mu(z)) = z - \frac{1}{G_\mu(z)}.
\end{align}
It is a consequence of the results of Nica and Speicher \cite{NS} that for each $\tau \in (0,1]$, there is a probability measure $[\mu]_\tau$ (which we call the normalised free compression) whose $R$-transform is precisely 
\begin{align} \label{eq:Rcomp}
R_{[\mu]_\tau}(s) = R_\mu( \tau s).
\end{align} 
% The existence of a measure with this $R$-transform is guaranteed by the Cauchy transform existence theorem (see, e.g.,\,\cite[Chapter 3.1]{MS}). (See \cite{ST} for an alternative construction in terms of fractional free additive convolution.)
We then simply define free compression by multiplying the probability measure $[\mu]_\tau$ by $\tau$ to construct a measure of total mass $\tau$, i.e.,
\begin{align} \label{eq:Rcomp2}
\mu_\tau := \tau [ \mu]_\tau. 
\end{align}

The measure $\mu_\tau$ constructed here is precisely the measure $\mu_\tau$ appearing in the statement of Theorem \ref{thm:compression}.
Through the remainder of this section $\{ F^*(x,\tau) : x \in \mathbb{R} , \tau \in [0,1] \}$ and the associated measures $(\mu_\tau)_{\tau \in [0,1]}$ will denote the free compression associated with $\mu$.

The following lemma is a variant on well-known calculations (see e.g. \cite{ST}). We would like to emphasise that our choice to renormalise the free compression so that it is a subprobability measure of total mass $\tau$ rather than a genuine probability measure appears to significantly simplify the form of the Burgers equation.

\begin{lemma}[\cite{ST}] \label{lem:cburgers}
The mixed Cauchy transforms $\{ G(z,\tau) : z \in \mathbb{C}^+, \tau \in [0,1] \}$
\begin{align*}
G(z,\tau) := \int_{-\infty}^\infty \frac{1}{z-x} \mu_\tau(\mathrm{d}x)
\end{align*}
of the free compression measures $(\mu_\tau)_{\tau \in [0,1]}$ defined via \eqref{eq:Rtran}, \eqref{eq:Rcomp} and \eqref{eq:Rcomp2} satisfy the Burgers equation
\begin{align} \label{eq:cburgers}
\partial_\tau G+\frac{1}{G}\partial_z G=0, \qquad z \in \mathbb{C}^+, \tau \in (0,1).
\end{align} 
\end{lemma}

\begin{proof}
Note that $G(z,\tau) = \tau G_{[\mu]_\tau}(z)$, and hence
$G_{[\mu]_\tau}(z)= \tau^{-1}G(z,\tau)$. Applying \eqref{eq:Rtran} to the probability measure $[\mu]_\tau$ gives
\[
R_{[\mu]_\tau}\!\left(\frac{G(z,\tau)}{\tau}\right)=z-\frac{\tau}{G(z,\tau)}.
\]
Using \eqref{eq:Rcomp}, we obtain
\begin{align}
\label{eq:compressed-implicit}
z=R_\mu(G(z,\tau))+\frac{\tau}{G(z,\tau)},
\end{align}
which implicitly characterises $G(z,\tau)$. By the inverse function theorem, $G(z,\tau)$ varies smoothly in $z$ and $\tau$. Taking respective derivatives of \eqref{eq:compressed-implicit} with respect to $z$ and $\tau$ gives 
\begin{align*}
1=\left(R_\mu'(G)-\frac{\tau}{G^2}\right)\partial_z G \qquad \text{and} \qquad 0=\left(R_\mu'(G)-\frac{\tau}{G^2}\right)\partial_\tau G+\frac{1}{G}.
\end{align*}
Eliminating the common factor $(R_\mu'(G)-\tau G^{-2})$ from these two identities yields \eqref{eq:cburgers}, as required.
\end{proof}

Given the collection of free compression measures $(\mu_\tau)_{\tau \in [0,1]}$, let us define functions $\{ G(z,\tau), H(z,\tau) : z \in \mathbb{C}^+, \tau \in [0,1]\}$ by setting
\begin{align*}
G(z,\tau) := \int_{-\infty}^\infty \frac{1}{z-x} \mu_\tau(\mathrm{d}x) \quad \text{and} \quad H(z,\tau) := \int_{-\infty}^\infty \log(z-x) \mu_\tau(\mathrm{d}x).
\end{align*}
We write 
\begin{align} \label{eq:GH}
G(z,\tau) = \pi( u(z,\tau) - iv(z,\tau)) \quad \text{and} \quad H(z,\tau) = \pi( U(z,\tau) - iV(z,\tau)),
\end{align}
and define $G(x,\tau), H(x,\tau), u(x,\tau)$ etc., for $x \in \mathbb{R}$ by taking limits $z = x+i\varepsilon$ as $\varepsilon \downarrow 0$.

By \eqref{eq:plemelj2} and \eqref{eq:VUdef}, the functions $u(x,\tau), v(x,\tau), U(x,\tau), V(x,\tau)$ can be constructed directly from the measures $\mu_\tau$; in particular, $v(x,\tau)$ is the density of $\mu_\tau$ at $x$. 

We define the \textbf{liquid region} associated with the free compression by 
\begin{align} \label{eq:liquid}
\mathcal{L} := \{ (x,\tau) \in \mathbb{R} \times (0,1) : v(x,\tau) > 0 \}.
\end{align}
On this liquid region, the boundary value function $\{G(x,\tau) : x \in \mathbb{R}, \tau \in [0,1]\}$ is real analytic and satisfies the real Burgers equation
\begin{align} \label{eq:realburgers}
\partial_\tau G + \frac{1}{G} \partial_x G = 0.
\end{align}

\begin{lemma} \label{lem:cata}
With $u = u(x,\tau)$ and $v = v(x,\tau)$, for $(x,\tau)$ in the liquid region $\mathcal{L}$ we have
\begin{align}
U_\tau(x,\tau) = - \frac{1}{2\pi} \log( u^2 + v^2) + \frac{1}{\pi} ( \log \tau - \log \pi),\label{eq:cata1c}
\end{align}
and
\begin{align} 
V_\tau(x,\tau) = \frac{1}{\pi} \mathrm{arctan}( u/v) - 1/2.\label{eq:cata2c}
\end{align}
\end{lemma}
\begin{proof}
Note that $\partial_z H = G$. In particular, the Burgers equation \eqref{eq:cburgers} gives $\partial_z ( H_\tau + \log G) = 0$, where we are taking the principal branch of the logarithm. It follows that for some function $c(\tau)$ we have $H_\tau = - \log G + c(\tau)$. Using the asymptotics $H(z,\tau) \sim \tau \log z$ and $G(z,\tau) \sim \tau/z$, we conclude that $c(\tau) = \log \tau$. Thus 
\begin{equation} \label{eq:crab}
    H_\tau = -\log G+ \log \tau, \qquad z \in \mathbb{C}^+, \tau \in (0,1),
\end{equation}
which by \eqref{eq:GH} reads 
\begin{align} \label{eq:clam}
\pi(U_\tau - i V_\tau) = -\log(\pi(u-iv)) + \log \tau, \qquad z \in \mathbb{C}^+, \tau \in (0,1).
\end{align}
The equations \eqref{eq:cata1c} and \eqref{eq:cata2c} now follow from taking real and imaginary parts of \eqref{eq:clam}, and then taking the boundary values in the liquid region.
\end{proof}

We now note that free compression $F^* = \{ F^*(x,\tau) : x \in \mathbb{R}, \tau \in [0,1] \}$ is indeed a compression. 
Clearly each $\mu_\tau$ is a measure of total mass $\tau$. 
It remains to establish the continuum interlacing inequalities.
In this direction, note from \eqref{eq:VUdef} and the fact that $F^*(\cdot,\tau)$ is the distribution function of a measure with total mass $\tau$, we have $V(x,\tau) = - (\tau - F^*(x,\tau))$. In particular, noting from \eqref{eq:cata2c} that $V_\tau(x,\tau)$ takes values in $[-1,0]$, we see that $F^*_\tau$ takes values in $[0,1]$, thereby establishing that $F^*_\tau$ satisfies the continuum interlacing inequalities \eqref{eq:continuum}.

Let us also note that by \eqref{eq:cata2c} we have $\frac{1}{\pi} \arctan(u/v) = F^*_\tau - 1/2$. Rearranging and using $v = F^*_x$ we obtain expressions for $u$ and $v$ in terms of $F^*_x$ and $F^*_\tau$:
\begin{align} \label{eq:trueform}
v = F^*_x \quad \text{and} \quad u = -F^*_x\frac{\cos(\pi F^*_\tau)}{\sin(\pi F^*_\tau)} \qquad (x,\tau) \in \mathcal{L}.
\end{align}

We would now like to relate the Burgers equation above to the Euler--Lagrange equations for the compression entropy functional, which we recall may be written 
\begin{align} \label{eq:ent2}
\mathcal{H}[F] := - \int_{-\infty}^\infty \int_0^1\Sigma( F_x, F_\tau)  \mathrm{d}\tau \mathrm{d}x,
\end{align}
where 
\begin{align} \label{eq:Sigma2}
\Sigma( F_x, F_\tau) = F_x \left\{ \log F_x - \log \sin (\pi F_\tau) + \log \pi - 1 \right\}.
\end{align}
The associated Euler--Lagrange equation is given by 
\begin{align} \label{eq:EL2}
0 =\frac{\partial}{\partial x} \frac{\partial \Sigma}{\partial F_x} + \frac{\partial}{\partial \tau} \frac{\partial \Sigma}{\partial F_\tau} =\frac{\partial}{\partial x} \left( \log F_x - \log \sin (\pi F_\tau) \right)   -  \frac{\partial}{\partial \tau}\left( \pi F_x \frac{ \cos(\pi F_\tau) }{ \sin (\pi F_\tau) }  \right).
\end{align}

We are now ready to prove the first part of Theorem \ref{thm:A}.

\begin{thm} \label{thm:freeEL}
Free compression $F^*$ satisfies the Euler--Lagrange equations \eqref{eq:EL2} associated with the compression entropy.
\end{thm}
\begin{proof}
In light of the fact that the Cauchy transforms associated with free compression satisfy the Burgers equation \eqref{eq:cburgers}, it is sufficient to establish the identity 
\begin{align} \label{eq:elb}
\mathrm{Re}\left( \partial_\tau G + \frac{\partial_x G}{G} \right) = \frac{\partial}{\partial x} \frac{\partial \Sigma}{\partial F_x}(F_x^*,F_\tau^*) 
+
\frac{\partial}{\partial \tau} \frac{\partial \Sigma}{\partial F_\tau}(F_x^*,F_\tau^*).
\end{align}

Using \eqref{eq:trueform} to obtain the second equality below we have  
\begin{align} \label{eq:tt1}
\mathrm{Re}\left( \partial_\tau G \right) = \pi u_\tau = - \pi \frac{\partial}{\partial \tau} \left( F^*_x\frac{\cos(\pi F^*_\tau)}{\sin(\pi F^*_\tau)} \right).
\end{align}
Conversely, by \eqref{eq:trueform}, $u^2 + v^2 = (F^*_x)^2 /\sin(\pi F^*_\tau)^2$, so that
\begin{align} \label{eq:tt2}
\mathrm{Re} \left( \frac{ \partial_x G}{ G} \right) = \mathrm{Re}( \partial_x \log G) = \partial_x \log |G|=  \frac{1}{2}\partial_x\log(u^2+v^2) = \partial_x\log F^*_x-\partial_x\log\sin(\pi F^*_\tau). 
\end{align}
Using \eqref{eq:tt1} and \eqref{eq:tt2} and comparing with \eqref{eq:EL2} we obtain \eqref{eq:elb}. 
\end{proof}

\subsection{The compression entropy of free compression}
Our first step to calculating $\mathcal{H}[F^*]$ is a local expression for $\Sigma( F^*_x, F^*_\tau)$ in terms of the logarithmic potentials $U(x,\tau)$.

\begin{lemma}
We have
\begin{align} \label{eq:Sigmanew}
\Sigma( F^*_x, F^*_\tau) =  v ( - \pi U_\tau + \log \tau - 1).
\end{align}
\end{lemma}

\begin{proof}
By \eqref{eq:trueform} we have $u^2+v^2 = (F^*_x)^2/\sin(\pi F^*_\tau)^2$, so that $\log F^*_x - \log \sin \pi F^*_\tau = \frac{1}{2}\log(u^2 + v^2)$. Using this fact and $F^*_x = v$ in the definition \eqref{eq:Sigma2} of $\Sigma(F^*_x,F^*_\tau)$ we obtain 
\begin{align*}
\Sigma(F^*_x,F^*_\tau) &= F^*_x \left\{  \frac{1}{2} \log(u^2 + v^2) + \log \pi - 1 \right\}.
\end{align*}
Using the identity \eqref{eq:cata1c} for $U_\tau$ completes the proof.
\end{proof}

In particular, since the integrand in $\mathcal{H}[F^*]$ is nonzero only when $v = F^*_x$ is nonzero, we can in fact write the compression entropy as an integral over the liquid region
\begin{align*}
\mathcal{H}[F^*] = \iint_{\mathcal{L}} v (  \pi U_\tau - \log \tau + 1) \mathrm{d}x\mathrm{d}\tau.
\end{align*}
We bear this fact in mind in the upcoming proof, in the knowledge that the integral is supported in precisely the region where the Burgers equation is satisfied. 

We conclude with our calculation of the compression entropy of free compression:

\begin{thm} \label{thm:compressionachieve}
Let $F^*$ be the free compression of $\mu$. Then $\mathcal{H}[F^*] = \chi[\mu]$.
\end{thm}
\begin{proof}
Using \eqref{eq:ent2} and \eqref{eq:Sigmanew} we have
\begin{align} \label{eq:fresh}
\mathcal{H}[F^*] =  \pi \int_0^1 \int_{-\infty}^\infty v(x,\tau) U_\tau(x,\tau) \mathrm{d}x ~ \mathrm{d}\tau - \int_0^1 (\log \tau - 1) \left( \int_{-\infty}^\infty v(x,\tau)\mathrm{d}x \right) \mathrm{d}\tau.
\end{align}
Now on the one hand, looking at the latter term on the right-hand side of \eqref{eq:fresh}, since $v(x,\tau)$ is the density of a measure with total mass $\tau$, we have 
\begin{align} \label{eq:fresh2}
\int_0^1 (\log \tau - 1) \left( \int_{-\infty}^\infty v(x,\tau)\mathrm{d}x \right) \mathrm{d}\tau = \int_0^1 \tau (\log \tau - 1)\mathrm{d}\tau = - \frac{3}{4}.
\end{align} 
As for the former term on the right-hand side of \eqref{eq:fresh}, define
\[
I(\tau) := \int_{-\infty}^\infty \int_{-\infty}^\infty  \log|y-x|\, v(x,\tau) v(y,\tau)\mathrm{d}x\mathrm{d}y.
\]
By symmetry and the definition of $U(x,\tau)$ we have
\begin{align*}
I'(\tau) = 2 \int_{-\infty}^\infty \int_{-\infty}^\infty  \log |y-x|\, v(x,\tau) v_\tau(y,\tau) \mathrm{d}x\mathrm{d}y = 2 \pi \int_{-\infty}^\infty v(x,\tau) U_\tau(x,\tau)\mathrm{d}x. 
\end{align*}
It follows that 
\begin{align} \label{eq:fresh3}
 \pi \int_0^1 \int_{-\infty}^\infty v(x,\tau) U_\tau(x,\tau) \mathrm{d}x ~ \mathrm{d}\tau
=  \frac{1}{2} \int_0^1 I'(\tau) \mathrm{d}\tau
=  \frac{1}{2}I(1),
\end{align}
where we have used the fact that $I(0) = 0$. 

Using \eqref{eq:fresh2} and \eqref{eq:fresh3} in \eqref{eq:fresh} we obtain
\[
\mathcal{H}[F^*] 
=
\frac{1}{2}\int_{-\infty}^\infty \int_{-\infty}^\infty  \log|y-x|\, v(x,1) v(y,1)\,\mathrm{d}x\mathrm{d}y
+\frac{3}{4} = \chi[\mu],
\]
where the final equality above follows from the fact that $v(\cdot,1)$ is the density function of $\mu$. 

That completes the proof.
\end{proof}

%%%%%%%%%%%%%%%%%%%%%%%%%%%%%%%%%%%%%%%%%%%%%%
%%%%%%%%%%%%%%%%%%%%%%%%%%%%%%%%%%%%%%%%%%%%%%
%%%%%%%%%%%%%%%%%%%%%%%%%%%%%%%%%%%%%%%%%%%%%%
\section{Basic properties of random surfaces and Gelfand--Tsetlin surfaces} \label{sec:randomsurfaces}
%%%%%%%%%%%%%%%%%%%%%%%%%%%%%%%%%%%%%%%%%%%%%%
%%%%%%%%%%%%%%%%%%%%%%%%%%%%%%%%%%%%%%%%%%%%%%
%%%%%%%%%%%%%%%%%%%%%%%%%%%%%%%%%%%%%%%%%%%%%%

In this section we study further basic properties of random surfaces and Gelfand--Tsetlin surfaces.

%%%%%%%%%%%%%%%%%%%%%%%%%%%%%%%%%%%%%%%%%%%%%%
\subsection{Pr\'ekopa--Leindler inequalities for random surfaces with convex potentials} \label{sec:PL}
%%%%%%%%%%%%%%%%%%%%%%%%%%%%%%%%%%%%%%%%%%%%%%

Before specializing to Gelfand--Tsetlin surfaces in Section \ref{sec:gtearnest}, in this section we prove the Pr\'ekopa--Leindler inequalities for random surfaces with general convex potentials given in Proposition \ref{prop:pl}. 

The classical Pr\'ekopa--Leindler inequality \cite{AGA} states that if $f,g,h:\mathbb{R}^n \to [0,\infty)$ are measurable functions such that for some $\lambda \in (0,1)$ and all $s,s' \in \mathbb{R}^n$
\[
h(\lambda s + (1 - \lambda) s') \geq f(s)^\lambda g(s')^{1 - \lambda},
\]  
then we have
\begin{align} \label{eq:PLclassical}
\int_{\mathbb{R}^n} h(s) \mathrm{d}s \geq \left( \int_{\mathbb{R}^n} f(s)\mathrm{d}s \right)^\lambda \left( \int_{\mathbb{R}^n} g(s) \mathrm{d}s \right)^{1 - \lambda};
\end{align}
see, e.g.,\ \cite[Section 1.4]{AGA}.

We say that a probability density function $f:\mathbb{R}^n \to [0,\infty)$ is log-concave if it can be written $f(s) = e^{ -W(s)}$, where $W:\mathbb R^n\to (-\infty,\infty]$ is a convex function. It is a well-known consequence of the Pr\'ekopa--Leindler inequality that if $s = (s_1,\ldots,s_n)$ is distributed according to a log-concave probability density function, then the marginal distribution of the first $k$ coordinates $(s_1,\ldots,s_k)$ of $s$ has a log-concave density on $\mathbb{R}^k$ \cite{AGA}. We now show that the standard Pr\'ekopa--Leindler inequality may be used to prove Proposition \ref{prop:pl}:

\begin{proof}[Proof of Proposition \ref{prop:pl}]
\textbf{Proof of (1)}. 
Define a functional $J:\mathbb{R}^C \to \mathbb{R}$ on functions $\phi:C \to \mathbb{R}$ by setting
\begin{align} \label{eq:Jdef}
J(\phi) := \exp \left\{ - \sum_{ \langle x,y \rangle \in C^* } V(\phi_y - \phi_x) \right\}.
\end{align} 
By the convexity of $V$, if $\phi,\phi':C \to \mathbb{R}$, we have
\begin{align} \label{eq:Jconv}
J(\lambda \phi + (1-\lambda)\phi') \geq J(\phi)^\lambda J(\phi')^{1-\lambda}.
\end{align}
Now given $\varphi:B \to \mathbb{R}$, define a functional $I_\varphi:\mathbb{R}^A \to [0,\infty)$ by setting 
\begin{align*}
I_\varphi(s) := J( \varphi \cup s),
\end{align*}
where we view $s \in \mathbb{R}^A$ as a function $s:A \to \mathbb{R}$, and $\varphi \cup s$ is the function on $C = A \sqcup B$ taking the value $s_z$ for $z \in A$ and $\varphi_z$ for $z \in B$. 

Note that with $Z_A(\varphi)$ as in \eqref{eq:pf} we have
\begin{align} \label{eq:Znote}
Z_A(\varphi) = \int_{\mathbb{R}^A} I_\varphi(s) \mathrm{d}s.
\end{align}

We also note that \eqref{eq:Jconv} implies
\begin{align} \label{eq:PLnew}
I_{\lambda \varphi + (1 - \lambda) \tilde{\varphi} } (\lambda s+ (1 - \lambda) s') \geq I_{\varphi}(s)^\lambda I_{\tilde{\varphi}}(s')^{1 - \lambda}.
\end{align} 

Applying the Pr\'ekopa--Leindler inequality \eqref{eq:PLclassical} with 
\[
h(s) = I_{\lambda \varphi + (1 - \lambda) \tilde{\varphi}}(s), 
\quad 
f(s) = I_{\varphi}(s), 
\quad 
g(s) = I_{\tilde{\varphi}}(s),
\]
and using \eqref{eq:PLnew} and \eqref{eq:Znote}, we obtain
\begin{align*}
Z_A(\lambda \varphi + (1-\lambda)\tilde{\varphi}) \geq Z_A(\varphi)^\lambda Z_A(\tilde{\varphi})^{1-\lambda},
\end{align*}
proving the first part of Proposition \ref{prop:pl}.

\textbf{Proof of (2)}. If $V$ is convex, then the negative of the exponent in \eqref{eq:rs} is a convex function of $(\phi_x)_{x \in C}$, and hence the associated density on $\mathbb{R}^A$ is log-concave. It follows that its marginal density on $\mathbb{R}^{|D|}$ associated with any subset $D$ of $A$ is also log-concave, as required.
\end{proof}

In the majority of the article we have considered the partition function $Z_A(\varphi)$ associated with integrating the interactions across a subset $A$ subject to boundary conditions from a function $\varphi:B \to \mathbb{R}$. A natural alternative form of partition function measures the integral of the interactions across the entire subset $C$ subject to staying within an $L_1$ distance $\delta > 0$ of a certain function. More specifically, given $\varphi:C \to \mathbb{R}$ and $\delta > 0$ consider the partition function of an $L_1$-type ball:
\begin{align} \label{eq:pfball}
\tilde{Z}_C( \varphi, \delta ) := \int_{\mathbb{R}^C} \mathrm{1} \left\{ \sum_{x \in C} | \phi_x - \varphi_x | \leq \delta \right\} \exp \left\{ - \sum_{ \langle x,y \rangle \in C^*} V(\phi_y - \phi_x) \right\} \prod_{x \in C} \mathrm{d}\phi_x.
\end{align}
The following result states that the partition functions for $L_1$ balls also satisfy the Pr\'ekopa--Leindler inequality. 

\begin{proposition} \label{prop:PLPF}
Given any functions $\varphi,\tilde{\varphi}:C \to \mathbb{R}$ and $\lambda \in (0,1)$ we have 
\begin{align*}
\tilde{Z}_C( \lambda \varphi+ (1-\lambda)\tilde{\varphi}, \delta ) \geq \tilde{Z}_C( \varphi, \delta )^\lambda \tilde{Z}_C( \tilde{\varphi}, \delta )^{1-\lambda}.
\end{align*}
\end{proposition}

\begin{proof}
With $J(\phi)$ as in \eqref{eq:Jdef}, set
\[
f_\varphi(\phi) := \mathrm{1} \left\{ \sum_{x \in C} | \phi_x - \varphi_x | \leq \delta \right\} J(\phi).
\]
Then $\tilde{Z}_C(\varphi,\delta) = \int_{\mathbb{R}^C} f_\varphi(\phi)\mathrm{d}\phi$. 
By the triangle inequality, if $
\sum_{x \in C} |\phi_x - \varphi_x | \leq \delta 
$ and $\sum_{x \in C} |\phi'_x - \tilde{\varphi}_x | \leq \delta$ 
then
\[
\sum_{x \in C} | \lambda \phi_x + (1-\lambda)\phi'_x - (\lambda \varphi_x + (1-\lambda)\tilde{\varphi}_x) | \leq \delta.
\]
Combining this with the convexity inequality for $J$, we obtain
\begin{align*}
f_{\lambda \varphi + (1-\lambda)\tilde{\varphi}}(\lambda \phi + (1-\lambda)\phi') 
\geq f_\varphi(\phi)^\lambda f_{\tilde{\varphi}}(\phi')^{1-\lambda}.
\end{align*}
Now apply the Pr\'ekopa--Leindler inequality.
\end{proof}

By taking logs and iterating the previous result, we see that whenever $(\varphi_j)_{j \in J}$ is a finite collection of functions $\varphi_j: C \to \mathbb{R}$ indexed by a finite set $J$, then if we define their average $\bar{\varphi} := |J|^{-1}\sum_{j \in J} \varphi_j$, we have 
\begin{align} \label{eq:smootha}
\log \tilde{Z}_C( \bar{\varphi}, \delta) \geq 
|J|^{-1} \sum_{j \in J} \log \tilde{Z}_C( \varphi_j,\delta).
\end{align}
The inequality \eqref{eq:smootha} will play a pivotal role in our proof of the large deviation upper bound in Theorem \ref{thm:C}. 
%%%%%%%%%%%%%%%%%%%%%%%%%%%%%%%%%%%%%%%%%%%%%%
%%%%%%%%%%%%%%%%%%%%%%%%%%%%%%%%%%%%%%%%%%%%%%
%%%%%%%%%%%%%%%%%%%%%%%%%%%%%%%%%%%%%%%%%%%%%%
\subsection{Basic properties of Gelfand--Tsetlin integrals} \label{sec:gtearnest}
%%%%%%%%%%%%%%%%%%%%%%%%%%%%%%%%%%%%%%%%%%%%%%
%%%%%%%%%%%%%%%%%%%%%%%%%%%%%%%%%%%%%%%%%%%%%%
%%%%%%%%%%%%%%%%%%%%%%%%%%%%%%%%%%%%%%%%%%%%%%

In this section we specialise to the hard-core interaction potential $V_{\mathrm{GT}}(\phi) := +\infty \mathrm{1}\{ \phi < 0 \}$. We emphasise that $V_{\mathrm{GT}}$ in \eqref{eq:hc} is a convex potential, and hence the conclusions of Proposition \ref{prop:pl} and Proposition \ref{prop:PLPF} hold for partition functions of Gelfand--Tsetlin surfaces on path-comparable sets $C$. It is also possible to show that these conclusions hold for non-path-comparable sets. 

Recall that for a general finite subset $C = A \sqcup B$ of $\mathbb{Z}^2$ with boundary conditions $\varphi:B \to \mathbb{R}$, the Gelfand--Tsetlin partition function is given by 
 \begin{align} \label{eq:pf2}
 Z_A(\varphi) = Z_A( \varphi_x : x \in B ) := \int_{\mathbb{R}^C}\mathrm{1}_{\mathrm{GT}}^C(\phi) \prod_{x \in B} \delta_{\varphi_x}(\mathrm{d}\phi_x) \prod_{x \in A} \mathrm{d}\phi_x,
 \end{align}
 where $ \mathrm{1}_{\mathrm{GT}}^C(\phi) = \prod_{x \leq y \in C} \mathrm{1}\{ \phi_x \leq \phi_y\}$ is the indicator function that $\phi:C \to \mathbb{R}$ is increasing. 
% In the majority of cases in the sequel $C$ will be path-comparable, and we will substitute $\mathrm{1}_{\mathrm{GT}}^C(\phi)$ for $\hat{\mathrm{1}}_{\mathrm{GT}}^C(\phi)$ in the integral without comment.
% Suppose that $C$ is a disjoint union $C = A \sqcup B$. Specialising the definition \eqref{eq:rs} to the interaction potential $V_{\mathrm{GT}}(x) := +\infty \mathrm{1} \{ x < 0\}$, a \textbf{random global Gelfand--Tsetlin surface} on $C$ with boundary conditions $\{ \varphi_x \,:\, x \in B \}$ is a random real-valued function $\{ \Phi_x : x \in C \}$ distributed according to a probability measure on $\mathbb{R}^C$ taking the form
% \begin{align} \label{eq:rs2}
% Z_A( \varphi)^{-1} \mathrm{1}_{\mathrm{GT}}^C(\phi) \prod_{x \in B} \delta_{\varphi_x}(\mathrm{d}\phi_x) \prod_{x \in A} \mathrm{d}\phi_x, 
% \end{align}
% provided the normalizing \textbf{global Gelfand--Tsetlin integral} $Z_A(\varphi)$, defined by

% is in $(0,\infty)$. If $C$ is path-comparable, then this definition is a special case of a Ginzburg-Landau $\nabla\phi$-interface model as defined in Section \ref{sec:randomsurfaces}.

We would like to identify precisely when $Z_A(\varphi)$ takes values in $(0,\infty)$. On the one hand, $Z_A(\varphi)$ is zero whenever $\varphi:B \to \mathbb{R}$ is not increasing on $B$. On the other hand, we say that the set $B$ bounds $C = A \sqcup B$ if for every $c \in C$, there exist $b, b' \in B$ such that $b \leq c \leq b'$. The integral $Z_A(\varphi)$ is finite whenever $B$ bounds $C$. A finite subset of $\mathbb{Z}^2$ is always bounded by its boundary (the set of points connected to its complement). However, a bounding set need not be a boundary: e.g., the diagonal set $D_n = \{(x_1,x_1) : 0 \leq x_1 \leq n \}$ bounds the triangle $\rt_n = \{ (x_1,x_2) : 0 \leq x_2 \leq x_1 \leq n\}$. 

In summary, $Z_A(\varphi) < \infty$ whenever $B$ bounds $C$. The precise conditions under which $Z_A(\varphi) > 0$ occurs are a bit more delicate. 
Given subsets $B$ and $C$ of $\mathbb{Z}^2$ such that $B$ bounds $C$, we say that a point $b$ in $B$ is a \textbf{partially lower} point if there exists $a \in A:= C-B$ such that $b \leq a$ and a \textbf{partially upper} point if there exists $a \in A$ such that $a \leq b$. We say that $b$ is a lower point if it is partially lower but not partially upper, and an upper point if it is partially upper but not partially lower. We say that $b$ is a mixed point if it is both partially lower and partially upper. We say that $b$ is an irrelevant point if it is neither partially lower nor partially upper.

\begin{example} \label{ex:upperlower}
Consider letting $C = \rt_n := \{ (x_1,x_2) : 0 \leq x_2 \leq x_1 \leq n \}$, and let $B = \partial \rt_n$ be its boundary. Then we can partition the boundary as follows:
\begin{itemize}
\item $S_n' := \{ (x_1,0) : 0 \leq x_1 \leq n-1\} \cup \{(1,1) \}$ are lower points.
\item $E'_n := \{ (n,x_2) : 1 \leq x_2 \leq n \} \cup \{(n-1,n-1)\}$ are upper points.
\item $D_n := \{ (x,x) : 2 \leq x \leq n-2 \}$ are mixed points.
\item $\{(n,0)\}$ is an irrelevant point.
\end{itemize}
\end{example}

\begin{lemma} \label{lem:pointtype}
Let $\varphi,\varphi':B \to \mathbb{R}$ be increasing. Suppose that $\varphi_x = \varphi_x'$ for every mixed point, $\varphi'_x \leq \varphi_x$ for every lower point, and  $\varphi'_x \geq \varphi_x$ for every upper point, then
\begin{align*}
Z_A(\varphi') \geq Z_A(\varphi).
\end{align*}
\end{lemma}
\begin{proof}
If the extension of a function $\phi:A \to \mathbb{R}$ to $C$ satisfies the Gelfand--Tsetlin inequalities with boundary conditions imposed by $\varphi:B \to \mathbb{R}$, it certainly satisfies them with the boundary conditions imposed by $\varphi':B \to \mathbb{R}$, which are easier to satisfy.
\end{proof}
Note that Lemma \ref{lem:pointtype} makes no stipulation about the relative values of $\varphi,\varphi'$ at any irrelevant points.
\medskip

We take a look at some of the basic scaling and translation properties of Gelfand--Tsetlin integrals $Z_A(\varphi)$.

\begin{lemma} \label{lem:GTprops} Let $C\subseteq \mathbb Z^2$ be finite, and assume that $C=A\sqcup B$ with $A\cap B=\emptyset$. Moreover, let $\varphi:B\to \mathbb R$ be a function. Then the following statements hold:
\begin{enumerate}
\item For any $\lambda > 0$, we have the homogeneity 
\begin{align} \label{eq:GTscal}
Z_A(\lambda \varphi) = \lambda^{ \# A} Z_A(\varphi).
\end{align}
\item Invariance under horizontal translations: for each $x\in \mathbb Z^2$, letting $x + C := \{ x + c \,:\, c \in C \}$ and defining the translated function $\varphi_{\to x}:(x+B) \to \mathbb{R}$ by $\varphi_{\to x}(x+b) := \varphi(b)$, we have
\begin{align*}
Z_{x+A}(\varphi_{\to x}) = Z_A(\varphi).
\end{align*}
\item Invariance under vertical translations: for each $c \in \mathbb{R}$, we have
\begin{align*}
Z_A( \varphi + c ) = Z_A(\varphi),
\end{align*}
where $(\varphi + c)(x) := \varphi(x) + c$. 
\end{enumerate}
\end{lemma}

\begin{proof}
Each statement follows from a simple change of variables.
\end{proof}

%%%%%%%%%%%%%%%%%%%%%%%%%%%%%%%%%%%%%%%%%
\subsection{Lower bounds for Gelfand--Tsetlin integrals via guaranteed space}
%%%%%%%%%%%%%%%%%%%%%%%%%%%%%%%%%%%%%%%%%

For $c \in \mathbb{R}$, we say that a function $\varphi:C \to \mathbb{R}$ defined on a subset $C$ of $\mathbb{Z}^2$ is $c$-increasing if and only if
\begin{align*}
\varphi_y - \varphi_x \geq c ( (y_1+y_2) - (x_1+x_2)) \qquad \text{for all $y \geq x \in C$}.
\end{align*}
For $\varphi:C \to \mathbb{R}$ we have the equivalent definitions
\begin{align*}
\text{$\varphi$ is increasing } \iff \text{$\varphi$ is $0$-increasing} \iff \mathrm{1}_{\mathrm{GT}}^C(\varphi) = 1.
\end{align*}
For $(u_1,u_2)\in\mathbb R^2$ with $u_1,u_2 > 0$, we let $\varphi^{u_1,u_2}:\mathbb{Z}^2 \to \mathbb{R}$ denote the linear function $\varphi^{u_1,u_2}(x_1,x_2) = u_1 x_1 + u_2 x_2$. The function $\varphi^{c,c}$ is $c$-increasing. The property that $\phi:C \to \mathbb{R}$ is $c$-increasing is equivalent to $\phi - \varphi^{c,c}$ being increasing. 
If $\varphi,\varphi':B \to \mathbb{R}$ are such that $\varphi$ is $c$-increasing, and $|\varphi'_x - \varphi_x| \leq r$ for all $x \in B$, then $\varphi'$ is $(c - 2r)$-increasing. 

If $B$ is a subset of $C$, we say that a function $\bar{\varphi}:C \to \mathbb{R}$ is an extension of $\varphi:B \to \mathbb{R}$ if $\bar{\varphi}(b) = \varphi(b)$ for all $b \in B$. 

\begin{lemma}[Guaranteed space lemma] \label{lem:space}
Let $C\subseteq \mathbb Z^2$ be finite, $B \subseteq C$ and $\varphi:B \to \mathbb{R}$ be a function. Set $A = C-B$. Suppose that for $c \in [0,\infty)$, $\varphi$ has an extension $\bar{\varphi}:C\to \mathbb R$ to all of $C$ that is $c$-increasing. Then 
\begin{align*}
Z_A(\varphi) \geq c^{ \# A}.
\end{align*}
\end{lemma}

\begin{proof}
Let $\bar{\varphi}:C\to \mathbb R$ be the extension of the function $\varphi:B\to \mathbb R$. Consider the set
\[
  E := \big\{ \phi: C \to\mathbb R\,:\, \forall x \in C:\, | \phi_x - \bar{\varphi}_x| \leq c/2 \big\} \subseteq \mathbb R^C.
\]   
Then every $\phi$ in $E$ is increasing, since for $x \leq y$ in $C$ with $y \neq x$ we have 
  \[
    \phi_y - \phi_x \geq \Big(\bar{\varphi}_y - \frac{c}{2}\Big) - \Big(\bar{\varphi}_x + \frac{c}{2}\Big) \geq 0;
  \]
where in the last step we used $\bar{\varphi}_y - \bar{\varphi}_x \geq c((y_1 +y_2) - (x_1+x_2)) \geq c$. It follows that
\begin{align*}
Z_A(\varphi)  = \int_{\mathbb{R}^C} \mathrm{1}_{\mathrm{GT}}^C(\phi) \prod_{x \in B} \delta_{\varphi_x}(\mathrm{d}\phi_x) \prod_{x \in A} \mathrm{d} \phi_x  \geq \int_E \prod_{x \in B} \delta_{\varphi_x}(\mathrm{d}\phi_x) \prod_{x \in A} \mathrm{d} \phi_x = c^{\# A},
\end{align*}
as required.
\end{proof}

Our next result says that $c$-increasing functions on boundaries of sets may always be extended to their interiors; in particular, in this situation we always have the lower bound on the Gelfand--Tsetlin integral established in the previous lemma.

\begin{lemma}[Extension lemma] \label{lem:extension}
Let $C$ be a finite subset of $\mathbb{Z}^2$, let $B$ be a subset of $C$ that bounds $C$, and let $\varphi:B \to \mathbb{R}$ be a $c$-increasing function. Then $\varphi$ has a $c$-increasing extension to all of $C$. In particular, if $A := C\setminus B$ then
\begin{align*}
Z_A(\varphi) \geq c^{ \# A}.
\end{align*}
\end{lemma}

\begin{proof}
Recall that a function $\varphi$ is $c$-increasing if and only if $\varphi-\varphi^{c,c}$ is increasing. Thus, it is sufficient to establish that every increasing $\varphi:B \to \mathbb{R}$ has an extension to $C$ that is also increasing. To this end, we simply define an extension $\bar{\varphi}:C \to \mathbb{R}$ by setting, for $x \in C$, 
\begin{align*}
\bar{\varphi}_x := \max_{y \leq x,\, y \in B} \varphi_y;
\end{align*}
since $B$ bounds $C$, there is certainly some element $y \in B$ less than $x \in C$. 

Then $\bar{\varphi}$ is increasing and agrees with $\varphi$ on $B$.
% , since increasing $x$ increases the set over which the maximum is taken on 
% since for $x,y\in C$ with $x\leq y$, we have
%   \[
%     \bar{\varphi}_y - \bar{\varphi}_x = \max_{v \leq y,\, v \in B} \varphi_v - \max_{w \leq x,\, w \in B} \varphi_w \geq \max_{v \leq y,\, v \in B} \varphi_v - \max_{w \leq y,\, w \in B} \varphi_w = 0.
%   \]
% Moreover $\bar{\varphi}_x = \varphi_x$ for $x \in B$, because $\varphi$ is increasing on $B$. 
The result now follows from Lemma \ref{lem:space}.
\end{proof}

%%%%%%%%%%%%%%%%%%%%%%%%%%%%%%%%%%%%%%%%
\subsection{Upper bounds for Gelfand--Tsetlin integrals} \label{sec:GTupper}
%%%%%%%%%%%%%%%%%%%%%%%%%%%%%%%%%%%%%%%%

Lemma \ref{lem:extension} provides a lower bound on Gelfand--Tsetlin integrals. Here we develop some upper bounds. In this direction, note that if $E = E_1 \cup \cdots \cup E_k$, we have the union inequality
\begin{align} \label{eq:unioneq}
\mathrm{1}_{\mathrm{GT}}^E(\phi) \leq \prod_{i=1}^k \mathrm{1}_{\mathrm{GT}}^{E_i}(\phi).
\end{align} 

The following result is an upper bound on certain Gelfand--Tsetlin integrals.

\begin{lemma} \label{lem:volupper}
Let $C$ be a subset of $\Box_n := \{0,1,\ldots,n\}^2$ with cardinality $\# C = m$. Let $a < b$ in $\mathbb{R}$. Then we have the upper bound
\begin{align*}
\int_{[an,bn]^C} \mathrm{1}_{\mathrm{GT}}^C(\phi) \prod_{x \in C} \mathrm{d}\phi_x 
\leq \exp \left\{  m \left( 2 + \log(b-a) + \log \left( \frac{(n+1)^2}{m } \right) \right) \right\}.
\end{align*}
\end{lemma}

\begin{proof}
By making a change of variable, we may assume without loss of generality that $b = 1$ and $a=0$. 

Consider the diagonals $D_k := \{(i,j) \in \Box_n : i - j = k \}$. Then $\Box_n = \bigcup_{ -n \leq k \leq n } D_k$, and writing $C_k := C \cap D_k$, we have $C = \bigcup_{ - n\leq k \leq n } C_k$. In particular, using \eqref{eq:unioneq} we have
\begin{align*}
\int_{[0,n]^C} \mathrm{1}_{\mathrm{GT}}^C(\phi) \prod_{x \in C} \mathrm{d}\phi_x 
\leq \int_{[0,n]^C} \left\{ \prod_{- n \leq k \leq n } \mathrm{1}_{\mathrm{GT}}^{C_k}(\phi) \right\} \prod_{x \in C} \mathrm{d}\phi_x  = \prod_{ - n \leq k \leq n } \int_{[0,n]^{C_k}} \mathrm{1}_{\mathrm{GT}}^{C_k}(\phi) \prod_{x \in C_k} \mathrm{d}\phi_x.
\end{align*}

Let $m_k := \# C_k$. Since $C_k$ is totally ordered, we have 
\begin{align*}
\int_{[0,n]^{C_k}} \mathrm{1}_{\mathrm{GT}}^{C_k}(\phi) \prod_{x \in C_k} \mathrm{d}\phi_x 
= \int_{0 \leq s_1 \leq \cdots \leq s_{m_k} \leq n } \mathrm{d}s_1 \cdots \mathrm{d}s_{m_k} 
= \frac{n^{m_k}}{m_k!}.
\end{align*}

It follows that
\begin{align*}
\int_{[0,n]^C} \mathrm{1}_{\mathrm{GT}}^C(\phi) \prod_{x \in C} \mathrm{d}\phi_x 
\leq \prod_{-n \leq k \leq n} \frac{n^{m_k} }{m_k!} 
= n^m \prod_{-n \leq k \leq n} \Gamma(m_k + 1)^{-1}. 
\end{align*}

The function $a \mapsto \log \Gamma(a+1)$ is convex. It follows that given real numbers $(a_k)_{- n \leq k \leq n}$ summing to $m$, the product $\prod_{-n \leq k \leq n} \Gamma(a_k+1)$ is minimised when $a_k = m/(2n+1)$ for each $k$. 
Therefore,
\begin{align*}
\int_{[0,n]^C} \mathrm{1}_{\mathrm{GT}}^C(\phi) \prod_{x \in C} \mathrm{d}\phi_x 
\leq n^{m} \Gamma\left(  \frac{m}{2n+1} +1 \right)^{-(2n+1)}.
\end{align*}

Using the lower bound $\Gamma(x+1) \geq (x/e)^x$ for all $x \geq 0$, and being generous with constants, we obtain
\begin{align*}
\int_{[0,n]^C} \mathrm{1}_{\mathrm{GT}}^C(\phi) \prod_{x \in C} \mathrm{d}\phi_x 
\leq n^m \left( \frac{m}{e(2n+1)} \right)^{- m} 
% &\leq \left( \frac{2e(n+1)^2}{m} \right)^m 
\leq \exp \left\{ m \left( 2 + \log \left( \frac{(n+1)^2}{m} \right) \right)   \right\},
\end{align*}
completing the proof.
\end{proof}

%%%%%%%%%%%%%%%%%%%%%%%%%%%%%%%%%%%%%%%%
\subsection{The triangular version of Theorem \ref{thm:B}} \label{sec:triangular}
%%%%%%%%%%%%%%%%%%%%%%%%%%%%%%%%%%%%%%%%

We start preparing the proof of Theorem \ref{thm:B}, which describes the surface tension of Gelfand--Tsetlin random surfaces. 
For reasons that will become clear in the following sections, rather than proving Theorem \ref{thm:B} directly, it is useful to begin by proving a version of this result on triangular --- rather than square --- interfaces. 

For concreteness, if $C$ is a subset of $\mathbb{Z}^2$, we define its boundary and interior by 
\begin{align*}
\partial C := \{ x \in C : \exists y \in \mathbb{Z}^2 \setminus C : \text{$y$ and $x$ are neighbours in $\mathbb{Z}^2$} \}, 
\qquad 
\interior{C} := C \setminus \partial C.
\end{align*}

We will consider the square and right-angled triangle of side-length $n$, which are given by 
\begin{align*}
\Box_n &:=  \{ (x_1,x_2) \in \mathbb{Z}^2 : 0 \leq x_1,x_2 \leq n  \},\\
\righttriangle_n &:= \{ (x_1,x_2) \in \mathbb{Z}^2 : 0 \leq x_2 \leq x_1 \leq n  \}.
\end{align*}

The associated square and triangular Gelfand--Tsetlin surface partition functions at slope $(u_1,u_2) \in \mathbb{R}^2$ are given by  
\begin{align}
    S_n(u_1,u_2) &:= Z_{\interior{\Box}_n}( \varphi^u_x : x \in \partial \Box_n ) := \int_{\mathbb{R}^{\Box_n} } \mathrm{1}_{\mathrm{GT}}^{\Box_n}(\phi) \prod_{x \in \partial \Box_n} \delta_{u_1 x_1 + u_2 x_2 }(\mathrm{d}\phi_x) \prod_{x \in \interior{\Box}_n} \mathrm{d}\phi_x, \label{eq:Sn}\\
T_n(u_1,u_2) &:= Z_{\interior{\righttriangle}_n}( \varphi^u_x : x \in \partial \righttriangle_n ) := \int_{\mathbb{R}^{\righttriangle_n} }\mathrm{1}_{\mathrm{GT}}^{\righttriangle_n}(\phi)\prod_{x \in \partial \righttriangle_n} \delta_{u_1 x_1 + u_2 x_2 }(\mathrm{d}\phi_x) \prod_{x \in \interior{\righttriangle}_n} \mathrm{d}\phi_x.\label{eq:Tn}
\end{align}

Recall that Theorem \ref{thm:B} states that $\lim_{n \to \infty} \frac{1}{n^2} \log S_n(u_1,u_2) = - \sigma(u_1,u_2)$ where $\sigma(u_1,u_2)$ is as in \eqref{eq:beadtension}. Over the next few sections, we will work first to establish the following triangular version of Theorem \ref{thm:B} with absolute bounds.

\begin{thm} \label{thm:linear}
For all $u_1,u_2>0$ and all integers $n \geq 2$ satisfying $(\log n)^{-1/2}  \leq u_1 \wedge u_2$ we have 
\begin{align*}
\left| \frac{1}{n^2/2} \log T_n(u_1,u_2) + \sigma(u_1,u_2) \right| \leq C \frac{\log n  + K_u}{n},
\end{align*}
where, as in \eqref{eq:beadtension}, 
\begin{align} \label{eq:st2}
\sigma(u_1,u_2):=- \log(u_1+u_2) - \log \sin \frac{ \pi u_1}{u_1+u_2} -1 + \log \pi.
\end{align}
Here $K_u = |\log (u_1 \wedge u_2) | +|\log(u_1 + u_2)|$ is a $(u_1,u_2)$-dependent constant. 
\end{thm}
%%%%%%%%%%%%%%%%%%%%%%%%%%%%
\subsection{Further properties of $T_n(u_1,u_2)$}
%%%%%%%%%%%%%%%%%%%%%%%%%%%%

In this section we look at basic properties of triangular integrals; recall the definition in \eqref{eq:Tn}. First we note that it is possible to give an easy lower bound using the extension lemma, Lemma \ref{lem:extension}.

\begin{lemma} \label{lem:charac}
Let $n\in\mathbb N$, $u_1,u_2\in\mathbb R$ with $u_1,u_2 > 0$. Then we have 
  \[
    T_n(u_1,u_2) \geq e^{ - CK_u n^2}
  \]  
\end{lemma}

\begin{proof}
Note $\varphi^u$ is $c$-increasing with $c = u_1 \wedge u_2$.  
% Let $x:=(x_1,x_2),y:=(y_1,y_2)\in \righttriangle_n$ with $x\leq y$. Consider $\varphi^{u_1,u_2}:\righttriangle_n \to \mathbb R$. Then 
%   \[
%     \varphi^{u_1,u_2}_y - \varphi^{u_1,u_2}_x = u_1(y_1-x_1) + u_2(y_2-x_2) \geq (u_1 \wedge u_2)\big[(y_1+y_2)-(x_1+x_2) \big],
%   \]
% and so $\varphi^{u_1,u_2}$ is $c$-increasing for $c:=u_1 \wedge u_2$. 
We use Lemma \ref{lem:extension} together with the fact that $\# \interior{\righttriangle}_n = (n-1)(n-2)/2$ to obtain $T_n(u_1,u_2) \geq (u_1 \wedge u_2)^{(n-1)(n-2)/2}$. Now use the definition of $K_u$. 
\end{proof} 
Our next result describes further properties of $T_n(u_1,u_2)$.

\begin{lemma} \label{lem:st}
Let $n\in\mathbb N$ and $u_1,u_2\in\mathbb R$ with $u_1,u_2 > 0$. Then $T_n(u_1,u_2)$ has the following properties:
\begin{enumerate}
\item $T_n$ is symmetric in $u_1$ and $u_2$.
\item $T_n$ is increasing in both $u_1$ and $u_2$. 
\item For every $\lambda>0$, we have the homogeneity
\begin{align} \label{eq:scaling}
T_n(\lambda u_1,\lambda u_2) = \lambda^{(n-1)(n-2)/2} T_n(u_1,u_2).
\end{align}
\end{enumerate}
\end{lemma}

\begin{proof}
The fact that $T_n$ is symmetric may be established by considering the change of variable sending $\phi:\righttriangle_n \to \mathbb{R}$ to $\phi':\righttriangle_n \to \mathbb{R}$ given by
\[
\phi'_{x_1,x_2} = (u_1 + u_2)n - \phi_{n-x_2,n-x_1}.
\]
This change of variable has unit Jacobian, and transfers a function on $\righttriangle_n$ agreeing with $\varphi^{u_1,u_2}$ on the boundary to one agreeing with $\varphi^{u_2,u_1}$ on the boundary.

As for the second point, consider $h_1,h_2 \geq 0$. If $\phi:\righttriangle_n \to \mathbb{R}$ agrees with $\varphi^{u_1,u_2}$ on the boundary, then the function $\phi':\righttriangle_n \to \mathbb{R}$ defined by
\[
\phi'_x := \phi_x + h_1 x_1 + h_2 x_2
\]
agrees with $\varphi^{u_1+h_1,u_2+h_2}$ on the boundary. Moreover, since $x \leq y$ implies $x_1 \leq y_1$ and $x_2 \leq y_2$, we have
\[
\phi_y \geq \phi_x \quad \Longrightarrow \quad \phi'_y \geq \phi'_x,
\]
so that $\mathrm{1}_{\mathrm{GT}}^{\righttriangle_n}(\phi') \geq \mathrm{1}_{\mathrm{GT}}^{\righttriangle_n}(\phi)$. Integrating, it follows that
\[
T_n(u_1+h_1,u_2+h_2) \geq T_n(u_1,u_2),
\]
and hence $T_n$ is increasing in both $u_1$ and $u_2$.

The final point is a special case of Lemma \ref{lem:GTprops}(1).
\end{proof}

% Of course, by the homogeneity given in \eqref{eq:scaling}, in order to compute $T_n(u_1,u_2)$, or prove Theorem \ref{thm:linear} for that matter, it is sufficient to consider the special case $u_1 + u_2=1$. 

In the next section we introduce several explicit formulas from \cite{johnston} for the bead model on the torus. We will ultimately use a determinantal formula to prove an upper bound for $T_n(u_1,u_2)$, and then a concentration result for log-concave density functions to prove a lower bound.

%%%%%%%%%%%%%%%%%%%%%%%%%%%%%%%%%%%%%%%%%%%%%%
%%%%%%%%%%%%%%%%%%%%%%%%%%%%%%%%%%%%%%%%%%%%%%
%%%%%%%%%%%%%%%%%%%%%%%%%%%%%%%%%%%%%%%%%%%%%%
\section{Proof of Theorem \ref{thm:linear}} \label{sec:import}
%%%%%%%%%%%%%%%%%%%%%%%%%%%%%%%%%%%%%%%%%%%%%%
%%%%%%%%%%%%%%%%%%%%%%%%%%%%%%%%%%%%%%%%%%%%%%
%%%%%%%%%%%%%%%%%%%%%%%%%%%%%%%%%%%%%%%%%%%%%%

%%%%%%%%%%%%%%%%%%%%%%%%%%
\subsection{Explicit formulas from \cite{johnston}} \label{sec:volumes}
%%%%%%%%%%%%%%%%%%%%%%%%%%

In the present section we import several results from \cite{johnston} which we later use to prove Theorem \ref{thm:linear}. We begin by defining bead configurations on the torus. For $n\in\mathbb N$ let $\mathbb{T}_n := [0,1) \times \{0,1,\ldots,n-1\}$ denote the semi-discrete torus of $n$ unit-length parallel strings. For any $k \in \mathbb N$, a bead configuration on the semi-discrete torus $\mathbb{T}_n$ is a collection of $nk$ distinct points $(y_1,\ldots,y_{nk})$ on $\mathbb{T}_n$ (i.e., $y_i = (t_i,h_i)$ for each $1\leq i \leq nk$) such that there are $k$ points on each string $[0,1) \times \{h\}$, $h\in \{0,1,\ldots,n-1\}$, and the $k$ points on neighbouring strings interlace. More specifically, for each $h \in \mathbb{Z}_n:=\mathbb Z/n\mathbb Z$, if $t_1 < \ldots < t_k$ and $t_1' < \cdots < t_k'$ are the horizontal coordinates of the points on strings $h$ and $h+1$ (mod $n$) respectively, then we have either 
\begin{align} \label{eq:doubleinterlace}
t_1 \le t_1' < \cdots < t_k \le t_k' \qquad \text{or} \qquad t_1' < t_1 \le \cdots \le t_k' < t_k.
\end{align}
In other words, between every two beads on one string there is a bead on the neighbouring string. See Figure \ref{fig:test0} for an example of a bead configuration on $\mathbb{T}_n$.

\begin{figure}[h!]
\centering\begin{tikzpicture}[scale=0.8]
\draw (0,0) -- (6,0);
\draw (0,1) -- (6,1);
\draw (0,2) -- (6,2);
\draw (0,3) -- (6,3);
\draw (0,4) -- (6,4);
\draw[lightgray] (0,0) -- (0,4);
\draw[lightgray] (6,0) -- (6,4);

\draw [very thick, CadetBlue] (0.2,0) circle [radius=0.1];
\draw [very thick, CadetBlue] (2.1,0) circle [radius=0.1];
\draw [very thick, CadetBlue] (4.31,0) circle [radius=0.1];
\draw [very thick, CadetBlue] (0.95,1) circle [radius=0.1];
\draw [very thick, CadetBlue] (2.7,1) circle [radius=0.1];
\draw [very thick, CadetBlue] (5.1,1) circle [radius=0.1];
\draw [very thick, CadetBlue] (1.75,2) circle [radius=0.1];
\node at (1.75,1.6) {$y_i$};
\draw [very thick, CadetBlue] (3.85,2) circle [radius=0.1];
\draw [very thick, CadetBlue] (5.6,2) circle [radius=0.1];
\draw [<->, fill, CadetBlue] (1.75,2.2) -- (3.85,2.2);
\node at (2.8,2.4) {$p_i$};
\draw[thick, dotted, CadetBlue] (1.75, 2) -- (1.75, 3.5);
\draw [very thick, CadetBlue] (0.6,3) circle [radius=0.1];
\draw [<->, fill, CadetBlue] (1.75,3.2) -- (2.7,3.2);
\node at (2.25,3.4) {$q_i$};
\draw [very thick, CadetBlue] (2.7,3) circle [radius=0.1];
\draw [very thick, CadetBlue] (4.75,3) circle [radius=0.1];
\draw [very thick, CadetBlue] (1.0,4) circle [radius=0.1];
\draw [very thick, CadetBlue] (3.3,4) circle [radius=0.1];
\draw [very thick, CadetBlue] (5.5,4) circle [radius=0.1];
\end{tikzpicture} 
\caption{A bead configuration on $n = 5$ strings with $k = 3$ beads per string. A particular bead $y_i$ is highlighted, as are the respective distances $p_i$ and $q_i$ to the next bead on the same string, and the next bead on the next string. For topological reasons, the average of the $q_i$ divided by the average of the $p_i$ takes the form $\ell/n$ for some integer $ 1 \leq \ell \leq n-1$.}
\label{fig:test0}
\end{figure}

To every bead configuration on $\mathbb{T}_n$ with $nk$ beads, we can associate a parameter $\gamma$ of the form $\gamma := \ell/n$ with $1 \leq \ell \leq n-1$ which we call the tilt of the configuration. (A notational warning: the tilt is denoted by $\tau$ in \cite{johnston}.) To define the tilt of a bead configuration $(y_1,\ldots,y_{nk})$, given a bead at $y_i$, let $p_i$ denote the horizontal distance to the next bead on the same string, and let $q_i$ denote the horizontal distance to the next bead on the string above (see Figure \ref{fig:test0}); if $y_i$ is the rightmost bead on a string, take these distances torically. 
The \textbf{tilt} $\gamma \in (0,1)$ of a bead configuration is the quotient of the averages of the $q_i$ and of the $p_i$:
\begin{align} \label{eq:tiltdef}
\gamma := \frac{ \frac{1}{nk} \sum_{i=1}^{nk} q_i  }{\frac{1}{nk} \sum_{i=1}^{nk} p_i  }.
\end{align}
In fact, since there are $k$ beads on each string and because distances are taken torically, the denominator $\frac{1}{nk} \sum_{i=1}^{nk} p_i$ in \eqref{eq:tiltdef} is equal to $1/k$. The fact that $0 \leq q_i \leq p_i$ for each $i$ guarantees $\gamma \in [0,1]$. It is possible to show that the tilt takes the form $\gamma = \ell/n$ for some integer $\ell \in \{1,\ldots,n-1\}$ \cite[Proof of (1.5)]{johnston} and we call $\ell$ the \textbf{occupation number} of the bead configuration. 

For $k,n\in\mathbb N$ and $\ell\in\{1,\dots,n-1\}$, the set of $(n,k,\ell)$ configurations, i.e., bead configurations on $\mathbb{T}_n$ with $nk$ beads and tilt $\gamma=\ell/n$, can be associated with a subset 
\[
\mathcal{W}^{(n)}_{k,\ell} := \Big\{(t_{h,j})_{0 \leq h \leq n-1,\,1 \leq j \leq k} \in [0,1)^{nk} \,:\, \text{$(n,k,\ell)$ configuration} \Big\}
\]
of $[0,1)^{nk}$ by letting $t_{h,j}$ denote the position of the $j^{\text{th}}$ bead on string $h$. Consequently, we can speak of the $nk$-dimensional Lebesgue measure of $\mathcal{W}_{k,\ell}^{(n)}$ and write 
\begin{align} \label{eq:voldef}
\mathrm{Vol}_{k,\ell}^{(n)} :=  \mathrm{Leb}_{nk}( \mathcal{W}_{k,\ell}^{(n)} ), \qquad k \in\mathbb N,\, 1 \leq \ell \leq n-1,
\end{align}
for the volume of the set of $(n,k,\ell)$ configurations. Also set $\mathrm{Vol}_{0,\ell}^{(n)} = \binom{n}{\ell}$, and $\mathrm{Vol}_{k,0}^{(n)} = \mathrm{Vol}_{k,n}^{(n)} = 0$.

The rescaled quantity $n^{n^2}\mathrm{Vol}^{(n)}_{n,\ell}$ denotes the volume of bead configurations on $n$ strings of length $n$ with $n$ beads per string and tilt $\gamma = \ell/n$. Of course, then in such a configuration the average distance between consecutive beads on the same string is $1$. 

An exact formula for $\mathrm{Vol}_{k,\ell}^{(n)}$ is given in \cite{johnston}. The asymptotics of this formula when $k = n$ are given by the following result.

\begin{thm}[\cite{johnston}] \label{thm:bead}
Let $n^{-1/4} \leq \gamma = \ell/n \leq 1 - n^{-1/4}$. Then the $n^2$-dimensional Lebesgue measure $\mathrm{Vol}_{n,\ell}^{(n)}$ satisfies
\begin{align} \label{eq:ordered9}
 \left( \frac{ e \sin( \pi \ell/n) }{ \pi } \right)^{-n^2}n^{n^2} \mathrm{Vol}^{(n)}_{n,\ell} = \frac{ e^{ \pi^2/6}}{ \sqrt{2 \pi }}\left( P( e^{-q^-_\gamma}) P( e^{-q_\gamma^+} ) + O\left(e^{ - c\gamma \log(n) }\right) \right),
\end{align}
where $P(s) := \prod_{i \geq 1}(1-s^i)^{-1}$ and $q_\gamma^\pm := 2 \pi^2 \pm 2 \pi^2 \iota \frac{\cos(\pi\gamma )}{ \sin(\pi \gamma)}$ and $c > 0$ is a universal constant.
\end{thm}

Recall from \eqref{eq:sigma} that $\sigma(u_1,u_2) = - \log(u_1+u_2) - \log \sin ( \pi u_1/(u_1+u_2)) - 1 + \log \pi$, so that in particular $\sigma(\gamma,1-\gamma) = - \log ( e \sin(\pi \gamma)/\pi)$. 

We will extract the following bound as a corollary:
\begin{cor}
Whenever $\log(n)^{-1/2} \leq \gamma = \ell/n \leq 1 - \log(n)^{-1/2}$ we have
\begin{align} \label{eq:ordered10}
\frac{1}{n^2} \log \left\{ n^{n^2} \mathrm{Vol}^{(n)}_{n,\ell} \right\} = - \sigma(\gamma,1-\gamma) + O(1/n^2).
\end{align}
\end{cor}

\begin{proof}
Note that for any $\gamma \in (0,1)$ we have $|e^{-q_\gamma^\pm}| = e^{-2\pi^2}$, and hence $P(e^{-q_\gamma^\pm})$ is bounded above and below uniformly in $\gamma$. In particular, there are constants $0 < C_1 < C_2$ such that provided $\log(n)^{-1/2} \leq \gamma \leq 1 - \log(n)^{-1/2}$ we have
\begin{align*}
0 < C_1 \leq \left( P( e^{-q^-_\gamma}) P( e^{-q_\gamma^+} ) + O\left(e^{ - c\gamma \log(n) }\right) \right) \leq C_2.
\end{align*}

It follows that taking logarithms of both sides of \eqref{eq:ordered9} and dividing by $n^2$, we obtain the result. 
\end{proof}

We turn to describing the correlation structure of \emph{random} bead configurations on the torus. For $\lambda \in \mathbb{R}$ and $T \geq 0$, define the partition function
\begin{align} \label{eq:pf000}
Z(\lambda,T) := \sum_{k \geq 0} \sum_{ 0 \leq \ell \leq n} T^{nk} e^{-\lambda \ell}\mathrm{Vol}^{(n)}_{k,\ell}.
\end{align}
Theorem 1.2 of \cite{johnston} states that we have the explicit expression
\begin{align} \label{eq:DA1}
    Z(\lambda,T) = \sum_{\theta =(\theta_1,\theta_2) \in \{0,1\}^2 } Z_\theta(\lambda,T), \qquad Z_\theta(\lambda,T) := \frac{1}{2}(-1)^{(\theta_1+1)(\theta_2+n+1)} \prod_{ z^n = (-1)^{\theta_2} } (e^{Tz} - (-1)^{\theta_1} e^{ - \lambda} ),
\end{align}
where the product is over all $n$ complex solutions of the equation $z^n = (-1)^{\theta_2}$.

We will also require a result describing the correlations of random bead configurations. 
Under a probability measure $\mathbf{P}^{(n),k,\ell}$, consider a uniformly chosen bead configuration with $k$ beads per string and tilt $\ell/n$. Following Section 1.4 of \cite{johnston}, let us also define a probability measure $\mathbf{P}^{(n),\lambda,T}$ governing bead configurations with a random number of beads per string and a random tilt by taking a weighted sum of the probability measures $\mathbf{P}^{(n),k,\ell}$:
\begin{align} \label{eq:weightsum}
\mathbf{P}^{(n),\lambda,T}(\cdot) = \sum_{k \geq 0} \sum_{ 0 \leq \ell \leq n} \frac{T^{nk}e^{-\lambda \ell} \mathrm{Vol}_{k,\ell}^{(n)}}{ Z(\lambda,T)} \mathbf{P}^{(n),k,\ell}(\cdot).
\end{align}
For $m \in\mathbb N$ let $\mathbb{T}_n^{[m]}$ denote the collection of $m$-tuples of distinct points in $\mathbb{T}_n$. 
Given a probability measure $\mathbf{P}$ on random bead configurations on $\mathbb{T}_n$, we define its $m^{\text{th}}$ correlation function to be the function $\rho_m:\mathbb{T}_n^{[m]} \to [0,\infty]$ defined by the limit
\begin{align} \label{eq:corrdef}
\rho_m(y_1,\ldots,y_m) := \lim_{\varepsilon \downarrow 0} \frac{1}{\varepsilon^m} \mathbf{P}\left( \text{For each $i=1,\ldots,m$, the interval $[y_i,y_i +\varepsilon \mathbf{e}_1)$ contains a bead}  \right),
\end{align}
(should the limit exist), 
where if $y = (t,h) \in \mathbb{T}_n = [0,1) \times \mathbb{Z}_n$, the interval $[y,y+\varepsilon \mathbf{e}_1) := \{ (s,h) : t \leq s < t+\varepsilon \}$ is an interval along the $h^{\text{th}}$ string of $\mathbb{T}_n$. 

Let $\rho_m^{(n),\lambda,T}$ and $\rho_m^{(n),k,\ell}$ denote the respective correlation functions of $\mathbf{P}^{(n),\lambda,T}$ and $\mathbf{P}^{(n),k,\ell}$. Then \eqref{eq:weightsum} tells us that for every $m \geq 1$ we have the relation 
\begin{align} \label{eq:weightsum2}
\rho^{(n),\lambda,T}_m(y_1,\ldots,y_m) = \sum_{k \geq 0} \sum_{ 0 \leq \ell \leq n} \frac{T^{nk}e^{-\lambda \ell} \mathrm{Vol}_{k,\ell}^{(n)}}{ Z(\lambda,T)} \rho^{(n),k,\ell}_m(y_1,\ldots,y_m)
\end{align}
between correlation functions. 

We now give a slightly abbreviated account of \cite[Theorem 1.7]{johnston}, which states that the correlation functions $\rho_m^{(n),\lambda,T}:\mathbb{T}_n^{[m]} \to [0,\infty]$ associated with $\mathbf{P}^{(n),\lambda,T}$ have a determinantal structure:

\begin{thm}[\cite{johnston}] \label{thm:detnew}
Write $d_\theta(n,\lambda,T) := Z_\theta(\lambda,T)/Z(\lambda,T)$ (as in \eqref{eq:DA1}). 
Then for all points $y_1,\ldots,y_m \in \mathbb{T}_n = [0,1) \times \mathbb{Z}_n$ with $y_i = (t_i,h_i)$, we have
%$\mathbf{P}^{(n),\lambda,T}$ 
\begin{align} \label{eq:Psum0}
\rho_m^{(n),\lambda,T}(y_1,\ldots,y_m)  = \sum_{ (\theta_1,\theta_2) \in \{0,1\}^2 } d_\theta(n,\lambda,T) \det_{i,j=1}^m H^{\lambda+\theta_1\pi \iota , \theta_2, T} ( y_j - y_i),
\end{align}
where for $t \in (-1,1), h \in \mathbb{Z}$, we have
\begin{align} \label{eq:detnewkernel}
H^{\beta,\theta_2,T}(t,h) =  \frac{T}{n} \sum_{ z^n = (-1)^{\theta_2} } z^{1-h} \frac{ e^{ - (\beta+Tz)[t] }}{ 1 - e^{ - (\beta+Tz)} } \qquad [t] := t + \mathrm{1}_{\{t<0\}},
\end{align}
and where $H^{\beta,\theta_2,T}$ is defined for $\beta \in \mathbb{C}$, $T\in(0,\infty)$ such that $\beta+Tz \notin 2\pi \iota \mathbb{Z}$ for $z \in \mathbb{C}$ where $z^n = (-1)^{\theta_2}$.
\end{thm}

That completes our import of formulas from \cite{johnston}. 

%%%%%%%%%%%%%%%%%%%%%%%%%%%%
\subsection{Unravelling the torus} \label{sec:unravel}
%%%%%%%%%%%%%%%%%%%%%%%%%%%%

% In order to prove Theorem \ref{thm:linear} we need to find a way of relating bead configurations on the semi-discrete torus $\mathbb{T}_n$ to triangular Gelfand--Tsetlin surfaces. 
The main idea of this section is to present a way of relating each bead configuration on the torus $\mathbb{T}_n$ with $n$ beads per string and tilt $\gamma=\ell/n$ to a \emph{pair} of Gelfand--Tsetlin surfaces defined on triangles in $\mathbb{Z}^2$ with linear boundaries and tilt $\ell/n$.
%\textcolor{red}{(We have not defined a 'tilt' for GT patterns, i.e., this part remains unclear at this stage. Probably we should write it differently. I'd just leave out the tilt part of the pattern and stop at 'linear boundaries'.)}. 
%% SJ: Agreed!

\begin{figure}[h!]
\centering
\begin{tikzpicture}[scale=0.8]
\draw (0,0) -- (8,0);
\draw (0,1) -- (8,1);
\draw (0,2) -- (8,2);
\draw (0,3) -- (8,3);
\draw (0,4) -- (8,4);
\draw[lightgray] (0,0) -- (0,4);
\draw[lightgray] (8,0) -- (8,4);

\draw [very thick, BrickRed] (0.2,0) circle [radius=0.1];
\node at (0.2,-0.3)   (a) {$X_{0,0}$};

\draw [very thick, BrickRed] (1.8,0) circle [radius=0.1];
\node at (1.8,-0.3)   (a) {$X_{1,1}$};
\draw [very thick, BrickRed] (2.7871,0) circle [radius=0.1];
\node at (2.7871,-0.3)   (a) {$X_{2,2}$};

\draw [very thick, FadedBrickRed]
(2.7871,5) circle [radius=0.1];
\node at (2.7871,4.7)   (a) {$X_{0,5}$};

\draw [very thick, BrickRed] (2.7871,0) circle [radius=0.1];

\draw [very thick, BrickRed] (4.38,0) circle [radius=0.1];
\node at (4.38,-0.3)   (a) {$X_{3,3}$};

\draw [very thick, BrickRed] (7.1,0) circle [radius=0.1];
\node at (7.1,-0.3)   (a) {$X_{4,4}$};

\draw [very thick, FadedBrickRed] (8.2,0) circle [radius=0.1];
\node at (8.2,-0.3)   (a) {$X_{5,5}$};

\draw [very thick, BrickRed] (0.55,1) circle [radius=0.1];
\node at (0.55,0.7)   (a) {$X_{0,1}$};
\draw [very thick, CadetBlue] (2.25,1) circle [radius=0.1];
\node [dgrey] at (2.25,0.7)   (a) {$X_{1,2}$};
\draw [very thick, CadetBlue] (3.78,1) circle [radius=0.1];
\node [dgrey] at (3.78,0.7)   (a) {$X_{2,3}$};
\draw [very thick, CadetBlue] (5.746,1) circle [radius=0.1];
\node [dgrey] at (5.746,0.7)   (a) {$X_{3,4}$};
\draw [very thick, BrickRed] (7.8,1) circle [radius=0.1];
\node at (7.8,0.7)   (a) {$X_{4,5}$};

\draw [very thick, BrickRed] (1.34,2) circle [radius=0.1];
\node at (1.34,1.7)   (a) {$X_{0,2}$};
\draw [very thick, CadetBlue] (2.685,2) circle [radius=0.1];
\node [dgrey] at (2.685,1.7)   (a) {$X_{1,3}$};
\draw [very thick, CadetBlue] (4.811,2) circle [radius=0.1];
\node [dgrey] at (4.811,1.7)   (a) {$X_{2,4}$};
\draw [very thick, BrickRed] (7.1,2) circle [radius=0.1];
\node at (7.1,1.7)   (a) {$X_{3,5}$};
\draw [very thick, CadetBlue] (0.2,2) circle [radius=0.1];
\node [dgrey] at (0.2,1.7)   (a) {};

\draw [very thick, CadetBlue] (0.6,3) circle [radius=0.1];
\node [dgrey] at (0.6,2.7)   (a) {};
\draw [very thick, BrickRed] (1.767,3) circle [radius=0.1];
\node at (1.767,2.7)   (a) {$X_{0,3}$};
\draw [very thick, CadetBlue] (3.683,3) circle [radius=0.1];
\node [dgrey] at (3.683,2.7)   (a) {$X_{1,4}$};
\draw [very thick, BrickRed] (6.03,3) circle [radius=0.1];
\node at (6.03,2.7)   (a) {$X_{2,5}$};

\draw [very thick, CadetBlue] (7.27,3) circle [radius=0.1];
\node [dgrey] at (7.27,2.7)   (a) {};

\draw [very thick, CadetBlue] (1.0,4) circle [radius=0.1];
\node [dgrey] at (1.0,3.7)   (a) {};
\draw [very thick, BrickRed] (2.157,4) circle [radius=0.1];
\node at (2.157,3.7)   (a) {$X_{0,4}$};
\draw [very thick, BrickRed] (4.178,4) circle [radius=0.1];
\node at (4.178,3.7)   (a) {$X_{1,5}$};
\draw [very thick, CadetBlue] (6.84,4) circle [radius=0.1];
\node [dgrey] at (6.84,3.7)   (a) {};
\draw [very thick, CadetBlue] (7.69,4) circle [radius=0.1];
\node [dgrey] at (7.69,3.7)   (a) {};

%\draw [very thick, BrickRed] (2.7871,5) circle [radius=0.1];
%\node at (2.7871,4.7)   (a) {$X_{0,5}$};

%\draw [very thick, BrickRed] (8.2,0) circle [radius=0.1];
%\node at (8.2,-0.3)   (a) {$X_{5,5}$};

%\draw [white] (2,-0.5) circle [radius=0.01];
%\draw [white] (2,4.3) circle [radius=0.1];

\end{tikzpicture}
\caption{A bead configuration on $n = 5$ strings of length $n$ with $n = 5$ beads per string and tilt $2/5$. We can extract a triangular subconfiguration $\{ X_{i,j} : 0 \leq i \leq j \leq n \}$ from this configuration by following consecutive beads, and also by setting $X_{0,n} = X_{\ell,\ell}$ and $X_{n,n} = n + X_{0,0}$ to obtain the corner values.}
  \label{fig:labelling}
\end{figure}%

Consider a bead configuration on $n$ strings of length $n$ with $n$ beads per string and tilt $\ell/n$. Note that these strings have length $n$ rather than length $1$ (compare with the rescaled quantity $n^{n^2}\mathrm{Vol}^{(n)}_{n,\ell}$ above). Such a configuration can be obtained by choosing an element $(t_{h,j})_{h \in \mathbb{Z}_n, 1 \leq j \leq n}$ uniformly from the subset $\mathcal{W}^{(n)}_{n,\ell}$ of $[0,1)^{n^2}$ (see Section \ref{sec:volumes}) and then considering the scaled points $nt_{h,j}$.

In an attempt to trace out a triangular bead configuration, we now label certain beads in the configuration with indices in the boundary of $\urt_n := \{ (i,j) : 0 \leq i \leq j \leq n \}$. Before describing this procedure, the reader may want to consult the diagram in Figure \ref{fig:labelling}. Let $X_{0,0}$, $X_{1,1}$,\ldots, $X_{n-1,n-1}$ denote the positions of the beads on the bottom string. Also set $X_{n,n} = X_{0,0}+n$. We will also set $X_{0,n} = X_{\ell,\ell}$, where $\ell/n$ is the tilt of the configuration. For $1 \leq i \leq n-1$, let $X_{0,i}$ denote the position of the first bead on string $i$ that occurs to the right of $X_{0,i-1}$. 
Let $X_{1,n}$ denote the first bead after $X_{0,n-1}$ on string $n-1$, and thereafter for $i = 2,\ldots,n-1$ let $X_{i,n}$ denote the first bead on string $n-i$ that lies to the right of $X_{i-1,n}$. It might be the case that $X_{n-1,n}$ continues past the end of the torus; in this case add $n$ to make the values increasing. 

In summary, to every bead configuration on $n$ strings of length $n$ with $n$ beads per string, we can associate a vector in $[0,\infty)^{3n-2}$ given by
\begin{align} \label{eq:Xdef}
\mathbf{X} := (X_{0,0},X_{0,1},\ldots,X_{0,n-1};X_{1,n},\ldots,X_{n-1,n};X_{1,1},\ldots,X_{n-1,n-1}),
\end{align}
and we can also associate with this vector the values $X_{0,n} := X_{\ell,\ell}$ and $X_{n,n} := X_{0,0} + n$. If we choose a bead configuration uniformly from $\mathcal{W}_{n,\ell}^{(n)}$ and rescale the $t_{h,j}$ by $n$, then $\mathbf{X}$ is a random vector. For suitable $i,j$, define 
\begin{align*}
Y_{i,j} = X_{i,j} - X_{0,0}.
\end{align*}
We will be most interested in the recentered random vector in $[0,n]^{3n-3}$ given by  
\begin{align} \label{eq:dd}
\mathbf{Y} = (Y_{0,1},\ldots,Y_{0,n-1};Y_{1,n},\ldots,Y_{n-1,n};Y_{1,1},\ldots,Y_{n-1,n-1}).
\end{align}
Since there are $n$ beads per string, and each string has length $n$, the average distance between consecutive beads on the same string is $1$. Also, since the tilt of the configuration is $\ell/n$, the average distance between a bead and the next bead on the next string up is $\ell/n$. In particular,
\begin{align} \label{eq:Ymean}
\mathbf{E}[Y_{0,i}] = i \ell/n, \qquad \mathbf{E}[Y_{i,i}] = i \qquad \text{and} \qquad \mathbf{E}[ Y_{i,n} ] = \ell + i(n-\ell)/n.
\end{align}
We write
\begin{align} \label{eq:yexp}
\mathbf{y}_0 := (\ell/n,\ldots,(n-1)\ell/n;\,\ell+(n-\ell)/n,\ldots,\ell+(n-1)(n-\ell)/n;\,1,\ldots,n-1) \in [0,n]^{3n-3}
\end{align}
for this expected value of $\mathbf{Y}$. 

Consider now conditioning on the event $\Gamma := \{ X_{0,0} = 0 , \mathbf{Y} = \mathbf{y}_0\}$ (of measure zero) that there is a bead located at $(0,0)$, and that the random vector $\mathbf{Y}$ is equal to its expectation. On the event $\Gamma$, we can associate with our bead configuration on $\mathbb{T}_n$ a pair of Gelfand--Tsetlin surfaces $(X_{i,j})_{0 \leq i \leq j \leq n}$ and $(X_{i,j})_{ 0 \leq j \leq i \leq n }$ defined on the triangular subsets
\begin{align*}
\urt_n := \{ (i,j) : 0 \leq i \leq j \leq n \} \quad \text{and} \quad \rt_n := \{ (i,j) : 0 \leq j \leq i \leq n \}
\end{align*}
of $\mathbb{Z}^2$; see Figure \ref{fig:labelling2}. These Gelfand--Tsetlin surfaces both have linear boundary conditions in the sense that
\begin{align*}
X_{i,j} := \varphi^{u_1,u_2}_{i,j} \qquad (i,j) \in \partial ~\urt_n \cup \partial \rt_n,
\end{align*}
with
\begin{align*}
u_1 = 1 - \ell/n \quad \text{and} \quad u_2 = \ell/n,
\end{align*}
and where we recall that $\varphi^{u_1,u_2}_x := u_1 x_1 + u_2 x_2$. 

\begin{figure}[h!]
\centering
\begin{subfigure}{.5\textwidth}
\centering
\begin{tikzpicture}[scale=0.8]

\draw[fill=Melon!20]    (0,0) -- (3.2,5) -- (0,5) --(0,0);
\draw[fill=YellowGreen!20]    (8,0) -- (3.2,5) -- (8,5) --(8,0);

\draw (0,0) -- (8,0);
\draw (0,1) -- (8,1);
\draw (0,2) -- (8,2);
\draw (0,3) -- (8,3);
\draw (0,4) -- (8,4);
\draw[lightgray] (0,0) -- (0,4);
\draw[lightgray] (8,0) -- (8,4);

\filldraw [very thick, BrickRed] (0,0) circle [radius=0.1];
\node at (0,-0.3)   (a) {$X_{0,0}$};

\filldraw [very thick, BrickRed] (1.6,0) circle [radius=0.1];
\node at (1.6,-0.3)   (a) {$X_{1,1}$};
\filldraw [very thick, BrickRed] (3.2,0) circle [radius=0.1];
\node at (3.2,-0.3)   (a) {$X_{2,2}$};
\filldraw [very thick, BrickRed] (4.8,0) circle [radius=0.1];
\node at (4.8,-0.3)   (a) {$X_{3,3}$};
\filldraw [very thick, BrickRed] (6.4,0) circle [radius=0.1];
\node at (6.4,-0.3)   (a) {$X_{4,4}$};

\filldraw [very thick, BrickRed] (0.64,1) circle [radius=0.1];
\node at (0.64,0.7)   (a) {$X_{0,1}$};
\filldraw [very thick, CadetBlue] (2.25,1) circle [radius=0.1];
\node [dgrey] at (2.25,0.7)   (a) {$X_{1,2}$};
\filldraw [very thick, CadetBlue] (3.78,1) circle [radius=0.1];
\node [dgrey] at (3.78,0.7)   (a) {$X_{2,3}$};
\filldraw [very thick, CadetBlue] (5.746,1) circle [radius=0.1];
\node [dgrey] at (5.746,0.7)   (a) {$X_{3,4}$};
\filldraw [very thick, BrickRed] (7.04,1) circle [radius=0.1];
\node at (7.04,0.7)   (a) {$X_{4,5}$};

\filldraw [very thick, BrickRed] (1.28,2) circle [radius=0.1];
\node at (1.28,1.7)   (a) {$X_{0,2}$};
\filldraw [very thick, CadetBlue] (2.685,2) circle [radius=0.1];
\node [dgrey] at (2.685,1.7)   (a) {$X_{1,3}$};
\filldraw [very thick, CadetBlue] (4.811,2) circle [radius=0.1];
\node [dgrey] at (4.811,1.7)   (a) {$X_{2,4}$};
\filldraw [very thick, BrickRed] (6.08,2) circle [radius=0.1];
\node at (6.08,1.7)   (a) {$X_{3,5}$};
\filldraw [very thick, Apricot] (0.2,2) circle [radius=0.1];
\node [dgrey] at (0.2,1.7)   (a) {};

\filldraw [very thick, Apricot] (0.3,3) circle [radius=0.1];
\node [dgrey] at (0.3,2.7)   (a) {};
\filldraw [very thick, BrickRed] (1.92,3) circle [radius=0.1];
\node at (1.92,2.7)   (a) {$X_{0,3}$};
\filldraw [very thick, CadetBlue] (3.683,3) circle [radius=0.1];
\node [dgrey] at (3.683,2.7)   (a) {$X_{1,4}$};
\filldraw [very thick, BrickRed] (5.12,3) circle [radius=0.1];
\node at (5.12,2.7)   (a) {$X_{2,5}$};
\filldraw [very thick, Bittersweet] (6.77,3) circle [radius=0.1];
\node [dgrey] at (6.77,2.7)   (a) {};

\draw[lightgray] (0,0) -- (3.2,5);
\draw[lightgray] (8,0) -- (3.2,5);

\filldraw [very thick, Apricot] (1.0,4) circle [radius=0.1];
\node [dgrey] at (1.0,3.7)   (a) {};
\filldraw [very thick, BrickRed] (2.56,4) circle [radius=0.1];
\node at (2.56,3.7)   (a) {$X_{0,4}$};
\filldraw [very thick, BrickRed] (4.16,4) circle [radius=0.1];
\node at (4.16,3.7)   (a) {$X_{1,5}$};
\filldraw [very thick, Bittersweet] (5.94,4) circle [radius=0.1];
\node [dgrey] at (5.94,3.7)   (a) {};
\filldraw [very thick, Bittersweet] (7.69,4) circle [radius=0.1];
\node [dgrey] at (7.69,3.7)   (a) {};

\draw [white] (4.8,-5.7) circle [radius=0.01];
\draw [white] (4.8,5) circle [radius=0.01];

\end{tikzpicture}
\caption{A bead configuration on $\mathbb{T}_n$ with some fixed\\
points in red along a triangle shape.}
  \label{fig:sub10}
\end{subfigure}%
\begin{subfigure}{.5\textwidth}
\centering
\begin{tikzpicture}[scale=0.8]

\draw[fill=Melon!20]    (8,0) -- (4.8,0) -- (4.8,-5) -- (8,0);
\draw[fill=YellowGreen!20]    (0,0) -- (4.8,0) -- (4.8,-5) --(0,0);

\draw (0,0) -- (8,0);
\draw (0,1) -- (8,1);
\draw (0,2) -- (8,2);
\draw (0,3) -- (8,3);
\draw (0,4) -- (8,4);
\draw (0,-1) -- (8,-1);
\draw (0,-2) -- (8,-2);
\draw (0,-3) -- (8,-3);
\draw (0,-4) -- (8,-4);
\draw[lightgray] (0,-4) -- (0,4);
\draw[lightgray] (8,-4) -- (8,4);

\filldraw [very thick, BrickRed] (0.96,-1) circle [radius=0.1];
\node at (0.96,-1.3)   (a) {$X_{1,0}$};
\filldraw [very thick, BrickRed] (1.92,-2) circle [radius=0.1];
\node at (1.92,-2.3)   (a) {$X_{2,0}$};
\filldraw [very thick, BrickRed] (2.88,-3) circle [radius=0.1];
\node at (2.88,-3.3)   (a) {$X_{3,0}$};
\filldraw [very thick, BrickRed] (3.84,-4) circle [radius=0.1];
\node at (3.84,-4.3)   (a) {$X_{4,0}$};

\filldraw [very thick, BrickRed] (0,0) circle [radius=0.1];
\node at (0,-0.3)   (a) {$X_{0,0}$};

\filldraw [very thick, BrickRed] (1.6,0) circle [radius=0.1];
\node at (1.6,-0.3)   (a) {$X_{1,1}$};
\filldraw [very thick, BrickRed] (3.2,0) circle [radius=0.1];
\node at (3.2,-0.3)   (a) {$X_{2,2}$};
\filldraw [very thick, BrickRed] (4.8,0) circle [radius=0.1];
\node at (4.8,-0.3)   (a) {$X_{3,3}$};
\filldraw [very thick, BrickRed] (6.4,0) circle [radius=0.1];
\node at (6.4,-0.3)   (a) {$X_{4,4}$};

\filldraw [very thick, BrickRed] (0.64,1) circle [radius=0.1];
\node at (0.64,0.7)   (a) {$X_{0,1}$};
\filldraw [very thick, CadetBlue] (2.25,1) circle [radius=0.1];
\node [dgrey] at (2.25,0.7)   (a) {$X_{1,2}$};
\filldraw [very thick, CadetBlue] (3.78,1) circle [radius=0.1];
\node [dgrey] at (3.78,0.7)   (a) {$X_{2,3}$};
\filldraw [very thick, CadetBlue] (5.746,1) circle [radius=0.1];
\node [dgrey] at (5.746,0.7)   (a) {$X_{3,4}$};
\filldraw [very thick, BrickRed] (7.04,1) circle [radius=0.1];
\node at (7.04,0.7)   (a) {$X_{4,5}$};

\filldraw [very thick, BrickRed] (1.28,2) circle [radius=0.1];
\node at (1.28,1.7)   (a) {$X_{0,2}$};
\filldraw [very thick, CadetBlue] (2.685,2) circle [radius=0.1];
\node [dgrey] at (2.685,1.7)   (a) {$X_{1,3}$};
\filldraw [very thick, CadetBlue] (4.811,2) circle [radius=0.1];
\node [dgrey] at (4.811,1.7)   (a) {$X_{2,4}$};
\filldraw [very thick, BrickRed] (6.08,2) circle [radius=0.1];
\node at (6.08,1.7)   (a) {$X_{3,5}$};

\filldraw [very thick, BrickRed] (1.92,3) circle [radius=0.1];
\node at (1.92,2.7)   (a) {$X_{0,3}$};
\filldraw [very thick, CadetBlue] (3.683,3) circle [radius=0.1];
\node [dgrey] at (3.683,2.7)   (a) {$X_{1,4}$};
\filldraw [very thick, BrickRed] (5.12,3) circle [radius=0.1];
\node at (5.12,2.7)   (a) {$X_{2,5}$};

\filldraw [very thick, Bittersweet] (3.61,-2) circle [radius=0.1];
\node [dgrey] at (3.61,-2.3)   (a) {$X_{3,1}$};

\draw[lightgray] (0,0) -- (3.2,5);
\draw[lightgray] (8,0) -- (3.2,5);
\draw[lightgray] (0,0) -- (4.8,-5);
\draw[lightgray] (8,0) -- (4.8,-5);

\filldraw [very thick, BrickRed] (2.56,4) circle [radius=0.1];
\node at (2.56,3.7)   (a) {$X_{0,4}$};
\filldraw [very thick, BrickRed] (4.16,4) circle [radius=0.1];
\node at (4.16,3.7)   (a) {$X_{1,5}$};

\filldraw [very thick, BrickRed] (3.2,5) circle [radius=0.1];
\node at (3.2,4.7)   (a) {$X_{0,5}$};
\filldraw [very thick, BrickRed] (4.8,-5) circle [radius=0.1];
\node at (4.8,-5.3)   (a) {$X_{5,0}$};
\filldraw [very thick, BrickRed] (8,0) circle [radius=0.1];
\node at (8,-0.3)   (a) {$X_{5,5}$};
\filldraw [very thick, BrickRed] (7.36,-1) circle [radius=0.1];
\node at (7.36,-1.3)   (a) {$X_{5,4}$};
\filldraw [very thick, BrickRed] (6.72,-2) circle [radius=0.1];
\node at (6.72,-2.3)   (a) {$X_{5,3}$};
\filldraw [very thick, BrickRed] (6.08,-3) circle [radius=0.1];
\node at (6.08,-3.3)   (a) {$X_{5,2}$};
\filldraw [very thick, BrickRed] (5.56,-4) circle [radius=0.1];
\node at (5.56,-4.3)   (a) {$X_{5,1}$};

\filldraw [very thick, Apricot] (5,-3) circle [radius=0.1];
\node [dgrey] at (5,-3.3)   (a) {$X_{4,1}$};
\filldraw [very thick, Apricot] (5.2,-2) circle [radius=0.1];
\node [dgrey] at (5.2,-2.3)   (a) {$X_{4,2}$};
\filldraw [very thick, Apricot] (5.9,-1) circle [radius=0.1];
\node [dgrey] at (5.9,-1.3)   (a) {$X_{4,3}$};

\filldraw [very thick, Bittersweet] (2.74,-1) circle [radius=0.1];
\node [dgrey] at (2.74,-1.3)   (a) {$X_{2,1}$};
\filldraw [very thick, Bittersweet] (4.49,-1) circle [radius=0.1];
\node [dgrey] at (4.49,-1.3)   (a) {$X_{3,2}$};

\end{tikzpicture}
\caption{A pair of generalized Gelfand--Tsetlin patterns on triangles with linear boundaries.}
  \label{fig:sub102}
\end{subfigure}%
\caption{If we condition on the event that $X_{0,0}$ lies at zero, and thereafter that $X_{0,i}$, $X_{i,n}$ and $X_{i,i}$ all lie at their 
conditionally expected positions given $\{X_{0,0} = 0\}$, 
we can rearrange the configuration to create a bead configuration on a square with fixed points along the diagonal.}
  \label{fig:labelling2}
\end{figure}%
We have the following lemma.

\begin{lemma} \label{lem:cuba}
Choose a bead configuration uniformly at random from the set of bead configurations on $n$ strings with $n$ beads per string and tilt $\gamma =\ell/n$. Let $\psi_{n,\ell}:[0,n]^{3n-3} \to [0,\infty)$ denote the probability density function of the random vector $\mathbf{Y}$ defined in \eqref{eq:dd}, and let $\mathbf{y}_0$ be the expectation vector of $\mathbf{Y}$, as in \eqref{eq:yexp}. Then
\begin{align*}
n^{n^2}\mathrm{Vol}^{(n)}_{n,\ell} \psi_{n,\ell}(\mathbf{y}_0) = T_n \left( \frac{\ell}{n}, \frac{n-\ell}{n} \right)^2,
\end{align*}
where $T_n(u_1,u_2)$ is the triangular partition function for Gelfand--Tsetlin surfaces defined in \eqref{eq:Tn}. 
\end{lemma}
\begin{proof}
Recall that $n^{n^2} \mathrm{Vol}^{(n)}_{n,\ell}$ is the volume of the set of bead configurations on $n$ toroidal strings of length $n$ with $n$ beads per string and tilt $\ell/n$. Given a value $\mathbf{y} \in [0,n]^{3n-3}$, 
\begin{align} \label{eq:nvol}
n^{n^2} \mathrm{Vol}^{(n)}_{n,\ell} \psi_{n,\ell}(\mathbf{y}) =  \text{Vol.\ arrangements.\ remaining $n^2-(3n-2)$ beads around $\mathbf{y} \cup \{0\}$},
\end{align}
where we now take a moment to explain the right-hand side of \eqref{eq:nvol}. 
With $\mathbf{X}$ defined in \eqref{eq:Xdef}, suppose we consider the set of bead configurations on $n$ strings of length $n$ with $n$ beads per string such that the random vector $\mathbf{X}$ takes the value $(0,\mathbf{Y})$ (i.e., we force $X_{0,0} = 0$, and set $X_{i,j} = Y_{i,j}$ for $(i,j)\neq (0,0)$). Then as in Section \ref{sec:volumes}, the set of ways of arranging the remaining $n^2-(3n-2)$ beads to form a bead configuration may be associated with a subset of $[0,n]^{n^2-3n+2}$; the right-hand side of \eqref{eq:nvol} refers to the $(n^2-3n+2)$-dimensional Lebesgue measure of this subset. 

Now note that any bead configuration with $\mathbf{Y} = \mathbf{y}_0$, i.e., with each of the points $Y_{0,i}$, $Y_{i,i}$, and $Y_{i,n}$ lying at its expected position, may be rearranged to create two triangular bead configurations, as in Figure \ref{fig:labelling2}. The result follows. 
\end{proof}

%%%%%%%%%%%%%%%%%%%%%%%%%%%
\subsection{Partition function upper bounds} \label{sec:51}
%%%%%%%%%%%%%%%%%%%%%%%%%%%
We now begin working towards our proof of Theorem \ref{thm:linear}, which is split into an upper and lower bound for $T_n(u_1,u_2)$. The more difficult proof of the upper bound for $T_n(u_1,u_2)$ appeals to the determinantal correlation formulas in Theorem \ref{thm:detnew} and subsequent linear algebra arguments. The proof of the complementary lower bound for $T_n(u_1,u_2)$ is more straightforward, and uses tools from the theory of log-concave functions. As we work towards our upper bound for $T_n(u_1,u_2)$, we begin with an upper bound on $Z_\theta(\lambda,T)$ defined in \eqref{eq:DA1}.

\begin{thm}\label{thm:Aound on Ztheta}
Let $p \in (-1,1)$, $n\in\mathbb N$, $T > 0$ and consider $\theta\in\{0,1\}^2$. Then we have
\begin{align} \label{eq:upper}
|Z_\theta(pT,T)| \leq  2^{n-1} e^{CT} \exp \left\{ Tn \left( \frac{\sqrt{1-p^2}}{\pi} - p (1 - \mathrm{arccos}(-p)/\pi)  \right) \right\},
\end{align}
for some universal constant $C > 0$ and the inverse cosine taken as $\arccos:(-1,1) \to (0,\pi)$. 
\end{thm}

\begin{proof}
Using \eqref{eq:DA1} we have
\begin{align} \label{eq:hy1}
|Z_\theta(pT,T)| \leq \frac{1}{2} \prod_{z^n = (-1)^{\theta_2}}( e^{T\mathrm{Re}(z)} + e^{ - pT} ) \leq 2^{n-1} \exp \left\{ T\sum_{ z^n = (-1)^{\theta_2} } \max \{ \mathrm{Re}(z), - p \} \right\}. 
\end{align}
Since the function $\phi\mapsto \max\{\cos\phi,-p\}$ is uniformly Lipschitz, its Riemann sums over the roots of $z^n=(-1)^{\theta_2}$ differ from the corresponding integral by at most a universal constant. Thus, controlling the sum with an integral, if $-p = \cos(\phi_0)$ then for some universal $C > 0$ we have 
\begin{align} \label{eq:hy2}
\sum_{ z^n = (-1)^{\theta_2} } \max \{ \mathrm{Re}(z), - p \}  &\leq \frac{n}{2\pi} \int_0^{2 \pi} \max \{ \cos(\phi) , - p \} \mathrm{d}\phi + C \nonumber \\
&= n(\sin(\phi_0)/\pi - p ( 1 - \phi_0/\pi)) + C.
\end{align} 
Using \eqref{eq:hy2} in \eqref{eq:hy1}, and the definition of $\phi_0$, the result follows. 
\end{proof}
%%%%%%%%%%%%%%%%%%%%%%%%%%%%%%%%%%%%%%%%%%%%%%
\subsection{Correlation function upper bounds}
%%%%%%%%%%%%%%%%%%%%%%%%%%%%%%%%%%%%%%%%%%%%%%

We would like to develop upper bounds for the correlation functions $\rho_m^{(n),k,\ell}$. 
We begin by observing that by bounding below the left-hand side of \eqref{eq:weightsum2} by a single summand on the right-hand side, and rearranging, we have the inequality
\begin{align} \label{eq:corrbound}
\mathrm{Vol}_{k,\ell}^{(n)} \rho^{(n),k,\ell}_m(y_1,\ldots,y_m) \leq 
 \frac{ Z(\lambda,T)}{T^{nk}e^{-\lambda \ell} } \rho^{(n),\lambda,T}_m(y_1,\ldots,y_m) .
\end{align}
% Our next two results bound the two quantities $Z(\lambda,T)/(T^{nk} e^{ - \ell \lambda})$ and $ \rho^{(n),\lambda,T}_m(y_1,\ldots,y_m)$ occurring on the right-hand side of \eqref{eq:corrbound}.
We now develop an upper bound for the correlation functions $\rho^{(n),\lambda,T}_m(y_1,\ldots,y_m)$. 
A classical corollary of Hadamard's inequality states that if $(A_{i,j})_{1\leq i,j \leq m}$ is a matrix with complex entries satisfying $|A_{i,j}| \leq B$ (where $B>0$), then we have $|\det_{i,j=1}^m (A_{i,j})| \leq B^mm^{m/2}$. 
In particular, using \eqref{eq:corrbound} and \eqref{eq:Psum0} and the definition of $d_\theta(\lambda,T)$ to obtain the first inequality below, then the determinant inequality to obtain the second, for any distinct points $y_1,\ldots,y_m$ in $\mathbb{T}_n$ we have 
\begin{align}
\mathrm{Vol}_{k,\ell}^{(n)} \rho^{(n),k,\ell}_m(y_1,\ldots,y_m) 
&\leq \frac{1}{T^{nk}e^{-\lambda \ell}}
\sum_{ (\theta_1,\theta_2)\in \{0,1\}^2 } \left| Z_\theta(\lambda,T)\right|
\left| \det_{i,j=1}^m H^{\lambda+\theta_1\pi \iota ,\theta_2,T}( y_j-y_i)\right| \nonumber \\
&\leq \frac{m^{m/2}}{T^{nk}e^{-\lambda \ell}}
\sum_{ (\theta_1,\theta_2)\in \{0,1\}^2 } \left| Z_\theta(\lambda,T)\right|
\left( \sup_{t \in (-1,1),\, h \in \mathbb{Z} } |H^{\lambda+\theta_1\pi \iota ,\theta_2,T}(t,h) | \right)^m. \label{eq:hada}
\end{align}

With a view to controlling the right-hand side of \eqref{eq:hada} we have the following definition:

\begin{df} \label{df:good}
We say that the pair $(p,T)$ is \textbf{good} if and only if for all integers $k \in \mathbb{Z}$,
\begin{align*}
\left| \frac{T \sqrt{1-p^2}}{\pi} -  k \right| \geq 1/4.
\end{align*} 
\end{df}

We make the following observation.

\begin{rem} \label{rem:good}
For any $p \in (-1,1)$ and $T > 0$, set $T ' = T + \pi/(2\sqrt{1-p^2})$. Then at least one of the pairs $(p,T)$ or $(p,T')$ is good.
\end{rem}

In Appendix \ref{sec:Hboundproof}, we prove the following bound.

\begin{lemma} \label{lem:Hbound}
If $(p,T)$ is good, then with $H^{\beta,\theta_2,T}$ defined in \eqref{eq:detnewkernel}, for all $(\theta_1,\theta_2) \in \{0,1\}^2$ we have
\begin{align} \label{eq:predet4}
  \sup_{t \in (-1,1),\,h \in \mathbb{Z} }\Big| H^{pT+\theta_1\pi \iota , \theta_2, T} (t,h)\Big| \leq C \frac{ T}{\sqrt{1-p^2}}.
\end{align}
\end{lemma}

We now complete our upper bound for the correlation functions under $\mathbf{P}^{(n),n,\ell}$. 

\begin{lemma} \label{lem:upa}
Let $m\in\mathbb N$ and assume that $\log(n)^{-1/2} \leq \gamma = \ell/n \leq 1/2$ for $1 \leq \ell \leq n-1$. Then, for any distinct points $y_1,\ldots,y_m$ in $\mathbb{T}_n$, we have
\begin{align*}
\mathrm{Vol}_{n,\ell}^{(n)} \rho^{(n),n,\ell}_m(y_1,\ldots,y_m) \leq \exp \left\{ - n^2 \log n - n^2 \sigma(\gamma,1-\gamma) + C (m+n) \log(m+n) \right\}
\end{align*}
for some universal $C \in(0,\infty)$ and $\sigma(u_1,u_2) = -\log(u_1+u_2) - \log \sin \frac{ \pi u_1}{u_1+u_2} -1 +  \log \pi$.
\end{lemma}
\begin{proof}
Setting $k=n$, $\lambda :=pT$ in \eqref{eq:hada}, and using \eqref{eq:predet4}, for any $p \in (-1,1)$ and $T > 0$ such that $(p,T)$ is good (Definition \ref{df:good}) we have
\begin{align} \label{eq:corrbound2}
\mathrm{Vol}_{n,\ell}^{(n)} \rho^{(n),n,\ell}_m(y_1,\ldots,y_m) 
\leq \frac{m^{m/2}}{ T^{n^2} e^{- \ell p T} } \sum_{\theta \in \{0,1\}^2} |Z_\theta(pT,T)| \left( C\frac{T}{\sqrt{1-p^2}} \right)^m.
\end{align} 

It remains to bound the factor $ \sum_{\theta \in \{0,1\}^2}|Z_\theta(pT,T)| /(T^{n^2} e^{ - \ell pT})$ occurring on the right-hand side of \eqref{eq:corrbound2}. Appealing now to the bound \eqref{eq:upper} on $Z_{\theta}(pT,T)$ from Theorem \ref{thm:Aound on Ztheta} we have 
\begin{align} \label{eq:geta}
\frac{  \sum_{\theta \in \{0,1\}^2}|Z_\theta(pT,T)|}{T^{n^2} e^{ - \ell pT} }\leq   2^{n+1} e^{CT} T^{ - n^2} e^{ Tn f_\gamma(p) },
\end{align}
where given $\ell=\gamma n$, we put 
\begin{align*}
f_\gamma:(-1,1)\to\mathbb R, \quad f_\gamma(p) := \frac{\sqrt{1-p^2}}{\pi} - p \bigl(1 - \mathrm{arccos}(-p)/\pi\bigr) + \gamma p.
\end{align*}

Note that \eqref{eq:geta} holds for any good pair $(p,T)$ with $p \in (-1,1)$ and $T > 0$, and we are at liberty to choose the best values possible. We begin by optimising over $p$. One readily checks (recalling that $\frac{\mathrm{d}}{\mathrm{d}x} \arccos(x) = - (1-x^2)^{-1/2}$) that 
   \[
     f_\gamma'(p) = \gamma + \frac{\mathrm{arccos}(-p)}{\pi} - 1,
   \]
and that the value $p := \cos(\pi \gamma)$ minimizes $f_\gamma$. Noting that $\mathrm{arccos}(-p) = \pi - \mathrm{arccos}(p)$ and using a standard trigonometric identity, we obtain that
\begin{align} \label{eq:apple1}
f_\gamma(\cos(\pi \gamma)) =  \sin(\pi \gamma)/\pi. 
\end{align}
Using \eqref{eq:corrbound2}, \eqref{eq:geta} and \eqref{eq:apple1}, with $p = \cos(\pi \gamma)$, provided $(p,T)$ is a good pair we have 
\begin{align} \label{eq:lapiz}
\mathrm{Vol}_{n,\ell}^{(n)} \rho^{(n),n,\ell}_m(y_1,\ldots,y_m) 
\leq 2^{n+1} e^{CT} T^{-n^2} e^{ Tn \frac{ \sin(\pi \gamma)}{ \pi} } m^{m/2} \left( \frac{CT}{ \sin(\pi \gamma)} \right)^m.
\end{align}

We now look to optimise over $T > 0$. As a function of $T$, the quantity 
\begin{align*}
T^{-n^2} e^{ Tn \frac{ \sin(\pi \gamma)}{ \pi} }= \exp \left( - n^2 \log T + T n \sin(\pi \gamma)/\pi \right)
\end{align*}
is minimized by setting $T = T_* := \pi n/\sin(\pi \gamma)$. 
Now recall from Remark \ref{rem:good} that either $T_*$ itself or $T' = T_* + \pi/(2\sqrt{1-p^2})$ makes $(p,T)$ a good pair. Since $p = \cos(\pi \gamma)$, it follows that there exists $\kappa \in \{0,1\}$ such that the pair $(p,T)$ with $ p = \cos( \pi \gamma)$ and 
\begin{align} \label{eq:Tchoice}
T = \frac{ \pi (n+\kappa/2) }{ \sin( \pi \gamma) } 
\end{align}
is good. Either way, with $T$ as in \eqref{eq:Tchoice}, we have 
 \begin{align} \label{eq:apple2}
T^{-n^2}e^{Tnf_\gamma(\cos(\pi \gamma))} &= \exp \left( - n^2 ( \log n + \log \pi - \log \sin (\pi \gamma) - 1 +O(\log^{1/2}(n)/n) ) \right) \nonumber \\
&= \exp \left( - n^2 \sigma( \gamma,1-\gamma) - n^2 \log n + O(n \log n)\right),
\end{align}
where we have used the condition in the statement of the lemma that $\sin(\pi\gamma) \geq \frac{1}{2} \gamma \geq \frac{1}{2} \log(n)^{-1/2}$.

Using \eqref{eq:apple2} in \eqref{eq:lapiz} with the good pair $(p,T) = (\cos(\pi \gamma), \pi(n+\kappa/2)/\sin(\pi \gamma))$, we obtain
\begin{align} \label{eq:lapiz2}
\mathrm{Vol}_{n,\ell}^{(n)} \rho^{(n),n,\ell}_m(y_1,\ldots,y_m) \leq \exp \left\{ - n^2 \sigma(\gamma,1-\gamma) - n^2 \log(n) + C (m+n) \log(m+n)  \right\}
\end{align}
for some sufficiently large universal $C\in(0,\infty)$. That completes the proof. 
\end{proof}

%%%%%%%%%%%%%%%%%%%%%%%%%%%%%%%%%%%%%%%%%%%%%%%%%%%%%%%%%%%%%%%%%%%%%%%%%%%
\subsection{Proof of Theorem \ref{thm:linear}} 
%%%%%%%%%%%%%%%%%%%%%%%%%%%%%%%%%%%%%%%%%%%%%%%%%%%%%%%%%%%%%%%%%%%%%%%%%%%

Before completing the proof of Theorem \ref{thm:linear}, we require an inequality from convex geometry. This inequality states that if $f:\mathbb{R}^n \to [0,\infty)$ is a log-concave probability density function of a random vector with expectation $x_0 \in \mathbb{R}^n$, then 
\begin{align} \label{eq:LV}
f(x_0) \leq \sup_{s \in \mathbb{R}^n} f(s) \leq C^n f(x_0);
\end{align}
for a universal constant $C > 0$. See \cite[Theorem 5.14b)]{LV}. 

We now claim that the density $\psi_{n,\ell}:[0,n]^{3n-3} \to [0,\infty)$ defined in Lemma \ref{lem:cuba} is log-concave. To see this, we note that to every bead configuration, by taking a labelling as in Figure \ref{fig:sub102}, we can associate the bead configuration on $n$ strings, having $n$ beads per string and tilt $\ell/n$ with a collection of numbers $\{ X_{i,j} : 0 \leq i,j \leq n-1\}$ taking values in $[0,n]$. A collection of numbers $\{ X_{i,j} : 0 \leq i,j \leq n-1 \}$ taking values in $[0,n)$ must satisfy a certain collection of linear inequalities in order to ensure they can be associated with a bead configuration. These inequalities are closed under convex combinations, so it follows that the collection of such $(X_{i,j})$ is a convex body in $[0,n)^{n^2}$. 
By the Pr\'ekopa--Leindler inequality, it follows that $\psi_{n,\ell}$, as the law of a projection of a convex body in $n^2$-dimensional space onto a $3n-3$-dimensional subspace, must be log-concave. 

We now complete the proof of Theorem \ref{thm:linear}.

\begin{proof}[Proof of Theorem \ref{thm:linear}]
We begin by proving Theorem \ref{thm:linear} in the special case $(u_1,u_2) = (\ell/n,1-\ell/n)$ for some $\ell$ with $ \log(n)^{-1/2} \leq \ell/n \leq 1/2$. 

By Lemma \ref{lem:cuba} we have 
\begin{align} \label{eq:ocean}
T_n \left( \frac{\ell}{n}, \frac{n-\ell}{n} \right) = n^{n^2/2} \left(\mathrm{Vol}^{(n)}_{n,\ell}\right)^{1/2} \psi_{n,\ell}(\mathbf{y}_0)^{1/2},
\end{align}
where $\psi_{n,\ell}$ is the probability density function on $\mathbb{R}^{3n-3}$ of the random vector $\mathbf{Y}$ and $\mathbf{y}_0$ is its expected value as introduced in Section \ref{sec:unravel}. 

Note that in a configuration with tilt $\ell/n$, any beads located at locations $\mathbf{y}_0$ must be precisely the beads defined in Section \ref{sec:unravel}. In particular,
\begin{align*}
\psi_{n,\ell}(\mathbf{y}_0) = \rho_{3n-2}^{(n),n,\ell}(\mathbf{y}_0, 0),
\end{align*}
where we note that $\mathbf{y}_0$ has $3n-3$ coordinates.

We can now prove an upper bound for $T_n \left( \frac{\ell}{n}, \frac{n-\ell}{n} \right)$ using Lemma \ref{lem:upa}. Indeed, setting $m = 3n-2$ in Lemma \ref{lem:upa}, provided $\log(n)^{-1/2} \leq \ell/n \leq 1/2$ we obtain
\begin{align} \label{eq:rome2}
\frac{1}{n^2/2} \log T_n \left( \frac{\ell}{n}, \frac{n-\ell}{n} \right) \leq -\sigma(\ell/n,1-\ell/n)  + C \log(n)/n.
\end{align}

We now develop a lower bound complementing \eqref{eq:rome2}. Since $\psi_{n,\ell}$ is log-concave, using the Lovász--Vempala inequality \eqref{eq:LV} in \eqref{eq:ocean} we obtain 
\begin{align} \label{eq:prerough}
T_n\left( \frac{\ell}{n} , \frac{n-\ell}{n} \right) 
\geq n^{n^2/2} \left( \mathrm{Vol}^{(n)}_{n,\ell} \right)^{1/2} C^{-n}  \sup_{ \mathbf{y} \in [0,n]^{3n-3} } \psi_{n,\ell}( \mathbf{y})^{1/2}.
\end{align}
Taking a very rough bound, since $\psi_{n,\ell}:\mathbb{R}^{3n-3} \to [0,\infty)$ is a probability density function supported on $[0,n]^{3n-3}$, we have
\begin{align} \label{eq:rough}
1 = \int_{[0,n]^{3n-3}}  \psi_{n,\ell}( \mathbf{y}) \mathrm{d}\mathbf{y} \leq n^{3n-3}   \sup_{ \mathbf{y} \in [0,n]^{3n-3} } \psi_{n,\ell}( \mathbf{y}). 
\end{align}
Using \eqref{eq:rough} in \eqref{eq:prerough} we obtain
\begin{align} \label{eq:tik1}
T_n\left( \frac{\ell}{n} , \frac{n-\ell}{n} \right) \geq n^{n^2/2} \left( \mathrm{Vol}^{(n)}_{n,\ell} \right)^{1/2} e^{ - C n \log n }.
\end{align}

Using the implied lower bound in \eqref{eq:ordered10} in \eqref{eq:tik1} it follows that whenever $\log(n)^{-1/2} \leq \ell/n \leq 1/2$ we have 
\begin{align*}
\frac{1}{n^2/2} \log T_n\left( \frac{\ell}{n} , \frac{n-\ell}{n} \right) \geq - \sigma(\ell/n, 1-\ell/n) - C\log(n)/n,
\end{align*}
for some universal $C\in(0,\infty)$. 

In summary, whenever $\log(n)^{-1/2} \leq \ell/n \leq 1/2$ we have
\begin{align*}
\frac{1}{n^2/2} \log T_n\left( \frac{\ell}{n} , \frac{n-\ell}{n} \right) =- \sigma(\ell/n,1-\ell/n) + O( \log n/n),
\end{align*}
where the implied constants are universal. That proves Theorem \ref{thm:linear} in the case where $(u_1,u_2)$ take the form $(\ell/n,1-\ell/n)$.

To generalise this result to general $(u_1,u_2)$ satisfying the condition $\log(n)^{-1/2}  \leq u_1 \wedge u_2$ of the statement, first note that since 
\begin{align*}
\sigma(\lambda u_1 , \lambda u_2) = - \log \lambda + \sigma(u_1,u_2), 
\end{align*}
by \eqref{eq:GTscal} and the fact that the interior of $\rt_n$ contains $(n-1)(n-2)/2$ points, we may assume without loss of generality that $u_1 + u_2 =1 $. From here, one may use parts (2) and (3) of Lemma \ref{lem:st} to obtain
\begin{align*}
\sigma(u_1+\varepsilon,u_2-\varepsilon) = \sigma(u_1,u_2) + O\left( \frac{ \varepsilon}{u_1 \wedge u_2} \right) = \sigma(u_1,u_2) + O \left( \frac{ \log(n) }{n} \right). 
\end{align*}
The proof may now be completed by rescaling $u_1$ and $u_2$ by $(u_1+u_2)$, and then choosing a value $\ell$ such that $(\ell/n,1-\ell/n) = (u_1/(u_1+u_2) + O(1/n), u_2/(u_1+u_2) + O(1/n))$. 
\end{proof}

%%%%%%%%%%%%%%%%%%%%%%%%%%%%%%%%%%%%%%%%%%%%%%
%%%%%%%%%%%%%%%%%%%%%%%%%%%%%%%%%%%%%%%%%%%%%%
%%%%%%%%%%%%%%%%%%%%%%%%%%%%%%%%%%%%%%%%%%%%%%
\section{Bounds for triangular partition functions} \label{sec:ST2proof}
%%%%%%%%%%%%%%%%%%%%%%%%%%%%%%%%%%%%%%%%%%%%%%
%%%%%%%%%%%%%%%%%%%%%%%%%%%%%%%%%%%%%%%%%%%%%%
%%%%%%%%%%%%%%%%%%%%%%%%%%%%%%%%%%%%%%%%%%%%%%
%%%%%%%%%%%%%%%%%%%%%%%%%%%%%%%%%%%%%%%%%%%%%%%%%%%%%%%%%%%%%%%%%%%%%%%%%%%%%%%
\subsection{Recapitulation and main statement}
%%%%%%%%%%%%%%%%%%%%%%%%%%%%%%%%%%%%%%%%%%%%%%%%%%%%%%%%%%%%%%%%%%%%%%%%%%%%%%%

Let us take stock of what we have achieved so far. 
Let $A$ and $B$ be disjoint subsets of $\mathbb{Z}^2$ and write $C = A \sqcup B$ for their union. 
Given a boundary function $\varphi:B \to \mathbb{R}$, we write $Z_A(\varphi)$ for the Gelfand--Tsetlin integral over $\mathbb{R}^A$ subject to the boundary conditions $\{ \varphi_x : x \in B \}$ (see \eqref{eq:pf2} for the precise formula). Recall that $\rt_n := \{ (x_1,x_2) \in \mathbb{Z}^2: 0 \leq x_2 \leq x_1 \leq n \}$. In the special case where $A = \interior{\rt_n}$ and $B = \partial \rt_n$, given a function $\varphi:\partial \rt_n \to \mathbb{R}$ let us introduce the shorthand
\begin{align} \label{eq:Tndef}
T_n(\varphi) := Z_{\interior{\rt_n}}(\varphi) := \int_{\mathbb{R}^{\rt_n}} \mathrm{1}_{\mathrm{GT}}^{ \rt_n}( \phi) \prod_{x \in \interior{\rt_n}} \mathrm{d}\phi_x \prod_{x \in \partial \rt_n } \delta_{\varphi_x}(\mathrm{d}\phi_x ).
\end{align}

When $\varphi$ is the linear function $\varphi^u_x := u_1x_1 + u_2x_2$ we write $T_n(u_1,u_2) := T_n(\varphi^u)$ for short. Theorem \ref{thm:linear} states that provided $(\log n)^{-1/2} \leq u_1 \wedge u_2$ we have 
\begin{align*}
\frac{1}{n^2/2} \log T_n(u_1,u_2) = -\sigma(u_1,u_2) + O \left( \frac{ \log n + K_u }{ n } \right),
\end{align*}
where $K_u := |\log(u_1 \wedge u_2)|+|\log(u_1+u_2)|$. 
In this section we will be concerned with extending this result to describe the logarithmic asymptotics of $T_n(\varphi)$ where $\varphi$ is a function close to $\varphi^u$. Let us recall that a function $\varphi:B \to \mathbb{R}$ defined on a subset $B$ of $\mathbb{Z}^2$ is $c$-increasing for some $c \geq 0$ if 
\begin{align*}
\varphi_y - \varphi_x \geq c( (y_1 - x_1) + (y_2 - x_2) ) = c\|y - x\|_1\qquad \forall x,y \in B \text{ with } x \leq y,
\end{align*}
where $\|y - x\|_1 = |y_1 - x_1| + |y_2 - x_2|$ is the $\ell_1$ norm on the lattice.
% Of course $\varphi:B \to \mathbb{R}$ is increasing if and only if it is $0$-increasing. We note that $\varphi^u$ itself is $(u_1 \wedge u_2)$-increasing. 
For $u = (u_1,u_2) \in (0,\infty)^2$, $r > 0$ and $c \geq 0$, define
\begin{align*}
B_n(u,r,c) := \{ \varphi:\partial \rt_n \to \mathbb{R} : |\varphi_x - \varphi^u_x| \leq r ~ \forall x \in \partial \rt_n, \text{ and $\varphi$ is $c$-increasing} \}
\end{align*}
to be the set of $c$-increasing functions in a uniform ball of radius $r$ around the linear function $\varphi^u:\partial\rt_n \to \mathbb{R}$. For short, we will write
\begin{align*}
B_n(u,r) = B_n(u,r,0) := \{ \varphi:\partial \rt_n \to \mathbb{R} : |\varphi_x - \varphi^u_x| \leq r ~ \forall x \in \partial \rt_n, \text{ and $\varphi$ is increasing} \}.
\end{align*}

The main result of this section is the following:

\begin{thm} \label{thm:upperlower}
Let $(u_1,u_2) \in (0,\infty)^2$ and $n \in \mathbb{N}$ satisfy $(\log n)^{-1/2} \leq u_1 \wedge u_2$ and $n \geq 10$. Given $r > 0$, define the error term 
\begin{align*}
\rho_r := \frac{ r + u_1 + u_2}{ (u_1 \wedge u_2) n}.
\end{align*}

Then we have the upper bound
\begin{align} \label{eq:upper1}
\frac{1}{n^2/2} \log T_n(\varphi) \leq - \sigma(u_1,u_2) + C \rho_r  \left( \log_+ \frac{1}{\rho_r }+ K_u \right),  \qquad \varphi \in B_n(u, r).
\end{align}

If additionally $r \leq n(u_1 \wedge u_2)/100$, we have the lower bound
\begin{align} \label{eq:lower1}
\frac{1}{n^2/2} \log T_n(\varphi) \geq -\sigma(u_1,u_2) -  C \rho_r \left( \log_+ \frac{1}{\rho_r }+ K_u \right), \qquad \varphi \in B_n(u,r,10r/n).
\end{align}
Here $C > 0$ is a universal constant and $K_u = |\log(u_1 \wedge u_2)| + |\log(u_1+u_2)|$, and $\log_+(x)$ is the maximum of $\log(x)$ and $0$.
\end{thm}

We note that our lower bound requires the function $\varphi:\partial\rt_n \to \mathbb{R}$ to be $10r/n$-increasing. Since $r \leq n(u_1 \wedge u_2)/100$, we have $10r/n \leq (u_1 \wedge u_2)/10$.

We note here that an error term of the order $r \log (1/r)$ also appears in Cohn, Kenyon and Propp \cite[Lemma 3.5]{CKP}.
%%%%%%%%%%%%%%%%%%%%%%%%%%%%%%%%%%%%%%%%%%%%%%%%%%%%%%%%%%%%%%%%%%%%%%%%%%%%%%%
\subsection{Upper and lower bounds for the surface tension function}
%%%%%%%%%%%%%%%%%%%%%%%%%%%%%%%%%%%%%%%%%%%%%%%%%%%%%%%%%%%%%%%%%%%%%%%%%%%%%%%

In this brief section we prove the following lemma.

\begin{lemma} \label{lem:sigmabound}
There are universal constants $c_1 < c_2 \in\mathbb{R}$ such that for all $u_1,u_2 > 0$ we have  
\begin{align} \label{eq:sigmabound}
c_1 + \log(u_1 \wedge u_2) \leq - \sigma(u_1,u_2) \leq c_2 + \log(u_1 \wedge u_2).
\end{align}
\end{lemma}
\begin{proof}
Since $- \sigma(\lambda u_1 ,\lambda u_2) = \log \lambda - \sigma(u_1,u_2)$ and $\sigma(u_2,u_1) = \sigma(u_1,u_2)$, to prove the result, it is no loss of generality to assume that $u_1 + u_2 = 1$ and that $u_1 \leq u_2$, which implies $u_1 \leq 1/2$. Then
\begin{align*}
\sigma(u_1,u_2) = - \log\sin \pi u_1 -1 + \log \pi.
\end{align*}
The result now follows from the inequality $2x/\pi \leq \sin(x) \leq x$ for $x \in [0,\pi/2]$. 
\end{proof}

As a corollary, let us note that
\begin{align} \label{eq:able}
|\sigma(u_1,u_2) | \leq C K_u.
\end{align}

Throughout this section we will regularly use the fact that
\begin{align*}
K_u \geq 1/2 \qquad \text{for all $u = (u_1,u_2) \in (0,\infty)^2$},
\end{align*}
in order to absorb error terms and produce simpler expressions. 

%%%%%%%%%%%%%%%%%%%%%%%%%%%%%%%%%%%%%%%%
\subsection{An integral upper bound}
%%%%%%%%%%%%%%%%%%%%%%%%%%%%%%%%%%%%%%%%
In the following statement, the integral $\int_{B_n(u,r')} T_n(\phi)\mathrm{d}\phi $ refers to the integral of $T_n(\phi)$ against Lebesgue measure on $B_n(u,r')$ when considered a subset of $\mathbb{R}^{\partial \rt_n}$. 

\begin{lemma} \label{lem:newupper}
Let $r \geq 0$. 
If $r ' =r + 2(u_1+u_2)$, then for all $\varphi \in B_{n+3}(u,r)$ we have
\begin{align*}
\frac{1}{n^2/2} \log T_{n+3}(\varphi) \leq \frac{1}{n^2/2} \log \int_{B_n(u,r')} T_n(\phi)\mathrm{d}\phi + C \frac{ \log(n)}{ n}.
\end{align*}
\end{lemma}
\begin{proof}
Consider the triangle $\righttriangle_{n+3}$. We can translate this triangle by $(-2,-1)$ to form the triangle
\begin{align*}
\tilde{\righttriangle}_{n+3} := (-2,-1) + \righttriangle_{n+3} := \{ (i-2,j-1) : (i,j) \in \righttriangle_{n+3} \}.
\end{align*}
We can also translate the function $\varphi:\partial \rt_{n+3} \to \mathbb{R}$ by defining $\tilde{\varphi}:\partial \tilde{\righttriangle}_{n+3} \to \mathbb{R}$ by setting $\tilde{\varphi}_x := \varphi_{x+(2,1)} - 2u_1 - u_2$ for $x \in \partial \tilde{\rt}_{n+3}$. Then the condition $|\varphi_x - \varphi^u_x| \leq r$ for all $x \in \partial\righttriangle_{n+3}$ is equivalent to the condition $|\tilde{\varphi}_x - \varphi^u_x| \leq r$ for all $x \in \partial\tilde{\righttriangle}_{n+3}$. By simple translation properties of Gelfand--Tsetlin integrals (see Lemma \ref{lem:GTprops} parts (2) and (3)), we have
\begin{align*}
T_{n+3}(\varphi) = Z_{\rt_n}(\tilde{\varphi}),
\end{align*} 
where we note that $\righttriangle_n$ is simply the interior of $\tilde{\righttriangle}_{n+3}$. We may then write
\begin{align} \label{eq:jungle}
T_{n+3}(\varphi) = \int_{\mathcal{C}}  T_n(\phi) \mathrm{d}\phi,
\end{align}
where $\mathcal{C}$ denotes the subset of $\mathbb{R}^{\partial \rt_n}$ consisting of increasing functions $\phi$ for which all of the Gelfand--Tsetlin inequalities are satisfied by the functions $\tilde{\varphi}:\partial \tilde{\rt}_{n+3}$ and $\phi:\partial \rt_n \to \mathbb{R}$ on every edge connecting elements of the neighbouring sets $\partial \tilde{\rt}_{n+3}$ and $\partial \rt_n$. Write
\begin{align*}
\partial \rt_n = D_n \cup E_n \cup S_n,
\end{align*}
where $D_n$ is the diagonal points, $E_n$ is the eastern side, and $S_n$ is the southern side, given respectively by 
\begin{align*}
D_n := \{ (i,i) : 0 \leq i \leq n \},  \quad 
E_n := \{ (n,i) : 0 \leq i \leq n \}, \quad \text{and} \quad 
S_n := \{ (i,0) : 0 \leq i \leq n \};
\end{align*}
these sets overlap at the corners but are otherwise disjoint.

\begin{figure}
\centering
\begin{tikzpicture}[scale=0.3]
  % Optional grid
  \draw[step=1,black!30,very thin] (-3,-2) grid (18,18);

  % ---- Triangles (black outlines) ----
  \draw[very thick,black] (0,0) -- (15,0) -- (15,15) -- cycle;
  \draw[very thick,black] (-2,-1) -- (16,-1) -- (16,17) -- cycle;

  % ---- Colored unit edges (30% transparent) ----
  % Yellow: (x,-1) -- (x,0), x=0..15
  \foreach \x in {0,...,15} {
    \draw[line width=2pt,yellow,opacity=0.3] (\x,-1) -- (\x,0);
  }

  % Green: (15,y) -- (16,y), y=0..15
  \foreach \y in {0,...,15} {
    \draw[line width=2pt,green,opacity=0.3] (15,\y) -- (16,\y);
  }

  % Blue: (x,x) -- (x,x+1), x=0..15
  \foreach \x in {0,...,15} {
    \draw[line width=2pt,blue,opacity=0.3] (\x,\x) -- (\x,\x+1);
  }

  % Red: (x,x) -- (x-1,x), x=1..15
  \foreach \x in {1,...,15} {
    \draw[line width=2pt,red,opacity=0.3] (\x,\x) -- (\x-1,\x);
  }

  % ---- Colored vertices ----
  % Yellow: (x,0), 0<x<15
  \foreach \x in {1,...,14} {
    \fill[yellow] (\x,0) circle (4pt);
  }

  % Green: (15,y), 0<y<15
  \foreach \y in {1,...,14} {
    \fill[green] (15,\y) circle (4pt);
  }

  % Purple: (15,0) and (x,x), 0<=x<=15
  \fill[purple] (15,0) circle (4pt);
  \foreach \x in {0,...,15} {
    \fill[purple] (\x,\x) circle (4pt);
  }
\end{tikzpicture}
\caption{The two black triangles denote $\partial \rt_n$ sitting immediately inside $\partial \tilde\rt_{n+3}$. The red/blue, yellow and green edges induce the inequalities in \eqref{eq:o1}, \eqref{eq:o2} and \eqref{eq:o3} respectively.}
\label{fig:between}
\end{figure}

Set $r ' =r  + 2(u_1+u_2)$. Then for $\phi \in \mathcal{C}$ we have
\begin{align}
\phi_x &\in [ \tilde{\varphi}_{x-(1,0)} , \tilde{\varphi}_{x + (0,1) } ] \subseteq [ \varphi^u_x - r' , \varphi^u_x + r '] \qquad &x \in D_n, \label{eq:o1}\\
\phi_x &\geq \tilde\varphi_{x - (0,1)} \geq \varphi^u_x - r' \qquad &x \in S_n, \label{eq:o2}\\
\phi_x &\leq \tilde\varphi_{x + (1,0)} \leq \varphi^u_x + r ' \qquad &x \in E_n.   \label{eq:o3}
\end{align}
For those $\phi$ in $\mathcal{C}$ we would like to have both upper and lower bounds for $\phi_x$ for all $x \in \partial \rt_n$. Taking a very generous bound, and using the fact that every value $\phi_x$ is squeezed between boundary values of $\tilde\varphi$ by monotonicity, for all $x \in \partial \rt_n$, if $\phi \in \mathcal{C}$ we have
\begin{align}  \label{eq:o4}
-r' \leq \min_{y \in \partial \tilde \rt_{n+3} } \tilde \varphi_y \leq \phi_x \leq \max_{y \in \partial \tilde \rt_{n+3} } \tilde \varphi_y \leq (u_1+u_2) n + r'.
\end{align}
Using \eqref{eq:o4} we can adapt \eqref{eq:o1}, \eqref{eq:o2} and \eqref{eq:o3} to provide (at times rather generous) upper and lower bounds in each coordinate:
\begin{align}
\phi_x &\in U_x := [ \varphi^u_x - r' , \varphi^u_x + r '] \qquad &x \in D_n, \label{eq:p1}\\
\phi_x &\in U_x := [\varphi_x^u-r', \varphi_x^u + r'']  \qquad &x \in S_n - D_n, \label{eq:p2}\\
\phi_x &\in U_x := [\varphi_x^u - r'', \varphi_x^u+ r'] \qquad &x \in E_n - ( D_n \cup S_n),  \label{eq:p3}
\end{align}
where $r'' := (u_1+u_2)n+r'$. In particular, if we let  $\mathcal{D}$ denote the collection of increasing functions $\phi:\partial\rt_n \to \mathbb{R}$ satisfying \eqref{eq:p1}, \eqref{eq:p2} and \eqref{eq:p3}, then $\mathcal{C} \subseteq \mathcal{D}$. In particular, by \eqref{eq:jungle} 
 we have
\begin{align} \label{eq:Ubound}
T_{n+3}(\varphi) \leq \int_{\mathcal{D}} T_n(\phi) \mathrm{d}\phi.
\end{align}

Recall Example \ref{ex:upperlower}, Lemma \ref{lem:pointtype}, and the surrounding discussion. In particular, there we identified monotonicity properties of Gelfand--Tsetlin partition functions $Z_A(\varphi)$ as the boundary values $\varphi:B \to \mathbb{R}$ change. In our case, so long as the function $\phi: \partial\rt_n \to \mathbb{R}$ is increasing, the functional $T_n(\phi)$ is increasing in $\phi_x$ for each $x \in E_n$, and decreasing in $\phi_x$ for each $x \in S_n$. (The functional $T_n(\phi)$ does not depend on the value of $\phi_x$ at $x = (n,0) \in S_n \cap E_n$, which in the language of that earlier section is an irrelevant point.) With this picture in mind, given a function $\phi:\partial\rt_n \to \mathbb{R}$, let us define a second function $\phi':\partial\rt_n \to \mathbb{R}$ by 
\begin{align*}
\phi'_x &:= \phi_x \qquad &x \in D_n,\\
\phi'_x &:= \varphi^u_x - r' \qquad &x \in S_n - D_n,\\
\phi'_x &:= \varphi^u_x + r' \qquad &x \in E_n - (D_n \cup S_n).
\end{align*}
Then, noting that $S_n - D_n$ consists only of lower points, and $E_n - (D_n \cup S_n)$ consists only of upper points, by Lemma \ref{lem:pointtype} we have $T_n(\phi') \geq T_n(\phi)$. In fact, more generally, for any $\lambda \in [0,1]$, we have 
\begin{align} \label{eq:xasd}
T_n( \lambda \phi + (1-\lambda)\phi' ) \geq T_n(\phi).
\end{align}

Let $\lambda_0 := 2r'/(r''+r')$. Now note that the map $\phi \mapsto \lambda_0 \phi + (1-\lambda_0)\phi' $ is an injection from $\mathcal{D}$ into $B_n(u,r')$. In particular, using \eqref{eq:xasd} to obtain the first inequality below, then using the fact that this map is an injection from $\mathcal{D}$ to $B_n(u,r')$ which rescales Lebesgue measure on $\mathbb{R}^{\partial \rt_n}$ by a factor of $\lambda_0^{ - \# ( E_n \cup S_n - D_n) } = \lambda_0^{ - (2n-1)}$, we have 
\begin{align*}
\int_{\mathcal{D}} T_n(\phi) \mathrm{d}\phi \leq \int_{\mathcal{D}} T_n(\lambda_0 \phi + (1-\lambda_0)\phi')  \mathrm{d}\phi \leq \left( \frac{r''+r'}{2r'} \right)^{ 2n-1} \int_{B_n(u,r')} T_n(\phi)\mathrm{d}\phi.
\end{align*}
Since $r''=(u_1+u_2)n+r'$ and $r'=r+2(u_1+u_2)\geq 2(u_1+u_2)$, we have
\[
\frac{r''+r'}{2r'} = 1 + \frac{(u_1+u_2)n}{2r'} \leq 1+\frac n4 \leq n
\]
for all $n \geq 2$. The result now follows from \eqref{eq:Ubound}.
\end{proof}

%%%%%%%%%%%%%%%%%%%%%%%%%%%%%%%%%%%%%%%%%%%%%%%%%%%%%%%%%%%%%%%%%%%%%%%%%%%%%%%%%%%%%%%%%
\subsection{An integral lower bound}
%%%%%%%%%%%%%%%%%%%%%%%%%%%%%%%%%%%%%%%%%%%%%%%%%%%%%%%%%%%%%%%%%%%%%%%%%%%%%%%%%%%%%%%%%

For integers $m,n \in \mathbb{N}$ with $m \geq 1$, we make a slight abuse of notation and write $\tilde{\righttriangle}_{n+3m} := \righttriangle_{n+3m} + (-2m,-m)$. We can embed the triangle $\righttriangle_n$ inside $\tilde{\righttriangle}_{n+3m}$ so that there is a boundary of width $m$ on each side; see Figure \ref{fig:trianglemouth}. 

Now let $\varphi:\partial \righttriangle_{n+3m} \to \mathbb{R}$ be a $c$-increasing function in $B_{n+3m}(u,r)$. Then we can translate it to create a function $\tilde{\varphi}: \partial \tilde{\righttriangle}_{n+3m} \to \mathbb{R}$ defined by setting
\begin{align*}
\tilde{\varphi}_x := \varphi_{x+(2m,m) } - 2m u_1 - m u_2, \qquad x \in \partial \tilde{\righttriangle}_{n+3m}.
\end{align*}
Note that $\tilde{\varphi}$ is also $c$-increasing, and satisfies $|\tilde{\varphi}_x - \varphi^u_x| \leq r$ for all $x \in \partial \tilde{\righttriangle}_{n+3m}$. 

By the invariance of Gelfand--Tsetlin integrals under horizontal and vertical translation, we have
\begin{align} \label{eq:partial0}
T_{n+3m}(\varphi) = Z_{\interior{\tilde{\righttriangle}}_{n+3m}}( \tilde{\varphi} ).
\end{align}
Moreover, by considering the contribution to the integral along the coordinates indexed by the subset $\partial \rt_n$ of $\interior{\tilde{\righttriangle}}_{n+3m}$ we may write 
\begin{align} \label{eq:partial1}
T_{n+3m}(\varphi) = Z_{\interior{\tilde{\righttriangle}}_{n+3m}}( \tilde{\varphi} ) = \int_{\mathbb{R}^{\partial \rt_n} }  T_n(\phi) F(\tilde{\varphi},\phi) \mathrm{d}\phi,
\end{align}
where for $\phi:\partial \righttriangle_n \to \mathbb{R}$, $F(\tilde{\varphi},\phi)$ is the total volume of the integral in the space $\tilde{\righttriangle}^\circ_{n+3m} - \righttriangle_n$ given boundaries $\tilde{\varphi}$ on $\partial \tilde{\righttriangle}_{n+3m}$ and $\phi$ on $\partial \righttriangle_n$. The volume $F(\tilde{\varphi}, \phi)$ can be defined more precisely as follows. Let 
\begin{align*}
C := \tilde{\righttriangle}_{n+3m} - \interior{\righttriangle}_n,
 \quad \text{and} \quad
B := \partial C = \partial \tilde{\righttriangle}_{n+3m } \cup \partial \rt_n.
\end{align*}
The set $B$ consists of the green and blue points in Figure \ref{fig:trianglemouth}. Let $A$ denote the interior of $C$, the collection of orange points in Figure \ref{fig:trianglemouth}. In other words, $A = \tilde{\righttriangle}^\circ_{n+3m} - \righttriangle_n$. 
Recall that if $f_i :B_i \to \mathbb{R}, i=1,2$ are two functions defined on disjoint sets $B_1$ and $B_2$, we define their union function $f_1 \cup f_2:B_1 \cup B_2 \to \mathbb{R}$ by letting $f_1 \cup f_2 = f_1$ on $B_1$ and $f_1 \cup f_2 = f_2$ on $B_2$. In particular, given $\phi:\partial \righttriangle_n \to \mathbb{R}$ and $\tilde{\varphi}:\partial \tilde{\righttriangle}_{n+3m } \to \mathbb{R}$, we may define their union function $\phi \cup \tilde{\varphi}$ on $\partial C$. Thus, the volume function $F(\tilde{\varphi},\phi)$ is defined precisely as the partition function
\begin{align} \label{eq:voly}
F(\tilde{\varphi},\phi) := Z_{ \interior{C} } (\phi \cup \tilde{\varphi}).
\end{align}

\begin{figure}
\centering
\begin{tikzpicture}[scale=0.2]
    % Draw grid lines
    \draw[step=1,black!30,very thin] (-11,-6) grid (21,26);
    % Draw triangle
 % Add OliveGreen circles at diagonal vertices
    \foreach \x in {0, ..., 15} {
        \fill[OliveGreen] (\x,\x) circle (2.1mm);
    }
 % Add OliveGreen circles at diagonal vertices
    \foreach \x in {-10, ..., 20} {
        \fill[blue] (\x,5+\x) circle (2.1mm);
    }
 % Add OliveGreen circles at diagonal vertices
    \foreach \x in {0, ..., 15} {
        \fill[OliveGreen] (\x,0) circle (2.1mm);
    }
 % Add OliveGreen circles at diagonal vertices
    \foreach \x in {-10, ..., 20} {
        \fill[blue] (\x,-5) circle (2.1mm);
    }

 % Add OliveGreen circles at diagonal vertices
    \foreach \x in {0, ..., 15} {
        \fill[OliveGreen] (15,\x) circle (2.1mm);
    }
 % Add OliveGreen circles at diagonal vertices
    \foreach \x in {-5, ..., 25} {
        \fill[blue] (20,\x) circle (2.1mm);
    }

 % Add OliveGreen circles at diagonal vertices
    \foreach \x in {-8, ..., 19} {
        \fill[Apricot] (\x,-4) circle (1.4mm);
    }
    \foreach \x in {-7, ..., 19} {
        \fill[Apricot] (\x,-3) circle (1.4mm);
    }

    \foreach \x in {-6, ..., 19} {
        \fill[Apricot] (\x,-2) circle (1.4mm);
    }
    \foreach \x in {-5, ..., 19} {
        \fill[Apricot] (\x,-1) circle (1.4mm);
    }

 % Add OliveGreen circles at diagonal vertices
    \foreach \x in {-8, ..., 19} {
        \fill[Apricot] (\x,-4) circle (1.4mm);
    }
    \foreach \x in {-7, ..., 19} {
        \fill[Apricot] (\x,-3) circle (1.4mm);
    }

    \foreach \x in {-6, ..., 19} {
        \fill[Apricot] (\x,-2) circle (1.4mm);
    }
    \foreach \x in {-5, ..., 19} {
        \fill[Apricot] (\x,-1) circle (1.4mm);
    }
    \foreach \x in {0, ..., 20} {
        \fill[Apricot] (16,\x) circle (1.4mm);
    }    
\foreach \x in {0, ..., 21} {
        \fill[Apricot] (17,\x) circle (1.4mm);
    }
    \foreach \x in {0, ..., 22} {
        \fill[Apricot] (18,\x) circle (1.4mm);
    }
    \foreach \x in {0, ..., 23} {
        \fill[Apricot] (19,\x) circle (1.4mm);
    }

   \foreach \x in {-1, ..., 15} {
        \fill[Apricot] (\x,\x+1) circle (1.4mm);
    }

 \foreach \x in {-2, ..., 15} {
        \fill[Apricot] (\x,\x+2) circle (1.4mm);
    }
 \foreach \x in {-3, ..., 15} {
        \fill[Apricot] (\x,\x+3) circle (1.4mm);
    }
 \foreach \x in {-4, ..., 15} {
        \fill[Apricot] (\x,\x+4) circle (1.4mm);
    }

    % Label non-diagonal vertices
    \node[below left, scale=0.7] at (0,0) {\small $(0,0)$};
    \node[below right, scale=0.7] at (15,0) {\small $(n,0)$};
    \node[above right,scale=0.7] at (15,15) {\small $(n,n)$};

    \node[below left, scale=0.7] at (-10,-5) {\small $(-2m,-m)$};
    \node[below right, scale=0.7] at (20,-5) {\small $(n+m,-m)$};
    \node[above right,scale=0.7] at (20,25) {\small $(n+m,n+m)$};
\end{tikzpicture}
\caption{A copy of the triangle $\righttriangle_n$ embedded within a translated copy $\tilde{\righttriangle}_{n+3m}$ of the triangle $\righttriangle_{n+3m}$.}
\label{fig:trianglemouth}
\end{figure}

%%%%%%%%%%%%%%%%%%%%%%% 

If $f_1,f_2 :B \to \mathbb{R}$ are such that $f_i$ is $c_i$-increasing for $i=1,2$, then for $\lambda \in [0,1]$ the convex combination $\lambda f_1 + (1 - \lambda )f_2$ is $(\lambda c_1 + (1-\lambda )c_2)$-increasing. We note in particular that the linear function $\varphi^u$ is $(u_1 \wedge u_2)$-increasing, and accordingly, if $\varphi^u$ is the linear function and $\phi:B \to \mathbb{R}$ is any increasing function, then for $\delta \in [0,1]$, $(1-\delta)\phi+ \delta\varphi^u$ is $\delta(u_1\wedge u_2)$-increasing.

Write 
\begin{align*}
 \tilde{B}_{n+3m}(u,r) := \{ \tilde\varphi:\partial\tilde{\righttriangle}_{n+3m} \to \mathbb{R} \,:\, |\tilde\varphi_x - \varphi^u_x| \leq r ~\forall x \in \partial\tilde{\righttriangle}_{n+3m} , \text{$\tilde\varphi$ is increasing}\},
 \end{align*}
and let $\tilde{B}_{n+3m}(u,r,c)$ be the subset of $\tilde{B}_{n+3m}(u,r) $ consisting of $c$-increasing functions. 

\begin{lemma} \label{lem:gspace}
Let $r > 0$, let $\lambda \in (0,1/2]$, and let $m \geq 4r/(u_1 \wedge u_2)$. Let $\phi \in B_n(u,r,\lambda(u_1 \wedge u_2))$ and $\tilde{\varphi} \in \tilde{B}_{n+3m}(u,r,\lambda(u_1 \wedge u_2))$.
Then their union function $\phi \cup \tilde{\varphi}$ is also $\lambda(u_1 \wedge u_2)$-increasing. 
\end{lemma}

\begin{proof}
Write $\zeta := \phi \cup \tilde{\varphi}$ for the union function. It is clear that 
\begin{align*}
\frac{ \zeta_y - \zeta_x }{ \|y-x\|_1} \geq \lambda(u_1 \wedge u_2) 
\end{align*}
whenever $x \leq y$ with either both $x,y \in \partial \tilde{\righttriangle}_{n+3m}$ or both $x,y \in \partial \righttriangle_n$.

It remains to examine the case where one of $x$ or $y$ is in $\partial \righttriangle_n$ and the other is in $\partial \tilde{\righttriangle}_{n+3m}$. To this end, take comparable points $x$ and $y$ with $x \leq y$, one lying in $\partial \righttriangle_n$ and the other lying in $\partial \tilde{\righttriangle}_{n+3m}$. Then since both $\phi$ and $\tilde{\varphi}$ are within $r$ of the linear function $\varphi^u$ on their respective domains, using the fact that $\varphi^u$ is $(u_1 \wedge u_2)$-increasing, the fact that $\|y-x\|_1 \geq m$, and then the fact that $m \geq 4r/(u_1 \wedge u_2)$, we have 
\begin{align*}
\frac{\zeta_y - \zeta_x}{\|y-x\|_1} 
\geq \frac{\varphi^u_y - \varphi^u_x - 2 r }{ \|y-x\|_1 } 
\geq (u_1 \wedge u_2) - \frac{2r}{\|y - x\|_1} 
\geq (u_1 \wedge u_2) - \frac{2r}{m} 
\geq \frac{1}{2}(u_1 \wedge u_2).
\end{align*}
Since $\lambda \leq 1/2$, the result follows.
\end{proof}

\begin{cor} \label{cor:space}
Under the conditions of Lemma \ref{lem:gspace}, together with the additional condition $m \leq n$, we have 
\begin{align*}
F(\tilde{\varphi},\phi) \geq \exp \left\{ - C m n \left( \log(1/\lambda)+ K_u \right) \right\}.
\end{align*}
\end{cor}
\begin{proof}
By Lemma \ref{lem:gspace}, $\zeta = \tilde{\varphi} \cup \phi$ is $\lambda(u_1 \wedge u_2)$-increasing. Applying the extension Lemma \ref{lem:extension} and using the fact that $F(\tilde{\varphi},\phi) = Z_{\interior{C}}(\zeta)$, we have $F(\tilde{\varphi},\phi) \geq (\lambda(u_1 \wedge u_2))^{ \# \interior{C} }$. Now note that $\# \interior{C} = O(mn)$ whenever $m \leq n$.
\end{proof}

Our final step in preparing the proof of the lower bound for triangular integrals is to observe that there is a natural injection $B_n(u,r) \to B_n(u,(1-\lambda)r,\lambda(u_1 \wedge u_2))$ defined by the map
\begin{align*}
\phi \mapsto (1-\lambda)\phi + \lambda \varphi^u.
\end{align*}
This map scales Lebesgue measure by a factor of $(1-\lambda)^{ \# \partial \righttriangle_n } = (1-\lambda)^{3n}$, and its image is contained in $B_n(u,(1-\lambda)r,\lambda(u_1 \wedge u_2))$. Thus, for any nonnegative functional $S:\mathbb{R}^{ \partial \righttriangle_n} \to [0,\infty)$ we have
\begin{align} \label{eq:Tmap}
\int_{B_n(u,r)} S( (1-\lambda)\phi + \lambda \varphi^u) \mathrm{d}\phi \leq (1-\lambda)^{-3n} \int_{B_n(u,(1-\lambda)r,\lambda(u_1 \wedge u_2))} S(\phi) \mathrm{d}\phi.
\end{align}

We are now ready to prove a lower bound for triangular integrals of well-increasing functions in terms of Lebesgue integrals. The Pr\'ekopa--Leindler inequality for partition functions plays an integral role in its proof.

\begin{lemma} \label{lem:newlower}
Let $\log(n)^{-1/2} \leq u_1 \wedge u_2$. 
Let $r >0, \lambda \in (0,1/2]$, and let $m \in \mathbb{N}$ be an integer satisfying $4r/(u_1 \wedge u_2) \leq m \leq \lambda n$. Let $\varphi \in B_{n+3m}(u,r)$ be $\lambda(u_1 \wedge u_2)$-increasing. 
Then
\begin{align*}
\frac{1}{n^2/2} \log T_{n+3m}(\varphi) \geq \frac{1}{n^2/2} \log  \int_{B_n(u,r)} T_n(\phi)\mathrm{d}\phi - C \lambda ( \log(1/\lambda) + K_u).
\end{align*}
\end{lemma}

\begin{proof}
Let $\varphi:\partial\rt_{n+3m} \to \mathbb{R}$ be an element of $B_{n+3m}(u,r)$ and be $\lambda (u_1 \wedge u_2)$-increasing. By restricting the domain of integration in \eqref{eq:partial1} to the subset $B_n(u,r,\lambda(u_1 \wedge u_2))$ of $\mathbb{R}^{ \partial \rt_n}$ to obtain the first inequality below, then using Corollary \ref{cor:space} to obtain the second, we have 
\begin{align*}
T_{n+3m}(\varphi) 
&\geq \int_{ B_n(u,r,\lambda(u_1 \wedge u_2)) }  T_n(\phi) F(\tilde{\varphi},\phi) \mathrm{d}\phi \\
&\geq \exp\left\{ - C m n \left( \log(1/\lambda)+ K_u \right) \right\} \int_{ B_n(u,r,\lambda(u_1 \wedge u_2)) }  T_n(\phi) \mathrm{d}\phi.
\end{align*}
Using the inclusion $B_n(u,r,\lambda(u_1 \wedge u_2)) \supseteq B_n(u,(1-\lambda)r,\lambda(u_1 \wedge u_2))$, \eqref{eq:Tmap}, the fact that $\lambda \leq 1/2$, and $m \leq \lambda n$, we obtain
\begin{align*}
T_{n+3m}(\varphi) 
\geq 2^{-3n} \exp \left\{ - C \lambda n^2 \left( \log(1/\lambda)+ K_u \right) \right\} \int_{B_n(u,r)  } T_n((1-\lambda)\phi+\lambda \varphi^u )\mathrm{d}\phi.
\end{align*}

Now by the Pr\'ekopa--Leindler inequality for partition functions of random surfaces, Proposition \ref{prop:pl} part (1), we have $T_n((1-\lambda)\phi+\lambda \varphi^u ) \geq T_n(\phi)^{1-\lambda}T_n(\varphi^u)^\lambda$. Accordingly, we have
\begin{align*}
T_{n+3m}(\varphi)
\geq \exp \left\{ - C \lambda n^2 \left( \log(1/\lambda)+ K_u \right) \right\} T_n(\varphi^u )^\lambda \int_{B_n(u,r)  } T_n(\phi)^{1-\lambda} \mathrm{d}\phi,
\end{align*}
where we have also swallowed the $2^{-3n}$ term using the fact that $\lambda \geq 1/n$.
Dividing through by the supremum of $T_n(\phi)$ over the ball $B_n(u,r)$, we obtain
\begin{align} \label{eq:nara}
T_{n+3m}(\varphi) \geq \exp \left\{ - C \lambda n^2 \left( \log(1/\lambda)+ K_u \right) \right\}\left( \frac{ T_n(\varphi^u )}{ \sup_{ \phi \in B_n(u,r) } T_n(\phi) } \right)^\lambda  \int_{B_n(u,r)  } T_n(\phi)\mathrm{d}\phi.
\end{align}

Using Lemma \ref{lem:volupper} with $a = -r/n$ and $b = u_1 + u_2 + r/n$ and $C = \rt_n$, we have 
\begin{align} \label{eq:supy}
\sup_{ \phi \in B_n(u,r) } T_n(\phi) \leq e^{ C n^2 ( 1 + \log(u_1 + u_2 + 2r/n)) } \leq e^{ C' n^2 K_u },
\end{align}
where to obtain the second inequality above, we have used the fact that $r \leq n(u_1 \wedge u_2)/8$ from the statement. Combining \eqref{eq:supy} with Lemma \ref{lem:charac} we obtain 
\begin{align} \label{eq:nara2}
 \frac{ T_n(\varphi^u )}{ \sup_{ \phi \in B_n(u,r) } T_n(\phi) } \geq e^{ -  C K_u n^2}.
\end{align}

Using \eqref{eq:nara2} in \eqref{eq:nara} and then taking logarithms of both sides, we obtain
\begin{align} \label{eq:nara3}
\frac{1}{n^2/2} \log T_{n+3m}(\varphi) \geq \frac{1}{n^2/2} \log \int_{B_n(u,r)} T_n(\phi) \mathrm{d}\phi - C \lambda \left( \log \frac{1}{\lambda} + K_u \right),
\end{align}
completing the proof.
\end{proof}

%%%%%%%%%%%%%%%%%%%%%%%%%%%%%%%%%%%%%%%%%%%%%%%%%%%%%%%%%%%%%%%%%%%%%%%%%%%%%%%%%%%%%%%%%
\subsection{Proof of Theorem \ref{thm:upperlower}}
%%%%%%%%%%%%%%%%%%%%%%%%%%%%%%%%%%%%%%%%%%%%%%%%%%%%%%%%%%%%%%%%%%%%%%%%%%%%%%%%%%%%%%%%%

We now use the results Lemma \ref{lem:newupper} and Lemma \ref{lem:newlower} to relate integrals of $T_n(\phi)$ to the surface tension $\sigma(u_1,u_2)$. 

\begin{lemma}
Let $n \geq 10$ and let $\log(n)^{-1/2} \leq u_1 \wedge u_2$.
For any $0 < r \leq \frac{n}{10}(u_1 \wedge u_2)$ we have the upper bound
\begin{align} \label{eq:prelow2}
\frac{1}{n^2/2}\log \int_{B_n(u,r)} T_n(\phi) \mathrm{d}\phi \leq - \sigma(u_1,u_2) + C \lambda_r (\log (1/\lambda_r) + K_u),
\end{align}
where $\lambda_{r} = \frac{1}{n} \lceil 4r/(u_1 \wedge u_2) \rceil$. 

Alternatively, for any $r \geq 2(u_1+u_2)$ we have the lower bound
\begin{align} \label{eq:prehigh}
\frac{1}{n^2/2}\log \int_{B_n(u,r)} T_n(\phi) \mathrm{d}\phi \geq - \sigma(u_1,u_2) - C (\log(n)+K_u)/n.
\end{align}
\end{lemma}

\begin{proof}
First, we prove \eqref{eq:prelow2}. Given some value $0 < r \leq \frac{n}{10}(u_1 \wedge u_2)$, consider setting $\varphi = \varphi^u$ in Lemma \ref{lem:newlower}, which is certainly an element of $B_{n+3m}(u,r)$ for any $r > 0$. Moreover, $\varphi^u$ is $\lambda(u_1\wedge u_2)$-increasing for any choice of $\lambda \in [0,1/2]$. Let 
\[
m = \left\lceil \frac{4r}{u_1\wedge u_2} \right\rceil,
\]
and note that $\lambda_r$ defined in the statement is the ratio $\lambda_r = m/n$.
The condition $r \leq \frac{n}{10} (u_1 \wedge u_2)$ guarantees $\lambda_r \leq 1/2$, so that by using Lemma \ref{lem:newlower} and rearranging we have 
\begin{align} \label{eq:prelow0}
\frac{1}{n^2/2} \log \int_{B_n(u,r)} T_n(\phi) \mathrm{d} \phi  \leq \frac{1}{n^2/2} \log T_{n+3m}(\varphi^u) + C\lambda_{r} \left( \log (1/\lambda_r) + K_u \right).
\end{align}
Using the upper bound for $T_{n+3m}(u_1,u_2) = T_{n+3m}(\varphi^u)$ in Theorem \ref{thm:linear} (and noting that $n+3m = (1 +3 \lambda_r)n$) we obtain 
\begin{align} \label{eq:prelow}
\frac{1}{n^2/2} \log \int_{B_n(u,r)} T_n(\phi) \mathrm{d} \phi  \leq -(1+3\lambda_r)^2 \sigma(u_1,u_2)  + C\lambda_{r} \left( \log (1/\lambda_r)  + K_u \right),
\end{align}
where we have used the fact that $\lambda_r \geq 1/n$ to absorb the additional error terms into the right-hand side. Using the bound in \eqref{eq:able}, together with the fact that $(1+3\lambda_r)^2 = 1 + O(\lambda_r)$, and changing the constant if necessary, we can absorb the $(1+3\lambda_{r})^2$ factor into the error term to obtain \eqref{eq:prelow2} as written.

To prove \eqref{eq:prehigh} we may assume without loss of generality that $r = 2(u_1 + u_2)$, since $\int_{B_n(u,r)} T_n(\phi) \mathrm{d}\phi$ increases with $r$. In this direction, setting $\varphi = \varphi^u$ and taking $r = 0$ in Lemma \ref{lem:newupper} (so that $r' = 2(u_1+u_2)$), and then using Theorem \ref{thm:linear} we obtain \eqref{eq:prehigh}.
\end{proof}
%%%%%%%%%%%%%%%%%%%%
We are now ready to complete the proof of Theorem \ref{thm:upperlower}.

\begin{proof}[Proof of Theorem \ref{thm:upperlower}]
First we prove the upper bound \eqref{eq:upper1}. Let us first observe that we may assume without loss of generality that $(u_1 + u_2) \leq \frac{n}{80}(u_1 \wedge u_2)$. If this condition is not satisfied, then $\rho_r \geq 1/80$, and the bound \eqref{eq:upper1} follows (with a sufficiently large universal constant $C>0$) as a consequence of the rough upper bound in \eqref{eq:supy}.

Under the assumption $(u_1 + u_2) \leq \frac{n}{80}(u_1 \wedge u_2)$, we will initially prove \eqref{eq:upper1} under the additional assumption $r \leq n(u_1 \wedge u_2)/40$, and thereafter show we can lift this additional assumption. Let $\varphi \in B_n(u,r)$. Replacing $n$ with $n-3$ in Lemma \ref{lem:newupper} (but using $n^2/2$ in place of $(n-3)^2/2$ as a denominator), we obtain 
\begin{align*}
\frac{1}{n^2/2} \log T_n(\varphi) \leq \frac{1}{n^2/2} \log \int_{B_{n-3}(u,r')} T_{n-3}(\phi)\mathrm{d}\phi + C \frac{\log(n)}{n},
\end{align*}
where $r' = r + 2(u_1 + u_2)$. These assumptions above entail that $r' \leq \frac{n}{20}(u_1 \wedge u_2) \leq \frac{n-3}{10}(u_1 \wedge u_2)$ (for all $n \geq 7$), so that by \eqref{eq:prelow2} we have
\begin{align}
\frac{1}{(n-3)^2/2}\log \int_{B_{n-3}(u,r')} T_{n-3}(\phi) \mathrm{d}\phi \leq - \sigma(u_1,u_2) + C \lambda_{r'} (\log (1/\lambda_{r'}) + K_u),
\end{align}
where $\lambda_{r'} = \frac{1}{n-3} \lceil 4r'/(u_1 \wedge u_2) \rceil$. Since $(n-3)^2 = n^2 + O(n)$, the inequality \eqref{eq:upper1}, in the case $r \leq n(u_1 \wedge u_2)/40$, now follows from noting that $\lambda_{r'} \leq 10 \rho_r$. 

Still working under the assumption $(u_1+u_2) \leq \frac{n}{80}(u_1 \wedge u_2)$, to extend the inequality to the case where $r > n(u_1 \wedge u_2)/40 =: r_0$, one may rearrange the Pr\'ekopa--Leindler inequality \eqref{eq:pfpl} to write
\begin{align} \label{eq:myriad}
\frac{1}{n^2/2} \log T_n(\varphi) \leq \frac{1}{\mu} \left[ \frac{1}{n^2/2} \log T_n(\mu \varphi + (1-\mu)\varphi^u ) - (1-\mu) \cdot \frac{1}{n^2/2} \log T_n(\varphi^u) \right].
\end{align} 
Now setting $\mu = r_0/r$, the function $\mu \varphi + (1-\mu)\varphi^u $ lies in $B_n(u,r_0)$, and one may apply the upper bound for $T_n(\mu \varphi + (1-\mu)\varphi^u )$, and combine it with the lower bound for $T_n(\varphi^u) $ in Theorem \ref{thm:linear}, to obtain \eqref{eq:upper1} for arbitrarily large $r > 0$.
\medskip

We start working towards \eqref{eq:lower1}. First, note that since $\rho_0 > 0$, it is no loss of generality to assume that $r \geq 2(u_1 + u_2)$, since, by the Pr\'ekopa-Leindler inequality, if the result holds for such $r$ it also holds for smaller $r$. Suppose now that $n,m,\lambda,r$ are as in the statement of Lemma \ref{lem:newlower}, and additionally suppose $r \geq 2(u_1+u_2)$. Combining Lemma \ref{lem:newlower} with the bound in \eqref{eq:prehigh} we obtain
\begin{align*}
\frac{1}{n^2/2} \log T_{n+3m}(\varphi) \geq - \sigma(u_1,u_2) - C \lambda ( \log (1/\lambda) + K_u),
\end{align*}
where we have used the fact that $\lambda \geq 1/n$ to absorb error terms. Using the fact that $m \leq \lambda n$ and $\lambda \leq 1/2$, we may replace the denominator $n^2/2$ by $(n+3m)^2/2$, and then use \eqref{eq:able} to obtain
\begin{align} \label{eq:prelower4a}
\frac{1}{(n+3m)^2/2} \log T_{n+3m}(\varphi) \geq - \sigma(u_1,u_2) -C \lambda ( \log (1/\lambda) + K_u),
\end{align}
for a possibly different universal $C > 0$. 

We now establish \eqref{eq:lower1} as written, but with $N$ in place of $n$. Now let $N >0$ be an integer, $r \leq N(u_1 \wedge u_2)/100$, and suppose that $\varphi \in B_N(u,r)$ is $\frac{10r}{N} := \lambda_r (u_1 \wedge u_2)$-increasing, with $\lambda_r = \frac{10 r }{N(u_1 \wedge u_2)} \leq 1/10$. We may write
\begin{align*}
N = n + 3m \qquad \text{where} \quad m: = \left\lceil \frac{4r}{u_1 \wedge u_2} \right\rceil, \qquad n := N - 3m.
\end{align*}
 Note that we have $m \leq \lambda n$ in this case, and hence \eqref{eq:prelower4a} holds with $N = n+3m$. The inequality \eqref{eq:lower1} now follows from noting that $\lambda_r \leq 10 \rho_r$. 
 
\end{proof}
%%%%%%%%%%%%%%%%%%%%%%%%%%%%%%%%%%%%%%%%%%%%%%%%%%%%%%%%%%%%%%%%%%%%%%%%%%%%%%%%%%%%%%%%%%%
\subsection{Proof of Theorem \ref{thm:B}}
%%%%%%%%%%%%%%%%%%%%%%%%%%%%%%%%%%%%%%%%%%%%%%%%%%%%%%%%%%%%%%%%%%%%%%%%%%%%%%%%%%%%%%%%%%%
We have established that partition functions of Gelfand--Tsetlin surfaces with triangular boundary conditions at an asymptotic tilt $(u_1,u_2)$ have their asymptotic behaviour governed by the surface tension function $\sigma(u_1,u_2)$. In this section we do the same with the square $\Box_n$, ultimately completing the proof of Theorem \ref{thm:B}. 

For a function $\varphi:\partial \Box_n \to \mathbb{R}$ define
\begin{align*}
S_n(\varphi) := \int_{\mathbb{R}^{\Box_n} } \mathrm{1}_{\mathrm{GT}}^{\Box_n}(\phi)   \prod_{x \in \partial \Box_n} \delta_{\varphi_x }(\mathrm{d}\phi_x) \prod_{x \in \interior{\Box}_n} \mathrm{d}\phi_x.
\end{align*}
In particular, set
\begin{align*}
S_n(u_1,u_2) := S_n(\varphi^u).
\end{align*}

Theorem \ref{thm:B} states that $\frac{1}{n^2} \log S_n(u_1,u_2)$ converges to $-\sigma(u_1,u_2)$ as $n \to \infty$. 

Define
\begin{align*}
B^{\mathrm{square}}_n(u,r) := \{ \varphi:\partial \Box_n \to \mathbb{R} : |\varphi_x - \varphi^u_x| \leq r ~\forall x \in \partial \Box_n, ~ \text{$\varphi$ is increasing} \}.
\end{align*}

\begin{proof}[Proof of Theorem \ref{thm:B}]
Assume the conditions of Theorem \ref{thm:upperlower}, so that in particular $n$ is sufficiently large that 
$\log(n)^{-1/2}  \leq u_1 \wedge u_2$. 
We begin by establishing a lower bound for $S_n(u_1,u_2)$. In fact, we will work more generally, and establish a lower bound for $S_n(\varphi)$ for all functions $\varphi:\partial \Box_n \to \mathbb{R}$ in $B^{\mathrm{square}}_n(u,(u_1 \wedge u_2)/4)$.

For short, let us write
\begin{align*}
B_n := \partial \Box_n \quad \text{and} \quad B_{n-2}' := \partial \Box_{n-2} + (1,1).
\end{align*}
Note that every element of $B_{n-2}'$ neighbours an element of $B_n$. Also write
\begin{align*}
\tilde{B}^{\mathrm{square}}_{n-2}(u,(u_1 \wedge u_2)/4) := \{ \varphi' : B_{n-2}' \to \mathbb{R} : |\varphi'_x - \varphi^u_x| \leq (u_1 \wedge u_2)/4 ~ \forall x \in B_{n-2}'\}.
\end{align*}
Note that the Gelfand--Tsetlin inequalities are automatically satisfied between any $\varphi\in B^{\mathrm{square}}_n(u,(u_1 \wedge u_2)/4)$ and $\varphi' \in \tilde{B}^{\mathrm{square}}_{n-2}(u,(u_1 \wedge u_2)/4)$. Let us thus record for later the lower bound
\begin{align*}
S_n(\varphi) \geq \int_{\tilde{B}^{\mathrm{square}}_{n-2}(u,(u_1 \wedge u_2)/4)} S_{n-2}(\varphi') \mathrm{d} \varphi',
\end{align*}
where $S_{n-2}(\varphi')$ is a natural abuse of notation denoting the partition function of the translated square $\Box_{n-2} + (1,1)$ with boundary conditions $\varphi':B_{n-2}' \to \mathbb{R}$. 

We now refine this bound further by placing two right-angled triangles $T^{(1)}$ and $T^{(2)}$ with respective side lengths $n-2$ and $n-3$ inside the square $\Box_n$. The triangle $T^{(1)}$ has vertices at $(1,1), (n-1,1), (n-1,n-1)$, and the triangle $T^{(2)}$ has vertices at $(1,2), (1,n-1), (n-2,n-1)$; see the left panel of Figure \ref{fig:combined}.

\begin{figure}[t]
\centering

\begin{minipage}{0.48\textwidth}
\centering

\begin{tikzpicture}[scale=0.28, line cap=round, line join=round]
  \def\n{20}

  % Square boundary
  \draw[very thick] (0,0) rectangle (\n,\n);

  % Blue lattice points in the square
  \foreach \i in {0,...,\n}{
    \foreach \j in {0,...,\n}{
      \fill[blue] (\i,\j) circle (0.14);
    }
  }

  % Green triangles
  \filldraw[
    fill=green!60!white,
    fill opacity=0.3,
    draw=green!80!black,
    thick
  ]
    (1,1) -- ({\n-1},1) -- ({\n-1},{\n-1}) -- cycle;

  \filldraw[
    fill=green!60!white,
    fill opacity=0.3,
    draw=green!80!black,
    thick
  ]
    (1,2) -- (1,{\n-1}) -- ({\n-2},{\n-1}) -- cycle;

  % Labels for square vertices
  \node[below left] at (0,0) {$(0,0)$};
  \node[below right] at (\n,0) {$(n,0)$};
  \node[above right] at (\n,\n) {$(n,n)$};
  \node[above left] at (0,\n) {$(0,n)$};

  % Labels
  \node[
    font=\bfseries\fontsize{20}{22}\selectfont,
    green!50!black
  ] at ({0.68*\n},{0.30*\n}) {$T^{(1)}$};

  \node[
    font=\bfseries\fontsize{20}{22}\selectfont,
    green!50!black
  ] at ({0.23*\n},{0.72*\n}) {$T^{(2)}$};

\end{tikzpicture}

\end{minipage}
\hfill
\begin{minipage}{0.48\textwidth}
\centering

\begin{tikzpicture}[scale=0.28, line cap=round, line join=round]
  \def\n{20}
  \def\m{8}

  % Outer right triangle boundary
  \draw[very thick] (0,0) -- (\n,0) -- (\n,\n) -- cycle;

  % Blue lattice points
  \foreach \i in {0,...,\n}{
    \foreach \j in {0,...,\i}{
      \fill[blue] (\i,\j) circle (0.14);
    }
  }

  % Inner green square
  \filldraw[
    fill=green!60!white,
    fill opacity=0.3,
    draw=green!80!black,
    thick
  ]
    ({\n-1},1)
    -- ({\n-1},{\m+1})
    -- ({\n-(\m+1)},{\m+1})
    -- ({\n-(\m+1)},1)
    -- cycle;

  % Triangle T^(1)
  \filldraw[
    fill=green!60!white,
    fill opacity=0.3,
    draw=green!80!black,
    thick
  ]
    (2,1) -- ({2+\m},1) -- ({2+\m},{1+\m}) -- cycle;

  % Triangle T^(2)
  \filldraw[
    fill=green!60!white,
    fill opacity=0.3,
    draw=green!80!black,
    thick
  ]
    ({\n-(\m+1)},{\m+2})
    -- ({\n-1},{\m+2})
    -- ({\n-1},{\n-2}) -- cycle;

  % Labels for vertices
  \node[below left] at (0,0) {$(0,0)$};
  \node[below right] at (\n,0) {$(n,0)$};
  \node[above right] at (\n,\n) {$(n,n)$};

  % Labels
  \node[
    font=\bfseries\fontsize{18}{20}\selectfont,
    green!50!black
  ] at ({0.28*\n},{0.18*\n}) {$T^{(1)}$};

  \node[
    font=\bfseries\fontsize{18}{20}\selectfont,
    green!50!black
  ] at ({0.83*\n},{0.63*\n}) {$T^{(2)}$};

  \node[
    font=\bfseries\fontsize{22}{24}\selectfont,
    green!50!black
  ] at ({0.70*\n},{0.28*\n}) {$S$};

\end{tikzpicture}

\end{minipage}

\caption{
Left: two right-angled triangles placed inside the square $\Box_n$. 
Right: a square and two right-angled triangles placed inside a larger triangular region.
}
\label{fig:combined}
\end{figure}

Consider any pair of functions $\phi^{(1)}:\partial T^{(1)} \to \mathbb{R}$ and $\phi^{(2)}:\partial T^{(2)} \to \mathbb{R}$ satisfying
\begin{align*}
|\phi^{(i)}_x - \varphi^u_x| \leq r := \frac{u_1 \wedge u_2}{4} \qquad x \in \partial T^{(i)}.
\end{align*}
Since $\varphi$, $\phi^{(1)}$ and $\phi^{(2)}$ all lie within distance $(u_1 \wedge u_2)/4$ of the linear function $\varphi^u$, their union function automatically satisfies the Gelfand--Tsetlin inequalities. More precisely, if $
\zeta := \varphi \cup \phi^{(1)} \cup \phi^{(2)}
$
is defined on 
$
J := \partial \Box_n \cup \partial T^{(1)} \cup \partial T^{(2)},
$
then $\zeta_x \leq \zeta_y$ whenever $x \leq y$ in $J$. Integrating over the respective interiors of $T^{(1)}$ and $T^{(2)}$, we obtain
\begin{align} \label{eq:citrus}
S_n(\varphi) \geq \left( \frac{u_1 \wedge u_2}{2} \right)^{ \# \partial T^{(1)} + \# \partial T^{(2)} } \inf_{ \phi \in B_{n-2}(u,(u_1 \wedge u_2)/4) } T_{n-2}(\phi) \inf_{ \phi \in B_{n-3}(u,(u_1 \wedge u_2)/4 ) } T_{n-3}(\phi).
\end{align}
Any function $\phi$ in $B_{n-2}(u,(u_1 \wedge u_2)/4)$ or $B_{n-3}(u,(u_1 \wedge u_2)/4)$ is automatically $(u_1 \wedge u_2)/2$-increasing. Applying \eqref{eq:lower1} with $r = (u_1 \wedge u_2)/4$, and using \eqref{eq:citrus}, we obtain, for $\varphi \in B^{\mathrm{square}}_n( u, (u_1 \wedge u_2)/4)$,
\begin{align} \label{eq:squarelower}
\frac{1}{n^2} \log S_n(\varphi) \geq - \sigma(u_1, u_2) - C \rho ( \log(1/\rho)+K_u),
\end{align}
where
\begin{align*}
\rho := \frac{u_1+u_2+(u_1 \wedge u_2)/4}{(u_1\wedge u_2)n}.
\end{align*}
In particular, setting $\varphi := \varphi^u$ and taking $\liminf$s, we obtain
\begin{align*}
\liminf_{n \to \infty} \frac{1}{n^2} \log S_n(u_1,u_2) \geq - \sigma(u_1,u_2).
\end{align*}

We now work towards the corresponding upper bound. Note that for $\varphi:B_n \to \mathbb{R}$, by integrating over the translated inner boundary $B_{n-2}'$, we have
\begin{align*}
S_n(\varphi) = \int_{\mathbb{R}^{B_{n-2}' }} \mathrm{1}_{\mathrm{GT}}^{B_n \cup B_{n-2}' } ( \varphi \cup \varphi' ) S_{n-2}(\varphi') \mathrm{d}\varphi',
\end{align*}
where $\mathrm{1}_{\mathrm{GT}}^{B_n \cup B_{n-2}' } ( \varphi \cup \varphi' )$ denotes the indicator function that all Gelfand--Tsetlin inequalities between $\varphi:B_n \to \mathbb{R}$ and $\varphi':B_{n-2}' \to \mathbb{R}$ are satisfied. Arguing as in the proof of Lemma \ref{lem:newupper}, one can show that
\begin{align} \label{eq:suite1}
S_n(\varphi) \leq \left( \frac{ C n(u_1 + u_2) }{ u_1 \wedge u_2} \right)^{O(n)} \int_{\tilde{B}_{n-2}^{\mathrm{square}}(u,(u_1 \wedge u_2)/4)} S_{n-2}(\varphi')\mathrm{d}\varphi'. 
\end{align}
It therefore remains to bound the integral on the right-hand side. To this end, note that we may place the translated square $\Box_{n-2} + (1,1)$ together with two translated triangles inside a larger triangle with side length $2n$; see the right panel of Figure \ref{fig:combined}. Integrating out over the complementary region, and arguing exactly as in the proof of the lower bound above, one obtains
\begin{align} \label{eq:suite2}
T_{2n}(u_1,u_2) \geq \left( \frac{u_1 \wedge u_2}{2} \right)^{ O(n) } \left( \int_{\tilde{B}_{n-2}^{\mathrm{square}}(u,(u_1 \wedge u_2)/4)} S_{n-2}(\varphi')\mathrm{d}\varphi' \right) \left( \int_{B_{n-2}(u,(u_1 \wedge u_2)/4)} T_{n-2}(\phi)\mathrm{d}\phi \right)^2.
\end{align}
Combining \eqref{eq:suite1} and \eqref{eq:suite2}, using the upper bound for $T_{2n}(u_1,u_2)$ from Theorem \ref{thm:linear}, and using the lower bound \eqref{eq:prehigh} for the triangular partition function integrals, we obtain
\begin{align*}
\frac{1}{n^2} \log S_n(\varphi) \leq -  \sigma(u_1, u_2) +C \rho ( \log(1/\rho)+K_u).
\end{align*}
In particular, setting $\varphi := \varphi^u$ and taking $\limsup$s, we obtain
\begin{align*}
\limsup_{n \to \infty} \frac{1}{n^2} \log S_n(u_1,u_2) \leq - \sigma(u_1,u_2),
\end{align*}
completing the proof of Theorem \ref{thm:B}.
\end{proof}

%%%%%%%%%%%%%%%%%%%%%%%%%%%%%%%%%%%%%%%%%%%%%%
%%%%%%%%%%%%%%%%%%%%%%%%%%%%%%%%%%%%%%%%%%%%%%
%%%%%%%%%%%%%%%%%%%%%%%%%%%%%%%%%%%%%%%%%%%%%%
%%%%%%%%%%%%%%%%%%%%%%%%%%%%%%%%%%%%%%%%%%%%%%
\section{The metric space of compressions and upper semicontinuity of entropy} \label{sec:metric}
%%%%%%%%%%%%%%%%%%%%%%%%%%%%%%%%%%%%%%%%%%%%%%
%%%%%%%%%%%%%%%%%%%%%%%%%%%%%%%%%%%%%%%%%%%%%%
%%%%%%%%%%%%%%%%%%%%%%%%%%%%%%%%%%%%%%%%%%%%%%

We begin by recalling some definitions. We call a function $\psi:U \subseteq \mathbb{R}^2 \to \mathbb{R}$ defined on a subset $U$ of $\mathbb{R}^2$ increasing if it is nondecreasing in both variables. An increasing surface is a function $\psi:\rt \to \mathbb{R}$ defined on $\rt := \{(s,t) \in \mathbb{R}^2 : 0 \leq t \leq s \leq 1 \}$
that is upper-semicontinuous and increasing. By Lebesgue's monotone differentiation theorem, the partial derivatives $\psi_s$ and $\psi_t$ of an increasing surface exist almost everywhere in $\rt$. We endow the set of increasing surfaces with the $L^1$ metric, given by 
\begin{align*}
\mathrm{d}(\psi,\tilde{\psi}) := \int_{\rt} |\tilde{\psi}(s,t) - \psi(s,t)| \mathrm{d}s\mathrm{d}t.
\end{align*}
The compression entropy $\mathcal{H}[\psi] \in [-\infty,+\infty]$ of an increasing surface is the integral
\begin{align*}
\mathcal{H}[\psi] := - \int_{\rt} \sigma( \nabla \psi(s,t)) \mathrm{d}s\mathrm{d}t,
\end{align*}
where
\begin{align} \label{eq:sigmaagain}
\sigma(u_1,u_2) := - \log(u_1+u_2) - \log \sin \left( \frac{\pi u_1}{u_1+u_2}\right) - 1 + \log \pi.
\end{align}

The main result of this section is the following:

\begin{thm} \label{thm:lowersemi}
Let $\mathcal{I}$ be the set of increasing surfaces $\psi:\rt \to [0,1]$ taking values in $[0,1]$. Then the entropy functional
\begin{align*}
\mathcal{H}:\mathcal{I} \to [-\infty,\infty)
\end{align*}
is upper-semicontinuous in the $L^1$ metric. In other words, if $(\psi_n)_{n\geq 1}$ is a sequence of elements of $\mathcal{I}$ such that $\mathrm{d}(\psi_n,\psi) \to 0$, then 
\begin{align*}
\limsup_{n \to \infty} \mathcal{H}[\psi_n ] \leq \mathcal{H}[\psi].
\end{align*}
\end{thm}

A crucial fact that will play an important role in this section is the following:

\begin{proposition}
The surface tension function $\sigma:\mathbb{R}_{>0}^2 \to (-\infty,\infty)$ is convex. 
\end{proposition}
\begin{proof}
This can be proved directly using the formula \eqref{eq:sigmaagain} for $\sigma(u_1,u_2)$. Alternatively, it follows immediately from Corollary \ref{cor:STconvex}. 
\end{proof}
%%%%%%%%%%%%%%%%%%%%%%%%%%%%%%%%%%%%%%%%%%%%%%%%%%%%%%%%%%%%%%%%%%%%%%%%%%%%%%%%
\subsection{Smoothed surfaces}
%%%%%%%%%%%%%%%%%%%%%%%%%%%%%%%%%%%%%%%%%%%%%%%%%%%%%%%%%%%%%%%%%%%%%%%%%%%%%%%%

Note that while $\mathrm{d}(\psi^n,\psi) \to 0$ guarantees the convergence of $\psi^n(s,t)$ to $\psi(s,t)$ for almost all $(s,t) \in \rt$, it makes no guarantee about the almost-everywhere convergence of the gradients $\nabla \psi^n(s,t)$ to $\nabla \psi(s,t)$. In order to obtain the convergence of gradients, we will see that it is necessary to take some sort of smoothing of the functions. 

For $0 < \varepsilon < 1/6$, let $\hollow^\varepsilon := \{(s,t) \in \rt : s \leq 1-\varepsilon, t \geq \varepsilon, t \leq s - 2 \varepsilon \}$. Then the square with side length $2\varepsilon$ and center at $(s,t)$ lies in $\rt$. 
The non-hypotenuse sides of $\hollow^\varepsilon$ have length $1-4\varepsilon$, so that its Lebesgue measure is given by 
\begin{align*}
|\hollow^\varepsilon| = \frac{1}{2}(1-4\varepsilon)^2.
\end{align*}

Throughout the remainder of the article, smoothed versions of increasing surfaces, obtained by averaging an increasing surface over a small square, will play an important role:

\begin{df}[The $\varepsilon$-smoothing and double $\varepsilon$-smoothing of $\psi:\rt \to \mathbb{R}$] \label{df:smoothing}
Let $0< \varepsilon < 1/6$. 
Given an increasing surface $\psi:\rt \to [0,1]$ define its $\varepsilon$-smoothing $\psi^\varepsilon:\hollow^\varepsilon \to \mathbb{R}$ by 
\begin{align*}
\psi^\varepsilon(s,t) := (2\varepsilon)^{-2} \int_{s-\varepsilon}^{s+\varepsilon} \int_{t- \varepsilon}^{t+\varepsilon} \psi(s',t')\mathrm{d}s'\mathrm{d}t', \qquad (s,t) \in \hollow^\varepsilon.
\end{align*} 
The double $\varepsilon$-smoothing of $\psi:\rt \to [0,1]$ is the function $\psi^{\varepsilon \varepsilon}:\hollow^{2\varepsilon} \to \mathbb{R}$ defined by 
\begin{align*}
\psi^{\varepsilon \varepsilon}(s,t) := (2\varepsilon)^{-2} \int_{s-\varepsilon}^{s+\varepsilon} \int_{t- \varepsilon}^{t+\varepsilon} \psi^\varepsilon(s',t')\mathrm{d}s'\mathrm{d}t', \qquad (s,t) \in \hollow^{2\varepsilon}.
\end{align*} 
\end{df}

We will use the double smoothing in Section \ref{sec:Cproof}.

While an increasing surface $\psi:\rt \to \mathbb{R}$ is differentiable almost everywhere in $\rt$, its $\varepsilon$-smoothing is differentiable on $\hollow^\varepsilon$, and its double $\varepsilon$-smoothing is twice differentiable on $\hollow^{2 \varepsilon}$. 
Note that the once- and twice-smoothings of an increasing surface are themselves increasing functions. In fact, if $\psi:\rt \to [0,1]$ takes values in $[0,1]$, the first-order derivatives of its $\varepsilon$-smoothing satisfy
\begin{align} \label{eq:scuba}
0 \leq \psi^\varepsilon_s(s,t) , \psi^\varepsilon_t(s,t)  \leq 1/(2\varepsilon) \qquad (s,t) \in \hollow^\varepsilon.
\end{align}
Likewise, 
\begin{align} \label{eq:scuba2}
0 \leq \psi^{\varepsilon \varepsilon}_s(s,t) , \psi^{\varepsilon \varepsilon}_t(s,t)  \leq 1/(2\varepsilon) \qquad (s,t) \in \hollow^{2\varepsilon}.
\end{align}
By the same argument, the second-order derivatives of the double $\varepsilon$-smoothing satisfy 
\begin{align} \label{eq:scuba3}
\left| \psi^{\varepsilon\varepsilon}_{ss}(s,t)\right|,
\left| \psi^{\varepsilon\varepsilon}_{tt}(s,t)\right|,
\left| \psi^{\varepsilon\varepsilon}_{st}(s,t)\right|
\leq \frac{1}{(2\varepsilon)^2},
\qquad (s,t)\in \hollow^{2\varepsilon}.
\end{align}

Introduce the notation
\begin{align*}
\mathcal{H}_\varepsilon[ \psi ] := - \int_{\hollow^\varepsilon} \sigma( \nabla \psi) \mathrm{d}s\mathrm{d}t.
\end{align*}
 
In a moment, we will relate $\mathcal{H}_\varepsilon[\psi^\varepsilon]$ to $\mathcal{H}[\psi]$. First, we have the following lemma, which may be regarded as a continuum analogue of Lemma \ref{lem:volupper}.
\begin{lemma}
Let $M \subseteq \rt~$ and let $f:M \to [0,1]$ be a function satisfying $\int_{M} f(s,t)\,\mathrm{d}s\mathrm{d}t = \varepsilon \leq 1/2$. Let $\psi:\rt \to [0,1]$ be an increasing surface. Then
\begin{align}\label{eq:evans2}
\int_{M} f \,\sigma(\nabla\psi)\,
\geq
- C \varepsilon \log(1/\varepsilon).
\end{align}

In particular, by letting $f$ be the indicator function of a set $M$, if $M \subseteq \rt$ has Lebesgue measure $|M|=\varepsilon$, then
\begin{align} \label{eq:evans}
\int_{M} \sigma(\nabla\psi)\,\mathrm{d}s\mathrm{d}t
\geq
- C \varepsilon \log(1/\varepsilon).
\end{align}
\end{lemma}

\begin{proof}
Since $\sigma(u_1,u_2)\geq -\log(u_1+u_2)+c$, we obtain
\[
\int_{M} f\,\sigma(\nabla\psi)
\geq
- \int_{M} f\log(\psi_s+\psi_t) + c\varepsilon.
\]
Since $\varepsilon^{-1}f$ is a probability density on $M$ and $\log$ is concave, Jensen's inequality gives
\[
- \int_{M} f\log(\psi_s+\psi_t)
\geq
-\varepsilon \log\!\left(\frac{1}{\varepsilon}\int_{M} f(\psi_s+\psi_t)\right).
\]
As $0\leq f\leq 1$, and $\psi:\rt \to [0,1]$ is increasing in both coordinates we have $\int_{M} f(\psi_s+\psi_t)\leq \int_{\rt}(\psi_s+\psi_t)\leq 2$. Combining these estimates yields the claim, after absorbing constants.
\end{proof}

\begin{lemma} \label{lem:energysmooth}
If $\psi:\rt \to [0,1]$ is an increasing surface and $0 < \varepsilon < 1/6$, then
\begin{align*}
\mathcal{H}_\varepsilon[ \psi^\varepsilon] \geq 
\mathcal{H}[\psi] - C \varepsilon \log(1/\varepsilon).
\end{align*}
\end{lemma}

\begin{proof}
Using the convexity of surface tension and Jensen's inequality, for almost every $(s,t) \in \hollow^\varepsilon$ we have
\begin{align*}
\sigma(\nabla \psi^\varepsilon(s,t)) \leq  (2\varepsilon)^{-2} \int_{s-\varepsilon}^{s+\varepsilon} \int_{t- \varepsilon}^{t+\varepsilon} \sigma ( \nabla \psi(s',t')) \mathrm{d}s'\mathrm{d}t'.
\end{align*}
Integrating over $(s,t) \in \hollow^\varepsilon$ we obtain
\begin{align*}
-\mathcal{H}_\varepsilon[ \psi^\varepsilon] \leq 
\int_{\rt} \sigma(\nabla \psi(s,t)) g_\varepsilon(s,t)\mathrm{d}s\mathrm{d}t,
\end{align*}
where
\begin{align*}
g_\varepsilon(s,t) := (2\varepsilon)^{-2} \int_{ \hollow^\varepsilon} \mathrm{1}\{ s \in [s'-\varepsilon,s'+\varepsilon], t \in [t'-\varepsilon,t' + \varepsilon] \} \mathrm{d}s'\mathrm{d}t'.
\end{align*}
The function $g_\varepsilon(s,t)$ takes values in $[0,1]$ and satisfies $\int_{\rt} g_\varepsilon(s,t) \mathrm{d}s\mathrm{d}t = |\hollow^\varepsilon| = \frac{1}{2}(1-4\varepsilon)^2$. 

We have 
\begin{align*}
-\mathcal{H}_\varepsilon[ \psi^\varepsilon] \leq 
-\mathcal{H}[\psi] - \int_{\rt} \sigma(\nabla \psi(s,t)) f_\varepsilon(s,t)\mathrm{d}s\mathrm{d}t,
\end{align*}
where $f_\varepsilon(s,t) = 1 - g_\varepsilon(s,t)$. Note that $0 \leq f_\varepsilon \leq 1$ and $\int_{\rt} f_{\varepsilon} \mathrm{d}s\mathrm{d}t = \frac{1}{2} (1 - (1-4\varepsilon)^2) \leq 4\varepsilon$. 
Using \eqref{eq:evans2} we have
\begin{align*}
\mathcal{H}_\varepsilon[ \psi^\varepsilon] \geq 
\mathcal{H}[\psi] - C\varepsilon \log(1/\varepsilon),
\end{align*}
completing the proof.
\end{proof}

%%%%%%%%%%%%%%%%%%%%%%%%%%%%%%%%%%%%%%%%%%%%%%%%%%%%%%%%%%%%%%%%%%%%%%%%%%%%%%%%%%%
\subsection{Convergence of smoothed entropy}
%%%%%%%%%%%%%%%%%%%%%%%%%%%%%%%%%%%%%%%%%%%%%%%%%%%%%%%%%%%%%%%%%%%%%%%%%%%%%%%%%%%

The following lemma records a truncation estimate for the surface tension integrand. 
\begin{lemma} \label{lem:trunc}
For $L \geq 2$ define $\sigma_L(u_1,u_2) := \max\{ \sigma(u_1,u_2), -L \}$. Let $M \subseteq \rt$. Then, for every increasing $\psi:M \to [0,1]$, we have 
\begin{align*}
\int_M \sigma_L(\nabla \psi) \,\mathrm{d}s\mathrm{d}t
-
\int_{M} \sigma(\nabla \psi) \,\mathrm{d}s\mathrm{d}t
\leq C e^{-L/2}.
\end{align*}
\end{lemma}

\begin{proof}
Since $\sigma(\nabla\psi)\geq -\log(\psi_s+\psi_t)-1+\log\pi$, we have the inclusion 
\begin{align*}
M_L := \{ (s,t) \in M: \sigma(\nabla\psi(s,t))\leq -L \} \subseteq \{(s,t) \in M: \psi_s(s,t)+\psi_t(s,t) \geq \tfrac{\pi}{e}e^L\} =: N_L.
\end{align*}
As $\psi$ is increasing on $M$ and takes values in $[0,1]$, and $M \subseteq \rt \subseteq [0,1]^2$, we have $\int_M \psi_s\leq 1$ and $\int_M \psi_t\leq 1$, and hence $\int_M (\psi_s + \psi_t) \leq 2$. With $|M_L|$ and $|N_L|$ denoting Lebesgue measure, using $M_L \subseteq N_L$ to obtain the first inequality below, and then Markov's inequality to obtain the second, we have 
\[
|M_L| \leq |N_L| \leq \varepsilon_L := \tfrac{2e}{\pi}e^{-L}.
\]
Since $\sigma_L(\nabla \psi)-\sigma(\nabla \psi)\geq 0$ and is supported on $M_L$,
\[
\int_{M_L} \sigma_L(\nabla \psi) \,\mathrm{d}s\mathrm{d}t
-
\int_{M_L} \sigma(\nabla \psi) \,\mathrm{d}s\mathrm{d}t
\leq  \int_{M_L} (- L - \sigma(\nabla \psi)) \leq - \int_{M_L}  \sigma(\nabla \psi),
\]
where we used the fact that $-L \leq 0$. Note $L \geq 2$ implies $\varepsilon_L \leq 1/2$. Using \eqref{eq:evans} and $|M_L|\leq \varepsilon_L$, 
\[
- \int_{M_L} \sigma(\nabla \psi)
\leq  C  \varepsilon_L \log(1/ \varepsilon_L),
\]
which is $\leq C e^{-L/2}$ for a sufficiently large constant $C > 0$. 
\end{proof}

%%%%%%%%%%%%%%%%%%%%%%%%%%%%%%%%

\begin{lemma} \label{lem:epsconv}
Let $\psi:\righttriangle \to [0,1]$ be an increasing surface such that $\mathcal{H}[\psi] > -\infty$. Then
\begin{align*}
\lim_{\varepsilon \to 0} \mathcal{H}_\varepsilon[\psi^\varepsilon] = \mathcal{H}[\psi].
\end{align*}
\end{lemma}

\begin{proof}
It is clearly a consequence of Lemma \ref{lem:energysmooth} that $\liminf_{\varepsilon \to 0} \mathcal{H}_\varepsilon[\psi^\varepsilon ] \geq \mathcal{H}[\psi]$. It remains to show that $\limsup_{\varepsilon \to 0} \mathcal{H}_\varepsilon[\psi^\varepsilon ] \leq \mathcal{H}[\psi]$.

Since $\psi$ is increasing, Lebesgue's monotone differentiation theorem guarantees that $\psi_s$ and $\psi_t$ exist almost everywhere. Moreover, by standard properties of convolution (see, e.g., Theorem 8.1.5 of \cite{folland}), for almost all $(s,t)$ in the interior of $\righttriangle$, the derivatives satisfy
\begin{align}
\label{eq:gradconv}
\lim_{ \varepsilon \to 0} \psi^\varepsilon_s(s,t) = \psi_s(s,t),
\quad
\lim_{ \varepsilon \to 0} \psi^\varepsilon_t(s,t) = \psi_t(s,t).
\end{align}
By continuity of $\sigma$ and $\sigma_L$, we therefore have
\begin{align}
\label{eq:sigconv}
\lim_{\varepsilon \to 0} \sigma(\nabla \psi^\varepsilon(s,t)) = \sigma( \nabla \psi(s,t)) \quad \text{and} \quad \lim_{\varepsilon \to 0} \sigma_L(\nabla \psi^\varepsilon(s,t)) = \sigma_L( \nabla \psi(s,t))
\end{align}
for almost all $(s,t) \in \righttriangle$. Applying Fatou's lemma to the nonnegative functions $\sigma_L(\nabla\psi^\varepsilon)+L$, we obtain
\begin{align}
\label{eq:fatou}
\liminf_{\varepsilon \to 0} \int_{\hollow^\varepsilon} \sigma_L \left( \nabla \psi^\varepsilon \right) \mathrm{d}s\mathrm{d}t 
\geq 
\int_{\righttriangle} \sigma_L \left( \nabla \psi\right) \mathrm{d}s\mathrm{d}t,
\end{align}
where we can treat the integrals inside the limit on the left-hand side of \eqref{eq:fatou} as integrals over $\rt$ rather than $\hollow^\varepsilon$ by simply setting $\sigma_L(\nabla \psi^\varepsilon) = 0$ for $(s,t) \in \rt - \hollow^\varepsilon$. 

On the other hand, by Lemma \ref{lem:trunc}, for each $\varepsilon>0$,
\begin{align}
\label{eq:truncstep}
\int_{\hollow^\varepsilon} \sigma \left( \nabla \psi^\varepsilon \right) \mathrm{d}s\mathrm{d}t
\geq
\int_{\hollow^\varepsilon} \sigma_L \left( \nabla \psi^\varepsilon \right) \mathrm{d}s\mathrm{d}t - C e^{-L/2}.
\end{align}
Combining \eqref{eq:fatou} and \eqref{eq:truncstep} to obtain the first inequality below, then simply using $\sigma_L \geq \sigma$ to obtain the second, we have 
\begin{align*}
\liminf_{\varepsilon \to 0} \int_{\hollow^\varepsilon} \sigma(\nabla \psi^\varepsilon) \mathrm{d}s \mathrm{d}t
\geq 
\int_{\righttriangle} \sigma_L(\nabla \psi) \mathrm{d}s \mathrm{d}t - C e^{-L/2}
 \geq 
-\mathcal{H}[\psi] - C e^{-L/2}.
\end{align*}
Equivalently, $\limsup_{\varepsilon \to 0} \mathcal{H}_\varepsilon[\psi^\varepsilon]
\leq \mathcal{H}[\psi] + C e^{-L/2}$. 
Letting $L \to \infty$ completes the proof.
\end{proof}

%%%%%%%%%%%%%%%%%%%%%%%%%%%%%%%%%%%%%%%%%%%%%%%%%%%%%%%%%%%%%%%%%%%%%%%%%
\subsection{Convergence under smoothing}
%%%%%%%%%%%%%%%%%%%%%%%%%%%%%%%%%%%%%%%%%%%%%%%%%%%%%%%%%%%%%%%%%%%%%%%%%

The following lemma is a smoothed version of Theorem \ref{thm:lowersemi}.

\begin{lemma} \label{lem:semicont}
If $(\psi_n)_{n\geq 1}$ is a sequence of increasing surfaces taking values in $[0,1]$, and such that $\mathrm{d}(\psi_n,\psi) \to 0$, then for every $0< \varepsilon < 1/6$ we have 
\begin{align*}
\limsup_{n \to \infty} \mathcal{H}_\varepsilon[\psi_n^\varepsilon ] \leq \mathcal{H}_\varepsilon[\psi^\varepsilon].
\end{align*}
\end{lemma}
\begin{proof}
Let $0< \varepsilon < 1/6$. 
Note that $\mathrm{d}(\psi_n,\psi) \to 0$ implies that $\psi_n^\varepsilon(s,t)$ converges pointwise to $\psi^\varepsilon(s,t)$ for all $(s,t) \in \hollow^\varepsilon$. In fact, by using dominated convergence and \eqref{eq:scuba}, we see that the first-order derivatives of $\psi_n^\varepsilon$ converge to those of $\psi^\varepsilon$ for almost all $(s,t) \in \hollow^\varepsilon$. It follows from Fatou's lemma that for every $L > 0$ we have 
\begin{align}
\label{eq:fatou_semicont}
\liminf_{n \to \infty} 
\int_{\hollow^\varepsilon} \sigma_L (\nabla \psi_n^{\varepsilon} ) \mathrm{d}s\mathrm{d}t 
\geq 
\int_{\hollow^\varepsilon} \sigma_L (\nabla \psi^{\varepsilon}) \mathrm{d}s\mathrm{d}t .
\end{align} 
Using Lemma \ref{lem:trunc} and $\sigma_L \geq \sigma$, we obtain
\begin{align}
\label{eq:final_semicont}
\liminf_{n \to \infty} \int_{\hollow^\varepsilon} \sigma(\nabla \psi_n^\varepsilon)\mathrm{d}s\mathrm{d}t
\geq 
\int_{\hollow^\varepsilon} \sigma(\nabla \psi^\varepsilon)\mathrm{d}s\mathrm{d}t - C e^{ - L/2}.
\end{align}
Equivalently, $\limsup_{n \to \infty} \mathcal{H}_\varepsilon[\psi_n^\varepsilon]
\leq \mathcal{H}_\varepsilon[\psi^\varepsilon] + C e^{-L/2}$. 
Since $L$ is arbitrary, the result follows. 
\end{proof}

%%%%%%%%%%%%%%%%%%%%%%%%%%%%%%%%%%%%%%%%%%%%%%%%%%%%%%%%%%%%%%%%%%%%%%%%%
\subsection{Proof of the upper semicontinuity of the entropy functional}
%%%%%%%%%%%%%%%%%%%%%%%%%%%%%%%%%%%%%%%%%%%%%%%%%%%%%%%%%%%%%%%%%%%%%%%%%

\begin{proof}[Proof of Theorem \ref{thm:lowersemi}]
Let $\delta > 0$ be arbitrary. Suppose that $(\psi_n)_{n \geq 1}$ is a sequence of increasing surfaces $\psi_n: \rt \to [0,1]$ that converges to $\psi:\rt \to [0,1]$ in the $L^1$ metric, i.e., $\mathrm{d}(\psi_n,\psi) \to 0$. Write 
\begin{align*}
\mathcal{H}[\psi_n] - \mathcal{H}[\psi] = A_{n,\varepsilon} + B_{n,\varepsilon} + C_\varepsilon,
\end{align*}
where
\begin{align*}
A_{n,\varepsilon} = \mathcal{H}[\psi_n] - \mathcal{H}_\varepsilon[\psi_n^\varepsilon], \quad 
B_{n,\varepsilon} = \mathcal{H}_\varepsilon[\psi_n^\varepsilon] - \mathcal{H}_\varepsilon[\psi^\varepsilon], \quad \text{and} \quad 
C_\varepsilon = \mathcal{H}_\varepsilon[\psi^\varepsilon] - \mathcal{H}[\psi].
\end{align*}
By Lemma \ref{lem:energysmooth} we have $A_{n,\varepsilon} \leq C\varepsilon \log(1/\varepsilon)
$
for all $n \geq 1$, uniformly in $n$. Thus, there exists $\varepsilon_1(\delta)$ such that $0 < \varepsilon \leq \varepsilon_1(\delta)$ implies
$
A_{n,\varepsilon} \leq \delta/3$.
Likewise, by Lemma \ref{lem:epsconv} there exists $\varepsilon_2(\delta)$ such that $0 < \varepsilon \leq \varepsilon_2(\delta)$ implies $C_\varepsilon \leq \delta/3$. Finally, set $\varepsilon = \varepsilon_1(\delta) \wedge \varepsilon_2(\delta)$. Then by Lemma \ref{lem:semicont} there exists $n(\delta)$ such that for $n \geq n(\delta)$,
$
B_{n,\varepsilon} \leq \delta/3.
$
Hence, for all $n \geq n(\delta)$,
$
\mathcal{H}[\psi_n] - \mathcal{H}[\psi] \leq \delta,
$
and since $\delta$ was arbitrary, we conclude that
$
\limsup_{n \to \infty} \mathcal{H}[\psi_n] \leq \mathcal{H}[\psi],
$
completing the proof.
\end{proof}

\subsection{Compressions and closures of Gelfand--Tsetlin patterns}

We close this section by proving that the set of compressions, or equivalently, the set of increasing surfaces $\psi:\rt \to \mathbb{R}$, are the closure of the Gelfand--Tsetlin patterns. 

\begin{proposition} \label{prop:closure}
Let $\psi:\rt \to \mathbb{R}$ be an increasing surface with finite $L^1$-norm and let $\delta > 0$. Then there exists $n \in \mathbb{N}$ and an increasing surface $\psi'$ that occurs as the surface associated with a $(n+1)$-dimensional Gelfand--Tsetlin pattern $(t_{k,j})_{0 \leq j \leq k \leq n}$ such that 
\begin{align*}
\mathrm{d}(\psi,\psi') := \int_{\rt} |\psi'(s,t) - \psi(s,t)|\mathrm{d}s\mathrm{d}t \leq \delta.
\end{align*}
\end{proposition}

\begin{proof}
If $\psi$ is bounded, we may simply set $\psi'(i/n,j/n) := \psi(i/n,j/n)$ and linearly interpolate $\psi'$ on triangles, so that $\psi'$ is the surface of a Gelfand--Tsetlin pattern. It is then an easy exercise using the increasingness of $\psi$ and $\psi'$ to verify that $\mathrm{d}(\psi,\psi') \leq CL/n$, where $L := \sup_{\rt} |\psi|$ and $C > 0$ is universal.

If $\psi$ is not bounded, simply approximate $\psi$ with a bounded function and use the triangle inequality.

% Alternatively, if $\psi$ is merely of finite $L^1$-norm but unbounded, for $L > 0$ define the truncation $\psi^L$ by
% \begin{align*}
% \psi^L(s,t) := \begin{cases}
% \psi(s,t) \qquad &\text{if $\psi(s,t) \in [-L,L]$}\\
% L \qquad &\text{if $\psi(s,t) > L$}\\
% -L \qquad &\text{if $\psi(s,t) < - L$}.
% \end{cases}
% \end{align*}
% Note that $\psi^L$ is also an increasing surface. 
% By the integrability of $\psi$, one can choose $L$ sufficiently large so that $\mathrm{d}(\psi,\psi^L) \leq \delta/2$, and by the previous part, one can choose an increasing surface $\psi'$ coming from a Gelfand--Tsetlin pattern such that $\mathrm{d}(\psi^L,\psi') \leq \delta/2$ also. By the triangle inequality it follows that $\mathrm{d}(\psi,\psi') \leq \delta$, completing the proof.
\end{proof}

%%%%%%%%%%%%%%%%%%%%%%%%%%%%%%%%%%%%%%%%%%%%%%%%%%%%%%%%%%%%
%%%%%%%%%%%%%%%%%%%%%%%%%%%%%%%%%%%%%%%%%%%%%%%%%%%%%%%%%%%%
%%%%%%%%%%%%%%%%%%%%%%%%%%%%%%%%%%%%%%%%%%%%%%%%%%%%%%%%%%%%
\section{Proof of Theorem \ref{thm:C}} \label{sec:Cproof}
%%%%%%%%%%%%%%%%%%%%%%%%%%%%%%%%%%%%%%%%%%%%%%%%%%%%%%%%%%%%
%%%%%%%%%%%%%%%%%%%%%%%%%%%%%%%%%%%%%%%%%%%%%%%%%%%%%%%%%%%%
%%%%%%%%%%%%%%%%%%%%%%%%%%%%%%%%%%%%%%%%%%%%%%%%%%%%%%%%%%%%
\subsection{Preparation}
Our proof of Theorem \ref{thm:C} is largely split into two parts: a large deviation upper bound and a large deviation lower bound. In preparation for both of these parts, let us recall the setup of Theorem \ref{thm:C}. We recall from Definition \ref{df:GTcomp} that the increasing surface associated with a Gelfand--Tsetlin pattern $(t_{k,j})_{ 0 \leq j \leq k \leq n}$ is defined by setting
\begin{align*}
\xi(x_1/n,x_2/n) := t_{n-x_1+x_2,x_2}, \qquad 0 \leq x_2 \leq x_1 \leq n,
\end{align*}
and linearly interpolating on triangles. As in the statement of Theorem \ref{thm:C}, under a probability measure $\mathbf{P}_\mu$, for each $n \geq 1$ let $\xi_n:\rt \to \mathbb{R}$ denote the increasing surface associated with the uniform random Gelfand--Tsetlin pattern $t^n := (t^n_{k,j})_{ 0 \leq j \leq k \leq n }$ whose bottom row is given by $t^n_{n,j} := Q_\mu(j/n)$. Here $Q_\mu:[0,1] \to \mathbb{R}$ is the quantile function of a compactly supported measure $\mu$ with $\chi[\mu] > - \infty$. The condition $\chi[\mu] > - \infty$ implies that $Q_\mu$ is strictly increasing (otherwise, $\mu$ has an atom somewhere, and $\chi[\mu] = -\infty$). 

Recall that $\rt_n' := \{ (x_1,x_2) \,:\,  0 \leq x_2 < x_1 \leq n \}$ is $\rt_n$ minus its diagonal entries. 
Altogether, given a measurable subset $\Gamma$ of the space of integrable increasing surfaces we may write
\begin{align*}
\mathbf{P}_\mu( \xi_n \in \Gamma) = D_n(\mu)^{-1} \int_{ \mathbb{R}^{\rt_n} } \mathrm{1} \left\{ \xi_n \in \Gamma \right\} \mathrm{1}^{\righttriangle_n}_{\mathrm{GT}}(\phi) \prod_{ x \in \rt_n' } \mathrm{d} \phi_x \prod_{k = 0}^{n} \delta_{Q_\mu(k/n)}( \mathrm{d}\phi_{k,k} ),
\end{align*}
where $\xi_n:\rt \to \mathbb{R}$ depends on the integrated variables $(\phi_x)_{x \in \rt_n}$, and is constructed from these variables by setting $\xi_n(x_1/n,x_2/n) := \phi_{x_1,x_2}$ and linearly interpolating. The normalisation factor is given by
\begin{align} \label{eq:Dnv}
D_n(\mu) := \mathrm{Leb}_{(n+1)n/2} \mathrm{GT}( Q_\mu(0/n),\ldots,Q_\mu(n/n)).
\end{align}
Note that $\rt_n'$ contains $n(n+1)/2$ points. In particular, by making a change of variables and changing the boundary conditions from $Q_\mu(k/n)$ to $n Q_\mu(k/n)$ for each $0 \leq k \leq n$, we can write
\begin{align} \label{eq:Dnv2}
\mathbf{P}_\mu( \xi_n \in \Gamma) = (n^{n(n+1)/2} D_n(\mu))^{-1} V_{n,\mu} ( \Gamma),
\end{align}
where 
\begin{align} \label{eq:Vdef}
V_{n,\mu} ( \Gamma) := \int_{ \mathbb{R}^{\rt_n} } \mathrm{1} \left\{ \frac{1}{n} \xi_n \in \Gamma \right\} \mathrm{1}^{\righttriangle_n}_{\mathrm{GT}}(\phi) \prod_{ x \in \rt_n' } \mathrm{d} \phi_x \prod_{k = 0}^{n} \delta_{n Q_\mu(k/n)}( \mathrm{d}\phi_{k,k} ).
\end{align}
Now using \eqref{eq:Dnv} and the Weyl dimension formula \eqref{eq:wdf} to obtain the first equality below, then using \eqref{eq:chi2} to obtain the second, we have
\begin{align} \label{eq:anka}
\lim_{n \to \infty} \frac{1}{n^2} \log \left( n^{n(n+1)/2} D_n(\mu) \right)
=
\lim_{n \to \infty} \frac{1}{n^2} \log \mathrm{Leb}_{(n+1)n/2} \mathrm{GT}( n Q_\mu(0/n),\ldots,n Q_\mu(n/n))
= \chi[\mu].
\end{align}

Thus, by \eqref{eq:Dnv2} and \eqref{eq:anka}, in order to prove the statement of Theorem \ref{thm:C}, it is sufficient to prove the large deviation upper bound
\begin{align} \label{eq:UPPER}
\limsup_{n \to \infty} \frac{1}{n^2} \log V_{n,\mu}(\Gamma) \leq \sup_{ \psi \in \overline{\Gamma} \cap \mathcal{A}_\mu } \int_{\rt} (-\sigma)(\nabla \psi),
\end{align}
and the large deviation lower bound 
\begin{align} \label{eq:LOWER}
\liminf_{n \to \infty} \frac{1}{n^2} \log V_{n,\mu}(\Gamma) \geq \sup_{ \psi \in \interior{\Gamma} \cap \mathcal{A}_\mu } \int_{\rt} (-\sigma)(\nabla \psi),
\end{align}
where in both cases
\begin{align} \label{eq:Enu} 
\mathcal{A}_\mu := \{ \psi: \rt \to \mathbb{R} \text{ increasing surface} \,:\, \psi(s,s) = Q_\mu(s) ~\text{ for all $s \in [0,1]$} \}
\end{align}
is the set of increasing surfaces agreeing with the quantile function of $\mu$ on the diagonal. Equivalently, $\mathcal{A}_\mu$ is the set of increasing surfaces associated with compressions of $\mu$. In both \eqref{eq:UPPER} and \eqref{eq:LOWER}, we use the convention that the supremum over an empty set is $-\infty$.

In establishing \eqref{eq:UPPER} and \eqref{eq:LOWER}, we will reduce the problem to considering the case where $\Gamma$ is a ball of radius $\delta$, that is,
\begin{align*}
\Gamma = B(\psi,\delta) := \{ \psi':\rt \to \mathbb{R} \text{ increasing and such that } \mathrm{d}(\psi,\psi') \leq \delta \}.
\end{align*}
We will find it useful to compare this integral to a sum over lattice points. In this direction, let $f,g:\rt \to [0,1]$ be increasing surfaces taking values in $[0,1]$, and set
\begin{align}
\mathrm{d}(f,g) := \int_{\rt} |g-f| \,\mathrm{d}s\mathrm{d}t \quad \text{and} \quad \mathrm{d}_n(f,g) := \frac{1}{n^2} \sum_{ (x_1,x_2) \in \rt_n} |g(x_1/n,x_2/n) - f(x_1/n,x_2/n)|.
\end{align}
Using the fact that $f$ and $g$ are increasing in both variables, it is a brief calculation to verify that  
\begin{align} \label{eq:dapprox}
|\mathrm{d}_n(f,g) - \mathrm{d}(f,g)| \leq C/n,
\end{align}  
where $C > 0$ is uniform on all increasing surfaces $f,g:\rt \to [0,1]$.

%%%%%%%%%%%%%%%%%%%%%%%%%%%%%%%%%%%%%%%%%%%%%%%%%%%%%%%%%%%%
%%%%%%%%%%%%%%%%%%%%%%%%%%%%%%%%%%%%%%%%%%%%%%%%%%%%%%%%%%%%
%%%%%%%%%%%%%%%%%%%%%%%%%%%%%%%%%%%%%%%%%%%%%%%%%%%%%%%%%%%%
\subsection{Proof of the large deviation upper bound \eqref{eq:UPPER}} \label{sec:mainupper}
%%%%%%%%%%%%%%%%%%%%%%%%%%%%%%%%%%%%%%%%%%%%%%%%%%%%%%%%%%%%
%%%%%%%%%%%%%%%%%%%%%%%%%%%%%%%%%%%%%%%%%%%%%%%%%%%%%%%%%%%%
%%%%%%%%%%%%%%%%%%%%%%%%%%%%%%%%%%%%%%%%%%%%%%%%%%%%%%%%%%%%
We now begin to prove the upper bound \eqref{eq:UPPER}. 

For starters, let $\psi:\righttriangle \to \mathbb{R}$ be an increasing function. Recall from Definition \ref{df:smoothing} that for $0 < \varepsilon < 1/6$, we write $\psi^\varepsilon$ and $\psi^{\varepsilon \varepsilon}$ for its once- and double-smoothing. These new functions are defined on the subsets $\hollow^\varepsilon$ and $\hollow^{2\varepsilon}$ of $\rt$. We will also consider
\begin{align*}
\rt' := \rt - \{(s,s) : 0 \leq s \leq 1 \} = \{ (s,t) \in \mathbb{R}^2 : 0\leq t < s \leq 1 \},
\end{align*}
% which is simply $\rt$ with its diagonal removed. Correspondingly, recall that $\rt_n'$ denotes $\rt_n$ with its diagonal removed.

Given a subset $E$ of $\mathbb{R}^2$, we define $E_n := \{ (x_1,x_2) \in \mathbb{Z}^2 : (x_1/n,x_2/n) \in E \}$. Note that this notation is consistent with our work so far in the case $E = \rt, \rt'$. If $\psi:E \to \mathbb{R}$ is a function, we define a new function $\psi^n:E_n \to \mathbb{R}$ by 
\begin{align*}
\psi_x^n := n \psi(x_1/n,x_2/n), \qquad x \in E_n.
\end{align*}
We will often consider $\psi^{\varepsilon n} :=(\psi^\varepsilon)^n:\hollow^\varepsilon_n \to \mathbb{R}$, defined on $\hollow^\varepsilon_n  := ( \hollow^\varepsilon)_n$. Given a function $\varphi:C \to \mathbb{R}$ defined on a subset $C$ of $\mathbb{Z}^2$ and $\delta > 0$, we define
\begin{align}
    K_n(\varphi,C,\delta) := \int_{[0,n]^{C}}  \mathrm{1} \left\{ \frac{1}{n^2} \sum_{ x \in C} |\phi_x - \varphi_x | \leq \delta n  \right\} \mathrm{1}_{\mathrm{GT}}^{C}(\phi)\prod_{x \in C} \mathrm{d}\phi_x, \label{eq:avball}\\
        K^\infty_n(\varphi,C,\delta) := \int_{[0,n]^{C}}  \mathrm{1} \left\{ \sup_{ x \in C} |\phi_x - \varphi_x | \leq \delta n  \right\} \mathrm{1}_{\mathrm{GT}}^{C}(\phi)\prod_{x \in C} \mathrm{d}\phi_x, \label{eq:absball}
\end{align}
to be the volumes of Gelfand--Tsetlin surfaces in scaled $\delta$-balls (in $L^1$ and $L^\infty$ respectively) around $\varphi:C \to \mathbb{R}$. More generally, given a function $\psi:E \to \mathbb{R}$ defined on a subset $E$ of $\mathbb{R}^2$, $\delta > 0$ and $n \in \mathbb{N}$, we will abuse notation and define 
\begin{align*}
    K_n(\psi,E,\delta) := K_n(\psi^n,E_n,\delta) \quad \text{and} \quad     K^\infty_n(\psi,E,\delta) := K^\infty_n(\psi^n,E_n,\delta).
\end{align*}
% We recall from Section \ref{sec:gtearnest} that $\mathrm{1}_{\mathrm{GT}}^{C}(\phi)$ is simply the indicator function that $\phi:C \to \mathbb{R}$ is an increasing function on the lattice.
Our first result uses the Pr\'ekopa--Leindler inequality for partition functions over balls to give an upper bound for $K_n(\psi,\rt',\delta)$ in terms of a smoothing of $\psi$. 
%%%%%%%%%%%%%%%%%%%%%%%%%%%%%%%%%%%%%%%%%%%%%%%
\begin{lemma} \label{lem:smoother}
Let $\delta > 0$, $0 < \varepsilon < 1/6$ and $n \in \mathbb{N}$ be such that $1/\varepsilon n \leq \delta/2$. Let $\psi:\righttriangle \to [0,1]$ be an increasing surface. Then
\begin{align*}
\frac{1}{n^2} \log K_n( \psi, \rt' , \delta) \leq \frac{1}{n^2}\log K_n(  \psi^{\varepsilon\varepsilon}, \hollow^{2\varepsilon},4\delta) + C\varepsilon \log(1/\varepsilon),
\end{align*}
where $\psi^{\varepsilon \varepsilon}:\hollow^{2\varepsilon} \to \mathbb{R}$ is the double $\varepsilon$-smoothing of $\psi$. 
\end{lemma}

\begin{proof}
Let $r = \lceil \varepsilon n \rceil$. 
Define
\begin{align*}
\righttriangle_n^{\circ r} := \{ (x_1,x_2) \in \righttriangle_n' : x_2 \geq r, x_1 \leq n-r, x_1-x_2 \geq 2r \},
\end{align*}
to be the right-angled triangle with corners at $(3r,r), (n-r,r), (n-r,n-3r)$ and non-hypotenuse side length $n-4r$. 
For any $\mathbf{j} = (j_1,j_2) \in \mathbb{Z}^2$ we define the translation
\begin{align*}
\righttriangle_n^{\circ r,\rightarrow \mathbf{j}} := \{ (x_1+j_1,x_2+j_2) : (x_1,x_2) \in \righttriangle_n^{\circ r} \}.
\end{align*}
From the definition, whenever $|j_1|,|j_2| \leq r-1$ we have $\righttriangle_n^{\circ r, \rightarrow \mathbf{j}} \subseteq \righttriangle_n'.$ 

Introducing the shorthand $A_{\mathbf{j}} := \righttriangle_n^{\circ r,\rightarrow \mathbf{j}}$ and setting $B_{\mathbf{j}} := \rt_n' - A_{\mathbf{j}}$, by the union inequality \eqref{eq:unioneq} for each $\mathbf{j}$ we have $\mathrm{1}_{\mathrm{GT}}^{\rt'_n}(\phi) \leq \mathrm{1}_{\mathrm{GT}}^{A_{\mathbf{j}}}(\phi) \mathrm{1}_{\mathrm{GT}}^{B_{\mathbf{j}}}(\phi)$. Ignoring the summands outside of $A_{\mathbf{j}}$, we have the following $\mathbf{j}$-dependent upper bound on $K_n( \psi, \rt' , \delta)$:
\begin{align} \label{eq:splitprod}
 K_n( \psi, \rt' , \delta)  \leq 
\int_{ [0, n]^{\rt_n'} } \mathrm{1} \left\{ \frac{1}{n^2} \sum_{x \in A_{\mathbf{j}}} | \phi_x - \psi^n_x | \leq \delta n  \right\} \mathrm{1}_{\mathrm{GT}}^{A_{\mathbf{j}}}(\phi) \mathrm{1}_{\mathrm{GT}}^{B_{\mathbf{j}}}(\phi) \prod_{ x \in \righttriangle_n' } \mathrm{d} \phi_x \leq  K_n( \psi^n|_{A_\mathbf{j}}, A_\mathbf{j} , \delta)  \mathcal{B}_\mathbf{j},
\end{align}
where $\mathcal{B}_\mathbf{j} := \int_{[0,n]^{B_{\mathbf{j}}}} \mathrm{1}_{\mathrm{GT}}^{B_{\mathbf{j}}}(\phi)\mathrm{d}\phi$, and $\psi^n|_{A_\mathbf{j}}$ denotes the restriction of $\psi^n$ to $A_{\mathbf{j}}$. 
Since $B_{\mathbf{j}}$ contains $O(\varepsilon n^2)$ points, by Lemma \ref{lem:volupper}, $\mathcal{B}_\mathbf{j} \leq e^{ C \varepsilon  (\log(1/\varepsilon)) n^2 }$ for some universal $C > 0$. Accordingly, taking logs and averaging over all $(2r-1)^2$ pairs of integers $\mathbf{j} = (j_1,j_2)$ with $|j_1|, |j_2| \leq r-1$, we obtain
\begin{align*}
\frac{1}{n^2}\log K_n( \psi, \rt' , \delta)  \leq C \varepsilon \log(1/\varepsilon) + \frac{1}{n^2} \frac{1}{(2r-1)^2} \sum_{|j_1|,|j_2| \leq r-1}  \log  K_n( \psi^n|_{A_\mathbf{j}}, A_\mathbf{j} , \delta).
\end{align*}
Using translation properties of Gelfand--Tsetlin partition functions together with the form of the Pr\'ekopa--Leindler inequality given in \eqref{eq:smootha}, we have 
\begin{align*}
\frac{1}{n^2}\log K_n( \psi, \rt' , \delta)  \leq C \varepsilon  \log(1/\varepsilon) + \frac{1}{n^2} \log K_n(\tilde{\psi},A_0,\delta),
\end{align*}
where $\tilde{\psi}:A_0 \to \mathbb{R}$ is defined by
\begin{align*}
\tilde{\psi}_x := \frac{1}{(2r-1)^2} \sum_{|j_1|,|j_2| \leq r-1} n \psi\!\left(\frac{x_1+j_1}{n},\frac{x_2 +j_2}{n}\right), \qquad x \in A_0.
\end{align*}
Note that $A_0 = \hollow^\varepsilon_n$. Also note that $\tilde{\psi}$ is $n$ times a discrete average of $\psi$ on $1/n$ mesh points over a square of side-length $2\varepsilon + O(1/n)$, whereas $\psi^\varepsilon$ is a continuous average of $\psi$ over a square of side-length $2\varepsilon$. Since $\psi$ is increasing and takes values in $[0,1]$, it follows that for all $x \in A_0$ we have 
\begin{align} \label{eq:riemann}
\left| \frac{1}{n} \tilde{\psi}_x - \psi^\varepsilon\!\left(\frac{x_1}{n},\frac{x_2}{n}\right) \right| \leq \frac{1}{\varepsilon n} \leq \delta. 
\end{align}
Indeed, \eqref{eq:riemann} follows by comparing the continuous average defining $\psi^\varepsilon$ with the corresponding discrete Riemann sum over mesh size $1/n$. Since $\psi$ is increasing and takes values in $[0,1]$, the discrepancy is bounded by the contribution from a boundary strip of width $1/n$ around the averaging square, which is $O(1/(\varepsilon n))$.

By \eqref{eq:riemann} and the triangle inequality, we have 
\begin{align*}
K_n(\tilde{\psi},A_0,\delta) \leq K_n(\psi^\varepsilon, \hollow^\varepsilon, 2\delta),
\end{align*}
so that
\begin{align*}
\frac{1}{n^2}\log K_n( \psi, \rt' , \delta)  \leq C \varepsilon  \log(1/\varepsilon) + \frac{1}{n^2}\log K_n(\psi^\varepsilon, \hollow^\varepsilon, 2\delta).
\end{align*}

To obtain the result as written, undertake a similar argument with $\psi^\varepsilon$ in place of $\psi$.
\end{proof}
%%%%%%%%%%%%%%%%%%%%%%%%%%%%%%%%%%%%%%%%%%%
The following lemma provides a means of converting $L^1$ bounds to $L^\infty$ bounds.

\begin{lemma} \label{lem:1inf}
Let $\varepsilon \geq 1/n$. Let $\Psi:\hollow^{2\varepsilon} \to [0,1]$ be an increasing function satisfying the inequality
\begin{align} \label{eq:oo}
0 \leq \Psi(s+s',t+t') - \Psi(s,t) \leq \frac{1}{\varepsilon} (s'+t')
\end{align}
for all $(s,t),(s+s',t+t') \in \hollow^{2\varepsilon}$, with $s',t' \geq 0$. Let $\phi:\hollow^{2\varepsilon}_n \to [0,n]$ be any increasing function satisfying the $L^1$-type bound
\begin{align*}
\frac{1}{n^2} \sum_{ x \in \hollow^{2\varepsilon}_n} |\phi(x_1,x_2) - n \Psi(x_1/n,x_2/n) | \leq \delta n.
\end{align*}
Then on $\hollow_n^{3\varepsilon}$, $\phi$ also satisfies the $L^\infty$-type bound
\begin{align*}
|\phi(x_1,x_2) - n\Psi(x_1/n,x_2/n) | \leq C \delta^{1/3} \varepsilon^{-2/3}n \qquad \text{for all $x \in \hollow^{3\varepsilon}_n$},
\end{align*}
where $C > 0$ is universal.
\end{lemma}
\begin{proof}
If $C\delta^{1/3}\varepsilon^{-2/3} \geq 1$, the claim is trivial, since both $\phi$ and $n \Psi$ take values in $[0,n]$. Thus suppose that for some $x \in \hollow^{3\varepsilon}_n$ we have $\phi(x_1,x_2) - n \Psi(x_1/n,x_2/n) \geq an$ with $a \in (0,1)$. Then using the fact that $\phi$ is increasing and $\Psi$ satisfies \eqref{eq:oo}, for all $y_1 \geq x_1$ and $y_2 \geq x_2$ with $(y_1,y_2) \in \hollow^{2\varepsilon}_n$ we have 
\begin{align} \label{eq:aa}
\phi(y_1,y_2) - n \Psi(y_1/n,y_2/n) \geq an - \frac{1}{\varepsilon}(y_1+y_2-x_1-x_2).
\end{align}
Note that if $(y_1,y_2) \in \mathbb{Z}^2$ with $y_1-x_1, y_2-x_2 \leq an\varepsilon/4$, then since $a \in (0,1)$, $(x_1,x_2) \in \hollow_n^{3\varepsilon}$ implies $(y_1,y_2) \in \hollow_n^{2 \varepsilon}$. Moreover, by \eqref{eq:aa}, any such $(y_1,y_2)$ satisfy $\phi(y_1,y_2) - n \Psi(y_1/n,y_2/n) \geq an/2$. 

Accordingly, if $\phi(x_1,x_2) - n \Psi(x_1/n,x_2/n) \geq an$ for $a \in (0,1)$, then since there are $(an\varepsilon/4 + O(1))^2$ choices of $y_1,y_2$ such that $y_1-x_1,y_2-x_2 \leq (an\varepsilon/4)$, we have 
\begin{align*}
\frac{1}{n^2} \sum_{y \in \hollow^{2\varepsilon}_n} |\phi(y_1,y_2) - n \Psi(y_1/n,y_2/n)| \geq \frac{1}{n^2}\frac{an}{2} (an\varepsilon/4 + O(1))^2 \geq c n a^3\varepsilon^2,
\end{align*}
for some universal constant $c > 0$. Since the assumed $L^1$ bound implies that $cn a^3 \varepsilon^2 \leq \delta n$, rearranging we see that $a \leq C \delta^{1/3} \varepsilon^{-2/3}$. 

A similar argument can be used to handle the case where $\phi(x_1,x_2) - n \Psi(x_1/n,x_2/n) \leq -an$ for some $a \in (0,1)$. Combining these two bounds for the magnitude of $a$, we obtain the result.
\end{proof}
%%%%%%%%%%%%%%%%%%%%%%%%%%%%%%%%%%%%%
\begin{cor} \label{cor:busan} 
Let $\Psi:\hollow^{2\varepsilon} \to [0,1]$ be as in the statement of Lemma \ref{lem:1inf}. 
Then
\begin{align} \label{eq:busan}
\frac{1}{n^2} \log K_n(\Psi,\hollow^{2\varepsilon},\delta) 
\leq \frac{1}{n^2} \log K^\infty_n(\Psi,\hollow^{3\varepsilon},C \delta^{1/3}\varepsilon^{-2/3}) + C \varepsilon \log(1/\varepsilon).
\end{align}
\end{cor}
\begin{proof}
By the previous lemma, with $\rho = C \delta^{1/3}\varepsilon^{-2/3}$ we have
\begin{align*}
K_n(\Psi,\hollow^{2\varepsilon},\delta) 
\leq \int_{[0,n]^{\hollow_n^{2\varepsilon}}} 
\mathrm{1} \left\{ |\phi_x - \Psi_x^n| \leq \rho n ~ \text{for all $x \in \hollow_n^{3\varepsilon}$} \right\} 
\mathrm{1}_{\mathrm{GT}}^{\hollow_n^{2\varepsilon}}(\phi)
\prod_{x \in \hollow_n^{2\varepsilon} }\mathrm{d}\phi_x.
\end{align*}
Letting $B_n := \hollow^{2\varepsilon}_n - \hollow^{3\varepsilon}_n$, by the union inequality we have $\mathrm{1}_{\mathrm{GT}}^{\hollow^{2\varepsilon}_n}(\phi) 
\leq 
\mathrm{1}_{\mathrm{GT}}^{\hollow^{3\varepsilon}_n}(\phi) 
\mathrm{1}_{\mathrm{GT}}^{B_n}(\phi)$, so that 
\begin{align*}
K_n(\Psi,\hollow^{2\varepsilon},\delta) 
\leq K^\infty_n(\Psi,\hollow^{3\varepsilon},\rho)\mathcal{B}_n, \qquad \text{where} \quad \mathcal{B}_n :=  
\int_{[0,n]^{ B_n} } 
\mathrm{1}_{\mathrm{GT}}^{B_n}(\phi)
\prod_{x \in B_n  }\mathrm{d}\phi_x.
\end{align*}
Now $B_n$ contains $O(\varepsilon n^2)$ points; thus we can again use Lemma \ref{lem:volupper} to obtain
$
\mathcal{B}_n \leq \exp\{ C \varepsilon \log(1/\varepsilon) n^2\}.
$
Taking logarithms and dividing by $n^2$ gives the result.
\end{proof}
%%%%%%%%%%%%%%%%%%%%%%%%%%%%%%%%%%%%%%%
In our proofs of both the upper and lower bound in Theorem \ref{thm:C}, we will approximate functions using piecewise linear functions on a linear mesh, so that the bounds developed in Theorem \ref{thm:upperlower} can be applied on this mesh. In this direction, let $m \in \mathbb{N}$ be an integer. For integers $a,b$ consider the subtriangles
\begin{align*}
\righttriangle^{m,ab} &:= \{ ((a+u)/m,(b+v)/m) : 0 \leq v < u \leq 1 \}\\
\urt^{m,ab} &:= \{ ((a+u)/m,(b+v)/m) : 0 < u \leq v < 1 \}.
\end{align*}
Note each $\righttriangle^{m,ab}$ is disjoint from any $\urt^{m,ab}$, though the union of such sets tiles space. In particular, for any $m \in \mathbb{N}$ we have the disjoint union
\begin{align} \label{eq:disj}
\righttriangle' := \{(s,t) \in \mathbb{R}^2 : 0 \leq t < s \leq 1 \} = \left( \bigcup_{0 \leq b \leq a \leq m-1 }  \righttriangle^{m,ab} \right)  \cup \left( \bigcup_{0 \leq b < a \leq m-1 }  \urt^{m,ab} \right) .
\end{align}

\begin{df}
We say that a function $\psi:U \to \mathbb{R}$ is \textbf{$m$-linear} if its restriction to every intersection $\righttriangle^{m,ab} \cap U $ and $ \urt^{m,ab} \cap U$ is affine, that is, on each such intersection, $\psi(x_1,x_2) = c + u_1 x_1 + u_2 x_2$ for some $c,u_1,u_2$ that depend on the intersection. 

Suppose $U := \bigcup_{i=1}^r U^i$ can be written as a union of triangles, where for $i$ we have  $U^i= \rt^{m,ab}$ or $U^i = \urt^{m,ab}$ for some $a,b$. Then given a function $\psi:U \to \mathbb{R}$, its \textbf{$m$-linearization} is the unique continuous $m$-linear function $\psi^m:U \to \mathbb{R}$ with $\psi(x_1/m,x_2/m) = \psi^m(x_1/m,x_2/m)$ for each pair $(x_1,x_2)$ of integers satisfying $0 \leq x_2 \leq x_1 \leq m$. 
\end{df}

In the following lemma, recall that if $U$ is a subset of $\mathbb{R}^2$, we write $U_n := \{ (x_1 ,x_2) \in\mathbb{Z}^2 : (x_1/n,x_2/n) \in U\}$. 
%%%%%%%%%%%%%%%%%%%%%%%%%%%%%%%%%%%%%%%%%%%
\begin{lemma} \label{lem:forlinear}
Let $\psi:U \to \mathbb{R}$ be an $m$-linear function defined on some union
$U = \bigcup_{i=1}^r U^i$ of triangles where each $U^i$ is equal to some
$\rt^{m,ab}$ or $\urt^{m,ab}$. Write $u^i = (u^i_1,u^i_2)$ for the slope
of $\psi$ on $U^i$. Suppose for some $\kappa > 0 $, the $u^i$ satisfy
\begin{align} \label{eq:kappacontrol}
\kappa \leq u^i_1 \wedge u^i_2 \leq u^i_1 +u^i_2 \leq 1/\kappa.
\end{align} 
Then for all $n \geq n_{\kappa,m,\rho}$ we have 
\begin{align} \label{eq:forlinear}
\frac{1}{n^2} \log K_n^\infty ( \psi, U, \rho  )
\leq - \int_U \sigma( \nabla \psi) \mathrm{d}s\mathrm{d}t
+ C (m\rho/\kappa)^{1/2}.
\end{align}
\end{lemma}
\begin{proof}
Write $U_n = \bigcup_{i=1}^r U^i_n$. Using the union inequality $\mathrm{1}_{\mathrm{GT}}^{U_n}(\phi) \leq \prod_{i=1}^r \mathrm{1}_{\mathrm{GT}}^{U^i_n}(\phi)$, 
and the fact that $|\phi_x - \psi^n_x| \leq \rho n$ for $x \in U_n$
certainly implies that $|\phi_x - \psi^n_x| \leq \rho n$ for each $x$ on
the boundary $\partial U^i_n$ of each $U^i_n$, we have the upper bound 
\begin{align} \label{eq:4l1}
\frac{1}{n^2} \log K_n^\infty ( \psi, U, \rho  ) \leq \frac{1}{m^2} \sum_{ i =1}^r \frac{n_i^2}{(n/m)^2} J^{(i)},
\end{align}
where 
\begin{align*}
J^{(i)}:= \frac{1}{n_i^2}\log \int_{[0,n]^{U^i_n}} 
\mathrm{1} \{ |\phi_x - \psi^n_x| \leq \rho n  \text{  for } x \in \partial U^i_n \}
\mathrm{1}_{\mathrm{GT}}^{U^i_n}(\phi) \prod_{x \in U^i_n} \mathrm{d}\phi_x,
\end{align*}
and $n_i$ is the side length of the triangle $U^i_n$, so that
$n_i = n/m+O(1)$. Then by integrating around the boundary
$[0,n]^{\partial U^i_n}$ and taking a generous bound, with the notation of Section \ref{sec:ST2proof},
we have
\begin{align} \label{eq:4l2}
J^{(i)} \leq \frac{1}{n_i^2} \log \left( n^{O(n_i)} \sup_{ \varphi \in B_{n_i}(u^i,\rho n) } T_{n_i}(\varphi)  \right), 
\end{align} 
where the
supremum is taken over the set of 
$\varphi: \partial \rt_{n_i} \to \mathbb{R}$ satisfying
$|\varphi_x - \varphi^{u^i}_x| \leq \rho n$ for all
$x \in \partial \rt_{n_i}$. 

For all $n$ sufficiently large, since $n_i = n/m+O(1)$ we will have $u_1^i \wedge u_2^i \geq \kappa \geq \log(n_i)^{-1/2}$. Thus, by appealing to the upper bound \eqref{eq:upper1} in Theorem
\ref{thm:upperlower}, we have
\begin{align} \label{eq:upperhere}
\frac{1}{n_i^2/2}\log \sup_{ \varphi \in B_{n_i}(u^i,\rho n) } T_{n_i}(\varphi)
\leq - \sigma(u^i_1,u^i_2) + E_i,
\qquad
E_i := C \rho_i ( \log(1/\rho_i) + g(u^i)),
\end{align}
where we recall that 
\begin{equation*}
\rho_i = \frac{ \rho n + u_1^i + u_2^i }{ (u^i_1 \wedge u^i_2) n_i } \qquad \text{and} \qquad g(u^i) := |\log(u^i_1 \wedge u^i_2)| + |\log(u^i_1+u^i_2)|.
\end{equation*}
Using \eqref{eq:kappacontrol}, we have $g(u^i) \leq 2 \log (1/\kappa)$. Also, since $\rho \geq 1/(\kappa n)$ for $n$ sufficiently large we have $\rho_i \leq 4m\rho/\kappa$. Altogether, taking a generous bound we have
\begin{align*}
E_i \leq C (m\rho/\kappa)^{1/2} \qquad \text{for all $i =1,\ldots,r$}.
\end{align*}

Since $n_i = n/m+O(1)$, the contribution from the boundary integration term
$n^{O(n_i)}$ is also absorbed into $ C (m\rho/\kappa)^{1/2}$. Combining
\eqref{eq:4l1}, \eqref{eq:4l2}, and \eqref{eq:upperhere}, we obtain
\begin{align} \label{eq:brubeck1}
\frac{1}{n^2} \log K_n^\infty ( \psi, U, \rho  )  
\leq
-\frac{1}{m^2}\sum_{i=1}^r \frac{n_i^2}{(n/m)^2}\frac{1}{2}\sigma(u^i)
+ C (m\rho/\kappa)^{1/2}.
\end{align}
Finally, since $\psi$ is linear with slope $u^i$ on $U^i$ and $U^i$ has Lebesgue measure $1/(2m^2)$, we have
\begin{align} \label{eq:brubeck2}
\frac{1}{m^2}\sum_{i=1}^r \frac{n_i^2}{(n/m)^2}\frac{1}{2}\sigma(u^i)
\leq
\int_U \sigma(\nabla \psi)\,\mathrm{d}s\mathrm{d}t + C \frac{m}{n} \log(1/\kappa),
\end{align}
where we have used \eqref{eq:able} and $g(u^i) \leq 2 \log(1/\kappa)$ to account for the discrepancy in dividing between $1/n_i^2$ and $1/(n/m)^2$. Combining \eqref{eq:brubeck1} and \eqref{eq:brubeck2}, we obtain \eqref{eq:forlinear}.
\end{proof}
%%%%%%%%%%%%%%%%%%%%%%%%%%%%%%%%%%%%%%%%%%%%%%%%%
Functions with bounded second derivatives are well approximated in the $L^\infty$ norm by their $m$-linearizations:

\begin{lemma} \label{lem:twicesmooth}
Let $f:U \to \mathbb{R}$ be a function defined on some union $U$ of the triangles $\rt^{m,ab}$ and $\urt^{m,ab}$ defined above. Suppose that all second-order derivatives of $f$ satisfy the bound
\begin{align*}
|f_{ss}| \leq L, \quad |f_{st}| \leq L, \quad |f_{tt}| \leq L
\end{align*}
for $(s,t) \in U$. Then if $f^m$ is its $m$-linearization we have
\begin{align*}
|f^m(s,t) - f(s,t)| \leq C L/m^2, \qquad \text{for all $(s,t) \in U$.}
\end{align*}
\end{lemma}
\begin{proof}
Follows from a Taylor expansion.
\end{proof}
%%%%%%%%%%%%%%%%%%%%%%%%%%%%%%%%%%%%%%%%%%%%%%%%
Now given any increasing function $\psi:\rt \to \mathbb{R}$ and $\kappa > 0$, define $\psi_\kappa$ by
\begin{align*}
\psi_\kappa(s,t) = \psi(s,t) + \kappa(s+t).
\end{align*}
Note that the operation $\psi \mapsto \psi_\kappa$ commutes with smoothing and $m$-linearizing, in that $(\psi_\kappa)^{\varepsilon \varepsilon m } = (\psi^{\varepsilon \varepsilon m})_\kappa$. On occasion, we unambiguously write $\psi_\kappa^{\varepsilon \varepsilon m}$ for either of these (ultimately identical) expressions. Below we will use the fact that for any increasing function $\psi:\rt \to (0,1)$, provided $\sup_{(s,t) \in \rt} \psi(s,t) + 2\kappa + \delta \leq 1$ we have
\begin{align} \label{eq:simpler}
    K_n^\infty(\psi,U,\delta) \leq K_n^\infty(\psi_\kappa,U,\delta).
\end{align}
The inequality \eqref{eq:simpler} can be proved by considering the change of variables $\phi'_x := \phi_x + \kappa (x_1 + x_2)$. 

\begin{thm} \label{thm:neuer}
 Let $\psi:\rt \to [0,1/2]$ be an increasing function. There exists $\delta_0 > 0$ such that for each $0 < \delta < \delta_0$, we can choose $\varepsilon_\delta > 0$ and $m_\delta > 0$ such that both $m_\delta$ and $\varepsilon_\delta m_\delta $ are integers, such that $\varepsilon_\delta \downarrow 0, m_\delta \uparrow +\infty$ as $\delta \downarrow 0$, and such that, for all sufficiently large $n \geq n_\delta$,
\begin{align*}
\frac{1}{n^2} \log K_n(\psi,\rt',\delta) \leq \int_{\hollow^{3\varepsilon}} (-\sigma)(\nabla \psi^{\varepsilon \varepsilon m }_\varepsilon ) + C \delta^{1/30},
\end{align*}
where $C > 0$ is universal. Specifically, we can take
\begin{align} \label{eq:goodchoices}
\varepsilon = (1+o(1))\delta^{1/25} \quad \text{and} \quad m = (1+o(1))\delta^{-1/5},
\end{align}
where $\varepsilon,m$ are chosen so that $\varepsilon m$ and $m$ are both integers.  
\end{thm}

\begin{proof}
Choose $0<\varepsilon<1/6$ such that $(\varepsilon n)^{-1}\leq \delta$. Then using Lemma \ref{lem:smoother}, for any such $\varepsilon$ we have
\begin{align} \label{eq:bearsmall1}
\frac{1}{n^2} \log K_n(\psi,\rt',\delta) \leq  \frac{1}{n^2} \log K_n( \psi^{\varepsilon \varepsilon}, \hollow^{2\varepsilon} , 4 \delta ) +C \varepsilon \log(1/\varepsilon).
\end{align}
Now noting by \eqref{eq:scuba2} that $\Psi := \psi^{\varepsilon \varepsilon}$ satisfies the conditions of Corollary \ref{cor:busan}, using \eqref{eq:busan} in \eqref{eq:bearsmall1} we can obtain
\begin{align} \label{eq:nuevo}
\frac{1}{n^2} \log K_n(\psi,\rt',\delta) \leq \frac{1}{n^2} \log K^\infty_n(\psi^{\varepsilon \varepsilon},\hollow^{3\varepsilon},C \delta^{1/3}\varepsilon^{-2/3}) +C \varepsilon \log(1/\varepsilon).
\end{align}

Let $m$ be an integer, and let $\varepsilon$ be chosen now so that $\varepsilon$ is an integer multiple of $1/m$. 
Let $\psi^{\varepsilon \varepsilon m}:\hollow^{3\varepsilon} \to \mathbb{R}$ be the $m$-linearization of $\psi^{\varepsilon \varepsilon}:\hollow^{3\varepsilon} \to \mathbb{R}$. 
The function $\psi^{\varepsilon \varepsilon}$ has all of its second derivatives bounded in absolute value by $L = C/\varepsilon^2$ by \eqref{eq:scuba3}, so accordingly, by Lemma \ref{lem:twicesmooth} we have 
\begin{align*}
|\psi^{\varepsilon\varepsilon m}(s,t) - \psi^{\varepsilon \varepsilon}(s,t) | \leq C/(\varepsilon m)^2  \qquad (s,t) \in \hollow^{3\varepsilon},
\end{align*}
for some possibly different universal $C>  0$. It follows from the triangle inequality and \eqref{eq:nuevo} that
\begin{align} \label{eq:bearsmall1b}
\frac{1}{n^2} \log K_n(\psi,\rt',\delta) \leq \frac{1}{n^2} \log K^\infty_n(\psi^{\varepsilon \varepsilon m},\hollow^{3\varepsilon}, \rho_{\delta,\varepsilon,m} )  + C \varepsilon \log(1/\varepsilon)
\end{align}
where $\rho_{\delta,\varepsilon,m} = C (\delta^{1/3}\varepsilon^{-2/3} + 1/(\varepsilon m)^2)$.

Since $\psi:\rt \to [0,1/2]$ takes values in $[0,1/2]$, so does $\psi^{\varepsilon \varepsilon m}$. In particular, provided $\rho_{\delta,\varepsilon,m} \leq 1/4$ and $\kappa \leq 1/8$, by \eqref{eq:simpler} we may replace $\psi^{\varepsilon \varepsilon m}$ with $\psi^{\varepsilon \varepsilon m}_\kappa$ in \eqref{eq:bearsmall1b} to obtain
\begin{align} \label{eq:bearsmall2}
\frac{1}{n^2} \log K_n(\psi,\rt',\delta) \leq \frac{1}{n^2} \log K^\infty_n(\psi^{\varepsilon \varepsilon m}_\kappa,\hollow^{3\varepsilon}, \rho_{\delta,\varepsilon,m} )  + C \varepsilon \log(1/\varepsilon).
\end{align}

Now recall that the derivatives of $\psi^{\varepsilon \varepsilon m}_\kappa$ are bounded above by $1/\varepsilon + \kappa$ and below by $\kappa$. In particular, setting $\varepsilon = \kappa$, using \eqref{eq:forlinear} with $\rho = \rho_{\delta,\varepsilon,m}$ and $U = \hollow^{3 \varepsilon}$ in \eqref{eq:bearsmall2} we obtain 
\begin{align} \label{eq:bearsmall3}
\frac{1}{n^2} \log K_n(\psi,\rt',\delta) \leq \int_{\hollow^{3\varepsilon}} (-\sigma)(\nabla \psi^{\varepsilon \varepsilon m }_\varepsilon ) + C \varepsilon \log(1/\varepsilon) + C(m \rho_{\delta,\varepsilon,m}/\varepsilon)^{1/2}.
\end{align}
With $\varepsilon$ and $m$ as in \eqref{eq:goodchoices}, using the definition $\rho_{\delta,\varepsilon,m} = C (\delta^{1/3}\varepsilon^{-2/3} + 1/(\varepsilon m)^2)$, it is a brief calculation to establish that the error terms in \eqref{eq:bearsmall3} are $ O(\delta^{1/30})$, completing the proof.
\end{proof}

%%%%%%%%%%%%%%%%%%%%%%%%%%%%%%%%%%%%%%%%%%%%%%%%%%%%%%%%%%%%%%%
In the following lemma, recall that when $\varepsilon$ is a multiple of $1/m$, given an increasing surface $\psi:\rt \to \mathbb{R}$, $\psi^{\varepsilon \varepsilon m }:\hollow^{2 \varepsilon} \to \mathbb{R}$ denotes the $m$-linearisation of the $\varepsilon$-double smoothing, and $\psi^{\varepsilon \varepsilon m }_\varepsilon:\hollow^{2 \varepsilon} \to \mathbb{R}$ is defined by $\psi^{\varepsilon \varepsilon m }_\varepsilon(s,t) := \psi^{\varepsilon \varepsilon m }(s,t) + \varepsilon(s+t)$. 
\begin{lemma} \label{lem:impish}
Let $\Lambda$ be a subset of the space of increasing surfaces $\psi:\rt \to [0,1]$. 
Let $\varepsilon := \varepsilon_\delta \downarrow 0$ and $m := m_\delta \uparrow \infty$ as $\delta \downarrow 0$ in a way such that $\varepsilon m \to \infty$. Then
\begin{align} \label{eq:impish}
\limsup_{\delta \downarrow 0} \sup_{ \psi \in \Lambda} \int_{\hollow^{3\varepsilon}} (-\sigma)(\nabla \psi^{\varepsilon \varepsilon m }_\varepsilon ) \leq \sup_{ \psi \in \overline{\Lambda}} \int_{ \rt} (-\sigma)(\nabla \psi).
\end{align}
\end{lemma}

\begin{proof}
Consider a sequence of elements of $\Lambda$ attaining the quantity on the left-hand side of \eqref{eq:impish}. More specifically, let $(\delta_n)_{n\geq 1}$ be a sequence of reals tending to zero, and write $\varepsilon_n := \varepsilon_{\delta_n}$ and $m_n := m_{\delta_n}$. Given a sequence $(\psi_n)_{n \geq 1}$ of elements of $\Lambda$, write $\psi'_n := \psi_{\varepsilon_n}^{\varepsilon_n \varepsilon_n m_n }$ for the associated smoothed and linearised function defined on $\hollow^{3 \varepsilon_n}$. By definition, there exists some such sequence $(\psi'_n)_{n \geq 1}$ such that
\begin{align*}
\lim_{n \to \infty} \int_{\hollow^{3 \varepsilon_n}} (-\sigma)(\nabla \psi_n') 
=
\limsup_{\delta \downarrow 0} \sup_{ \psi \in \Lambda} \int_{\hollow^{3\varepsilon}} (-\sigma)(\nabla \psi^{\varepsilon \varepsilon m }_\varepsilon ).
\end{align*}

Now by the compactness in the $L^1$ topology of the set of increasing surfaces $\psi:\rt \to [0,1]$, $(\psi_n)_{n \geq 1}$ has some subsequence $(\psi_{n_k})_{k \geq 1}$ converging to an element $\psi$ of $\overline{\Lambda}$. Thus, to complete the proof of the statement, it is sufficient to establish that
\begin{align} \label{eq:archaic}
\limsup_{k \to \infty} \int_{\hollow^{3 \varepsilon_{n_k}}} (-\sigma)(\nabla \psi_{n_k}') \leq \int_{ \rt} (-\sigma)(\nabla \psi).
\end{align}

We now claim that for each fixed $\varepsilon_0 > 0$, we have
\begin{align} \label{eq:relish}
\lim_{k\to \infty} \int_{\hollow^{\varepsilon_0}} | \psi_{n_k}' - \psi | = 0,
\end{align}
where the limit starts with $k$ sufficiently large so that $\hollow^{\varepsilon_0}$ is a subset of $\hollow^{3\varepsilon_{n_k}}$, and hence $\psi_{n_k}'$ is defined on the former. Indeed, the convergence in \eqref{eq:relish} follows from $\int_{\rt} |\psi_{n_k} - \psi| \to 0$, the easily verified bound $\int_{\hollow^{\varepsilon_0}} |\psi_\varepsilon^{\varepsilon \varepsilon m} - \psi | \leq O(1/m) + O(\varepsilon)$ (uniformly over all $\psi:\rt \to [0,1]$ and $\varepsilon_0 > 0$) and the triangle inequality. 

    Now suppose that the integral on the right-hand side of \eqref{eq:archaic} is in $(-\infty,+\infty)$. Then for all $\varepsilon_0 > 0$ we have 
\begin{align} \label{eq:relish2}
\left| \int_{ \rt - \hollow^{\varepsilon_0} } (-\sigma)(\nabla \psi) \right|  \leq f(\varepsilon_0),
\end{align}
where $f(\varepsilon_0) \downarrow 0$ as $\varepsilon_0 \downarrow 0$. 

Now by \eqref{eq:evans2}, and the fact that the Lebesgue measure of $\hollow^{\varepsilon_0} \setminus \hollow^{3 \varepsilon_{n_k}} $ is $O(\varepsilon_0)$ for any $k$, we have 
\begin{align} \label{eq:relish3}
\int_{\hollow^{\varepsilon_0} \setminus \hollow^{3 \varepsilon_{n_k}}  } (-\sigma)(\nabla \psi_{n_k}') \leq O ( \varepsilon_0 \log(1/\varepsilon_0)).
\end{align}

Now by applying Theorem \ref{thm:lowersemi} to $\hollow^{\varepsilon_0}$ in place of $\rt$, and using the convergence in \eqref{eq:relish}, we have
\begin{align} \label{eq:relish4}
\limsup_{ k \to \infty} \int_{\hollow^{\varepsilon_0}} (-\sigma)(\nabla \psi_{n_k}') \leq \int_{\hollow^{\varepsilon_0}} (-\sigma)(\nabla \psi).
\end{align}

Combining \eqref{eq:relish2}, \eqref{eq:relish3}, and \eqref{eq:relish4}, we obtain 
\begin{align} \label{eq:archaic2}
\limsup_{k \to \infty} \int_{\hollow^{3 \varepsilon_{n_k}}} (-\sigma)(\nabla \psi_{n_k}') 
\leq 
\int_{ \rt} (-\sigma)(\nabla \psi) + f(\varepsilon_0) + C (\varepsilon_0 \log (1/\varepsilon_0)).
\end{align}

Since $\varepsilon_0$ is arbitrary, this completes the proof of \eqref{eq:archaic} in the case where $\int_{ \rt} (-\sigma)(\nabla \psi)$ is finite. 

The proof in the case where $\int_{ \rt} (-\sigma)(\nabla \psi) = -\infty$ is similar; we leave the details to the reader.
\end{proof}
%%%%%%%%%%%%%%%%%%%%%%%%%%%%%%%%%%%%%%%%%%%%%%%%%%%%%%%%%%%%%%%%%%%%%%%%%%%%%%%%%%%%%%%%%%%%%%
We are now ready to complete the proof of the large deviation upper bound \eqref{eq:UPPER}. 

\begin{proof}[Proof of \eqref{eq:UPPER}]

We may assume without loss of generality that $\Gamma$ is a subset of the set of increasing surfaces taking values in $[0,1]$.
Since the set of increasing surfaces taking values in $[0,1]$ is compact, for each $\delta > 0$ there exists an integer $N_\delta$ such that for any nonempty measurable set $\Gamma$ of increasing surfaces taking values in $[0,1]$, there are elements $\psi_1,\ldots,\psi_{N_\delta}$ in $\Gamma$ for which we have the inclusion
\[
\Gamma \subseteq \bigcup_{i=1}^{N_\delta} B(\psi_i , \delta).
\]
By the union bound we have 
$V_{n,\mu}(\Gamma) \leq \sum_{i=1}^{N_\delta} V_{n,\mu}(B(\psi_i,\delta))$. 
Taking limsups, we see that for any $\delta > 0$ we have the upper bound
\begin{align*}
\limsup_{n \to \infty}\frac{1}{n^2} \log V_{n,\mu}(\Gamma) 
&\leq 
\max_{i=1}^{N_\delta} \limsup_{n \to \infty} \frac{1}{n^2} \log V_{n,\mu}(B(\psi_i,\delta))  \\
&\leq 
\sup_{\psi \in \Gamma} \limsup_{n \to \infty} \frac{1}{n^2} \log V_{n,\mu}(B(\psi,\delta)).
\end{align*} 

Now since we assume that $\mu$ has compact support, by rescaling, in proving \eqref{eq:UPPER} we may assume without loss of generality that $\mu$ is supported in $[0,1/2]$, so that $Q_\mu:[0,1] \to \mathbb{R}$ takes values in $[0,1/2]$. 
The boundary conditions $\phi_{k,k} = n Q_\mu(k/n)$ impose $\phi_x \in [0,n/2] \subseteq [0,n]$ for all $x \in \rt_n$. Using this bound, but otherwise ignoring the boundary conditions, provided that $\delta>0$ is larger than the error term appearing in \eqref{eq:dapprox}, we have
\begin{align} \label{eq:neuer2}
V_{n,\mu}(B(\psi,\delta)) \leq K_n(\psi,\rt',2\delta).
\end{align}

Applying Theorem \ref{thm:neuer} and Lemma \ref{lem:impish} to \eqref{eq:neuer2}, and subsequently taking the infimum over all $\delta > 0$, we obtain 
\begin{align} \label{eq:cooper}
\limsup_{n \to \infty}\frac{1}{n^2} \log V_{n,\mu}(\Gamma) \leq \sup_{\psi \in \overline{\Gamma}} \int_{\rt} ( -\sigma)(\nabla \psi).
\end{align} 

We are not yet finished with the proof of \eqref{eq:UPPER}: it still remains to show that the supremum may be restricted to $\psi \in \overline{\Gamma} \cap \mathcal{A}_\mu$, where $\mathcal{A}_\mu$ is the set of increasing surfaces agreeing with $Q_\mu$ on the diagonal. In this direction, define $B(\mathcal{A}_\mu,\delta)$ to be the set of increasing surfaces $\psi:\rt \to \mathbb{R}$ lying within a distance $\delta$ of an element of $\mathcal{A}_\mu := \{ \psi:\rt \to \mathbb{R} : \psi(s,s) = Q_\mu(s) ~ \forall s \in [0,1]\}$.  For any $\delta > 0$, for all sufficiently large $n$ we have
\begin{align*}
V_{n,\mu}( \Gamma) = V_{n,\mu}( \Gamma \cap B(\mathcal{A}_\mu,\delta) ),
\end{align*}
since any function $\psi:\rt \to \mathbb{R}$ satisfying $\psi(k/n,k/n) =Q_\mu(k/n)$ for all $0 \leq k \leq n$ is necessarily within an $O(1/n)$ distance of an element of $\mathcal{A}_\mu$ in the $L^1(\rt)$ metric. 

In particular, applying \eqref{eq:cooper} with $\Gamma \cap B(\mathcal{A}_\mu,\delta)$ in place of $\Gamma$, for arbitrary $\delta > 0$ we have
\begin{align*}
\limsup_{n \to \infty}\frac{1}{n^2} \log V_{n,\mu}(\Gamma) 
=
\limsup_{n \to \infty}\frac{1}{n^2} \log V_{n,\mu}(\Gamma \cap B(\mathcal{A}_\mu,\delta))  \leq 
\sup_{\psi \in \overline{\Gamma \cap B(\mathcal{A}_\mu,\delta)}} \int_{\rt} ( -\sigma)(\nabla \psi).
\end{align*} 

Recall that $\mathcal{H}[\psi] := \int_{\rt} ( -\sigma)(\nabla \psi)$.  By the upper semicontinuity of $\psi \mapsto \mathcal{H}[\psi]$, if $(F_\delta)_{\delta > 0}$ is any collection of closed sets for which $\delta_1 \leq \delta_2$ implies $F_{\delta_1} \subseteq F_{\delta_2}$ and we set $F_0 := \cap_\delta F_\delta$, we have
\begin{align*}
\limsup_{\delta \downarrow 0} \sup_{\psi \in F_\delta} \mathcal{H}[\psi] = \sup_{ \psi \in F_0} \mathcal{H}[\psi]. 
\end{align*}

Setting $F_\delta := \overline{\Gamma \cap B(\mathcal{A}_\mu,\delta)}$, since $\mathcal{A}_\mu$ is closed we have $F_0 = \bigcap_{\delta>0} \overline{\Gamma \cap B(\mathcal{A}_\mu,\delta)} = \overline{\Gamma} \cap \mathcal{A}_\mu$, which completes the proof of \eqref{eq:UPPER}.
\end{proof}
%%%%%%%%%%%%%%%%%%%%%%%%%%%%%%
\subsection{Proof of the large deviation lower bound \eqref{eq:LOWER}} \label{sec:mainlower}
%%%%%%%%%%%%%%%%%%%%%%%%%%%%%%

In this section, we prove the large deviation lower bound \eqref{eq:LOWER}, whose proof --- as is often the case in the theory of large deviations --- is easier than the corresponding upper bound \eqref{eq:UPPER}. In its proof, we will appeal to the following standard argument. 

\begin{lemma} \label{lem:finder} 
Let $\psi:\rt \to \mathbb{R}$ be an increasing surface satisfying
\begin{align*}
\int_{\rt} \sigma( \nabla \psi) \mathrm{d}s\mathrm{d}t  < \infty.
\end{align*}
Then for every $\delta > 0$, there exist arbitrarily small $\varepsilon > 0$ and arbitrarily large $m > 0$ (which may be chosen so that $\varepsilon m$ is an integer), for which there exists another increasing surface $\tilde{\psi}:\rt \to \mathbb{R}$ such that:
\begin{enumerate}
\item $\tilde{\psi}$ agrees with $\psi$ on the diagonal:
\begin{align*}
\psi(s,s) = \tilde{\psi}(s,s) \qquad \text{for } s \in [0,1].
\end{align*}
\item $\tilde{\psi}$ is close to $\psi$ in the $L^1(\rt)$ metric:
\begin{align*}
\mathrm{d}(\psi,\tilde{\psi}) \leq \delta/2.
\end{align*}
\item $\tilde{\psi}$ has surface tension not too much greater than that of $\psi$:
\begin{align} \label{eq:nottoo}
\int_{\rt} \sigma( \nabla \tilde\psi) \mathrm{d}s\mathrm{d}t  \leq \int_{\rt} \sigma( \nabla \psi) \mathrm{d}s\mathrm{d}t  + \varepsilon.
\end{align}
\item $\tilde{\psi}$ has the following favourable properties:
\begin{enumerate}
\item $\tilde{\psi}$ is $m$-linear on $\tre^\varepsilon:= \{ (s,t) \in \mathbb{R}^2 : 0 \leq t \leq t+ \varepsilon \leq s \leq 1 \}$.
\item There exists $\kappa \in (0,1)$ such that the first derivatives of $\tilde{\psi}$ satisfy
\begin{align} \label{eq:firstder}
0 < \kappa \leq \tilde{\psi}_s(s,t) \wedge \tilde{\psi}_t(s,t) \leq \tilde{\psi}_s(s,t) + \tilde{\psi}_t(s,t) \leq 1/\kappa  \qquad \text{for all $(s,t) \in \tre^\varepsilon$}.
\end{align}
\end{enumerate}
\end{enumerate}
\end{lemma}

\begin{proof}
Since the proof is technical, and not particularly insightful, we omit it. The key idea is that $\psi$ may first be smoothed away from the diagonal, and thereafter linearised away from the diagonal. 
\end{proof}

\begin{proof}[Proof of \eqref{eq:LOWER}]
Let $\psi$ be an element of $\interior{\Gamma} \cap \mathcal{A}_\mu$, where we recall that $\mathcal{A}_\mu$ is the set of increasing surfaces agreeing with the quantile function $Q_\mu$ of $\mu$ on the main diagonal. 
Since $\psi$ lies in the interior of $\Gamma$, there exists $\delta > 0$ such that $B(\psi,\delta) \subseteq \Gamma$. In particular,
\[
V_{n,\mu}(\Gamma) \geq V_{n,\mu}( B(\psi,\delta)),
\]
so that in order to establish \eqref{eq:LOWER} it is sufficient to show that whenever $\psi \in \mathcal{A}_\mu$, for all $\delta > 0$ we have
\begin{align} \label{eq:berlioz}
\liminf_{n \to \infty} \frac{1}{n^2} \log V_{n,\mu}( B(\psi,\delta)) \geq \int_{\rt} (-\sigma)(\nabla \psi).
\end{align}

For sufficiently small $\varepsilon > 0$, chosen so that $\varepsilon m$ is an integer, let $\tilde{\psi}:\rt \to \mathbb{R}$ be as in the statement of Lemma \ref{lem:finder}. Then by property (2), $B(\tilde{\psi},\delta/2) \subseteq B(\psi,\delta)$, so that 
\begin{align} \label{eq:volu0}
V_{n,\mu}( B(\psi,\delta)) \geq V_{n,\mu}( B(\tilde{\psi},\delta/2)  ).
\end{align}
We now look for a lower bound on $V_{n,\mu}( B(\tilde{\psi},\delta/2)  )$. 

Write $\tre^\varepsilon = \bigcup_{i=1}^r U^i$ as a union of right-angled triangles of side-length $1/m$, so that $(\tre^\varepsilon)_n = \bigcup_{i=1}^r U^i_n$. Recall that $\rt' := \{ (s,t) \in \mathbb{R}^2 : 0 \leq t < s \leq 1 \}$, so that $\rt'_n := \{ 0 \leq x_2 < x_1 \leq n \}$ is $\rt_n$ minus the diagonal.  
Define $W := \rt' - \tre^\varepsilon$, and let $W_n = (\rt' - \tre^\varepsilon)_n$. Recall that given $\tilde{\psi}:\rt \to \mathbb{R}$ we define $\tilde{\psi}^n:\rt_n \to [0,n]$ by
\begin{align*}
\tilde{\psi}^n_x := n \tilde{\psi}(x_1/n,x_2/n).
\end{align*}

With $\kappa$ as in \eqref{eq:firstder}, let 
\[
E_n(\tilde{\psi},\kappa) := \left\{
\phi:\rt'_n \to [0,n] \,:\,
\begin{array}{ll}
|\phi_x-\tilde{\psi}^n_x| \leq \kappa/4,
& \text{for every } x \text{ in some } \partial U^i_n, \\[2mm]
\phi_x \in
\big[
(\tilde\psi^n_{x-\mathbf e_1} \vee \tilde\psi^n_{x-\mathbf e_2}),
\tilde\psi^n_x
\big],
& \text{for every } x \in W_n 
\end{array}
\right\},
\]
where, in case $x \in W_n$ but $x-\mathbf{e}_2 \notin \rt_n'$ (i.e., when $x_2 = 0$), we simply set $\tilde{\psi}^n_{x-\mathbf{e}_2} = - \infty$. 

Note that $\tilde{\psi}(s,s) = \psi(s,s) = Q_\mu(s)$ for all $s \in [0,1]$. 

We now show that any $\phi \in E_n(\tilde{\psi},\kappa)$, while defined on $\rt_n'$, is compatible with the Gelfand--Tsetlin inequalities imposed by $\phi_{k,k} := nQ_\mu(k/n)$, $k=0,\ldots,n$ along the diagonal. Indeed, consider the conditions
\[
\phi_x \in
\big[
(\tilde\psi^n_{x-\mathbf e_1} \vee \tilde\psi^n_{x-\mathbf e_2}),
\tilde\psi^n_x
\big]
\]
corresponding to $x \in Q_n$ in the first subdiagonal $Q_n := \{ (k+1,k) : 0 \leq k \leq n-1\}$. Here, for $\phi \in E_n(\tilde{\psi},\kappa)$, for $0 \leq k \leq n-1$ we have 
\begin{align*}
\phi_{k+1,k} \in [( \tilde{\psi}^n_{k,k} \vee \tilde{\psi}^n_{k+1,k-1} ) , \tilde{\psi}^n_{k+1,k}  ] \subseteq [\tilde{\psi}^n_{k,k}, \tilde{\psi}^n_{k+1,k+1}] = [n Q_\mu(k/n), n Q_\mu((k+1)/n)].
\end{align*}
That is, any function $\phi \in E_n(\tilde{\psi},\kappa )$, by definition a function on $\rt_n'$, can be extended to $\rt_n$ with $\phi_{k,k} = nQ_\mu(k/n)$ in a way that the Gelfand--Tsetlin inequalities associated with edges between the main diagonal $D_n := \{(k,k) : 0 \leq k \leq n\}$ and the first subdiagonal $Q_n$ are satisfied. 

Moreover, it is easily verified that if $\phi \in E_n(\tilde{\psi},\kappa)$, and $\xi_n:\rt \to \mathbb{R}$ is the associated increasing surface, then
\begin{align*}
\mathrm{d}(n^{-1}\xi_n,\tilde{\psi}) := \int_{\rt} |n^{-1}\xi_n - \tilde{\psi}| \leq O( \kappa/n + 1/m ).
\end{align*}
In particular, provided $m$ is chosen sufficiently large, any $n^{-1}\xi_n$ associated with $\phi \in E_n(\tilde{\psi},\kappa)$ lies in $B(\tilde{\psi},\delta/2)$, and therefore, using the definition of $V_{n,\mu}(\Gamma)$, for all sufficiently large $n$ we have
\begin{align} \label{eq:volu1}
V_{n,\mu} ( B(\tilde{\psi},\delta/2) ) \geq  \mathcal{E}_n := \int_{ \mathbb{R}^{\righttriangle_n'} } \mathrm{1} \left\{ \phi \in E_n(\tilde{\psi},\kappa) \right\} \mathrm{1}^{\rt'_n}_{\mathrm{GT}}(\phi) \prod_{ x \in \righttriangle_n' } \mathrm{d} \phi_x .
\end{align}

Note that every $\phi \in E_n(\tilde{\psi},\kappa)$ automatically satisfies the Gelfand--Tsetlin inequalities on $W_n$, and since $\kappa$ is smaller than the derivative of $\tilde{\psi}$ anywhere on $\rt_n' - W_n$, the Gelfand--Tsetlin inequalities between neighbouring $U^i_n$ and $U^j_n$, or between any $U^i_n$ neighbouring $W_n$, are also automatically satisfied. In other words, for $\phi \in E_n(\tilde{\psi},\kappa)$, in place of the union inequality \eqref{eq:unioneq} for Gelfand--Tsetlin indicators we in fact have the equality
\begin{align*}
\mathrm{1}_{\mathrm{GT}}^{\rt_n}(\phi) = \mathrm{1}_{\mathrm{GT}}^{W_n}(\phi) \prod_{i=1}^r \mathrm{1}_{\mathrm{GT}}^{U^i_n}(\phi) \qquad \text{for }\phi \in E_n(\tilde{\psi},\kappa).
\end{align*}
It follows that with $\mathcal{E}_n$ as in \eqref{eq:volu1} we have
\begin{align} \label{eq:ekwal}
\mathcal{E}_n = \mathcal{L}_n \prod_{i=1}^r T_i,
\end{align}
where 
\begin{align*}
\mathcal{L}_n := \prod_{x \in W_n} L_x, \qquad L_x = \tilde{\psi}_x^n - (\tilde{\psi}_{x-\mathbf{e}_1}^n \vee \tilde{\psi}_{x-\mathbf{e}_2}^n) 
\end{align*}
is the product over $x \in W_n$ of the lengths $L_x$ of the intervals restricting the position of $\phi_x$, and where
\begin{align*}
T_i := \int_{ \mathbb{R}^{U^i_n} } \mathrm{1} \left\{ |\phi_x - \tilde{\psi}_x^n | \leq \kappa/4 ~ \forall x \in \partial U^i_n \right\} \mathrm{1}_{\mathrm{GT}}^{U_n^i}(\phi) \prod_{x \in U^i_n} \mathrm{d}\phi_x.
\end{align*}

We now look for lower bounds on the terms $\mathcal{L}_n$ and $T_1,\ldots,T_r$ in \eqref{eq:ekwal}. In the direction of the former, we note that since $\tilde \psi^n_x := n \tilde\psi(x_1/n,x_2/n)$, by using Fatou's lemma one can verify that
\begin{align*}
\liminf_{ n \to \infty} \frac{1}{n^2} \log \mathcal{L}_n \geq \int_{ \rt - \tre^\varepsilon} \log ( \tilde{\psi}_s \wedge \tilde{\psi}_t ),
\end{align*}
where here, $\tilde{\psi}_s$ and $\tilde{\psi}_t$ denote derivatives with respect to $s$ and $t$. 

Now, using Lemma \ref{lem:sigmabound} (which we recall states that $- \sigma(u_1,u_2)$ differs from $\log(u_1 \wedge u_2)$ by at most a constant for all $u_1,u_2 > 0$), and the fact that $\rt-\tre^\varepsilon$ has Lebesgue measure at most $O(\varepsilon)$, we may alternatively write 
\begin{align} \label{eq:ekwal2}
\liminf_{ n \to \infty} \frac{1}{n^2} \log \mathcal{L}_n \geq - O(\varepsilon) +  \int_{ \rt - \tre^\varepsilon} (-\sigma)(\nabla \tilde{\psi}).
\end{align}
We emphasise the $O(\varepsilon)$ term is universal.

We turn to finding lower bounds on the $T_i$. 
Note that $\tilde{\psi}^n:\partial U^i_n \to \mathbb{R}$ is $\kappa$-increasing since the derivatives of $\tilde{\psi}$ are bounded below by $\kappa > 0$. It follows that any function $\phi:\partial U^i_n \to \mathbb{R}$ satisfying $|\phi_x - \tilde{\psi}_x^n | \leq \kappa/4$ for $x \in \partial U^i_n$ is certainly $\kappa/2$-increasing. Now, letting $n_i$ be the side-length of the triangle $U^i_n$, we can bound $T_i$ below by 
\begin{align} \label{eq:Tilower}
T_i \geq (\kappa/4)^{ \# \partial U^i_n } \inf_{ \varphi \in B_{n_i}(u^i,\kappa/4,\kappa/2) } T_{n_i}( \varphi),
\end{align}
where we recall our notation from Section \ref{sec:ST2proof} that $B_{n_i}(u^i,\kappa/4,\kappa/2)$ denotes the set of increasing functions $\varphi:\partial \rt_{n_i} \to \mathbb{R}$ that lie within $\kappa/4$ of the linear function $\varphi^{u^i}$, and are $\kappa/2$-increasing. 
Now using the lower bound \eqref{eq:lower1} from Theorem \ref{thm:upperlower}, for all sufficiently large $n$ we have 
\begin{align*}
\frac{1}{n_i^2/2} \log \inf_{ \varphi \in B_{n_i}(u^i,\kappa/4,\kappa/2) } T_{n_i}( \varphi) \geq (-\sigma)(u^i_1,u^i_2) - \rho_i ( \log(1/\rho_i) + g(u^i)),
\end{align*}
where $\rho_i = \frac{ \kappa/4 + u_1^i + u_2^i }{ (u_1^i \wedge u_2^i) n_i }$. Note by \eqref{eq:firstder} that $\rho_i \leq \frac{2m}{\kappa^2 n}$ and $g(u^i) \leq 2 \log (1/\kappa)$. 

Now using the fact that $n_i = (1+o(1))n/m$, the fact that $g(u^i) \leq C \log (1/\kappa)$, and the fact that $\tre^\varepsilon$ is covered by $r = O(m^2)$ triangles, by \eqref{eq:Tilower} we have
\begin{align} \label{eq:ekwal3}
\frac{1}{n^2} \log \prod_{i=1}^r T_i &\geq \frac{m}{n} \log (\kappa/4) + \frac{1}{2 m^2}\sum_{i=1}^r (-\sigma)(u_1^i, u_2^i) - \frac{2m}{\kappa^2 n}( \log(n/m) + \log (1/\kappa)) \nonumber \\
&= \int_{\tre^\varepsilon} (- \sigma)(\nabla \tilde{\psi})  - \frac{2m}{\kappa^2 n}( \log(n/m) + \log (1/\kappa)).
\end{align}

Combining \eqref{eq:volu0}, \eqref{eq:volu1}, \eqref{eq:ekwal}, \eqref{eq:ekwal2}, and \eqref{eq:ekwal3} to obtain the first inequality below, and then using \eqref{eq:nottoo} to obtain the second, we have 
\begin{align*}
\liminf_{n \to \infty} \frac{1}{n^2} \log V_{n,\mu}(B(\psi,\delta)) \geq \int_{\rt} (- \sigma)(\nabla \tilde{\psi}) - O(\varepsilon) \geq \int_{\rt} (- \sigma)(\nabla \psi) - O(\varepsilon).
\end{align*}
Now since the $O(\varepsilon)$ term is universal and $\varepsilon > 0$ is arbitrary, that completes the proof of \eqref{eq:berlioz}, and by extension, \eqref{eq:LOWER}.
\end{proof}

We emphasise that with \eqref{eq:UPPER} and \eqref{eq:LOWER} now established, the proof of Theorem \ref{thm:C} is complete. 

\subsection{Completion of the proof of Theorem \ref{thm:A}}

We now complete the proof of Theorem \ref{thm:A}. Recall that in Section \ref{sec:free} we established that free compression $F^*$ satisfies the Euler--Lagrange equations associated with the compression entropy functional, and that $\mathcal{H}[F^*] = \chi[\mu]$. The final step is the following:

\begin{thm} \label{thm:loop}
We have $\sup \{ \mathcal{H}[F] : \text{$F$ is a compression of $\mu$} \} = \chi[\mu]$.
\end{thm}

\begin{proof}
Recall that if $\psi$ is the increasing surface associated with a compression $F$, we write $\mathcal{H}[\psi] = \mathcal{H}[F]$. 
Setting $\Gamma$ to be the space of \emph{all} integrable increasing surfaces in the setting of Theorem \ref{thm:C}, we see that $\inf_{ \psi \in \mathcal{A}_\mu } I_\mu[\psi] = 0,$
which is equivalent to $
\sup_{ F \in \mathcal{A}_\mu} \mathcal{H}[F] = \chi[\mu]$. 
% However, we saw in Theorem \ref{thm:compressionachieve} that if $F^*$ is the free compression of $\mu$, then $\mathcal{H}[F^*] = \chi[\mu]$. The result follows. 
\end{proof}

Theorem \ref{thm:A} now follows from Theorem \ref{thm:compressionachieve}, Theorem \ref{thm:loop} and the strict concavity of the compression entropy functional.

We assure the reader that there is no issue of circularity: at no point in establishing Theorem \ref{thm:B} or Theorem \ref{thm:C} did we assume any part of Theorem \ref{thm:A}. 

%%%%%%%%%%%%%%%%%%%%%%%%%%%%%%%%%%%%%%%%%%%%%%%%%%%%%%%%
%%%%%%%%%%%%%%%%%%%%%%%%%%%%%%%%%%%%%%%%%%%%%%%%%%%%%%%%
%%%%%%%%%%%%%%%%%%%%%%%%%%%%%%%%%%%%%%%%%%%%%%%%%%%%%%%%
\appendix

\section{Proof of Lemma \ref{lem:Hbound}} \label{sec:Hboundproof}

In this section, we prove Lemma \ref{lem:Hbound}. Recall that $H^{\beta,\theta_2,T}:(-1,1) \times \mathbb{Z} \to \mathbb{C}$ is given by 
\begin{align} \label{eq:detnewkernelagain}
H^{\beta,\theta_2,T}(t,h) =  \frac{T}{n} \sum_{ z^n = (-1)^{\theta_2} } z^{1-h} \frac{ e^{ - (\beta+Tz)[t] }}{ 1 - e^{ - (\beta+Tz)} } \qquad [t] := t + \mathrm{1}_{\{t<0\}}.
\end{align}
Let $p \in (-1,1)$ and $T > 0$. 
Lemma \ref{lem:Hbound} states that whenever $(p,T)$ is good (meaning that $T\sqrt{1-p^2}/\pi$ lies at least a distance $1/4$ away from an integer) we have the uniform bound 
\begin{align} \label{eq:predet4again}
  \sup_{t \in (-1,1),\,h \in \mathbb{Z} }\Big| H^{pT+\theta_1\pi \iota , \theta_2, T} (t,h)\Big| \leq C \frac{ T}{\sqrt{1-p^2}}.
\end{align}

%%%%%%%%%%%%%%%%%%%%%%%%%%%%%%%%%%%%
We now prove \eqref{eq:predet4again}. Setting $\lambda := pT$ and using \eqref{eq:detnewkernel} together with the triangle inequality and the fact that all $n$ of the $n^{\text{th}}$ roots of $(-1)^{\theta_2}$ lie on the circle $\{z \in \mathbb{C} : |z|=1\}$, we have
\begin{align} \label{eq:predet}
  \sup_{t \in (-1,1),\,h \in \mathbb{Z} }\Big| H^{pT+\theta_1\pi \iota , \theta_2, T} (t,h)\Big| 
  & = \sup_{t \in (-1,1),\,h \in \mathbb{Z} }\Bigg| \frac{T}{n} \sum_{ z^n = (-1)^{\theta_2} } z^{1-h} \frac{ e^{ - (pT+\theta_1\pi \iota+Tz)[t] }}{ 1 - e^{ - (pT+\theta_1\pi \iota+Tz)} }\Bigg| \nonumber \\
  & \leq T\sup_{|z|=1}  \sup_{t \in [0,1)}  \frac{ |e^{ - (T(p+z)+\theta_1\pi \iota)} |^t}{ |1 - e^{ - (T(p+z)+\theta_1\pi \iota)}| } \nonumber \\
  & \leq T \sup_{|z|=1} \max_{\kappa \in \{-1,1\} } \frac{ 1}{ |1 - e^{ -  \kappa (T(p+z)+\theta_1\pi \iota)}| },
\end{align}
where to obtain the final inequality above we used the fact that the supremum in $t$ is either attained at $t = 0$ or as $t \uparrow 1$. Now for $\omega \in \mathbb{C}$ and $\kappa \in \{-1,1\}$ we have the bound 
\begin{align*}
\frac{1}{|1 -e^{- \kappa \omega}|} \leq C + \frac{C}{ \inf_{ k \in \mathbb{Z} } | \omega - 2 \pi k \iota |},
\end{align*}
for some sufficiently large universal constant $C > 0$.

In particular, after accounting for $\kappa = \pm 1$ and either $\theta_1 = 0$ or $\theta_1 = 1$ we have
\begin{align} \label{eq:predet2}
\max_{\kappa \in \{-1,1\} } \frac{ 1}{ |1 - e^{ -  \kappa (T(p+z)+\theta_1\pi \iota)}| } \leq C + \frac{C}{\inf_{k \in \mathbb{Z}} | T(p+z) - \pi k \iota | }.
\end{align}
Combining \eqref{eq:predet} and \eqref{eq:predet2} we obtain

\begin{align} \label{eq:predet90}
  \sup_{t \in (-1,1),\,h \in \mathbb{Z} }\Big| H^{pT+\theta_1\pi \iota , \theta_2, T} (t,h)\Big| 
\leq CT \left( 1 + \frac{1}{ \inf_{|z| = 1 } \inf_{k \in \mathbb{Z}} | T(p+z) - \pi k \iota| } \right).
\end{align}
Our next lemma controls the denominator in \eqref{eq:predet90} for good $(p,T)$. 

\begin{lemma}
If $(p,T)$ is good, then 
\begin{align} \label{eq:predet3} 
\inf_{|z| = 1 } \inf_{k \in \mathbb{Z}} | T(p+z) - \pi k \iota| \geq c \sqrt{1-p^2}
\end{align}
for some universal constant $c>0$.
\end{lemma}
\begin{proof}
We will sketch the proof. The set $\{ T(p+z) : |z| = 1 \}$ intersects the imaginary axis at the point $\pm \iota T \sqrt{1-p^2}$. By the goodness of $(p,T)$, this point lies at least a distance $\pi/4$ from any points of the form $\pi \iota k$ where $k \in \mathbb{Z}$. The gradient of the curve at this intersection is $\sqrt{1-p^2}$. In particular, the closest this curve comes to any element of $\{ \pi \iota k: k \in \mathbb{Z} \}$ is of the order at least $\sqrt{1-p^2}$. 
\end{proof}

Using \eqref{eq:predet3} in \eqref{eq:predet90} we obtain \eqref{eq:predet4again}, as required.

\section*{Acknowledgements}
JP has been supported by
the German Research Foundation under DFG project 516672205. Part of this work was carried out during SJ’s visit to JP at the University of Passau. The authors gratefully acknowledge the excellent working conditions provided by the University of Passau.

% The authors are grateful to Octavio Arizmendi, Colin McSwiggen, Hariharan Narayanan and Jenny Pi for their valuable comments.

%%%%%%%%%%%%%%%%%%%%%%%%%%%%%%%%%%%%%%%%%%%%%%%%%%%%%%%%
%%%%%%%%%%%%%%%%%%%%%%%%%%%%%%%%%%%%%%%%%%%%%%%%%%%%%%%%
%%%%%%%%%%%%%%%%%%%%%%%%%%%%%%%%%%%%%%%%%%%%%%%%%%%%%%%%


\begin{thebibliography}{99}

 \bibitem{AGZ}
 \textsc{Anderson, G. W., Guionnet, A., \& Zeitouni, O.} (2010).
 \emph{An introduction to random matrices}.
 Cambridge University Press.

\bibitem{AJ}
\textsc{Arizmendi, O. \& Johnston, S. G. G.} (2023). Free probability via entropic optimal transport. arXiv preprint arXiv:2309.12196.

\bibitem{AGA}
\textsc{Artstein-Avidan, S., Giannopoulos, A., \& Milman, V. D.} (2021).
\emph{Asymptotic geometric analysis. Part II}.
Mathematical Surveys and Monographs, Vol. 261, American Mathematical Society.

% \bibitem{BCEHK}
% \textsc{Bao, Z., Cipolloni, G., Erdős, L., Henheik, J., \& Kolupaiev, O.} (2025).
% Decorrelation transition in the Wigner minor process.
% \emph{Probability Theory and Related Fields}, 1--52.

\bibitem{baryshnikov}
\textsc{Baryshnikov, Y.} (2001).
GUEs and queues.
\emph{Probability Theory and Related Fields} 119(2), 256--274.

\bibitem{BBatoms}
\textsc{Belinschi, S. T. \& Bercovici, H.} (2004).
Atoms and regularity for measures in a partially defined free convolution semigroup.
\emph{Mathematische Zeitschrift} 248(4), 665--674.

\bibitem{BGH}
\textsc{Belinschi, S., Guionnet, A., \& Huang, J.} (2022).
Large deviation principles via spherical integrals.
\emph{Probability and Mathematical Physics} 3(3), 543--625.

\bibitem{BAG}
\textsc{Ben Arous, G. \& Guionnet, A.} (1997).
Large deviations for Wigner's law and Voiculescu's non-commutative entropy.
\emph{Probability Theory and Related Fields} 108(4), 517--542.

\bibitem{BP}
 \textsc{Bercovici, H., Pata, V., \& Biane, P.} (1999).
 Stable laws and domains of attraction in free probability theory.
 \emph{Annals of Mathematics} 149(3), 1023--1060.

\bibitem{BV}
\textsc{Bercovici, H. \& Voiculescu, D.} (1993). Free convolution of measures with unbounded support. Indiana University Mathematics Journal, 42(3), 733-773.

\bibitem{BBS}
\textsc{Berestycki, J., Brunet, É., \& Shi, Z.} (2017).
Accessibility percolation with backsteps.
\emph{ALEA} 14, 45--62.

\bibitem{berestycki}
\textsc{Berestycki, N. \& Powell, E.} (2025).
\emph{Gaussian free field and Liouville quantum gravity}.
Cambridge University Press.


\bibitem{BG}
\textsc{Borodin, A. \& Gorin, V.} (2015).
General $\beta$-Jacobi corners process and the Gaussian free field.
\emph{Comm. Pure Appl. Math.} \textbf{68}(10), 1774--1844.

\bibitem{boutillier}
\textsc{Boutillier, C.} (2009).
The bead model and limit behaviors of dimer models.
\emph{Annals of Probability} 37(1), 107--142.

\bibitem{BuG}
\textsc{Bufetov, A. \& Gorin, V.} (2018).
Fluctuations of particle systems determined by Schur generating functions.
\emph{Adv. Math.} \textbf{338}, 702--781.

\bibitem{CL}
\textsc{Carlen, E. \& Loss, M.} (1992). Competing symmetries, the logarithmic HLS inequality and Onofri's inequality on $S^n$. \emph{Geometric \& Functional Analysis GAFA}, 2(1), 90-104.

% \bibitem{BGVV}
% \textsc{Brazitikos, S., Giannopoulos, A., Valettas, P., \& Vritsiou, B.-H.} (2014).
% \emph{Geometry of isotropic convex bodies}.
% Mathematical Surveys and Monographs, Vol. 196, American Mathematical Society.

\bibitem{CGW}
\textsc{Chandra, A., Gunaratnam, T. S., \& Weber, H.} (2022).
Phase transitions for $\phi^4_3$.
\emph{Communications in Mathematical Physics} 392(2), 691--782.

\bibitem{CKP}
\textsc{Cohn, H., Kenyon, R., \& Propp, J.} (2001).
A variational principle for domino tilings.
\emph{Journal of the American Mathematical Society} 14(2), 297--346.

% \bibitem{CM}
% \textsc{Collins, B., \& Metcalfe, A.} (2019).
% Gelfand--Tsetlin polytopes and random contractions away from the limiting shapes.
% \emph{arXiv:1911.00842}.

\bibitem{CMZ}
\textsc{Coquereaux, R., McSwiggen, C., \& Zuber, J.-B.} (2020).
On Horn's problem and its volume function.
\emph{Communications in Mathematical Physics} 376(3), 2409--2439.

\bibitem{CM}
\textsc{Collins, B. \& Metcalfe, A.} (2023).
Gelfand--Tsetlin polytopes and random contractions away from the limiting shape.
\emph{Ann. Fac. Sci. Toulouse Math.} \textbf{32}(3), 423--533.

\bibitem{sfree3}
\textsc{Collins, B., Mingo, J. A. \& Śniady, P., Speicher, R.} (2007). Second order freeness and fluctuations of random matrices. III: Higher order freeness and free cumulants. \emph{Documenta Mathematica}, 12, 1-70.

\bibitem{CZ}
\textsc{Coquereaux, R., \& Zuber, J.-B.} (2018).
From orbital measures to Littlewood--Richardson coefficients and hive polytopes.
\emph{Annales de l'Institut Henri Poincaré D} 5(3), 339--386.

\bibitem{mattdiana}
\textsc{De Armas Bellon, D., \& Roberts, M. I.} (2026).
Accessibility percolation with rough Mount Fuji labels.
\emph{arXiv:2603.29561}.

\bibitem{DZ}
\textsc{Dembo, A., \& Zeitouni, O.} (2010).
\emph{Large deviations techniques and applications}.
Springer.

\bibitem{DGI}
\textsc{Deuschel, J.-D., Giacomin, G., \& Ioffe, D.} (2000).
Large deviations and concentration properties for $\nabla \varphi$ interface models.
\emph{Probability Theory and Related Fields} 117(1), 49--111.

\bibitem{DM1}
\textsc{Duse, E. \& Metcalfe, A.} (2015).
Asymptotic geometry of discrete interlaced patterns: Part I.
\emph{Int. J. Math.} \textbf{26}(11), 1550093.

\bibitem{DM2}
\textsc{Duse, E. \& Metcalfe, A.} (2020).
Asymptotic geometry of discrete interlaced patterns: Part II.
\emph{Ann. Inst. Fourier (Grenoble)} \textbf{70}(1), 375--436.



\bibitem{folland}
\textsc{Folland, G. B.} (1999).
\emph{Real analysis: modern techniques and their applications}.
John Wiley \& Sons.

% \bibitem{FNforest}
% \textsc{Forrester, P. J., \& Nagao, T.} (2011).
% Determinantal correlations for classical projection processes.
% \emph{Journal of Statistical Mechanics: Theory and Experiment} 2011(08), P08011.

\bibitem{funaki}
\textsc{Funaki, T.} (2005).
Stochastic interface models.
In: \emph{Lectures on Probability Theory and Statistics},
Lecture Notes in Mathematics, Vol. 1869, Springer, 103--274.

\bibitem{FN}
\textsc{Funaki, T., \& Nishikawa, T.} (2001).
Large deviations for the Ginzburg--Landau $\nabla\phi$-interface model.
\emph{Probability Theory and Related Fields} 120(4), 535--568.

\bibitem{FS}
\textsc{Funaki, T., \& Spohn, H.} (1997).
Motion by mean curvature from the Ginzburg--Landau $\nabla\phi$ interface model.
\emph{Communications in Mathematical Physics} 185(1), 1--36.

\bibitem{hari1}
\textsc{Gangopadhyay, A. \& Narayanan, H.} (2025).
On the randomized Horn problem and the surface tension of hives.
\emph{Annales de l'Institut Henri Poincaré D} (published online first).

\bibitem{GSS}
\textsc{Garza-Vargas, J., Srivastava, N., \& Stier, Z.} (2026). Finite free information inequalities. arXiv preprint arXiv:2602.15822.

\bibitem{GT1950}
\textsc{Gel'fand, I. M. \& Tsetlin, M. L.} (1950).
Finite-dimensional representations of the group of unimodular matrices.
\emph{Doklady Akademii Nauk SSSR} 71, 825--828.

\bibitem{gordenko}
\textsc{Gordenko, A.} (2020).
Limit shapes of large skew Young tableaux and a modification of the TASEP process.
\emph{arXiv:2009.10480}.

\bibitem{gorin}
\textsc{Gorin, V.} (2021).
\emph{Lectures on random lozenge tilings}.
Cambridge University Press.

\bibitem{gustavsson}
\textsc{Gustavsson, J.} (2005).
Gaussian fluctuations of eigenvalues in the GUE.
\emph{Annales de l'Institut Henri Poincaré Probabilités et Statistiques} 41, 151--178.

\bibitem{HS}
\textsc{Helffer, B. \& Sjöstrand, J.} (1994). On the correlation for Kac-like models in the convex case. \emph{Journal of statistical physics,} 74 (1), 349-409.

\bibitem{horn}
\textsc{Horn, A.} (1962).
Eigenvalues of sums of Hermitian matrices.
\emph{Pacific Journal of Mathematics} 12(1), 225--241.

\bibitem{HK}
\textsc{Hoskins, J. \& Kabluchko, Z.} (2023).
Dynamics of zeroes under repeated differentiation.
\emph{Experimental Mathematics} 32(4), 573--599.

\bibitem{johansson}
\textsc{Johansson, K.} (1998). On fluctuations of eigenvalues of random Hermitian matrices,
\emph{Duke Mathematical Journal} 91(1), 151-204.

\bibitem{johnston}
\textsc{Johnston, S. G. G.} (2025).
Continuous Kasteleyn theory for the bead model.
\emph{Journal of the European Mathematical Society} (published online first).

\bibitem{johnstonGUE}
\textsc{Johnston, S. G. G.} (2026).
Log-concavity and concentration bounds for a single gap between GUE eigenvalues.
\emph{arXiv:2601.04869}.

% \bibitem{JM}
% \textsc{Johnston, S. G. G., \& McSwiggen, C.} (2024).
% On the limiting Horn inequalities.
% \emph{arXiv:2410.08907}.

\bibitem{JO}
\textsc{Johnston, S. G. G. \& O'Connell, N.} (2020).
Scaling limits for non-intersecting polymers and Whittaker measures.
\emph{Journal of Statistical Physics} 179(2), 354--407.

\bibitem{kabluchko}
\textsc{Kabluchko, Z.} (2015).
Critical points of random polynomials with independent identically distributed roots.
\emph{Proceedings of the American Mathematical Society} 143(2), 695--702.

\bibitem{KO}
\textsc{Kenyon, R. \& Okounkov, A.} (2007).
Limit shapes and the complex Burgers equation.
\emph{Acta Mathematica} 199(2), 263--302.

\bibitem{KP}
\textsc{Kenyon, R. \& Prause, I.} (2022).
Gradient variational problems in $\mathbb{R}^2$.
\emph{Duke Mathematical Journal} 171(14), 3003--3022.

\bibitem{KP2}
\textsc{Kenyon, R. \& Prause, I.} (2024).
Limit shapes from harmonicity: dominos and the five vertex model.
\emph{Journal of Physics A: Mathematical and Theoretical} 57(3), 035001.

\bibitem{klyachko}
\textsc{Klyachko, A. A.} (1998).
Stable bundles, representation theory and Hermitian operators.
\emph{Selecta Mathematica} 4, 419--445.

\bibitem{KT}
\textsc{Knutson, A. \& Tao, T.} (1999).
The honeycomb model of $\mathrm{GL}_n(\mathbb{C})$ tensor products I: Proof of the saturation conjecture.
\emph{Journal of the American Mathematical Society} 12(4), 1055--1090.

\bibitem{LV}
\textsc{Lovász, L. \& Vempala, S.} (2007).
The geometry of logconcave functions and sampling algorithms.
\emph{Random Structures \& Algorithms} 30(3), 307--358.

\bibitem{MP}
\textsc{Magazinov, A. \& Peled, R.} (2022).
Concentration inequalities for log-concave distributions with applications to random surface fluctuations.
\emph{The Annals of Probability} 50(2), 735--770.

\bibitem{MSS}
\textsc{Marcus, A. W., Spielman, D. A., \& Srivastava, N.} (2022). Finite free convolutions of polynomials. \emph{Probability Theory and Related Fields}, 182(3), 807-848.

\bibitem{metcalfe}
\textsc{Metcalfe, A. P.} (2013).
Universality properties of Gelfand--Tsetlin patterns.
\emph{Probability Theory and Related Fields} 155, 303--346.

\bibitem{miller}
\textsc{Miller, J.} (2011).
Fluctuations for the Ginzburg--Landau $\nabla\phi$ interface model on a bounded domain.
\emph{Communications in Mathematical Physics} 308(3), 591--639.

\bibitem{sfree2}
\textsc{Mingo, J. A., Śniady, P., \& Speicher, R.} (2007). Second order freeness and fluctuations of random matrices: II. Unitary random matrices. \emph{Advances in Mathematics}, 209(1), 212-240.

\bibitem{MS}
\textsc{Mingo, J. A. \& Speicher, R.} (2017).
\emph{Free probability and random matrices}.
Springer.

\bibitem{sfree1}
\textsc{Mingo, J. A. \& Speicher, R.} (2006). Second order freeness and fluctuations of random matrices: I. Gaussian and Wishart matrices and cyclic Fock spaces. \emph{Journal of Functional Analysis,} 235(1), 226-270.
% \bibitem{montgomery}
% \textsc{Montgomery, H. L.} (1973).
% The pair correlation of zeros of the zeta function.
% In \emph{Analytic Number Theory, Proc. Sympos. Pure Math.} 24, 181--193.

\bibitem{hari2}
\textsc{Narayanan, H.} (2025).
On the limit of random hives with GUE boundary conditions.
\emph{arXiv:2502.06414}.

\bibitem{hari3}
\textsc{Narayanan, H.} (2025).
Random discrete concave functions on an equilateral lattice with periodic Hessians.
\emph{Journal of Fourier Analysis and Applications} 31(4), 50.

\bibitem{hari4}
\textsc{Narayanan, H.} (2026). Hives from deformed GUE minor processes. arXiv preprint arXiv:2607.04138.

\bibitem{NarS}
\textsc{Narayanan, H. \& Sheffield, S.} (2024).
Large deviations for random hives and the spectrum of the sum of two random matrices.
\emph{The Annals of Probability} 52(3), 1093--1152.

\bibitem{NST}
\textsc{Narayanan, H., Sheffield, S., \& Tao, T.} (2024).
Sums of GUE matrices and concentration of hives from correlation decay of eigengaps.
\emph{Probability Theory and Related Fields} 190(3--4), 1121--1165.

\bibitem{NS}
\textsc{Nica, A. \& Speicher, R.} (1996).
On the multiplication of free $N$-tuples of noncommutative random variables.
\emph{American Journal of Mathematics} 118(4), 799--837.

% \bibitem{nishikawa}
% \textsc{Nishikawa, T.} (2003).
% Hydrodynamic limit for the Ginzburg--Landau $\nabla\varphi$ interface model with boundary conditions.
% \emph{Probability Theory and Related Fields} 127(2), 205--227.

\bibitem{oconnell}
\textsc{O'Connell, N.} (2012).
Directed polymers and the quantum Toda lattice.
\emph{Annals of Probability} 40(2), 437--458.

\bibitem{OSZ}
\textsc{O'Connell, N., Seppäläinen, T., \& Zygouras, N.} (2014).
Geometric RSK correspondence, Whittaker functions and symmetrized random polymers.
\emph{Inventiones Mathematicae} 197(2), 361--416.

% \bibitem{odlyzko}
% \textsc{Odlyzko, A. M.} (1987).
% On the distribution of spacings between zeros of the zeta function.
% \emph{Mathematics of Computation} 48, 273--308.

\bibitem{petrov}
\textsc{Petrov, L.} (2014).
Asymptotics of random lozenge tilings via Gelfand--Tsetlin schemes.
\emph{Probab. Theory Related Fields} \textbf{160}(3--4), 429--487.

\bibitem{santa}
\textsc{Santambrogio, F.} (2015). \emph{Optimal Transport for Applied Mathematicians: Calculus of Variations, PDEs, and Modeling}. Volume 87 of Progress in Nonlinear Differential Equations and Their Applications.

\bibitem{sheffield}
\textsc{Sheffield, S.} (2005).
\emph{Random surfaces}.
Astérisque 304.

\bibitem{ST}
\textsc{Shlyakhtenko, D. \& Tao, T.} (2022).
Fractional free convolution powers.
\emph{Indiana University Mathematics Journal} 71(6), 2551--2594.

 \bibitem{steinerberger}
\textsc{Steinerberger, S.} (2019).
A nonlocal transport equation describing roots of polynomials under differentiation.
 \emph{Proceedings of the American Mathematical Society} 147(11), 4733--4744.

\bibitem{steinerbergerfree}
\textsc{Steinerberger, S.} (2023).
Free convolution powers via roots of polynomials.
\emph{Experimental Mathematics} 32(4), 567--572.

\bibitem{sun}
\textsc{Sun, W.} (2018).
Dimer model, bead model and standard Young tableaux: finite cases and limit shapes.
\emph{arXiv:1804.03414}.
\bibitem{taoblog}

\textsc{Tao, T.} (2007).
\href{https://terrytao.wordpress.com/2007/04/19/open-question-what-is-a-quantum-honeycomb/}{What is a quantum honeycomb? (blog post)}.

% \bibitem{taoRMT}
% \textsc{Tao, T.} (2012).
% \emph{Topics in random matrix theory}.
% American Mathematical Society.

\bibitem{taoblog2}
\textsc{Tao, T.} (2023).
\href{https://terrytao.wordpress.com/2023/06/20/sums-of-gue-matrices-and-concentration-of-hives-from-correlation-decay-of-eigengaps/}{Sums of GUE matrices and concentration of hives from correlation decay of eigengaps (blog post)}.

\bibitem{taoGUE}
\textsc{Tao, T.} (2025).
On the distribution of eigenvalues of GUE and its minors at fixed index.
\emph{Random Matrices: Theory and Applications} 14(4), 2550023.

\bibitem{Vcon1}
\textsc{Voiculescu, D.} (1986).
Addition of certain non-commuting random variables.
\emph{Journal of Functional Analysis} 66, 323--346.

\bibitem{Vcon2}
\textsc{Voiculescu, D.} (1987).
Multiplication of certain non-commuting random variables.
\emph{Journal of Operator Theory} 18(2), 223--235.

\bibitem{Vcon3}
\textsc{Voiculescu, D.} (1991).
Limit laws for random matrices and free products.
\emph{Inventiones Mathematicae} 104(1), 201--220.

\bibitem{Vent1}
\textsc{Voiculescu, D.} (1993).
The analogues of entropy and of Fisher's information measure in free probability theory, I.
\emph{Communications in Mathematical Physics} 155(1), 71--92.


\bibitem{Vent2}
\textsc{Voiculescu, D.} (1994).
The analogues of entropy and of Fisher's information measure in free probability theory, II.
\emph{Inventiones Mathematicae} 118(3), 411--440.

\bibitem{Vent3}
\textsc{Voiculescu, D.} (1996).
The analogues of entropy and of Fisher's information measure in free probability theory III:
The absence of Cartan subalgebras.
\emph{Geometric and Functional Analysis} 6(1), 172--199.

\bibitem{VIMRN}
\textsc{Voiculescu, D.} (1998). A strengthened asymptotic freeness result for random matrices with applications to free entropy. 
\emph{IMRN: International Mathematics Research Notices}, 1998(1), 41–63.
\end{thebibliography}
\end{document}